\documentclass[a4paper,12pt,twoside,openright,pdftex]{report}
\pdfoutput=1
\usepackage[applemac]{inputenc} 
\usepackage[T1]{fontenc} 
\usepackage{lmodern}
\usepackage[frenchb,english]{babel}
\usepackage{xspace}
\usepackage[centertags]{amsmath}
\usepackage{amssymb}
\usepackage[all,2cell,pdf]{xy} 
\UseAllTwocells
\SelectTips{cm}{11}
\usepackage[amsmath,thmmarks,thref,hyperref]{ntheorem}
\usepackage{calrsfs}
\usepackage{dsfont} 
\usepackage{mbboard}
\usepackage{mathtools}
\usepackage{extarrows}
\usepackage{chemarrow}
\usepackage{yhmath}
\usepackage{enumerate}
\usepackage{leftidx}
\usepackage{cmll}
\usepackage{verbatim}

\usepackage{abstract}
\usepackage{index}
\usepackage{hyperref}

\usepackage{cite}

\usepackage[textheight=26.7cm,hmargin=2.5cm]{geometry} 
\theoremheaderfont{\medskip\normalfont\bfseries}
\theoremseparator{}
\theoremstyle{break}
\newtheorem{Prop}{Proposition}[section]
\newtheorem{Lem}[Prop]{Lemma}
\newtheorem{Thm}[Prop]{Theorem}
\newtheorem{Cor}[Prop]{Corollary}
\newtheorem{Axn}[Prop]{Axiom}
\newtheorem{Axns}[Prop]{Axioms}

\theoremstyle{nonumberbreak}
\newtheorem{Axiom}{Axiom}

\theoremstyle{plain}
\theoremseparator{ :}
\newtheorem{Def}[Prop]{Definition}
\newtheorem{Defp}[Prop]{Definitions}
\newtheorem{Not}[Prop]{Notation}

\theorembodyfont{\upshape}
\newtheorem{Rem}[Prop]{Remark}
\newtheorem{Remp}[Prop]{Remarks}
\newtheorem{Ex}[Prop]{Example}
\newtheorem{Exp}[Prop]{Examples}
\newtheorem{Term}[Prop]{Terminology}

\theoremstyle{nonumberplain}
\theoremseparator{ :}

\theoremheaderfont{\slshape}
\theoremsymbol{\fbox{}}
\newtheorem{Pf}{Proof}
\qedsymbol{\PfSymbol}

\newenvironment{Exs}[1][]{%
\begin{Exp}#1\begin{enumerate}}%
{\end{enumerate}\end{Exp}}

\newenvironment{Defs}[1][]{%
\begin{Defp}#1
\begin{enumerate}}%
{\end{enumerate}\end{Defp}}

\newenvironment{Rems}[1][]{%
\begin{Remp}#1\begin{enumerate}}%
{\end{enumerate}\end{Remp}}

\numberwithin{equation}{section} 

\newcommand\N{\ensuremath{\mathds{N}}\xspace}
\newcommand\Z{\ensuremath{\mathds{Z}}\xspace}
\newcommand\Q{\ensuremath{\mathds{Q}}\xspace}
\newcommand\R{\ensuremath{\mathds{R}}\xspace}

\newcommand{\Par}{\medskip\\}
\newcommand{\indexb}[1]{\index{#1|textbf}}
\newcommand{\resp}{resp.\ }

\newcommand{\ra}{\xrightarrow}
\let\power\wp
\renewcommand{\epsilon}{\varepsilon}
\newcommand{\tm}{\textnormal{($\tilde{\text{M}}$)}\xspace}
\newcommand{\tmp}{\textnormal{($\tilde{\text{M}}$')}\xspace}
\newcommand{\tl}{\textnormal{($\tilde{\text{L}}$)}\xspace}
\newcommand{\tc}{\textnormal{($\tilde{\text{C}}$)}\xspace}
\newcommand{\cqfd}{\hfill\PfSymbol}
\newcommand{\ssi}{if and only if\xspace}
\newcommand{\subscript}[1]{\ensuremath{_\mathrm{#1}}}
\newcommand{\op}{{\normalfont(}}
\newcommand{\fp}{{\normalfont\,)}\xspace}
\renewcommand{\frak}{\mathfrak}
\newcommand{\id}{\ensuremath{\mathrm{id}}\xspace}

\newcommand{\C}{\ensuremath{\mathcal C}\xspace}
\newcommand{\D}{\ensuremath{\mathcal D}\xspace}
\newcommand{\E}{\ensuremath{\mathcal E}\xspace}
\newcommand{\B}{\ensuremath{\mathcal B}\xspace}
\newcommand{\T}{\ensuremath{\mathcal T}\xspace}
\newcommand{\V}{\ensuremath{\mathcal V}\xspace}
\newcommand{\A}{\ensuremath{\mathcal A}\xspace}
\newcommand{\I}{\ensuremath{\mathcal I}\xspace}
\newcommand{\fL}{\ensuremath{\mathcal L}\xspace}
\newcommand{\fR}{\ensuremath{\mathcal R}\xspace}
\newcommand{\dA}{\ensuremath{\mathds A}\xspace}
\newcommand{\dB}{\ensuremath{\mathds B}\xspace}
\newcommand{\one}{\ensuremath{\mathds 1}\xspace}
\newcommand{\two}{\ensuremath{\mathbb 2}\xspace}

\newcommand{\FIB}{\ensuremath{\mathit{FIB}}\xspace}
\newcommand{\FIBc}{\ensuremath{\mathit{FIB}_c}\xspace}
\newcommand{\FIBcf}{\ensuremath{\mathit{FIB}_{c,f}}\xspace}
\newcommand{\FIBcgp}{\ensuremath{\mathit{FIB}_{c,gp}}\xspace}
\newcommand{\FIBgp}{\ensuremath{\mathit{FIB}_{gp}}\xspace}
\newcommand{\CFIB}{\ensuremath{\mathit{CFIB}}\xspace}
\newcommand{\BIFIB}{\ensuremath{\mathit{BIFIB}}\xspace}
\newcommand{\BIFIBADJ}{\ensuremath{\mathit{BIFIBADJ}}\xspace}
\newcommand{\OPFIB}{\ensuremath{\mathit{OPFIB}}\xspace}
\newcommand{\IND}[1]{\ensuremath{\mathit{IND(#1)}}\xspace}
\newcommand{\INDc}[1]{\ensuremath{\mathit{IND}_c(#1)}\xspace}
\newcommand{\OPIND}[1]{\ensuremath{\mathit{OPIND(#1)}}\xspace}
\newcommand{\OPINDoc}[1]{\ensuremath{\mathit{OPIND}_{oc}(#1)}\xspace}
\newcommand{\BIIND}[1]{\ensuremath{\mathit{BI\text{-}IND(#1)}}\xspace}
\newcommand{\BIINDc}[1]{\ensuremath{\mathit{BI\text{-}IND_c(#1)}}\xspace}
\newcommand{\BIINDoc}[1]{\ensuremath{\mathit{BI\text{-}IND_{oc}(#1)}}\xspace}
\newcommand{\BIINDbc}[1]{\ensuremath{\mathit{BI\text{-}IND_{bc}(#1)}}\xspace}

\newcommand{\MONCAT}{\ensuremath{\mathit{MONCAT}}\xspace}
\newcommand{\rcMONCAT}{\ensuremath{\mathit{RCMONCAT}}\xspace}
\newcommand{\SYMMON}{\ensuremath{\mathit{SYMMON}}\xspace}
\newcommand{\rcSYMMON}{\ensuremath{\mathit{RCSYMMON}}\xspace}
\newcommand{\CARTMON}{\ensuremath{\mathit{CARTMON}}\xspace}
\newcommand{\Mod}{\ensuremath{\mathit{Mod}}\xspace}
\newcommand{\Mon}{\ensuremath{\mathit{Mon}}\xspace}
\newcommand{\Comm}{\ensuremath{\mathit{Comm}}\xspace}
\newcommand{\cRing}{\ensuremath{\leftidx{_c}{\mathit{Ring}}{}}\xspace}

\newcommand{\MONFIB}{\ensuremath{\mathit{MONFIB}}\xspace}
\newcommand{\SMONFIB}{\ensuremath{\mathit{SMONFIB}}\xspace}
\newcommand{\MONIND}[1]{\ensuremath{\mathit{MONIND}(#1)}\xspace}
\newcommand{\MONINDs}[1]{\ensuremath{\mathit{MONIND}_s(#1)}\xspace}

\newcommand{\BICAT}{\ensuremath{\mathit{BICAT}}\xspace}
\newcommand{\XCAT}{\ensuremath{\mathcal{CAT}}\xspace}
\newcommand{\CAT}{\ensuremath{\mathit{CAT}}\xspace}
\newcommand{\Cat}{\ensuremath{\mathit{Cat}}\xspace}
\newcommand{\ADJ}{\ensuremath{\mathit{ADJ}}\xspace}
\newcommand{\DBL}{\ensuremath{\mathit{DBLCAT}}\xspace}

\newcommand{\F}{\ensuremath{\mathcal F}\xspace}

\newcommand{\SITE}{\ensuremath{\mathit{SITE}}\xspace}
\newcommand{\SITEs}{\ensuremath{\mathit{SITE}_s}\xspace}
\newcommand{\Site}{\ensuremath{\mathit{Site}}\xspace}
\newcommand{\PSITE}{\ensuremath{\mathit{PSITE}}\xspace}
\newcommand{\PSite}{\ensuremath{\mathit{PSite}}\xspace}
\newcommand{\PCSITE}{\ensuremath{\mathit{PCSITE}}\xspace}
\newcommand{\PCSite}{\ensuremath{\mathit{PCSite}}\xspace}

\newcommand{\Top}{\ensuremath{\mathit{Top}}\xspace}
\newcommand{\Diff}{\ensuremath{\mathit{Diff}}\xspace}
\newcommand{\Triv}{\ensuremath{\mathit{Triv}}\xspace}

\newcommand{\Loc}{\ensuremath{\mathit{Loc}}\xspace}

\newcommand{\Sch}{\ensuremath{\mathit{Sch}}\xspace}
\newcommand{\Oc}[1]{\ensuremath{\mathcal{O}_{#1}}\xspace}
\newcommand{\Set}{\ensuremath{\mathit{Set}}\xspace}
\newcommand{\Grp}{\ensuremath{\mathit{Grp}}\xspace}
\newcommand{\Ab}{\ensuremath{\mathit{Ab}}\xspace}
\newcommand{\Ring}{\ensuremath{\mathit{Ring}}\xspace}
\newcommand{\Ringed}{\ensuremath{\mathit{Ringed}}\xspace}
\newcommand{\LRinged}{\ensuremath{\mathit{LRinged}}\xspace}

\newcommand{\M}{\ensuremath{\mathcal M}\xspace}
\newcommand{\G}{\ensuremath{\mathcal G}\xspace}
\newcommand\Ob{\ensuremath{\mathrm{Ob}\,}\xspace}
\newcommand\Mor{\ensuremath{\mathrm{Mor}\,}\xspace}
\newcommand\cod{\ensuremath{\mathrm{cod}}\xspace}
\newcommand\dom{\ensuremath{\mathrm{dom}}\xspace}
\newcommand{\Sh}{\ensuremath{\mathit{Sh}}\xspace}
\newcommand{\res}{\ensuremath{\mathrm{res}}\xspace}
\newcommand{\PSh}{\ensuremath{\mathit{PSh}}\xspace}
\newcommand{\Spec}{\ensuremath{\mathrm{Spec}}\xspace}
\newcommand{\im}[1]{\ensuremath{\mathit{Im}\,#1}\xspace}
\newcommand{\rim}[1]{\ensuremath{\leftidx{_r}{\mathit{Im}}{}\,#1}\xspace}

\newcommand{\rep}[1]{\ensuremath{\leftidx{_r}{#1}{}}\xspace}
\newcommand{\rfim}[1]{\ensuremath{\leftidx{_{\mathit{rf}}}{\mathit{Im}}{}\,#1}\xspace}

\newcommand{\frep}[1]{\ensuremath{\leftidx{_{\mathit{rf}}}{#1}{}}\xspace}
\newcommand{\fim}[1]{\ensuremath{\leftidx{_f}{\mathit{Im}}{}\,#1}\xspace}
\newcommand{\full}[1]{\ensuremath{\leftidx{_{\mathit{f}}}{#1}{}}\xspace}

\newcommand{\U}{\ensuremath{\mathcal U}\xspace}
\newcommand{\ZFC}{\ensuremath{ZFC}\xspace}
\newcommand{\ZF}{\ensuremath{ZF}\xspace}
\newcommand{\ZFI}{\ensuremath{ZFI}\xspace}
\newcommand{\NBG}{\ensuremath{NBG}\xspace}
\newcommand{\NBI}{\ensuremath{NBG_I}\xspace}

\newcommand{\adjoint}[4]{$$\xymatrix{
	#3:#1\ar@<.9ex>[r]\ar@{}[r]|{\perp} & #2:#4 \ar@<.9ex>[l]
	}$$}

\makeindex\proofmodefalse
\newindex{not}{ndx}{nnd}{Glossary of notations}

\begin{document}

\pagenumbering{Roman}
\begin{titlepage}
  {\centering\huge \vspace*{9em}Categorical Foundations for $K$-Theory\\[5em]
  \large Nicolas MICHEL\\[2em]
  PhD Thesis\\\small
  Revised and expanded version\\[2em]
  \large Ecole Polytechnique Fédérale de Lausanne\\Switzerland\\[2em]
  November 7, 2011\\}
  \vfill
  \small
  \noindent\textbf{Address}\\ 
  Nicolas Michel\\
  Station 8\\ 
  CH-1015 Lausanne\\
  Switzerland\\
  \href{mailto:nicolas.michel@epfl.ch}{nicolas.michel@epfl.ch}
\end{titlepage}

\thispagestyle{empty} 
\cleardoublepage
\pagenumbering{roman}

\thispagestyle{empty}
\selectlanguage{frenchb}
\vspace*{2cm}
\begin{quotation}
A mes parents, pour leur enseignement et leur amour.
\end{quotation}

\thispagestyle{empty} 
\clearpage

\thispagestyle{empty} 
\selectlanguage{english}
\vspace*{3cm}
\begin{quote}
``This is our fate, to live with doubts, to pursue a subject whose absoluteness we are not certain of, in short to realize that the only `true' science is itself of the same mortal, perhaps empirical, nature as all other human undertakings.''

\medskip
\hfill Paul Cohen\footnote{P.J. Cohen: \textit{Comments on the foundations of set theory}. In D.S. Scott (ed.): \textit{Axiomatic Set Theory}, vol. XIII, Part 1 of \textit{Proceedings of Symposia in Pure Mathematics}, pp. 9–15. American Mathematical Society, 1971.}
\end{quote}

\vspace{5cm}

\begin{quote}
``For now we see in a mirror dimly, but then face to face. Now I know in part; but then shall I know even as also I am known.''

\medskip
\hfill Apostle Paul\footnote{Apostle Paul: \textit{First Epistle to the Corinthians}. Chapter 13. The Bible.}
\end{quote}

\thispagestyle{empty} 
\cleardoublepage


\currentpdfbookmark{Table of contents}{toc}
\tableofcontents
\cleardoublepage


\selectlanguage{frenchb}
\phantomsection\addcontentsline{toc}{chapter}{Acknowledgments}
\renewcommand*\abstractname{Remerciements}
\begin{abstract}
\renewcommand*{\abstracttextfont}{\normalfont\normalsize}
\setcounter{page}{7}
\thispagestyle{plain}
J'aimerais chaleureusement remercier ma directrice de thèse, Professeure Kathryn Hess Bellwald. Elle a su me donner entière liberté tout en s'intéressant de près à mon travail, montrer de l'enthousiasme pour ma recherche, tout en ayant des exigences élevées. Merci Kathryn pour ta présence stimulante et amicale tout au long de ces années.

Je remercie les experts, Professeurs Marino Gran, Tom Leinster et Jacques Thévenaz, pour avoir soigneusement lu mon manuscrit, ainsi que pour les nombreux commentaires constructifs qu'ils ont faits. Je remercie Professeur Tom Mountford pour avoir dirigé avec sympathie la défense privée.

J'aimerais aussi remercier certains mathématiciens qui ont pris du temps pour répondre à mes questions durant ma thèse. Je remercie d'abord Professeur Steve Lack pour l'attention qu'il a portée à mon travail tout au long de ma thèse et pour les partages amicaux. Je remercie les Professeurs Ross Street, André Joyal et Zoran Skoda pour leur accueil chaleureux devant mes questions. Je remercie enfin Professeur Charles Weibel pour avoir assidûment répondu à mes questions par email.

Je veux aussi remercier tous mes collègues que je ne saurais tous citer ici, avec qui nous avons échangé sur les mathématiques comme sur les tics des matheux. Je pense en particulier à Théophile longtemps mon collègue de bureau. Merci pour ton respect et ta disponibilité. Je pense aussi à Ilias, avec qui j'aime tant discuter de la vie et de Celui qui en est la source. Je pense à Patrick avec qui je partage le bureau. Merci pour les discussions personnelles que nous avons. Je pense à Jérôme. Merci pour ton amitié, ta disponibilité et pour tout le travail que tu as fait en plus pour m'épargner du temps en fin de thèse. Je remercie Maria, notre secrétaire. Merci pour votre disponibilité et votre efficacité dans le travail. Merci en particulier pour votre aide précieuse dans les dernières heures avant que je rende mon manuscrit.

Je remercie les amis pour tous les partages. Je ne saurais vous nommer tous ici, mais vous êtes chers à mon c\oe ur.

Je remercie ma belle-famille pour leur accueil ouvert et chaleureux, ainsi que pour leur soutien. Je remercie mon frère Luc et ma s\oe ur Juliette pour leur fidélité depuis toujours. Je remercie mes parents à qui je dédie cette thèse.

Je remercie mon fils Théodore qui fait bondir mon c\oe ur de joie. Je remercie mon épouse Priscille, toi qui partages le jardin de mon c\oe ur. Merci pour ta présence stimulante et douce à mes côtés.

Je remercie le Christ en qui se trouve la Paix qui surpasse toute intelligence.
\end{abstract}

\cleardoublepage

\selectlanguage{english}
\phantomsection\addcontentsline{toc}{chapter}{Abstract}
\renewcommand*\abstractname{Abstract}
\begin{abstract}
\setcounter{page}{9}\thispagestyle{plain}

$K$-Theory was originally defined by Grothendieck as a contravariant functor from a subcategory of schemes to abelian groups, known today as $K_0$. The same kind of construction was then applied to other fields of mathematics, like spaces and (not necessarily commutative) rings. In all these cases, it consists of some process applied, not directly to the object one wants to study, but to some category related to it: the category of vector bundles over a space, of finitely generated projective modules over a ring, of locally free modules over a scheme, for instance. 

Later, Quillen extracted axioms that all these categories satisfy and that allow the Grothendieck construction of $K_0$. The categorical structure he discovered is called today a \emph{Quillen-exact category}. It led him not only to broaden the domain of application of $K$-theory, but also to define a whole \emph{$K$-theory spectrum} associated to such a category. Waldhausen next generalized Quillen's notion of an exact category by introducing categories with weak equivalences and cofibrations, which one nowadays calls \emph{Waldhausen categories}. $K$-theory has since been studied as a functor from the category of suitably structured  small categories (Quillen-exact, Waldhausen, symmetric monoidal, and topologically enriched variants of these) to some category of spectra%
\footnote{Works of Lurie, Toën and Vezzosi have shown that $K$-theory really depends on the $(\infty,1)$-category associated to a Waldhausen category \cite{TV04}. Moreover, topological $K$-theory of spaces and Banach algebras takes the fact that the Waldhausen category is \emph{topological} in account \cite{Pal96,Mit01}.}.
This has given rise to a huge field of research, so much so that there is a whole journal devoted to the subject.

In this thesis, we want to take advantage of these tools to begin studying $K$-theory from another perspective. Indeed, we have the impression that, in the generalization of topological and algebraic $K$-theory that has been started by Quillen, something important has been left aside. $K$-theory was initiated as a (contravariant) functor from the various categories of spaces, rings, schemes, \ldots, not from the category of Waldhausen small categories. Of course, one obtains information about a ring by studying its Quillen-exact category of (finitely generated projective) modules, but still, the final goal is the study of the ring, and, more globally, of the category of rings. Thus, in a general theory, one should describe a way to associate not only a spectrum to a structured category, but also a structured category to an object. Moreover, this process should take the morphisms of these objects into account. This gives rise to two fundamental questions. 
\begin{enumerate}
\item
What kind of mathematical objects should $K$-theory be applied to?
\item
Given such an object, what category ``over it'' should one consider and how does it vary over morphisms? 
\end{enumerate}
Considering examples, we have made the following observations. Suppose $\C$ is the category that is to be investigated by means of $K$-theory, like the category of topological spaces or of schemes, for instance.
\begin{itemize}
\item
The category associated to an object of \C is a sub-category of the category of modules over some monoid in a monoidal category with additional structure (topological, symmetric, abelian, model).
\item
The situation is highly ``fibred'': not only morphisms of \C induce (structured) functors between these sub-categories of modules, but the monoidal category in which theses modules take place might vary from one object of \C to another.
\item
In important cases, the sub-categories of modules considered are full sub-categories of ``locally trivial'' modules with respect to some (possibly weakened notion of) Grothendieck topology on \C. That is, there are some specific modules that are considered sufficiently simple to be called trivial and locally trivial modules are those that are, locally over a covering of the Grothendieck topology, isomorphic to these.
\end{itemize}
In this thesis, we explore, \emph{with $K$-theory in view}, a categorical framework that encodes these kind of data. We also study these structures for their own sake, and give examples in other fields. We do not mention in this abstract set-theoretical issues, but they are handled with care in the discussion. Moreover, an appendix is devoted to the subject. 

After recalling classical facts of Grothendieck fibrations (and their associated indexed categories), we provide new insights into the concept of a bifibration. We prove that there is a 2-equivalence between the 2-category of bifibrations over a category \B and a 2-category of pseudo double functors from \B into the double category $\frak{ADJ}$ of adjunctions in \CAT. We next turn our attention to composable pairs of fibrations $\E\ra{Q}\D\ra{P}\B$, as they happen to be fundamental objects of the theory. We give a characterization of these objects in terms of pseudo-functors $\B^{op}\to \FIB_c$ into the 2-category of fibrations and Cartesian functors.

We next turn to a short survey about Grothendieck (pre-)topologies. We start with the basic notion of \emph{covering function}, that associate to each object of a category a family of coverings of the object. We study separately the saturation of a covering function with respect to sieves and to refinements. The Grothendieck topology generated by a pretopology is shown to be the result of these two steps.

Inspired by Street \cite{Str04}, we define the notion of (locally) trivial objects in a \emph{fibred site}, that is, in a fibred category $P\colon\E\to\B$ equipped with some notion of covering of objects of the base \B. The \emph{trivial objects} are objects chosen in some fibres. An object $E$ in the fibre over $B\in\B$ is \emph{locally trivial} if there exists a covering $\{f_i\colon B_i\to B\mid i\in I\}$ such the inverse image of $E$ along $f_i$ is isomorphic to a trivial object. Among examples are torsors, principal bundles, vector bundles, schemes, locally constant sheaves, quasi-coherent and locally free sheaves of modules, finitely generated projective modules over commutative rings, topological manifolds, … We give conditions under which locally trivial objects form a subfibration of $P$ and describe the relationship between locally trivial objects with respect to subordinated covering functions.

We then go into the algebraic part of the theory. We give a definition of \emph{monoidal fibred categories} and show a 2-equivalence with \emph{monoidal indexed categories}. We develop algebra (monoids and modules) in these two settings. Modules and monoids in a monoidal fibred category $\E\to\B$ happen to form a pair of fibrations $Mod(\E)\ra{Mon} Mon(\E)\ra{P}\B$. 

We end this thesis by explaining how to apply this categorical framework to $K$-theory and by proposing some prospects of research.
\bigskip
\begin{center}
-----------------------------------------------------\medskip

KEYWORDS\bigskip

$K$-theory -- Local triviality -- Grothendieck fibration -- Grothendieck topology -- Monoidal fibred category -- Module
\end{center}

\end{abstract}

\chapter*{Foreword for this version}
\addcontentsline{toc}{chapter}{Foreword} 
This is a revised and expanded version of my PhD thesis \cite{MyThesis}. Moreover, an article is in preparation that is going to summarize some parts of this thesis and answer some of the open questions mentioned in \autoref{cha:Prospects}. 

We have added in this version a section on \emph{fibred sites}, where we study different ways of mixing the notions of fibred categories and Grothendieck sites. Our notion of fibred site happens to be a particular case of an internal fibration in some 2-category of sites.

We have quite deeply remodeled the \autoref{cha:LocTriv} about locally trivial objects. In particular, we now start with a section about locally trivial in a \emph{site} before talking about such objects in a \emph{fibred site}. Indeed, we first want to introduce the subject of locally trivial objects by defining them in the more common context of sites, rather than directly in the context of fibred sites. Moreover, locally trivial objects in a site appear as the locally trivial objects \emph{in the base category} in the fibred context. We also clarify in this chapter the previously defined notions of image of a morphism of fibrations using the language of factorization systems.

The \autoref{cha:FibMod} has not much changed, except that in \autoref{ssec:fibred-algebra}, we have simplified the definition of monoids and modules in a monoidal fibred category.

Finally, we have revised \autoref{cha:Prospects} thanks to these new developments.


\pagenumbering{arabic}
\chapter{Introduction}
This thesis has started with a very simple question. Let $\C$ be a category and suppose one would like to define the $K$-theory of objects $C$ of \C, as the $K$-theory of spaces or of schemes. What sensible (Quillen-exact, Waldhausen, …) category $\A_C$ should be associated to an object $C$ of \C so that one can apply the $K$-theory functor to it and thus obtain a meaningful $K$-theory spectrum $K(C)\coloneqq K(\A_C)$ of $C$. It then leads to the attendant question: how should the categories $\A_C$ vary over a morphism $C\to C'$ in \C? 

Well, the question so asked is too general. One should keep it in mind for a while and first face the following question: in fact, what objects $K$-theory should be applied to? The observation of the different kinds of existing $K$-theories leads us to the conclusion that $K$-theory is designed to be applied to categories of monoids in a monoidal fibred category or to a category that has a functor into such a category. 

A monoidal fibred category over a category \B can be equivalently defined to be a monoidal object $(\E\to\B,\otimes)$ in the 2-category of fibrations or as a contravariant pseudo-functor from \B to the 2-category of monoidal categories, monoidal functors and monoidal natural transformations. A monoidal fibration $P\colon\E\to\B$ gives rise to a fibration of modules over monoids $\Mod(\E)\to \Mon(\E)$. A monoid $R$ in $P$ is just a monoid $R$ in a fibre of $P$, which is a monoidal category. An $R$-module $(R,M)$ in $P$ is an $R$-module in the fibre of $R$. For example, in algebraic geometry, a sheaf of $\Oc X$-module $\F$ is a module $(X,\Oc X,\F)$ in the monoidal fibration of sheaves of abelian groups over \Top, the category of topological spaces. Morphisms of monoids and modules in a monoidal fibration can however cross from one fibre to another.

Notice that if the category \E is trivially fibred, that is, fibred over the terminal object \one of \CAT, then $\E$ is just a monoidal category and one recovers the usual fibration of modules over monoids. In the case of a category \C together with a functor $F\colon\C\to\Mon(\E)$, one obtains a pullback fibration over \C from $\Mod(\E)\to \Mon(\E)$. In the latter, the category \C will always satisfy this condition. In particular, \C might be just $\Mon(\E)$.

The answer to the first question partly answers to the question we had kept in mind. Indeed, there is a category of modules over each monoid $R$. Moreover, these categories of modules form a fibration over the category of monoids and so are related by the inverse image functors of the fibration (or, in a more ``fibrational'' fashion, by Cartesian arrows). Thus, in the situation of a category $\C$ equipped with a functor $F\colon\C\to\Mon(\E)$, one has a candidate for the category $\A_{C}$ over $C\in\C$: the category of modules of the monoid $F(C)$. This candidate is good, but still not the right one. It is too big (very often a proper class) and too complicated. One thus first need to restrict its size by considering ``finitely presented'' objects in some sense (finitely presented modules, or vector bundles whose fibres are finite dimensional, for instance). 

The second restriction is more subtle and leaves some choices that lead to different $K$-theories. There are indeed among the modules of the fibrations some that are trivial, in the sense that they are easier to understand. They are for instance the trivial bundles, free modules, ``affine'' sheaves of modules $\tilde M$. The $K$-theory of an object is sometimes defined directly on the subcategory of trivial modules over it. Yet, one often takes an intermediate step. For example, one considers finitely generated projective modules (instead of free), vector bundles (instead of trivial), locally free sheaves of modules (insteed of free). In important cases, like the one just cited, these are modules that are ``locally trivial'' with respect to some weak notion of Grothendieck pretopology on \C. Moreover, we observe that in these examples, the fact that the trivial objects are finitely presented in some sense implies that the locally trivial ones also are.

This thesis is meant to give a precise categorical foundation for these ideas. $K$-theory is the motivation of this work, but we also study the categorical concepts for their own sake and with an interest for wider applications. We now give a short summary of the thesis.

\hyperref[cha:Locality]{Chapter \ref*{cha:Locality}} provides the necessary theory of defining a central notion in our work, that of \emph{(op-, bi-)fibred site}.  After recalling classical facts of Grothendieck fibrations (and their associated indexed categories), we provide new insights into the concept of a bifibration. We prove that there is a 2-equivalence between the 2-category of bifibrations and a 2-category of pseudo double functors into the double category $\frak{ADJ}$ of adjunctions in \CAT. 

We next turn our attention to pairs of composable fibrations $\E\ra{Q}\D\ra{P}\B$ over a category $\B$, as they happen to be fundamental objects of the theory. We call them, after Hermida \cite{Her99}, \emph{fibrations over fibrations}. It is a classical fact that the composite $P\circ Q$ and the fibres $Q_A$ for $A\in\B$, are fibrations. We show that Cartesian lifts of $Q$ can be decomposed as an \emph{horizontal} lift pre-composed with a \emph{vertical} lift, i.e., a lift in $P\circ Q$ pre-composed by a lift of a fibre fibration. We finally give a characterization of these objects in terms of pseudo-functors $\B^{op}\to \FIB_c$ into the 2-category of fibrations and Cartesian functors.

We then investigate some weak notions of Grothendieck (pre-)topologies. We start with the basic notion of \emph{covering function}, that associate to each object of a category a family of coverings of the object. We then consider different axioms on these covering functions. We study separately the saturation of a covering function with respect to sieves and to refinements. The Grothendieck topology generated by a pretopology is shown to be the result of these two steps. We then consider particular sites, sites with pullbacks, and sheaves with value in a category \A on these sites. Number of examples of covering functions are considered.

We end this chapter by exploring the combination of the notions of fibred categories and sites, what leads us to define \emph{fibred sites}. We compare this notion to the ones introduced by Grothendieck and by Jardine. We finally show that our fibred sites are examples of internal fibrations in some 2-category of sites.

\hyperref[cha:LocTriv]{Chapter \ref*{cha:LocTriv}} is devoted to the study of trivial and locally trivial objects, notions that are inspired by the foundational article \cite{Str04} of Street. We start by defining trivial and locally trivial objects in a site. Given a site $(\C,K)$, trivial objects in \C are just objects of a chosen replete subcategory $\Triv\subset\C$. An object in this setting is locally trivial if it can be $K$-covered by trivial objects. We then transpose these notions to the fibred world. We now have a fibred site $(P,K)$ and a replete subfunctor $\Triv\subset P$ determining trivial objects in the total and the base categories of $P$. We explore the very common situation where the subfunctor is determined by a morphism of fibrations into $P$. In order to define trivial objects, we adapt properties defined fo functors to morphisms of fibrations and study factorization systems in the category of fibrations. We then introduce the notion of locally trivial objects in a fibred site $(P\colon\E\to\B,K)$. An object $E$ in the fibre over $B\in\B$ is \emph{locally trivial} if there exists a covering $\{f_i\colon B_i\to B\mid i\in I\}$ such the inverse image of $E$ along $f_i$ is isomorphic to a trivial object. 

Among examples of locally trivial objects in a site or in a fibre site are torsors, principal bundles, vector bundles, schemes, locally constant sheaves, quasi-coherent and locally free sheaves of modules, finitely generated projective modules over commutative rings, differentiable manifolds, … We give conditions under which locally trivial objects form a subfibration of $P$ and describe the relationship between locally trivial objects with respect to subordinated covering functions. We finally consider the particular situation of locally trivial objects in a fibration over a fibration, situation that arises in algebraic geometry.

\hyperref[cha:FibMod]{Chapter \ref*{cha:FibMod}} treats the subject of modules and monoids in a monoidal fibred category. We start with the classical case of modules and monoids in a monoidal category. We then turn to a short introduction to the internal versions of this notions: internal abelian groups, rings and modules, but also tensor product of internal abelian groups. We next study the notion of \emph{monoidal fibred categories} and show a 2-equivalence with \emph{monoidal indexed categories}. We briefly study the situation of a monoidal bifibration in a double category theoretic fashion. 

We then develop algebra (monoids and modules) in these two settings. Modules and monoids in a monoidal fibred category $\E\to\B$ happen to form a fibration over a fibration $Mod(\E)\ra{Mon} Mon(\E)\ra{P}\B$. This implies that the Cartesian lift of a morphism of monoids factors as an horizontal lift and a vertical lift, which is given by the bifibration of modules over monoids in the monoidal category of the fibre.

\hyperref[cha:Prospects]{Chapter \ref*{cha:Prospects}} explains how to apply the categorical framework developed in the previous chapters to $K$-theory. We also state what are the main open questions to fulfill this goal and some ideas of application.

\hyperref[cha:Found]{Appendix \ref*{cha:Found}} contains a discussion about the set-theoretical issues of category theory and an exposition of the point of view adopted in our work.


\chapter{Fibred sites}\label{cha:Locality}
Consider a category \B. We would like to study this category by means of ``categories of structures over its objects''. This is actually a very common strategy. For example, one gets information about a ring by studying its category of modules, its category of chain complexes or its derived category. In another field, one gets information about a topological space by studying the categories of bundles of different types over it: principal $G$-bundles, vector bundles, fibrations and so on. So to each object $B$ of \B one associates a category $\E_{B}$ ``over'' it that is supposed to encode information about $B$.  Of course, these categories are not unrelated to each other; they also have to reflect the relationships between objects of \B, namely the morphisms in \B. If there is a morphism $f\colon A\to B$ in \B, one asks that there be a functor $f^*$ between the categories over $A$ and $B$, in one direction or the other, depending on the context. Considering the contravariant case, this gives rise to a correspondence
$$\xymatrix{
			&B\ar@{|->}[r] & \E_{B} 			&\\
	A\ar[r]^f	&B\ar@{|->}[r] & \E_{B}\ar[r]^{f^{*}} 	&\E_{A}.
}$$
This correspondence should somehow respect composition and identities of \B. It turns out that it is too strong to require strict preservation, but since the codomain is a (extra large) 2-category, the extra large 2-category \CAT of (possibly large) categories, we have sufficient structure to afford a looser condition: the correspondence has to be a pseudo-functor
$$\B^{op}\longrightarrow \CAT$$
(the category \B is considered as a discrete 2-category). The whole data of the category \B is now transported to the category of categories: objects, morphisms, composition and identities. Such pseudo-functors are called \emph{indexed categories over \B}.

Understanding the whole category of structures over an object is often too difficult a task. In order to simplify the theory, one usually restricts the size of the structures. One can impose set-theoretical conditions on the category, like being small or skeletally small. More specifically one can impose finiteness conditions on the objects, e.g., by restricting to finitely generated modules or finite dimensional vector bundles, which imply some set-theoretical size restriction.

Another simplification comes from trivial (or free) objects. There is often a ``trivial'' manner to obtain structures over an object (e.g., free $R$-modules or product bundles), and we can restrict our attention to them. An intermediate step is to consider objects that are only ``locally trivial'' over some notion of \emph{Grothendieck topology} on the category \B

Suppose that in each category $\E_{B}$ there is a set of objects considered as trivial. One would like to consider objects $E\in\E_{B}$ that are not trivial but that become so when ``restricted'' to some objects $A_{i}$ along arrows $f_{i}\colon A_{i}\to B$ of a \emph{covering} in \B. This means that $f_{i}^*(E)$ is a trivial object in $\E_{A_{i}}$. One sees that the local triviality condition depends on what collections of arrows are considered as coverings in the category \B, that is on the topology on \B.
\paragraph{}
The goal of this chapter is to define and study a categorical framework for describing locally trivial structures over a category. The two main axes are the studies of \emph{Grothendieck fibrations}, which are the ``intrinsic'' indexed categories,  and of \emph{covering functions}, which are generalized Grothendieck topology. This chapter takes of course its inspiration in the ideas developed by Grothendieck, in particular his excellent article on \emph{fibred categories}, \cite[Exposé VI]{SGA1}. We will also give specific references throughout the chapter.

\section{Grothendieck fibrations}\label{sec:GroFib}
Here are diverse references that I recommend on the subject \cite{BorII94,Shu08a,Gra69,Joh202,Str05,SGA1,JT94,Web07}.

As mentioned in the introduction of this chapter, we would like to study a category \B by means of categories of ``structures over its objects''. We have explained how one is naturally led to consider pseudo-functors $F\colon\B^{op}\to \CAT$. These objects are called (\B -)\emph{indexed categories}. There is another categorical framework that encodes the same information: \emph{Grothendieck fibrations}. These are functors $\E\to\B$ satisfying certain lifting conditions. The difference between the latter and indexed categories is like the difference between the data of a category \C having all limits with respect to a diagram small category \D and the data of the category \C with an explicit limit functor $\C^{\D}\ra{\lim}\C$.

Principal $G$-bundles give a good geometric feeling for the correspondence between indexed and fibred categories. Recall that a principal $G$-bundle is a continuous map $p\colon E\to B$ with a fibre-preserving, continuous $G$-action satisfying a condition of \emph{existence} of ``local trivializations'' over some open covering $\{U_{\alpha}\}$ of $B$. Such a principal bundle determines a set of continuous functions, called ``$G$-transition functions'', $$\{g_{\alpha\beta}\colon U_{a}\cap U_{\beta}\to G\},$$ subject to some \emph{coherence} conditions. Conversely, given such set of transition functions, one can glue together the trivial pieces $U_{\alpha}\times G$ using them and obtain a principal bundle. Up to equivalence, principal $G$-bundles and sets of $G$-transitions functions contain the same information and both views are useful depending on the context\footnote{One can actually make this sentence precise. There are a category of principal $G$-bundles and a category of sets of $G$-transition functions, and these two categories are equivalent.}.
\paragraph{}
Before going at the heart of the subject, we need to fix some notations and terminology about bicategories.

\subsection{Notations for bicategories}

In order to fix notation, we recall now shortly some definitions of bicategory theory.
\begin{Defs}
\item\index{Bicategory|textbf}
	A bicategory $\mathds{A}$ consists of
	\begin{itemize}
\item
A class $\mathds{A}_{0}$ of \emph{objects}.
\item
For all pair $(A,B)$ of objects, a category $\mathds{A}(A,B)$.
\begin{itemize}
\item
The objects of $\mathds{A}(A,B)$ are called \emph{arrows} (or \emph{morphisms}, or \emph{1-cells}) of $\mathds{A}$ and denoted by $A \xrightarrow{f}B$.
\item
The morphisms of $\mathds{A}(A,B)$ are called \emph{2-cells} of $\mathds{A}$ and denoted $\xymatrix@1{A \rtwocell^{f}_{g}{\alpha} & B }$.
\item
The composition of 2-cells in $\mathds{A}(A,B)$ is denoted by $\beta\bullet\alpha$.
\item
The identity 2-cell at $f\colon A\to B$ is denoted $\xymatrix@1{A\rtwocell^{f}_{f}{\ \iota_{f}} & B}$.
\end{itemize}
\item
For all triple $(A,B,C)$ of objects, a \emph{composition} functor 
$$\circ\colon\mathds{A}(A,B)\times \mathds{A}(B,C)\to \mathds{A}(A,C).$$
\item
For all object $A$, an identity morphism $1_{A}\colon A\to A$. We denotes $\xymatrix@1{A\rtwocell^{1_{A}}_{1_{A}}{\ \iota_{A}} & A}$ its identity 2-cell.
\item
For all triple of composable arrows $(f,g,h)$, an isomorphism natural in $f$, $g$ and $h$, called \emph{associator}, 
$$ \alpha_{f,g,h}\colon (h\circ g)\circ f\xRightarrow{\cong}h\circ(g\circ f).$$
\item
For all morphism $f\colon A\to B$, two isomorphisms natural in $f$, called \emph{left} and \emph{right unitors},
$$\lambda_{f}\colon 1_{B}\circ f\xRightarrow{\cong} f\quad\text{and}\quad \rho_f\colon f\circ 1_{A}\xRightarrow{\cong}f.$$
\end{itemize}
These data are subject to coherence axioms. A bicategory is a \emph{2-category}\index{2-category|see{Bicategory}}\index{Bicategory!category@2-category|textbf} if the associators and unitors are identities.\bigskip

\begin{Not}[whiskers]\index{Whisker|textbf}\index{Bicategory!Whisker|textbf}
	Let $\mathds{A}$ be a bicategory. In the following situations in $\mathds A$,
	$$\xymatrix{
	A\rtwocell^f_g{\alpha} & B\ar[r]^h & C &\text{and}& A\ar[r]^f & B\rtwocell^g_h{\alpha} & C
	}$$
	we denote $h\cdot\alpha$ and $\alpha\cdot f$ respectively the composite 2-cells $\iota_h\circ\alpha$ and $\alpha\circ\iota_f$.
\end{Not}

Let now $\mathds{A}$ and $\mathds{B}$ be bicategories.
\item
A \emph{lax functor}\indexb{Lax functor|see{Functor}}\indexb{Functor!Lax --} from $\mathds{A}$ to $\mathds B$ is a triple $(\Phi,\gamma,\delta)$ where:
\begin{itemize}
\item
\(\Phi\) consists of the following data.
\begin{itemize}
\item
A function $$\Phi_{0}\colon \mathds{A}_{0}\to \mathds{B}_{0}.$$
\item
For all pairs $(A,A')$ of objects of $\mathds{A}$, a functor 
$$\Phi_{A,A'}\colon \mathds{A}(A,A')\to \mathds{B}(\Phi_{0}(A),\Phi_{0}(A')).$$
\end{itemize}
\item
\(\gamma\) consists of a 2-cell $$ \gamma_{f,g}\colon \Phi(g)\circ \Phi(f)\Rightarrow \Phi(g\circ f),$$ for all composable pair of arrows $f$ and $g$ of $\mathds{A}$, that is natural in $f$ and $g$.
\item
\(\delta\) consists of a 2-cell $$ \delta_{A} \colon 1_{\Phi_{0}(A)}\Rightarrow\Phi_{A,A}(1_{A}),$$ for all object $A$ of $\mathds{A}$.
\end{itemize}
These data are subject to coherence axioms. The natural transformations $ \gamma$ and \(\delta\) are called the \emph{structure morphisms} of the lax functor.\\
A lax functor is a \emph{pseudo-functor}\indexb{Functor!Pseudo-functor}\index{Pseudo-functor|see{Functor}} if its structure morphisms are isomorphisms and a \emph{2-functor}\index{functor@2-functor|see{Functor}}\indexb{Functor!Functor@2-functor} if these are identities.
\item
Let $\Phi,\Psi \colon\mathds{A}\rightrightarrows\mathds{B}$ a pair of lax functors. A \emph{lax natural transformation}\indexb{Lax natural transformation|see{Natural transformation}}\indexb{Natural transformation!Lax --} from \(\Phi\) to \(\Psi\) is a pair $(\tau,\xi)$ where
\begin{itemize}
\item
\(\tau\) consists of an arrow $ \tau \colon\Phi(A)\to\Psi(A)$ for all objects $A$ of $\mathds{A}$.
\item
\(\xi\) consists of a 2-cell $\xi_{f}$ for each arrow $f \colon A\to B$ in $ \mathds{A}$, which fills in the following square
$$\xymatrix{
\Phi(A)\ar[r]^{\tau_{A}} \ar[d]_{\Phi(f)} \drtwocell<\omit>{\ \xi_{f}}& \Psi(A) \ar[d]^{\Psi(f)}\\
\Phi(B) \ar[r]_{\tau_{B}} & \Psi(B) 
}$$
and which is natural in $f$.
\end{itemize}
These data are subject to coherence axioms.\\ 
A \emph{pseudo-natural transformation}\indexb{Pseudo-natural transformation|see{Natural transformation}}\indexb{Natural transformation!Pseudo-natural transformation} is a lax natural transformation $(\tau,\xi)$ between pseudo-functors whose structure morphism \(\xi\) is an isomorphism. A \emph{2-natural transformation}\indexb{natural transformation@2-natural transformation|see{Natural transformation}}\indexb{Natural transformation!natural transformation@2-natural transformation} is a lax natural transformation between 2-functors whose structure morphism is an identity. 

An \emph{oplax natural transformation}\indexb{Oplax natural transformation|see{Natural transformation}}\indexb{Natural transformation!Oplax --} from \(\Phi\) to \(\Psi\) is defined in the same manner as a pair $(\tau,\xi)$, with the difference that the 2-cell $\xi_f$ goes the opposite direction.
\item
Let $\xymatrix@1@C=2.5cm{\mathds{A}\rtwocell<4>^{\Phi}_{\Psi}{\omit(\sigma,\xi)\big\Downarrow\ \big\Downarrow(\tau,\zeta)} &\mathds{B}}$ be a pair of lax natural transformations.\\
A \emph{modification} $\Xi$ from \(\sigma\) to \(\tau\) consists of a 2-cell 
$$
\xymatrix@1@C=1.5cm{\Phi(A) \rtwocell^{\sigma_{A}}_{\tau_{A}}{\quad\Xi_{A}} & \Psi(A)}
$$
for all objects $A$ of $\mathds{A}$. It is subject to the following axiom.\\
For each 2-cell $\xymatrix@1@C=1.5cm{A \rtwocell^{f}_{g}{\alpha}& B}$ in $ \mathds{A}$, the following equality holds.
$$\xymatrix@C=1.8cm@R=.8cm{
\Phi(A) \ar@/_1pc/[ddr]_{ \Phi(g)} \ar[r]^{\sigma_{A}} \rlowertwocell_{\tau_{A}}{\quad\Xi_{A}}& \Psi(A) \ar[r]^{\Psi(f)} \rlowertwocell_{\Psi(g)}{\quad\Psi(\alpha)}& \Psi(B)\\
{}\rrtwocell<\omit>{\quad\zeta_{g}}&&{}\\
& \Phi(B)\ar@/_1pc/[uur]_{\tau_{B}}&
}
=
\xymatrix@C=1.8cm@R=.8cm{
\Phi(A) \ddrtwocell<5>^{\quad\Phi(f)}_{\Phi(g)\quad}{\quad\Phi(\alpha)} \ar[r]^{\sigma_{A}}  & \Psi(A) \ar[r]^{ \Psi(f)} & \Psi(B)\\
{}\rrtwocell<\omit>{\quad\xi_{f}}&&{}\\
&\Phi(B) \uurtwocell<4.8>^{\sigma_{B}\quad}_{\quad\tau_{B}}{\quad\Xi_{B}}&}
$$
\end{Defs}
Bicategories, lax functors, lax natural transformations and modifications organize in an extra large tricategory $\BICAT$ \cite{GPS95}.\index[not]{BICAT@\BICAT} 

We often work the bi-XL-category $\mathds{B}^{\mathds{A}}$\index[not]{BA@$\mathds{B}^{\mathds{A}}$} of pseudo-functors from $\mathds{A}$ to $\mathds{B}$, oplax natural transformations and modifications for bicategories $\mathds{A}$ and $\mathds{B}$. We specify by a right index ${\mathds{B}^{\mathds{A}}}_{lax}$ the corresponding bi-XL-category whose morphisms are \emph{lax} natural transformations. We specify by a left index $\leftidx{_{lax}}{\mathds{B}}{^{\mathds{A}}}$ the 2-category of \emph{lax} functors, oplax natural transformations and modifications. The 0-cell- and 2-cell-full sub-bi-XL-category of $\mathds{B}^{\mathds{A}}$ with morphisms pseudo-natural transformations is denoted ${\mathds{B}^{\mathds{A}}}_p$. Note that in general we write left indices to specify objects and right indices to specify morphisms.

When $\mathds{B}$ is a 2-category, then all these bi-XL-categories are 2-XL-categories. In fact, the full sub-tri-XL-category of $\BICAT$ consisting of 2-categories is a 3-XL-category \cite{BorI94}. If one strictifies everything, that is, if one considers 2-categories, 2-functors, 2-natural transformations, and modifications, then one gets a 3-XL-category denoted $2$-\CAT. Its 2-XL-categories of 2-functors from $\mathds{A}$ to $\mathds{B}$ are denoted $\leftidx{_{2}}{\mathds{B}}{^\mathds{A}}$.

We denote \two the category with two objects and one non-identity morphism. As we work only with 2-categories of 2-functors and 2-natural transformations from \two to a 2-category $\mathds{B}$, we make an exception in our notation and write $\mathds{B}^{\two}$\index[not]{B2@$\mathds{B}^{\two}$} for $\leftidx{_{2}}{\mathds{B}}{^\two}$. It is called the \emph{strict arrow 2-category}\indexb{Arrow 2-category}\indexb{Bicategory!Arrow 2-category} of $\mathds{B}$. More concretely, its objects are arrows in $\mathds{B}$, its morphisms commuting squares in $\mathds{B}$, and its 2-cells commutative diagrams of 2-cells of the type:
$$\xymatrix@=1.5cm{A \ar[d]_f \rtwocell^{h_1}_{k_1}{\ \alpha_1} & B \ar[d]^g\\
A' \rtwocell^{h_0}_{k_0}{\ \alpha_0} & B'.
}$$
Composition and identities are levelwise those of $\mathds{B}$.

Finally, we will deal with duality for bicategories. There are three different ways of dualizing a bicategory $\mathds A$. The bicategory $\mathds A^{op}$\index[not]{op@$\mathds A^{op}$} has reversed arrows, i.e., $$\mathds A^{op}(A,B)=\mathds A(B,A).$$ The bicategory $\mathds A^{co}$\index[not]{Co@$\mathds A^{co}$} has reversed 2-cells, i.e., $$\mathds A^{co}(A,B)=\mathds A(A,B)^{op}.$$ The bicategory $\mathds A^{coop}$\index[not]{Coop@$\mathds A^{coop}$} has both arrows and 2-cells reversed.

\subsection{Basic notions}
\subsubsection{Cartesian arrows}

\begin{Defs}
[Let $P\colon\E\to\B$ be a functor.]

\item
 An object $E$ of \E is said to be (or sit) \emph{over} an object $B$ of \B if $P(E)=B$. \textit{Idem} for arrows.
\item
The \emph{fibre of $P$ at}\indexb{Fibration!Fibre of a --} $B$ is the subcategory, denoted $P^{-1}(B)$ or $\E_{B}$, defined as the preimage of the discrete subcategory $\{B\}\subset\B$. The morphisms of the fibres are called \emph{vertical}.
\item
	An arrow $h\colon D\to E$ in \E is \emph{Cartesian over}\indexb{Cartesian!-- arrow} $f$ if it sits over $f$ and for all $K\ra k E$ in \E and all $P(K)\ra g A$ in \B such that $f\circ g=P(k)$ there exists a unique $K\ra l D$ in $E$ such that $h\circ l = k$.
		$$\xymatrix@=1.5cm{
		K\ar@{-->}[r]_{\exists!\,l}\ar@/^2pc/[rr]^{\forall k} &D\ar[r]_h & E & \E\ar[d]^P\\
		P(K)\ar[r]^{\forall g}\ar@/_2pc/[rr]_{P(k)} & A\ar[r]^{f} & B & \B
		}$$
In this situation, one says that $h$ is a \emph{Cartesian lift}\indexb{Cartesian!-- lift} of $f$ (\emph{at} $E$) and that $D$ is an \emph{inverse image of $E$ over $f$}.\indexb{Inverse image!-- of an object}
\end{Defs}

Here is an overview of the behaviour of Cartesian arrows.
\begin{Prop}\label{Prop:Cartesian}
	Let $P\colon\E\to\B$ be a functor.
	\begin{enumerate}[(i)]
		\item Cartesian lifts of a given arrow at a given object are unique up to a unique vertical isomorphism.\\
Precisely, if $D\ra h E$ and $D'\ra{h'} E$ are two Cartesian arrows over the same arrow $A\ra f B$ in \B, then there exists a unique vertical isomorphism $D'\ra{\cong} D$ making the following triangle commute:
		$$\xymatrix@=1.2cm{
			D'\ar@{-->}[d]_{\exists!\,\cong}\ar[dr]^{h'}\\
			D\ar[r]_{h} & E\\
			A\ar[r]^f & B
		}$$
		\item\label{Prop:Cartesianii} Let $D \ra{h}E \ra{k} K$ be arrows in \E with $k$ Cartesian. Then the composite $k\circ h$ is Cartesian \ssi $h$ is Cartesian.

		\item Any isomorphism of \E is Cartesian.\\
		In particular, $1_{E}\colon E\to E$ is Cartesian for all object $E\in\E$.

		\item An arrow in \E is an isomorphism \ssi it is a Cartesian lift of an isomorphism.\label{Prop:Cartesian:liftiso}
	\end{enumerate}
\qed
\end{Prop}

\subsubsection{Grothendieck fibrations and Grothendieck construction}\label{sssec:GroCons}
\begin{Defs}
	\item
A functor $P\colon\E\to\B$ is a \emph{(Grothendieck) fibration}\indexb{Fibration} (\emph{over \B} if, for every arrow $f$ in \B and every object $E$ over its codomain, there exists a Cartesian arrow over $f$ with codomain $E$. In this situation, \E is said to be \emph{fibred (in categories) over \B}.\index{Fibred category|see{Fibration}}
	\item
A choice of a Cartesian arrow $\bar f_E\colon f^{*}(E)\to E$ for all $f$ in \B and all $E\in\E$ over $\cod f$ is called a \emph{cleavage}\indexb{Cleavage} of the fibration $P$. One sometimes denotes a cleavage by his choice function $\kappa$ so that $\kappa(f,E)=\bar f_E$. A cleavage $\kappa$ is \emph{normal} if its Cartesian lifts of identities are identities. A couple $(P,\kappa)$ of a fibration and a cleavage (resp. a normal cleavage)\indexb{Cleavage!Normal --} is called a \emph{cloven fibration}\indexb{Fibration!Cloven --} (resp. \emph{normalized cloven fibration})\indexb{Fibration!Normalized cloven --}\footnote{Unless there is a canonical choice of Cartesian arrows, defining a cleavage on a fibration requires \emph{global axiom of choice}. See the first paragraph of \autoref{fn:GAC} in \autoref*{cha:FibMod}, \autopageref*{fn:GAC}. Notice that there is always a canonical Cartesian lift at an object of an identity morphism: the identity of this object (see \thref{Prop:Cartesian} (iii)). Therefore, it is always possible to \emph{normalize} a cleavage.}.

\end{Defs}

\begin{Exs}\label{ex:fibration}
\item
Any category is a fibration over the terminal object $\one$ of \CAT.
	\item\label{ex:fibration:isos}
An identity functor $Id\colon \B\to\B$ is a trivial example of a fibration. More generally, any isomorphism of categories is a fibration.
\item
The projection of the product $\prod \A_i\to\A_i$ of categories onto one of them is a fibration.
\item (Canonical fibration over \C)\indexb{Fibration!Canonical --}
A more interesting example is given for any category $\C$ with pullbacks and is called the \emph{canonical fibration} over \C. It is the codomain functor $cod\colon\C^{\two}\to\C$ from the category $\C^{\two}$ of arrows in \C to \C. The fibre at $C$ of this fibration is isomorphic to the slice category $\C/C$. We sometimes call \emph{bundles over $C$} the objects of $\C/C$ and denote them by greek letters $\xi=(E\ra p C)$. We also write $Bun(\C)$ for $\C^{\two}$.

Given an arrow $f\colon A\to B$ in \C and an object $\xi=(E\ra p B)$ in $\C^{\two}$ over $B$, a morphism $(g,f)\colon\zeta\to\xi$ with $\zeta=(D\ra{q} A)$ is Cartesian over $f$ \ssi the following square is a pullback in \C:
$$\xymatrix{
D\ar[r]^{g}\ar[d]_{q} & E\ar[d]^p\\
A\ar[r]_f & B.
}$$
We usually write $f^*\xi=(f^*E\ra{\bar p}A)$ such a Cartesian arrow.

\item\label{ex:fibration:gbundles}\indexb{Bundles!G-bundles@$G$-bundles}
(Fibration of $G$-bundles) If $\C$ has in addition a terminal object, \C admits a structure of a Cartesian monoidal category. We can then enrich the previous example by considering \emph{$G$-bundles} for an internal group $G$ of \C (see \autoref{cha:FibMod}, if needed, for definitions of the concepts involved such as internal groups, modules and Cartesian monoidal categories). A $G$-bundle over $A\in\C$ is an object $E\ra p A$ of $\C/A$ together with a \emph{$G$-object} structure $E\times G\ra{\kappa}E$ on $E$, i.e., a $G$-module structure for the monoid $G$ in the Cartesian monoidal structure on \C, such that the action $\kappa$ \emph{preserves the fibre}, i.e., the following diagram
$$\xymatrix{
E\times G\ar[r]^{\kappa}\ar[d]_{pr} & E\ar[d]^p\\
E\ar[r]_p & B.
}$$
commutes. Morphisms of $G$-bundles are morphisms $$\phi=(g,f)\colon(D\to A)\to(E\to B)$$ in $\C^{\two}$ such that $g$ is \emph{equivariant}, i.e., a morphism of $G$-modules in \C.

These data form a category, which we denote $G\text{-}Bun(\C)$ and the codomain functor yields a fibration $G\text{-}Bun(\C)\to\C$. A Cartesian lift of an arrow $f\colon A\to B$ in $\C$ at a $G$-bundle $\xi$ over $B$ is given by a Cartesian lift $\bar f_{\xi}\colon f^*(\xi)\to\xi$ in $Bun(\C)$ together with the unique $G$-bundle structure on $f^*(\xi)$ such that $\bar f_{\xi}$ is a morphism of $G$-bundles. Thus, a morphism $(g,f)\colon\zeta\to\xi$ in $G$-Bun(\C) with $\zeta=(D\ra{q} A)$ is Cartesian over $f$ \ssi the following square
$$\xymatrix{
D\ar[r]^{g}\ar[d]_{q} & E\ar[d]^p\\
A\ar[r]_f & B.
}$$
is a pullback in \C and $D$ comes with the unique fibre preserving $G$-action such that $g$ is equivariant.
\item\label{ex:fibration:bunvect}
Let us now consider the category $Top$ of topological spaces. The two previous examples apply. We consider another type of bundles. Since the functor $$\theta_B\colon Top\to Top/B$$ with $\theta_B^F=(B\times F\ra{pr}B)$ is finite product preserving, it preserves internal rings. Thus, one has a natural internal ring structure on $\theta_B^\R$ in $Top/B$. Internal $\theta_B^\R$-modules in $Top/B$ are called \emph{bundles of vector spaces}. These bundles also form a fibration over \Top, which we denote by $\mathit{VBun}\to\Top$.
\item\label{ex:fibration:cat}
The object functor $\Ob\colon \Cat\to \Set$ from the category of small categories to the category of sets is a fibration. Note that its fibre above the terminal $*$ of $Set$ is isomorphic to the category of monoids (in \Set). It has a canonical cleavage with respect to which the Cartesian lift of an inclusion at a category \C is the inclusion of the full sub-category. For a general function $f\colon X\to Ob\C$, the Cartesian lift $\bar f_\C\colon f^*\C\to\C$ is defined this way. The category $f^*\C$ has objects $X$ and $f^*(\C)(x,y)=\{x\}\times\C(f(x),f(y))\times\{y\}$. Composition and identities are those of \C. The functor $\bar f_\C$ is the function $f$ on objects and projection on the second component on morphisms.
\item
Modules in a monoidal category form a fibration $Mod(\V)\to Mon(\V)$ over the category of monoids in \V (see part \ref{ssec:modmon}).
\end{Exs}

\paragraph{Indexed categories}\label{par:indexedcat}
The structure we have considered above is the ``intrinsic'', that is, choice-free definition of the notion of a fibration. As we said in the introduction, it is possible, by making choices, to capture the same information in a pseudo-functor into \CAT, and that's what we look at now. It is not difficult to define this pseudo-functor, but in the following we try to motivate in an (we hope) intuitive way not only its form but the fact that essentially no information is lost (we will later on give a precise 2-categorical meaning of the latter property). The following discussion happens to also lead quite directly to the opposite construction, from a pseudo-functor to a fibration.

Consider a cloven fibration $P\colon \E\to\B$. We encode this data in a new way, by means of its fibres and ``inverse image'' functors between them.\\
One starts with the class of the fibres of $P$, indexed by the objects of \B. What is missing? Let us study first how one can view all objects and morphisms of \E as objects and morphisms of fibres plus some additional data. It is a trivial fact that all objects of \E live in fibres. Moreover, every arrow in \E factors uniquely as a vertical arrow with the same domain composed with a Cartesian arrow of the cleavage. More concretely, consider the Cartesian lift $\bar f_E\colon f^*E \to E$ of a morphism $f\colon A\to B$ in \B at an object $E$ of $\E_B$. It determines a bijection:
\begin{equation}\label{eq:bijvert}
\{h\text{ in }\E\mid P(h)=f\text{ and } cod(h)= E\}\xlongleftrightarrow{1:1}\{h \text{ in }\E_A\mid cod(h)=f^*E\}
\end{equation}
The bijection is defined this way. Let $h\colon D\to E$ be a morphism in \E over $f$, then the factorization is given by
\begin{equation}\label{eq:factCart}\begin{aligned}
\xymatrix@C=1.5cm{
D \ar[dr]^{h} \ar@{-->}[d]_{\exists!\,\bar h} & \\
f^*E \ar[r]_-{\overline{f}_E} & E.}
\end{aligned}\end{equation}
Thus, up to bijection, all morphisms in \E are given by vertical morphisms with codomain $f^*E$ for some arrow $f$ in \B and some object $E$ of $\E_{\cod(f)}$. This leads us to the conclusion that in addition to the fibres, one has to remember that in each fibre, there are special objects: the $f^*E$. We have not yet looked at the whole structure of \E of course. What about identities and composition? The question of identities also resolves the ambiguity about vertical morphisms: they already are in the fibres, but seem to be modified by the preceding bijection (replacing $f$ by the identity that they are over). 

Identities of \E are in fibres, but the identity $1_E$ at an object $E$ over an object $B$ becomes under the bijection \hyperref[eq:bijvert]{(\ref*{eq:bijvert})} (replacing $f$ by $1_B$) the arrow $E\ra{(\delta_B)_E} 1_B^*E$ that makes the following diagram commute.
\begin{equation}\label{eq:identity}\begin{aligned}
\xymatrix@C=1.5cm{
E \ar[dr]^{1_E} \ar@{-->}[d]_{\exists!\,(\delta_B)_E} & \\
(1_B)^*E \ar[r]_-{\overline{(1_B)}_E} & E.}
\end{aligned}
\end{equation}
This is an isomorphism by Proposition \ref{Prop:Cartesian}. If the cleavage is normalized, then it is an identity and the bijection acts as an identity on fibres.\\
In general, one should keep track of these isomorphisms, because they record that among all morphisms in $\E_B(E,id_B^*E)$, there is one that is the ``same'' as $1_E$. And from this information, one automatically obtains the compatibility between vertical morphisms $h\colon D\to E$ and their counterparts under the bijection:
$$
\xymatrix@=1.5cm{
D \ar[dr]^{h} \ar@{-->}[d]_{\bar h} & \\
(1_B)^*E \ar[r]_-{\overline{(1_B)}_E} \ar@{-->}[d]_{{(\delta_B)_E}^{-1}} & E.\\
E\ar[ru]_{1_E}
}$$

What about composition? Consider a composable pair of arrows $D \ra{h} E \ra{k} K$ in \E over $A \ra{f}B \ra{g}C$. Under our new description of morphisms in \E, they become vertical arrows $D \ra{\bar h}f^*E$ and $E \ra{\bar k}g^*K$, which are not composable in general. In fact, our data so far is incomplete. The following diagram shows how one can recover composition of general arrows in \E inside fibres.
\begin{equation}\label{eq:vertcomp}\begin{aligned}\shorthandoff{;:!?}\xymatrix@!=1cm{
D \ar[r]^h \ar@{-->}[d]_{\exists!\,\bar h} & E\ar[r]^k \ar@{-->}[d]_{\exists!\,\bar k} & K\\
f^*E \ar[ur]_(.4){\bar f_E} \ar@{-->}[d]_{\exists!\,f^*(\bar g)}& g^*K \ar[ur]_{\bar g_K} &\\
f^*g^* K \ar@{-->}[d]_{\exists!\,(\gamma_{f,g})_K} \ar[ur]_{\bar f_{g^*K}} & &\\
(g\circ f)^* K\ar@/_2pc/[rruuu]_{\overline{(g\circ f)}_K} & &
}\end{aligned}\end{equation}
The outer triangle commutes and therefore, its vertical edge is the vertical arrow corresponding to the composite $k\circ h$ under the bijection \hyperref[eq:bijvert]{(\ref*{eq:bijvert})} (replacing $f$ by $g\circ f$ and and $E$ by $K$). Notice that the morphism $(\gamma_{f,g})_K$ is an isomorphism by Proposition \ref{Prop:Cartesian}. 

In conclusion, in order to keep track of the identities and the composition of \E, one needs to remember, in addition to the fibres, the vertical isomorphisms $(\delta_B)_E$ where $E\in\E$ is over $B\in\B$ and $(\gamma_{f,g})_K$ for all composable pairs $A \ra{f}B \ra{g}C$ and $K\in\E$ over $C$. Moreover, it is not sufficient to know the objects $f^*E$ for all $f$ in \B and all $E$ in $\E_{\cod(f)}$; one also needs to retain the morphisms $f^*(h)$ for all morphisms $h$ in $\E_{\cod(f)}$. The following commutative diagram defines them.
$$\xymatrix{
f^*D\ar@{-->}[d]_{\exists!f^*(h)} \ar[r]^{\bar f_D} & D\ar[d]^{h}\\
f^* E \ar[r]_{\bar f_E} & E.
}$$
\begin{Lem}\label{lem:pseudo}
These data together determine a pseudo-functor with structure isomorphisms given by \(\delta\) and \(\gamma\) defined above.
$$\begin{array}{rcl}
\Phi_P\colon\B^{op}	& \longrightarrow	& \CAT\\
A 		& \longmapsto		& \E_A\\
A \ra{f} B 	& \longmapsto		& f^*\colon \E_B\to\E_A
\end{array}
$$
\cqfd
\end{Lem}
\begin{Def}
	Let \B be a category. A \emph{\B-indexed category}\indexb{Indexed category} (or a \emph{category indexed over \B}) is a pseudo-functor $\Phi\colon\B^{op}\to \CAT$. The functors 
	$$\Phi(f)\colon\Phi(B)\to\Phi(A)$$
	for $f\colon A\to B$
	in \B are called \emph{inverse image functors}\indexb{Inverse image!-- functor}. Moreover, if it does not induce confusion, they are simply denoted $f^*$.
\end{Def}

\begin{Exs}\label{ex:indcat}
\item
Given a topological space $X$, its topology $\T(X)$ forms a preorder under inclusion and therefore a small category. Moreover, a continuous map $$f\colon X\to Y$$ induces a monotone function $f^{-1}\colon \T(Y)\to\T(X)$, and therefore a functor. One thus obtains a functor (thus an indexed category over \Top)
$$
\T\colon\Top^{op}\longrightarrow\Cat.
$$
\item\label{ex:incatrep}
Let $\mathds A$ be a small 2-category. Then, for an object $A\in \mathds A$, there is an obvious representable 2-functor with values in the 2-category \Cat of small categories
$$
\mathds A(-,A)\colon \mathds A^{op}\to\Cat
$$
For the same reason, given an 2-XL-category $\mathds A$ and an object $A\in \mathds A$, there is a representable 2-functor with values in the XXL-2-category $\XCAT$ of XL-categories
$$
\mathds A(-,A)\colon \mathds A^{op}\to\XCAT
$$
(remember that XL-2-categories are the small 2-categories in the ambient \NBG). In particular, this applies to the 2-XL-category $\CAT$ and a category \A. When you restrict the domain of the representable 2-functor $\CAT(-,\A)$ to the 2-category of small categories, it takes its images in \CAT. One therefore obtains a 2-functor
$$
\A^{(-)}\colon\Cat^{op}\to\CAT,
$$
which restricts to an indexed category over the category \Cat:
$$\xymatrix@R=.5cm{
	\Cat^{op}\ar[r]^{\A^{(-)}}	& \CAT\\
	\C\ar@{|->}[r] 			&  \A^{\C}\\
	\B\ra{F}\C \ar@{|->}[r] 	& \A^{\C} \ar[rr]^{F^*}		& 			& \A^{\B}\\
					   	& \C\ar[r]^{G}				& \A\ar@{|->}[r]	& \B\ar[r]^{F} &\C\ar[r]^{G} & \A\\
					   	& \C\rtwocell^G_H{\alpha}	& \A\ar@{|->}[r]	& \B\ar[r]^{F} & \C\rtwocell^G_H{\alpha} & \A
	}$$
\end{Exs}
\paragraph{Grothendieck construction}
We consider now the other direction: given an indexed category $\Phi\colon\B^{op}\to \CAT$, one can construct a cloven fibration $\E_{\Phi} \ra{P_{\Phi}}\B$ with no loss of information. Actually, if one considers the pseudo-functor $\Phi_P$ associated to a cloven fibration $P\colon\E\to\B$, then the new description of $P$ that we gave in the preceding paragraph \hyperref[par:indexedcat]{``Indexed categories''} is precisely the fibration constructed from $\Phi_P$. If one just has a pseudo-functor, then one mimics the latter description pretending one has Cartesian arrows between fibres, even though one just has the inverse image functors.
\begin{Lem}\label{lem:grothcons}
	Let $(\Phi,\gamma,\delta)\colon\B^{op}\to \CAT$ be an indexed category. The following data, called the \emph{Grothendieck construction of $\Phi$}\indexb{Grothendieck!-- construction}, is a cloven fibration, denoted $P_{\Phi}\colon\E_{\Phi}\to\B$.
\normalfont
\begin{itemize}
\item
Objects of $\E_{\Phi}$: Pairs $(B,E)$ where $B\in\B$ and $E\in\Phi(B)$.
\item
Morphisms of $\E_{\Phi}$: Pairs $(A,D) \ra{(f,h)} (B,E)$ where $f\colon A\to B$ is in \B and \linebreak${h\colon D\to f^*E}$ is in $\Phi(A)$ (compare diagram \ref{eq:factCart}).
\item
Composition in $\E_{\Phi}$: Let $(A,D) \ra{(f,h)} (B,E) \ra{(g,k)}(C,K)$ be a composable pair in $\E_{\Phi}$. Then its composite is the pair whose first component is $g\circ f$ and second component is the composite in $\E_A$:
$$\xymatrix@1@C=1.5cm{D \ar[r]^{h} & f^*E \ar[r]^-{f^*(k)} & f^*g^*(K) \ar[r]^{(\gamma_{f,g})_K}_{\cong} & (g\circ f)^*(K).}
$$
(Compare diagram \ref{eq:vertcomp})
\item
Identities in $\E_{\Phi}$: $(A,D) \ra{(1_A,(\delta_A)_D)}(A,D)$ (compare diagram \ref{eq:identity})
\item
Functor $P_{\Phi}$: The projection on the first component.
\item
Cleavage of $P_{\Phi}$: The Cartesian lift of $f\colon A\to B$ at $(B,E)$ is the morphism
$$
	(A,f^*E)\ra{(f,1_{f^*E})}(B,E).
$$
\end{itemize}
\cqfd
\end{Lem}
When constructing a \B-indexed category from a cloven fibration over \B, we took care not to lose information in the process. As noticed above, we actually directly construct an indexed category and its Grothendieck construction, and argued that these three objects contain the same information but organized in a different manner.  A more precise and more categorical way of stating this fact consists of defining a \emph{category} of cloven fibrations over \B and in showing that a cloven fibration is isomorphic to the Grothendieck construction of its corresponding indexed category. The hom-set bijections of this isomorphism are in fact determined by the bijection \hyperref[eq:bijvert]{(\ref*{eq:bijvert})}. It is then natural to ask whether this statement extends to the respective categories. It does in fact extend to an equivalence of categories, and the latter isomorphism provides the essential surjectivity of the functor. Moreover, fibrations over \B and indexed \B-categories naturally live in 2-categories and these are 2-equivalent.

\subsection{The Grothendieck construction as a 2-equivalence}
In this part we define the 2-XL-category of (cloven) fibrations and recall the well known 2-equivalence between fibrations and indexed categories. We end by a study of pullback of fibrations along functors and of the associated concept in the world of indexed categories.
\subsubsection{The 2-equivalence}
\begin{Defs}
\item
The \emph{2-XL-category of fibrations}\indexb{Fibration!XL-category of@2-XL-category of --s}, denoted \FIB\index[not]{FIB@\FIB}, is the full sub-2-XL-category of $\CAT^{\two}$ whose objects are fibrations. Its morphisms are called \emph{morphisms of fibrations}\indexb{Fibration!Morphisms of --s} and 2-cells, \emph{transformations}\indexb{Fibration!Transformations}. We also consider the \emph{2-category \CFIB\index[not]{CFIB@\CFIB} of cloven fibrations} whose morphisms and arrows are the morphisms and arrows of the underlying fibrations. The following definitions can be applied to both \FIB and \CFIB.

\item
A \emph{Cartesian morphism of fibrations}\indexb{Fibration!Cartesian morphisms of --s}\indexb{Cartesian!-- morphism of fibrations|see{Fibration}} is a morphism of fibrations that preserves Cartesian arrows. Explicitly, a morphism from the fibration $Q$ to the fibration $P$
$$\xymatrix{
\D \ar[r]^{F_1} \ar[d]_Q & \E \ar[d]^{P}\\
\A \ar[r]_{F_0} & \B,
}$$
is Cartesian if, given a Cartesian arrow $D \ra{h} D'$ in $\D$, its image $F_1(D) \ra{F_1(h)}F_1(D')$ is Cartesian in $\E$. The 2-cell full sub-2-XL-category of $\CAT^{\two}$ with objects fibrations and morphisms Cartesian morphisms of fibrations is denoted $\FIBc$\index[not]{FIBc@\FIBc}. 

\item
Finally, both \FIB and $\FIBc$ can be restricted to a fixed base category \B. These categories, denoted $\FIB(\B)$\index[not]{FIB(B)@\FIB(\B)} and $\FIBc(\B)$\index[not]{FIBc(B)@\FIBc(\B)} are the sub-2-XL-categories of \FIB and $\FIBc$ respectively that have fibrations over \B as objects, morphisms of fibrations (resp.\ Cartesian morphisms) $(F_0,F_1)$ with base functor $F_0=Id_{\B}$ as morphisms and transformations with base natural transformation $\alpha_0=\iota_B$ as 2-cells. Their morphisms are called (\emph{Cartesian}) \emph{functors over \B}\indexb{Functor!-- over a category} and 2-cells \emph{natural transformations over \B}\indexb{Natural transformation!-- over a category}.
\item
The \emph{2-XL-category of \B-indexed categories} is the 2-XL-category $\CAT^{\B^{op}}$ of pseudo-functors, oplax natural transformations and modifications. We sometimes denote it $\IND\B$. Its 2-cell full sub-XL-2-category $(\CAT^{\B^{op}})_p$ whose morphisms are pseudo-natural transformations is denoted $\INDc\B$.
\end{Defs}

\begin{Rem}\label{rem:CartMor}
	A morphism of fibrations from $\D \to \A$ to $\E\to\B$ is Cartesian if and only if it preserves a Cartesian lift for each pair $(f,D)$ of an arrow $f$ in \A and an object $D$ in $\D_{\cod(f)}$ (use Proposition \ref{Prop:Cartesian}).
\end{Rem}

There is a forgetful functor $\CFIB\to FIB$ that forgets cleavages and is the identity on 1- and 2-cells. If one assumes the global axiom of choice, then it determines a 2-equivalence. In the preceding part, we have announced a 2-equivalence between cloven fibrations and indexed categories, via the Grothendieck construction. Indeed, we have constructed an indexed category from a cloven fibration and conversely a cloven fibration from an indexed category. Nevertheless, we assume such an axiom of choice and, therefore, consider a cleavage as a mere choice, not an additional structure. Consequently, we state the result directly for non-cloven fibrations. This is the way the 2-equivalence is expressed in \cite{BorII94,Shu08a} for example, but Johnstone claims it only for cloven fibrations \cite{Joh102}. The theorem is proven in \cite{BorII94,Joh102}, but for the 2-category $\FIBc(\B)$ of fibrations and Cartesian functors over \B. The proof goes exactly the same in the more general setting of the 2-category $\FIB(\B)$ of fibrations and functors over \B, but a subtlety arises that was hidden in the former context. A functor over \B between two fibrations does indeed not yield a lax natural transformation, but an oplax one.

\begin{Thm}\label{thm:gcons}
	The Grothendieck construction determines a 2-equivalence\index{Grothendieck!-- construction}\index[not]{G@$\mathcal G$}
	$$
	\mathcal G\colon \CAT^{\B^{op}}\xlongrightarrow{\simeq} \FIB(\B)
	$$
	which restricts to a 2-equivalence
	$$
	(\CAT^{\B^{op}})_p\xlongrightarrow{\simeq} \FIBc(\B).
	$$
	\cqfd
\end{Thm}

This theorem shows that the apparently different concepts of fibrations and indexed categories are two views on the same objects. In this work we will go back and forth in this 2-equivalence.
When going from fibred categories to indexed categories, universal properties are replaced by ``functional'' correspondences, which one is much more accustomed to. Nevertheless, the indeterminacy of the universal arrows up to isomorphism is now reflected by the ``pseudo'' nature of the functor. One has arbitrarily fixed something indeterminate in nature and so what one gains in intuition and habit, one loses it in simplicity. In this connection, the reader might want to have a look at the remark at the end of the chapter ``Catégories fibrées et descente'' of \cite{SGA1}. 

One should be careful though with statements of the kind ``theses two concepts are the same because they live in (2-)equivalent (2-)categories''. Indeed, there is an elementary theory of fibrations that can be interpreted in any 2-category with suitable 2-categorical limits \cite{Str74a} and even in a bicategory \cite{Str80,Str87}. The 2-equivalence stands for the internalization of the notion of a fibration in the particular (extra large) 2-category \CAT. There is no internal notion of an indexed category in a general 2-category. They are notions particular of \CAT, that has the special property of being a model of the same theory as its objects%
\footnote{It is the theory of 2-categories, categories being seen as discrete 2-categories. Objects of \CAT are model of it in the Grothendieck universe, whereas \CAT is a model of it in the larger set theory that contains the universe.}.%
In a general 2-category $\mathds A$, it does not make sense to talk about $\mathds A^B$ for an object $B\in\mathds A$. The article \cite{Ben85} contains a deep (and virulent) discussion on that matter.

\subsubsection{Pullbacks of fibrations}
The following lemma is classical and easy to check \cite{Joh102,BorII94}.
\begin{Lem}\label{lem:pullfib}
	Let $P\colon\E\to \B$ be a fibration and $F\colon \A\to\B$ a functor. Then, pullbacks of $P$ along $F$ are also fibrations, i.e., in a pullback-square in \CAT
	$$\xymatrix{
	F^*\E\ar[d]_{\bar P}\ar[r]^{\bar F} & \E\ar[d]^P\\
	\A\ar[r]_{F} & \B,
	}$$
	the functor $\bar P\colon F^*\E\to\A$ is a fibration. Moreover, the morphism $(\bar F,F)$ is Cartesian.\\
	If one considers the canonical choice of pullback, then a Cartesian arrow in $F^*\E$ over $A'\ra{f}A$ in \A at $(A,E)$ is given by $(f,\overline{F(f)}_E)\colon(A',F(f)^*E)\to(A,E)$.
\end{Lem}
Note that the lemma also provides us with a canonical cleavage on the pullback of a cloven fibration. The following proposition is not difficult to check.
\begin{Prop}
	Once a choice of pullbacks in $CAT$ is made, pullbacks of fibrations over $\B$ along a functor $F\colon\A\to\B$ yield a 2-functor
	$$
	F^*\colon\FIB(\B)\longrightarrow\FIB(\A),
	$$
	which restricts to a 2-functor
	$$
	F^*\colon\FIB_c(\B)\longrightarrow\FIB_c(\A).
	$$
	\cqfd
\end{Prop}
It is natural to ask what is the effect of a pullback on the corresponding indexed category. This is in fact just precomposition with $F^{op}$. Given a functor $F\colon\A\to\B$, there is indeed a precomposition 2-functor $F^*\colon\CAT^{\B^{op}}\to\CAT^{\A^{op}}$ (because 2-XL-categories form a 3-XXL-category 2-$\XCAT$). The following proposition is proved in \cite{Joh202} at the level of objects. The proof of the 2-naturality is straightforward, but quite a juggle.
\begin{Prop}\label{prop:grothpull}
	Given a functor $F\colon\A\to\B$, there is a 2-natural isomorphism
	$$\xymatrix@=1.3cm{
	\IND\B\ar[r]^{\G}\ar[d]_{F^*} \ar@{}[dr]|{\cong}& \FIB(\B)\ar[d]^{F^*}\\
	\IND\A\ar[r]_\G & \FIB(\A).}$$
	\cqfd
\end{Prop}
\begin{Ex}\label{ex:presheafA}
	We use the notation of Example \hyperref[ex:incatrep]{\ref*{ex:indcat} (2)}. We have already described the representable indexed category $\A^{(-)}\colon\Cat^{op}\to\CAT$ determined by a category \A. One can precompose it with the (dual of the) \emph{opposite} automorphism of the category \Cat
$$
(-)^{op}\colon\Cat\xlongrightarrow{\cong}\Cat
$$
that takes a category \C to its dual $\C^{op}$. The composite functor $$\PSh(-;\A)\coloneqq\A^{(-)^{op}}\colon\Cat^{op}\to\CAT$$ then takes a small category \C to the category $\PSh(\C;\A)\coloneqq\A^{\C^{op}}$ of \emph{presheaves on \C with values in \A.} We therefore call it the \emph{indexed category of presheaves with values in \A.}\indexb{Sheaf!Indexed category of presheaves with values in \A}\index{Indexed category!-- of presheaves with values in \A|see{Sheaf}} By abuse of notation, we still write $F^*$ the functor $\PSh(F;A)$, even though, to be coherent with Example \hyperref[ex:incatrep]{\ref*{ex:indcat} (2)}, one should write $(F^{op})^*$.
\end{Ex}

We now use pullbacks of fibrations in order to define a (2-)product on them. Given two fibrations over \B, $\D\ra{Q}\B\xleftarrow{P}\E$, their pullback inherits a fibration structure over \B. One can check it directly, or by using the fact that a composition of fibrations is a fibration (see \thref{lem:compfib}). The specific form of its Cartesian lifts are also a consequence of this lemma.
\begin{Lem}
	Let $\D\ra{Q}\B\xleftarrow{P}\E$ be fibrations over \B. Then their pullback
	$$\xymatrix@=1.3cm{
	**[l]{\D\times_\B\E}\ar[d]\ar[r]\ar[dr]^R & \E\ar[d]^P\\
	\D\ar[r]_{Q} & \B
	}$$
	is a fibration over \B.\\
	A Cartesian lift of an arrow $f\colon A\to B$ at $(D,E)$ is given by the pair of Cartesian lifts $(f^*D,f^*E)\ra{(\bar f_D,\bar f_E)}(D,E)$.
\end{Lem}

We now get to the main result. This is probably well known.
\begin{Prop}\label{prop:2prodfib}
	The pullback-fibration defined above yields a binary 2-product in $\FIB(\B)$ and in $\FIB_c(\B)$. Therefore, $\FIB(\B)$ and $\FIB_c(\B)$ have all finite 2-products.
\end{Prop}
\begin{Pf}
	First remark that the binary product in $\CAT$ is actually a binary 2-product. It is then not difficult to prove that the binary product of $\CAT/\B$ is a binary 2-product. Now, we have just seen that fibrations are preserved by this 2-product. It is also a binary 2-product in $\FIB_c(\B)$. Indeed, by the previous lemma, the projections of the product are Cartesian and the functors into the product induced by a pair of Cartesian functors in the component of the product are Cartesian. Finally, note that the identity functor $Id_\B$ of the category \B is a 2-terminal object of both $\FIB(\B)$ and $\FIB_c(\B)$. Indeed, it is a 2-terminal object in $\CAT/\B$. Moreover, $Id_\B$ is a fibration over \B (see example \ref*{ex:fibration}~(\ref{ex:fibration:isos}). Finally, it is a 2-terminal object in $\FIB_c(\B)$ because, as noticed above, a morphism of fibrations that is itself a fibration is Cartesian.
\end{Pf}

\subsection{Opfibrations and bifibrations}\label{ssec:opfibbifib}
In this part, we first study the notion dual to that of a fibration, called an \emph{opfibration}. We state some important duals of the theorems for fibrations and explain how the Grothendieck construction behaves with respect to dualization. We then turn to \emph{bifibrations}, which are both fibrations and opfibrations. We first recall a classical characterization of bifibrations in term of adjunctions, that goes back to Grothendieck \cite{SGA1}. We finally give a new characterization of their pseudo-functor counterpart.
\subsubsection{Opfibrations}
The basic notions of Cartesian morphisms, cleavages and (cloven) fibrations can be dualized. They are distinguished by an \emph{op} prefix: \emph{OpCartesian morphisms}\indexb{OpCartesian arrow}\index{OpCartesian lift|see{opCartesian arrow}}, \emph{opcleavage}\index{Opcleavage|see{Opfibration}} and (\emph{opcloven}) \emph{opfibrations}\indexb{Opfibration}\indexb{Opfibration!Opcloven --}. Given a functor $P\colon\E\to\B$, the ``\emph{op}-notion'' is defined to be the original notion but with respect to the dual functor $P^{op}\colon\E^{op}\to \B^{op}$. So a morphism in \E is opCartesian if it is Cartesian in $\E^{op}$ for the dual functor $P^{op}\colon\E^{op}\to\B^{op}$, and the functor $P$ itself is an opfibration if its dual $P^{op}$ is a fibration%
\footnote{In category theory, one usually distinguish a notion from its dual by the prefix \emph{co}-, like in the word ``colimits'' and all its particular cases. Grothendieck, in fact, called ``cofibration'' the dual of a fibration. From a purely categorical point of view, this name would have been preferable since it is self-understandable: just reverse the arrows in the definition. There is even a deeper reason; indeed, opfibrations are internal fibrations not in $\CAT^{op}$, but in $\CAT^{co}$ \cite{Web07}. Nevertheless, we decided here to follow the now quite established terminology. It has been introduced in order not to conflict with the topologists' intuition of a cofibration. It is coherent in a certain sense; the term ``fibration'' was most probably chosen because of the lifting property of such functors, which is similar to the lifting property defining a topological fibration. But categorical opfibrations are also defined by a lifting property, not by an \emph{extension} property like topological cofibrations. So if the name of fibrations was meant to be intuitive to topologists, the name of their dual should at least not be misleading for them.}.

Given a theorem about Cartesian morphisms, Cartesian functors or fibrations, one obtains its dual following the usual procedure. For instance, given an opfibration, apply the theorem to its dual, which is a fibration, and then go back to the opfibration by taking the dual of the results. Most of the time, the dual theorem is obvious and we will not even state it. However, we point out some results for opfibration when they are of great importance or when the dualization of the fibration result seems tricky. A useful observation is that taking the opposite of a category is an automorphism of the XL-category \CAT:
$$
(-)^{op}\colon \CAT\to\CAT.
$$
It is indeed its self-inverse\footnote{It is not a 2-functor from \CAT to \CAT though, since it reverses natural transformations. It is in fact a 2-isomorphism $(-)^{op}\colon\CAT^{co}\to\CAT$.}. Therefore, it preserves limits and colimits in \CAT.

\indexb{Opfibration!XLcategory@2-XL-category of --s}As for fibrations, opfibrations form a 2-XL-category, the sub-2-XL-category $\OPFIB$\index[not]{OPFIB@\OPFIB} of the strict arrow-2-XL-category $\CAT^{\two}$ consisting of opfibrations. We also consider restricting the morphisms to opCartesian functors, which we denote $\OPFIB_{oc}$\index[not]{OPFIBoc@$\OPFIB_{oc}$}. Finally, there are versions over a particular category \B, which we denote $\OPFIB(\B)$\index[not]{OPFIB(B)@\OPFIB(\B)} and $\OPFIB_{oc}(\B)$\index[not]{OPFIBoc(B)@$\OPFIB_{oc}(\B)$}. 

A choice of an opcleavage on an opfibration $P\colon\E\to\B$ determines a pseudo-functor $\Phi\colon\B\to\CAT$. Such pseudo-functors are called \emph{opindexed categories \op\emph{over} \B\fp}\indexb{Opindexed category}. Given an arrow $f\colon A\to B$ in the base, the functor $\Phi(f)$ is written
$$
f_*\colon \E_A\to\E_B,
$$
and called \emph{direct image functor of $f$}. There is a Grothendieck op-construction\indexb{Grothendieck!-- op-construction} for opindexed categories over \B, giving rise to an opfibration over \B. We give it here in detail for applications.

\index{Grothendieck!-- fibration|see{Fibration}}\index{Grothendieck!-- opfibration|see{Opfibration}}

\begin{Lem}[Dual of \ref{lem:grothcons}]\label{lem:grothopcons}
	Let $(\Phi,\gamma,\delta)\colon\B\to \CAT$ be an opindexed category. The following data, called the \linebreak\emph{Grothendieck op-construction of $\Phi$}, is a cloven opfibration, denoted $$\G_{op}(\Phi)=(\E_{\Phi}\ra{P_{\Phi}}\B).$$\index[not]{Gop@$\G_{op}$}\normalfont
\begin{itemize}
\item
Objects of $\E_{\Phi}$: Pairs $(B,E)$ where $B\in\B$ and $E\in\Phi(B)$.
\item
Morphisms of $\E_{\Phi}$: Pairs $(A,D) \ra{(f,h)} (B,E)$ where $f\colon A\to B$ is in \B and \linebreak${h\colon f_*D\to E}$ is in $\Phi(B)$.
\item
Composition in $\E_{\Phi}$: Let $(A,D) \ra{(f,h)} (B,E) \ra{(g,k)}(C,K)$ be a composable pair in $\E_{\Phi}$. Then its composite is the pair whose first component is $g\circ f$ and second component is the composite in $\E_A$:
$$\xymatrix@1@C=1.5cm{(g\circ f)_*D \ar[r]^{(\gamma_{f,g}^{-1})_D}_{\cong} & g_*f_*D \ar[r]^-{g_*(h)} & g_*E \ar[r]^{k} & K.}
$$
\item
Identities in $\E_{\Phi}$: $(A,D) \ra{(1_A,(\delta_A^{-1})_D)}(A,D)$.
\item
Functor $P_{\Phi}$: The projection on the first component.
\item
Cleavage of $P_{\Phi}$: The Cartesian lift of $f\colon A\to B$ at $(B,E)$ is the morphism
$$
	(A,E)\ra{(f,1_{f_*E})}(B,f_*E).
$$
\end{itemize}
\cqfd
\end{Lem}

Let us shortly study the dual relation between the Grothendieck construction and op-construction.\label{page:grothconsopcons} First note that a pseudo-functor $\Phi\colon\B\to\CAT$, being seen as a lax functor, gives rise to an oplax functor by inverting its identity and composition structure isomorphisms. Together with this inverted structure, it is thus a lax functor $$\Phi\colon\B\to\CAT^{co}.$$ One can therefore post-compose such a $\Phi$ with the 2-functor $(-)^{op}\colon\CAT^{co}\to\CAT$. 

Now, $\Phi$ can be considered either as an opindexed category over \B, or as an indexed category over $\B^{op}$. Take the Grothendieck construction $\G(\Phi)$ of the latter. This is a fibration over $\B^{op}$, and therefore, its dual $\G(\Phi)^{op}$ is an opfibration over \B. On the other hand, consider $\Phi$ as an opindexed category over \B. It is also a pseudofunctor $\Phi\colon\B\to\CAT^{co}$ and thus one obtains an opindexed catgory $^{op}\circ\Phi$ by post-composing $\Phi$ with the 2-functor $(-)^{op}$ already mentioned. Its Grothendieck op-construction is an opfibration $\G_{op}(^{op}\circ\Phi)$ over \B. These two processes yield equal opfibrations:
$$
(\G(\Phi))^{op}=\G_{op}((-)^{op}\circ\Phi).
$$

\begin{Def}
	The \emph{2-XL-category of opindexed categories}\indexb{Opindexed category!XLcategory of@2-XL-category of opindexed categories} is the 2-XL-category $(\CAT^{\B})_{lax}$ of pseudo-functors from \B to \CAT, lax natural transformations between them and modifications. It is denoted $\OPIND\B$\index[not]{OPIND(B)@\OPIND\B}. Its restriction to pseudo-natural transformations is written $\OPINDoc\B$\index[not]{OPINDoc(B)@\OPINDoc\B}.
\end{Def}

\begin{Thm}[Dual of \ref{thm:gcons}]\label{thm:grothopcons}
	The Grothendieck op-construction determines a 2-equivalence\index{Grothendieck!-- op-construction}\index[not]{Gop@$\G_{op}$}
	$$
	\mathcal G_{op}\colon (\CAT^{\B})_{lax}\xlongrightarrow{\simeq} \OPFIB(\B)
	$$
	which restricts to a 2-equivalence
	$$
	(\CAT^{\B})_p\xlongrightarrow{\simeq} \OPFIB_{oc}(\B).
	$$
\end{Thm}

Opfibrations are also stable under pullbacks. Indeed, as noticed above, taking the opposite is an automorphism of \CAT. So the opposite of a pullback square is a pullback square. Therefore, one obtains the following dual result.

\begin{Lem}[Dual of \ref{lem:pullfib}]\label{lem:pullopfib}
	Let $P\colon\E\to \B$ be an opfibration and $F\colon \A\to\B$ a functor. Then, pullbacks of $P$ along $F$ are also opfibrations, i.e., in a pullback-square in \CAT
	$$\xymatrix@=1.3cm{
	F^*\E\ar[d]_{\bar P}\ar[r]^{\bar F} & \E\ar[d]^P\\
	\A\ar[r]_{F} & \B,
	}$$
	the functor $\bar P\colon F^*\E\to\A$ is an opfibration.\\
	If one considers the canonical choice of pullback, then an opCartesian arrow in $F^*\E$ over $A\ra{f}A'$ in \A at $(A,E)$ is given by $(f,\underline{F(f)}_E)\colon(A,E)\to(A',F(f)_*E)$.
	
	\cqfd
\end{Lem}

\begin{Ex}\label{ex:presheafTop}
Given a category \A, we have defined the indexed category $$\PSh(-;\A)\colon\Cat^{op}\to\CAT$$ over $\Cat$ of presheaves with values in \A (see \thref{ex:presheafA}). In algebraic geometry, one considers its Grothendieck construction $PSh(\A)\to\Cat$, the Grothendieck fibration of presheaves with values in \A, or a pullback of it. An important case arises when pulling back this fibration along the functor $\T\colon\Top^{op}\to\Cat$, defined in \thref{ex:indcat}, that associates to a space $X$ its topology $\T(X)$ viewed as a preorder under inclusion. This is a fibration over $\Top^{op}$, and one is naturally led to consider its corresponding opfibration $$PSh_{\Top}(\A)^{op}\to\Top$$ over \Top.

By combining \thref{prop:grothpull} and the discussion following \thref{lem:grothopcons}, one obtains the latter opfibration by the following procedure\footnote{This fact is already stated in Grothendieck's paper on fibred categories \cite{SGA1}. We hope our explanation is useful.}: precompose $\PSh(-;\A)$ with $$\T^{op}\colon\Top\to\Cat^{op}$$ and postcompose it with $(-)^{op}\colon\CAT^{co}\to\CAT$. This gives rise to the \emph{opindexed category of presheaves on topological spaces with values in \A}\indexb{Sheaf!Opindexed category of presheaves on topological spaces}\indexb{Opindexed category!-- of presheaves on topological spaces|see{Sheaf}}:
$$
\PSh_{\Top}(-;\A)^{op}\colon\Top\longrightarrow\CAT.
$$
Its direct image functor associated to a map $f\colon X\to Y$ is the dual of the functor $$f_*\colon\PSh(X;\A)\to\PSh(Y;\A),$$ where $f_*=(f^{-1})^*$ in the notation of the indexed category $\PSh(-;\A)$.
When the context is clear, we denote this indexed category simply by $\PSh(-;\A)^{op}$. Its Grothendieck op-construction is, up to isomorphism, the opfibration $$\PSh_{\Top}(\A)^{op}\to\Top$$ defined above. It has the following form.
\begin{itemize}
\item
\textbf{Objects}: Pairs $(X,P)$ where $X$ is a topological space and $P\in\PSh(X;\A)$ is a presheaf on $X$ with values in \A.
\item
\textbf{Morphisms}: Pairs $(X,P)\ra{(f,\alpha)}(Y,Q)$, where $f\colon X\to Y$ is a continuous map and $\alpha\colon Q\Rightarrow f_*P$ is a morphism in $\PSh(Y;\A)$ and $f_*P=P\circ (f^{-1})^{op}$.
\end{itemize}
\end{Ex}

\subsubsection{Bifibrations}\label{sssec:Bifib}
We will quite often encounter functors that are both fibrations and opfibrations. These objects are called \emph{bifibrations}\indexb{Bifibration} and bifibrations with both a cleavage and an opcleavage, \emph{bicloven bifibrations}\indexb{Bifibration!Bicloven --}. Morphisms between them can be either plain functors, Cartesian functors, opCartesian functors or both. The latter are called \emph{biCartesian morphisms of bifibrations}\indexb{BiCartesian morphism of bifibrations|see{Bifibration}}. We denote these respective 2-XL-categories \BIFIB, $\BIFIB_c$, $\BIFIB_{oc}$, $\BIFIB_{bc}$, and similarly for their version over \B.\indexb{Bifibration!XLcategory of@2-XL-category of --s}\index[not]{BIFIB@\BIFIB}\index[not]{BIFIBc@$\BIFIB_c$}\index[not]{BIFIBoc@$\BIFIB_{oc}$}\index[not]{BIFIBbc@$\BIFIB_{bc}$}

We now give a characterization of bifibrations in terms of adjunctions of their inverse and direct image functors. This is a classical result \cite{SGA1,Shu08a}, which we give here with some additional practical information that we found very useful.

\begin{Prop}\label{prop:bifibadjoint}
	Let $P\colon\E\to\B$ be a bicloven bifibration. Then, for each morphism $f\colon A\to B$ in \B, there is an adjunction with left adjoint the direct image functor and right adjoint the inverse image functor:
	
	$$\xymatrix{
	f_*:\E_A \ar@<.9ex>[r]\ar@{}[r]|{\perp} & \E_B:f^*\ar@<.9ex>[l]
	}$$
	Its unit $\eta^f$ and counit $\epsilon^f$ are given by:
	$$\xymatrix{
	f^*f_*D\ar[dr]^{\bar f_{f_*D}} &\\
	D \ar@{-->}[u]^{\exists!\eta^f_{_D}} \ar[r]_{\underline{f}_D} & f_*D\\
	A\ar[r]_f & B}\quad\quad
	\xymatrix{& f_*f^*E\ar@{-->}[d]^{\exists!\epsilon^f_{_E}}\\
	f^*E\ar[r]_{\bar f_E}\ar[ru]^{\underline{f}_{f^*E}} & E\\
	A\ar[r]_f & B.}
	$$
	Moreover, if $h\colon D\to f^*E$ is in $\E_A$ and $k\colon f_*D\to E$ is in $\E_B$, then one obtains the transpose $h_\sharp\colon f_*D\to E$ of $h$, resp $k^\sharp\colon D\to f^*E$ of $k$, this way:
	$$\xymatrix{
	D\ar[r]^{\underline{f}_D}\ar[d]_h & f_*D\ar@{-->}[d]^{\exists!\, h_\sharp}\\
	f^*E\ar[r]_{\bar f_E} & E\\
	A\ar[r]_f & B}\quad\text{and}\quad\xymatrix{
	D\ar[r]^{\underline{f}_D}\ar@{-->}[d]_{\exists!\, k^\sharp} & f_*D\ar[d]^{k}\\
	f^*E\ar[r]_{\bar f_E} & E\\
	A\ar[r]_f & B
	}$$ 
	
	Conversely, if $P\colon\E\to\B$ is a cloven fibration such that each inverse image functor $f^*$ admits a left adjoint $f_*$ under an adjunction with unit $\eta^f$, then $P$ is a bifibration with opCartesian morphisms given by:
	$$\xymatrix{
	f^*f_*D\ar[dr]^{\bar f_{f_*D}} &\\
	D \ar[u]^{\eta^f_{_D}} \ar@{-->}[r]_{\underline{f}_D} & f_*D\\
	A\ar[r]_f & B.}
	$$
	If $P\colon\E\to\B$ is an opcloven opfibration such that each direct image functor $f_*$ admits a right adjoint $f^*$ under an adjunction with counit $\epsilon^f$, then $P$ is a bifibration with Cartesian morphisms given by:
	$$\xymatrix{
	&f_*f^*E\ar[d]^{\epsilon^f_E}\\
	f^*E\ar@{-->}[r]_{\bar f_E}\ar[ur]^{\underline f_{f^*E}} & E\\
	A \ar[r]^f & B
	}$$
	\cqfd
\end{Prop}
This proposition gives us a hint on how to define the concept of a \emph{bi-indexed category}, the pseudo-functor version of a bifibration. One could define it either as an opindexed category all of whose direct image functors have a right adjoint or as an indexed category all of whose inverse image functors have a left adjoint. One have then to prove that these two definition are equivalent ``up to isomorphism''. 

We do not define bi-indexed categories this way for several reasons. Firstly, there is an arbitrary choice to make between these two options : have the right adjoints as part of the data, and existence of left adjoints as an axiom (they are then determined only up to isomorphism), or the contrary. In practice, one is then force to go back and forth between these two equivalent categories, depending on the example one has to deal with. Secondly, one goes from (op-)fibrations to (op-)indexed categories by choosing a cleavage. Given a bifibration with a bicleavage, one obtains a choice of both the left and the right adjoints. Thus, it seems to us coherent that bi-indexed categories have both an associated opindexed and an associated indexed category (without going through a choice). The structure we obtain, which, to the best of our knowledge, is new, is an interesting application of the theory of \emph{double categories}.

\paragraph{}
	We now define a notion of \emph{bi-indexed category} such that a bicleavage on a bifibration determines such an object. This first requires defining a new 2-XL-category built from \CAT. There is a \emph{2-XL-category of adjunctions}\indexb{Adjunction!XLcategory@2-XL-category of --s} in \CAT, denoted $\ADJ$\index[not]{ADJ@\ADJ} \cite[IV.7-8]{McL97}\footnote{This in fact remains true for adjunctions in an arbitrary 2-(XL)-category $\mathds A$ \cite{KS74}, \cite[VII, Exercices 26-31]{Tay99},\cite[XII.4]{McL97}.}. Its objects are all categories. Arrows from $\A$ to $\B$ are adjunctions
	$$\xymatrix{
	F:\A \ar@<.9ex>[r]\ar@{}[r]|{\perp} & \B:G.\ar@<.9ex>[l]
	}$$
	Two-cells \begin{equation}\label{eq:conjugatepairs}\xymatrix@C=1.3cm{\A\rtwocell^{F\dashv G}_{F'\dashv G'}{} & \B}\end{equation}
	are called \emph{transformations of adjoints}, but also \emph{conjugate pairs}\indexb{Adjunction!Conjugate (pair)}\index{Conjugate (pair)|see{Adjunction}}. They are pairs $(\alpha,\beta)$ of natural transformations ${\alpha\colon F'\Rightarrow F}$ and ${\beta\colon G\Rightarrow G'}$ (note the reverse direction), called \emph{conjugates}, satisfying the following axiom:
	\begin{Axn}[Conjugate pair]\label{axiom:conjugate}
	For each arrow $F(A)\to B$ in \B, the transpose of the composite $F'(A)\ra{\alpha_A}F(A)\to B$ under the adjunction $F'\dashv G'$ is equal to the composite $A\to G(B)\ra{\beta_B}G'(B)$ of the transpose of $F(A)\to B$ under the adjunction $F\dashv G$ with $\beta_B$. 
	\end{Axn}
	
	Vertical composition of 2-cells is just the componentwise vertical composition of natural transformations.
	
	Now, consider two adjunctions
	$$\xymatrix{
	\A \ar@<.9ex>[r]^F\ar@{}[r]|{\perp} & \B\ar@<.9ex>[l]^G \ar@<.9ex>[r]^{F'}\ar@{}[r]|{\perp} & \C\ar@<.9ex>[l]^{G'}
	}$$
	with units $\eta$ and $\eta'$, and counits $\epsilon$ and $\epsilon'$, respectively. Their composite is the adjunction
	$$\xymatrix{
	F'\circ F:\A \ar@<.9ex>[r]\ar@{}[r]|{\perp} & \C:G\circ G'\ar@<.9ex>[l]
	}$$
	with unit $(G\cdot\eta'\cdot F)\bullet \eta$ and counit $\epsilon'\bullet (F'\cdot\epsilon\cdot G')$. Horizontal composition of 2-cells is given by componentwise horizontal composition of natural transformations.
	
	There are forgetful 2-functors $L\colon \ADJ^{co}\to\CAT$ and $R\colon \ADJ^{op}\to\CAT$ that select the left (resp. right) adjoint functors and the natural transformations between these. They are trivially surjective on objects. More interestingly, their local functors $L_{\A,\B}\colon\ADJ(\A,\B)^{op}\to\B^{\A}$ and $R_{\A,\B}\colon\ADJ(\A,\B)\to\A^{\B}$ are full and faithful. In other words, given two adjunctions from \A to \B, a natural transformation $\alpha$ between the left adjoints determines a unique natural transformation $\beta$ between the right adjoints such that the pair $(\alpha,\beta)$ is a morphism of these adjunctions. The same holds when one starts with a natural transformation between the right adjoints. Thus, $L$ and $R$ are 2-equivalences between $\ADJ^{co}$ (resp. $\ADJ^{op}$) and the 2-cell full sub-2-XL-category of \CAT consisting of all categories and left adjoint (resp. right adjoint) functors. In particular, the bijection between conjugate pairs and their left or right component preserves vertical and horizontal compositions, as well as identities.	
	
\begin{Lem}\label{lem:biindexedbifib}
	Let $P\colon\E\to\B$ be a bicloven bifibration. Consider the associated opindexed category ${\Phi_P}_*\colon\B\to\CAT$ with structure isomorphisms ${\gamma_{f,g}}_*\colon(g\circ f)_*\Rightarrow g_*\circ f_*$, for each composable pair $A\ra f B\ra g C$ in \B and ${\delta_A}_*\colon {1_A}_*\Rightarrow Id_{\E_A}$, for each object $A\in\B$. Consider also the associated indexed category ${\Phi_P}^*\colon\B^{op}\to\CAT$ with structure isomorphisms ${{\gamma_{f,g}}^*\colon f^*\circ g^*\Rightarrow(g\circ f)^*}$ and ${\delta_A}^*\colon Id_{\E_A}\Rightarrow {1_A}^*$.\Par	
	Then, there is a pseudo-functor
	$$\begin{array}{rcl}
\Omega_P\colon\B& \longrightarrow	& \ADJ\\
A 		& \longmapsto		& \E_A\\
A \ra{f} B 	& \longmapsto		& \xymatrix{f_*\colon\E_A \ar@<.9ex>[r]\ar@{}[r]|{\perp} & \E_B:f^*,\ar@<.9ex>[l]}
\end{array}
$$
whose structure isomorphisms are the conjugate pairs $({\gamma_{f,g}}_*,{\gamma_{f,g}}^*)$ and $({\delta_A}_*,{\delta_A}^*)$.
\end{Lem}
\begin{Pf}[Sketch]
	One proves that the latter pairs are conjugate by means of Axiom \ref{axiom:conjugate} and by the characterization of the transpose morphisms for an adjunction $f_*\dashv f^*$ given in Proposition \ref{prop:bifibadjoint}. The coherence axioms of a pseudo-functor are automatically satisfied because the compositions of natural transformations in \ADJ are compositions in \CAT on each component of conjugate pairs. Again, pseudo-functor coherence axioms were already satisfied in ${\Phi_P}_*$ and ${\Phi_P}^*$.
\end{Pf}
The previous lemma leads to a definition of bi-indexed categories.

\begin{Def}[First definition]
	A \emph{bi-indexed category over \B}\indexb{Bi-indexed category} is a pseudofunctor
	$$\begin{array}{rcl}
\Omega\colon\B& \longrightarrow	& \ADJ\\
A 		& \longmapsto		& \Omega(A)\\
A \ra{f} B 	& \longmapsto		& \xymatrix{\Omega(f)_*\colon\Omega(A) \ar@<.9ex>[r]\ar@{}[r]|{\perp} & \Omega(B):\Omega(f)^*\ar@<.9ex>[l]}
\end{array}
$$
It comes thus for each object $B\in\B$ with a conjugate pair $\delta_B=({\delta_B}_*,{\delta_B}^*)$ of natural isomorphisms, where ${\delta_B}_*\colon\Omega(1_B)_*\Rightarrow Id_{\Omega(B)}$ and ${\delta}_B^*\colon Id_{\Omega(B)}\Rightarrow\Omega(1_B)^*$ and similarly for the composition isomorphisms $\gamma_{f,g}$.
\end{Def}
The lax functor $\Omega$ gives rise to an oplax functor (thus a lax functor $\B\to\CAT^{co}$) by inverting the structure isomorphisms $\delta$ and $\gamma$, and we also denote it $\Omega$. Any bi-indexed category determines thus both an opindexed and an indexed category by post-composition with $L$ and $R$ respectively:
\begin{equation}\label{eq:ADJLR}
\B\ra{\Omega}\ADJ^{co}\ra{L}\CAT\quad\text{and}\quad\B^{op}\ra{\Omega^{op}}\ADJ^{op}\ra R\CAT.
\end{equation}
It would be natural now to define morphisms of bi-indexed categories over \B as lax natural transformations of pseudo-functors from \B to \ADJ. This is not correct though, because a functor $F$ over \B between two bifibrations does not give rise to such a transformation. Indeed, the component at an object $B\in\B$ of the induced lax natural transformation between the associated opindexed categories is just the restriction of $F$ to the fibres over $B$. This is a plain functor, not an adjunction, and therefore not a morphism in \ADJ. One sees that two different kinds of morphisms are needed in the target of bi-indexed categories in order to capture both bi-indexed categories and their morphisms.

There is actually a \emph{double XL-category}\indexb{Adjunction!Double XL-category of --s} of adjunctions in \CAT, that we denote $\mathfrak{ADJ}$\index[not]{ADJ2@$\mathfrak{ADJ}$} \cite{GP99,KS74,Shu08a,Pal71}\footnote{The notion of double category goes back to Ehresman\cite{Ehr65}. A recent good introduction can be found in \cite{Fio07}. A deep treatment of the subject can be found in \cite{GP99,GP04}. See also \cite{DPP06}.}. Recall in short that a \emph{double category}\indexb{Double category}  $\frak A$ consists of a class of objects $\Ob\frak A$, a class of horizontal arrows Hor$\frak A$, a class of vertical arrows Ver$\frak A$ and a class of cells Cell$\frak A$. Arrows have objects as domain and codomain. We differentiate horizontal and vertical arrows by denoting the first ones by normal arrows and the second ones by dotted arrows $\xymatrix@1{A\ar[r]|{\bullet} & B}$. Objects and horizontal (resp. vertical) arrows form a category whose identities we denote $1_A^h\colon A\to A$ (resp. $\xymatrix@1{1^v_A\colon A\ar[r]|-{\bullet} & A}$). The cells are squares
$$
\xymatrix{A\ar[r]^{f}\ar[d]|{\bullet}_{h} & B\ar[d]|{\bullet}^{k}\\
C\ar[r]_{g} & D\ar@{}[ul]|{\alpha}}
$$
with horizontal domain and codomain $h$ and $k$ and vertical domain and codomain $f$ and $g$. Cells admit both horizontal and vertical composition. Horizontal (resp. vertical) arrows and cells with vertical (resp. horizontal) composition form a category. We denote
$$
\xymatrix{A\ar[r]^{f}\ar[d]|{\bullet}_{1^v_A} & B\ar[d]|{\bullet}^{1^v_B}\\
A\ar[r]_{f} & B\ar@{}[ul]|{\iota^v_f}}\quad\quad
\xymatrix{A\ar[r]^{1_{A}^h}\ar[d]|{\bullet}_{k} & A\ar[d]|{\bullet}^{k}\\
B\ar[r]_{1^h_B} & B\ar@{}[ul]|{\iota^h_k}}
$$
the identities of these respective categories. We denote double categories by Gothic letters $\frak A,\frak B,\dots,\frak H,\frak I,\ldots,\frak U,\frak V,\dots$. 

Now, the double category $\frak{ADJ}$ is defined as follows. Its objects are all categories, its vertical arrows are functors (with composition of \CAT), its horizontal arrows are adjunctions (with composition of $\ADJ$) and its cells are \emph{adjoint squares} or \emph{mate pairs}\indexb{Mate (pair)|see{Adjunction}}\indexb{Adjunction!mate (pair)}, which we define now. An \emph{adjoint square} consists of the data
$$\xymatrix@=1.5cm{
\A\ar[d]_H \ar@<.9ex>[r]^F\ar@{}[r]|{\perp}\ar@{}[rd]|{(\xi,\zeta)} & \B\ar@<.9ex>[l]^G\ar[d]^{K}\\
\A' \ar@<.9ex>[r]^{F'}\ar@{}[r]|{\perp} & \B'\ar@<.9ex>[l]^{G'}
}$$
where $\xi$ and $\zeta$ are natural transformations, called \emph{mates}, filling in the following diagrams
$$\xymatrix{
\A\ar[r]^F\drtwocell<\omit>{^\xi}\ar[d]_H & \B\ar[d]^K\\
\A'\ar[r]_{F'} & \B'}\hspace{1.5cm}
\xymatrix{
\A\ar[d]_H & \B\ar[d]^K\ar[l]_G\\
\A'\urtwocell<\omit>{\zeta} & \B'\ar[l]^{G'}
}$$
and satisfying the following axiom. 
\begin{Axn}[Mate pair]\label{axiom:matepair}
Given an arrow $F(A)\ra{f}B$ in \B, consider the composite $$F'H(A) \ra{\xi_A}KF(A) \ra{K(f)}K(B),$$ and let $H(A) \ra{h}G'K(B)$ be its transpose under $F'\dashv G'$. On the other hand, consider the transpose $A \ra{g}G(B)$ of $f$ under $F\dashv G$ and the composite $$H(A) \ra{H(g)}HG(B)\ra{\zeta_B}G'KB.$$ The latter must coincide with the arrow $h$ just defined.
\end{Axn}

See \cite{KS74} for the axiom (beautifully) expressed in terms of diagrams of 2-cells. Vertical and horizontal compositions of cells are given by componentwise \emph{pasting} in \CAT. Their respective identities are
$$\xymatrix@=1.3cm{
\A\ar[d]_{Id} \ar@<.9ex>[r]^F\ar@{}[r]|{\perp}\ar@{}[rd]|{(\iota_F,\iota_G)} & \B\ar@<.9ex>[l]^G\ar[d]^{Id}\\
\A \ar@<.9ex>[r]^{F}\ar@{}[r]|{\perp} & \B\ar@<.9ex>[l]^{G}
}\quad
\xymatrix@=1.3cm{
\A\ar[d]_H \ar@<.9ex>[r]^{Id}\ar@{}[r]|{\perp}\ar@{}[rd]|{(\iota_H,\iota_H)} & \A\ar@<.9ex>[l]^{Id}\ar[d]^{H}\\
\A' \ar@<.9ex>[r]^{Id}\ar@{}[r]|{\perp} & \A'\ar@<.9ex>[l]^{Id}
}$$

\begin{Rem}\label{rem:bijmates}
Notice that conjugate pairs defined above \hyperref[eq:conjugatepairs]{(\ref*{eq:conjugatepairs})} are precisely the mate pairs whose vertical functors are identities. Like conjugates, mates $\xi$ and $\zeta$ of a mate pair $(\xi,\zeta)$ uniquely determine each other \cite{KS74}. In fact, we have decided to have both natural transformations in the definition, in order not to make a choice and to have them both available in practice. However, one should remember that the bijection between mate natural transformations preserves all the structure of the double category. This means that it preserves vertical and horizontal composition, as well as vertical and horizontal identities. Therefore, when checking the commutativity of a diagram of cells in the double category $\mathfrak{ADJ}$, one only needs to check it on the 2-cells between the left adjoints, or on the 2-cells between the right adjoints. The commutativity is then ensured for their mates.
\end{Rem}

Given a 2-category $\mathds A$, there is a horizontal way of associating to it a double category, the \emph{horizontal double category of \dA}\indexb{Double category!Horizontal -- of a 2-category}, that is denoted $\mathfrak H\mathds A$. It has the same objects as $\dA$, its horizontal morphisms are the morphisms of \dA, its vertical morphisms are vertical identities and its cells are 2-cells in \dA. On the other hand, given a double category $\mathfrak A$, there is a horizontal way of associating to it a 2-category $\mathds H\mathfrak A$, called the \emph{horizontal 2-category of $\frak A$}\indexb{Double category!Horizontal 2-category of a --}. They have the same objects, the morphisms of $\mathds H\mathfrak A$ are the horizontal morphisms of $\mathfrak A$ and the 2-cells are cells of $\mathfrak A$ whose vertical arrows are identities. Symmetrically, one has a \emph{vertical double category}\indexb{Double category!Vertical -- of a 2-category} $\mathfrak V \mathds{A}$ associated to a 2-category $\mathds{A}$ and a \emph{vertical 2-category}\indexb{Double category!Vertical 2-category of a --} $\mathds V\mathfrak A$ associated to a double category $\frak A$.

\begin{Ex}
	The horizontal 2-category $\mathds{H}\frak{ADJ}$ is precisely the 2-category \ADJ of adjunctions, whereas its vertical 2-category $\mathds{V}\frak{ADJ}$ is isomorphic to \CAT.
\end{Ex}

\begin{Def}
	A \emph{lax double functor}\indexb{Double functor!Lax --} $F\colon\frak A\to\frak B$ between double categories consists of the following data. 
	\begin{itemize}
\item
It is a function that associate objects, horizontal arrows, vertical arrows and cells of $\frak A$ to the same type of structure in $\frak B$, respecting domains and codomains.
\item
It preserves strictly the vertical structure, i.e., composition and identities of vertical arrows, and vertical composition and vertical identities of cells.
\item
It preserves the horizontal structure up to cells. More precisely, there are, for each object $A\in\frak A$, a cell
$$
\xymatrix{FA\ar[r]^{1_{FA}^h}\ar[d]|{\bullet}_{1^v_{FA}} & FA\ar[d]|{\bullet}^{1^v_{FA}}\\
FA\ar[r]_{F(1^h_A)} & FA\ar@{}[ul]|{\delta_A}}
$$
and, for each pair of composable horizontal arrows $A \ra{f} B \ra{g} C$, a cell 
$$\xymatrix{
FA\ar[r]^{Ff}\ar[d]|{\bullet}_{1^v_{FA}} & FB\ar[r]^{Fg} & FC\ar[d]|{\bullet}^{1^v_{FC}}\\
FA\ar[rr]_{F(g\circ f)} && FC.\ar@{}[ull]|{\gamma_{f,g}}
}$$
These data are subject to axioms, two axioms of \emph{naturality} and two axioms of coherence (see the dual notion of oplax double functor in \cite[1.9]{DPP06}). A \emph{pseudo double functor}\indexb{Double functor!Pseudo --} is a lax functor whose structure cells $\delta$ and $\gamma$ are isomorphisms for the vertical composition of cells.
\end{itemize}
\end{Def}
Now, one easily verifies that there is a bijection between lax double functors from $\frak H \dA$ to $\frak B$ and lax functors from $\dA$ to $\mathds H\frak B$:
\begin{equation}\label{eq:laxdouble}
\leftidx{_{lax}}{\DBL(\frak H\dA,\frak B)}{}\cong \leftidx{_{lax}}{\BICAT(\dA,\mathds{H}\frak B)}{}.
\end{equation}
This bijection restricts to pseudo objects. We thus obtain our second definition of a bi-indexed category.

\begin{Def}[Second definition]
	A \emph{bi-indexed category over \B}\indexb{Bi-indexed category} is a pseudo double functor $$\Omega\colon\frak H \B\longrightarrow\frak{ADJ},$$ where $\B$ is seen as 2-category with only identity 2-cells.
\end{Def}

We now define the 2-XL-category of bi-indexed categories. This requires some definitions of double category theory.

\begin{Defs}
	\item A \emph{vertical transformation}\indexb{Vertical transformation} $\tau$ between lax double functors $$\tau\colon F\Rightarrow G\colon\frak A\to\frak B$$ is given by the following data
	\begin{itemize}
\item
For each object $A\in\frak A$, a \emph{vertical} arrow $\xymatrix@1{FA\ar[r]|{\bullet}^{\tau_A} & GA}$,
\item
For each \emph{horizontal} arrow $f\colon A\to B$, a cell
$$\xymatrix@=1.3cm{
FA\ar[d]|{\bullet}_{\tau_A} \ar[r]^{Ff} & FB\ar[d]|{\bullet}^{\tau_B}\\
GA \ar[r]_{Gf} & GB\ar@{}[ul]|{\tau_f}
}$$
\item
These data are subject to axioms: two axioms of \emph{naturality} and two axioms of \emph{coherence} with the structures of $F$ and $G$ (see \cite[1.13]{DPP06}).
\end{itemize}
	\item A \emph{modification}\indexb{Modification of vertical transformations} between two vertical transformations $\sigma,\tau\colon F\Rightarrow G\colon\frak A\to\frak B$ is given by the following data\footnote{I obtained this notion by adapting the concept of \emph{strong modification} as defined in \cite{GP99} in the case where the strong horizontal transformations are just identities (recall that the role of horizontal and vertical is reversed in our text regards to conventions adopted in \cite{GP99}).}:
	\begin{itemize}
		\item
		For each object $A\in\frak A$, a cell in $\frak B$
		$$\xymatrix@=1.3cm{
FA\ar[d]|{\bullet}_{\sigma_A} \ar[r]^{1^h_{FA}} & FA\ar[d]|{\bullet}^{\tau_A}\\
GA \ar[r]_{1^h_{GA}} & GA.\ar@{}[ul]|{\Xi_A}
}$$
		\item
		These data are subject to the following axioms:
		\begin{itemize}
\item[(M1)]
For each horizontal arrow $f\colon A\to B$ in $\frak A$,
$$\xymatrix@=1.3cm{
FA \ar[d]|{\bullet}_{\sigma_A} \ar[r]^{1^h_{FA}} & FA \ar[r]^{Ff} \ar[d]|{\bullet}^{\tau_A} & FB \ar[d]|{\bullet}^{\tau_B}\\
GA \ar[r]_{1^h_{GA}} & GA \ar[r]_{Gf} \ar@{}[ul]|{\Xi_A}& GB \ar@{}[lu]|{\tau_f}
}\quad = \quad 
\xymatrix@=1.3cm{
FA \ar[d]|{\bullet}_{\sigma_A} \ar[r]^{Ff} & FB \ar[r]^{1^h_{FB}} \ar[d]|{\bullet}^{\sigma_B} & FB \ar[d]|{\bullet}^{\tau_B}\\
GA \ar[r]_{Gf} & GB \ar[r]_{1^h_{GB}} \ar@{}[ul]|{\sigma_f}& GB. \ar@{}[lu]|{\Xi_B}
}$$
\item[(M2)]
For each vertical arrow $\xymatrix@1{f\colon A\ar[r]|-{\bullet} & B}$ in $\frak A$,
$$\xymatrix@=1.3cm{
FA \ar[r]^{1^h_{FA}} \ar@{}[dr]|{\Xi_A} \ar[d]|{\bullet}_{\sigma_A} & FA \ar[d]|{\bullet}^{\tau_A}\\
GA \ar[d]|{\bullet}_{Gf} \ar[r]^{1^h_{GA}} \ar@{}[rd]|{\iota^h_{Gf}} & GA\ar[d]|{\bullet}^{Gf}\\
GB \ar[r]_{1^h_{GB}} & GB
}\quad=\quad
\xymatrix@=1.3cm{
FA \ar[r]^{1^h_{FA}} \ar@{}[dr]|{\iota^h_{Ff}} \ar[d]|{\bullet}_{Ff} & FA \ar[d]|{\bullet}^{Ff}\\
FB \ar[d]_{\sigma_B} \ar[r]^{1^h_{FB}} \ar@{}[rd]|{\Xi_B} & FB\ar[d]|{\bullet}^{\tau_B}\\
GB \ar[r]_{1^h_{GB}} & GB.
}$$
\end{itemize}
\end{itemize}
\end{Defs}

\begin{Lem}
	Given double categories $\frak A$ and $\frak B$, the lax double functors from $\frak A$ to $\frak B$, their vertical transformations and modifications organize into a 2-XL-category $$\leftidx{_{lax}}{\DBL(\frak A,\frak B)}{}.$$
\end{Lem}
\begin{Pf}
	The proof follows quite directly from the axioms. One uses an important coherence property of double categories that is a consequence of the axioms and that is worth mentioning \cite{Fio07}: given an object $A$ of a double category $\frak A$, there is an equality of identity cells
$$\xymatrix@=1.3cm{
A\ar[r]^{1^h_A}\ar[d]_{1^v_A} & A\ar[d]^{1^v_A}\\
A\ar[r]_{1^h_A} & A\ar@{}[ul]|{\iota^v_{1^h_A}}
}\quad=\quad\xymatrix@=1.3cm{
A\ar[r]^{1^h_A}\ar[d]_{1^v_A} & A\ar[d]^{1^v_A}\\
A\ar[r]_{1^h_A} & A\ar@{}[ul]|{\iota^h_{1^v_A}}
}$$
\cqfd
\end{Pf}

\begin{Def}
	The \emph{2-XL-category of bi-indexed categories over \B}\indexb{Bi-indexed category!XLcategory of@2-XL-category of bi-indexed categories} is the 2-XL-category $\DBL(\frak H\B,\frak{ADJ})$ of pseudo double functors from $\frak H\B$ to $\frak{ADJ}$, their vertical transformations and modifications, where $\B$ is considered as a 2-category with only identity 2-cells. It is also denoted $\BIIND\B$.\index[not]{BI-IND(B)@\BIIND\B}
\end{Def}
We denote $\leftidx{_{\mathit{right}}}{\IND\B}{}$ the full sub-2-XL-category of $\IND\B$ of indexed categories whose inverse image functors are right adjoints. Similarly, we denote $\leftidx{_{\mathit{left}}}{\OPIND\B}{}$ the full sub-2-XL-category of $\OPIND\B$ of opindexed categories whose direct image functors are left adjoints. We define ``forgetful'' 2-functors:
$$
\leftidx{_{\mathit{left}}}{\OPIND\B}{}\xlongleftarrow{U}\BIIND\B\xlongrightarrow{V}\leftidx{_{\mathit{right}}}{\IND\B}{}
$$
We describe $U$, and $V$ is similar. Given a pseudo double functor $$(\Omega,\gamma,\delta)\colon\frak H\B\to\frak{ADJ},$$ there is, under bijection \hyperref[eq:laxdouble]{(\ref*{eq:laxdouble})}, a pseudo-functor $\B\to\ADJ$, which we denote the same way, $(\Omega,\gamma,\delta)$. The opindexed category considered is the composite $\B\ra{\Omega}\ADJ^{co}\ra{L}\CAT$, where $L$ is the 2-functor defined in \hyperref[eq:ADJLR]{(\ref*{eq:ADJLR})}. In detail, it has the same action as $\Omega$ on objects and it selects the left adjoint $\Omega(f)_*$ for a morphism $f$ in \B. Its structure isomorphisms are defined as follows. For each $B\in\B$ and for each composable pair $A\ra f B\ra g C$ in \B,
$$\xymatrix@C=1.5cm{
\Omega(A)\rtwocell<5>^{Id_{\Omega(A)}}_{\Omega(1_A)_*}{\quad\delta_{A*}^{-1}}&\Omega(A), & \Omega(A)\ar[r]^{\Omega(f)_*} \rrlowertwocell<-10>_{\Omega(g\circ f)_*}{\quad\gamma_{f,g*}^{-1}}& \Omega(B) \ar[r]^{\Omega(g)_*} & \Omega(C).
}$$
We denote $\Omega_*\coloneqq U(\Omega)$.

Next, if $\tau\colon\Omega\Rightarrow\Omega'$ is a vertical transformation, the it gives rise to a lax transformation $\tau_*\colon\Omega\Rightarrow\Omega'$ this way. For each $B\in\B$, take the same functor $\tau_{*A}=\tau_A\colon\Omega(A)\to\Omega'(A)$, and for each arrow $f\colon A\to B$ in \B, it select the left part of the mate pair $\tau_f$, i.e., ${(\tau_*)_f\coloneqq(\tau_f)_*}$.

Finally, given a modification $\Xi\colon\sigma\Rrightarrow\tau$ between vertical transformations $\Omega\Rightarrow\Omega'$, one obtains a modification between the corresponding lax transformations $\Xi\colon\sigma_*\Rrightarrow\tau_*$. Indeed, for each $B\in\B$, the component at $B$ of $\Xi\colon\sigma\Rrightarrow\tau$ is a cell in $\frak{ADJ}$
$$\xymatrix{
\Omega(A) \ar[d]_{\sigma_A}\ar@{=}[r]& \Omega(A)\ar[d]^{\tau_A}\\
\Omega'(A)\ar@{=}[r] & \Omega'(A).\ar@{}[lu]|{\Xi_A}
}$$
In this case, the mate pair $\Xi_A=(\Xi_{A*},\Xi_A^*)$ must have its mates equal $\Xi_{A*}=\Xi_A^*$, and that is why we do not write $\Xi_*$, but just $\Xi$ for the induced modification of lax transformations.
\begin{Thm}\label{thm:Bifib}\index{Bifibration!XLcategory of@2-XL-category of --s}\index{Grothendieck!-- construction}\index{Grothendieck!-- op-construction}
	The forgetful 2-functors $U$ and $V$ defined above are 2-equivalences. Moreover, there is a 2-natural isomorphism filling the following diagram of 2-equivalences.
	$$\xymatrix@=1.5cm{
	\BIIND\B \ar[d]_U^\sim \ar[r]^V_\sim & \leftidx{_{\mathit{right}}}{\IND\B}{}\ar[d]_\sim^{\G}\\
	\leftidx{_{\mathit{left}}}{\OPIND\B}{} \ar[r]^\sim_{\G_{op}} & \BIFIB(\B)\ar@{}[lu]|{\cong}
	}$$
\end{Thm}
\begin{Pf}
	The fact that $U$ and $V$ are 2-equivalences is based on remark \ref{rem:bijmates} (which is essentially Proposition 2.2 of \cite{KS74}): the bijection between mate natural transformations of adjunctions preserves all the double category structure of $\frak{ADJ}$. The fact that Grothendieck (op-)construction 2-functor restricts to a 2-functor as indicated in the diagram is a consequence of Proposition \ref{prop:bifibadjoint} and of the fact that Grothendieck construction is a 2-equivalence. Indeed, given, for instance, an indexed category $\Phi\colon\B^{op}\to\CAT$, its Grothendieck construction is a canonically cloven fibration $P_{\Phi}$. Moreover, the indexed category $\Phi_{P_{\Phi}}$ associated to this cloven fibration is strict-naturally isomorphic to the original pseudo-functor $\Phi$. Thus, given a morphism $f$ in \B, the left adjoint to $\Phi(f)$ induces a left adjoint to the inverse image functor $f^*$ induced by the cleavage on $\E_{\Phi}$. One can then apply the Proposition \ref{prop:bifibadjoint}. In addition, these restricted Grothendieck (op)-constructions are 2-equivalences because the non-restricted 2-functors are equivalences (Theorems \ref{thm:gcons} and \ref{thm:grothopcons}) and because of Proposition \ref{prop:bifibadjoint}.
	
	We now turn to the 2-natural isomorphism filling the diagram of the theorem. Consider a pseudo double functor $\Omega\colon\frak H\B\to\frak{ADJ}$. Let $\Omega_*=U(\Omega)$ and $\Omega^*=V(\Omega)$ and consider their Grothendieck op-construction, resp. construction. These are respectively bifibrations $P_{\Omega_*}\colon\E_{\Omega_*}\to\B$ and $P_{\Omega^*}\colon\E_{\Omega^*}\to\B$. The isomorphism over \B between them is defined as follows. It acts as identity on objects. Let $(f,h)\colon(A,D)\to(B,E)$ be a morphism in $\E_{\Omega_*}$. Thus, $h\colon f_*D\to E$. We turn it into a morphism $(f,h^\sharp)\colon(A,D)\to(B,E)$ of $\E_{\Omega^*}$ with $h^\sharp\colon D\to f^*E$ the adjoint morphism of $h$ under $\Omega(f)_*\dashv\Omega(f)^*$. In order to prove that this indeed is a functor, one uses the fact that the composition and identity structure isomorphisms of $\Omega_*$ and $\Omega^*$ are conjugates. Moreover, this functor is obviously an isomorphism. 
	
	One proves the naturality of this isomorphism with respect to vertical transformation $\tau$ by the fact that the functors $U$ and $V$ select the respective mates of the cell $\tau_f$ for every $f$ in \B. The naturality with respect to modifications uses the conjugate relation between the identity structure isomorphisms $\delta_*$ and $\delta^*$ of $\Omega_*$ and $\Omega^*$.
\end{Pf}

\begin{Rem}\label{rem:choiceadjpseudo}
As a consequence of this theorem, given, for instance, an op-indexed category whose direct image functors are left adjoints, any choice of right adjoints for each direct image functor gives rise to an indexed category. Moreover, the Grothendieck construction of the latter is isomorphic the Grothendieck op-construction of the former.	
\end{Rem}

We haven't addressed yet the question of the restriction of the diagram of this theorem to morphisms of bifibrations that are Cartesian, opCartesian or biCartesian functors. (Op-)Cartesian functors correspond to pseudo-natural transformations of (op-)indexed categories under the Grothendieck (op-)construction 2-equivalence. Yet, this condition of pseudo-naturality (instead of mere lax or oplax naturality) does not correspond to a double category concept, and that is why a natural transformation of a mate pair can be an isomorphism without its mate being an isomorphism. Let us explain this a little more. Consider an adjoint square, i.e., a cell of $\frak{ADJ}$:
	$$\xymatrix@=1.5cm{
\A\ar[d]_H \ar@<.9ex>[r]^F\ar@{}[r]|{\perp}\ar@{}[rd]|{\tau=(\tau_*,\tau^*)} & \B\ar@<.9ex>[l]^G\ar[d]^{K}\\
\A' \ar@<.9ex>[r]^{F'}\ar@{}[r]|{\perp} & \B'\ar@<.9ex>[l]^{G'}
}$$
One cannot express by a condition in the double category $\frak{ADJ}$ that $\tau_*$ and $\tau^*$ are natural isomorphisms, because the pair $((\tau_*)^{-1},(\tau^*)^{-1})$ formed by their inverses does not live in this double category (and the composition under which they are inverse from each other is neither the vertical, nor the horizontal composition $\frak{DBL}$, but some kind of ``diagonal composition''). 

Now, there are 2-cell and 0-cell full sub-2-XL-categories $\BIINDc\B$, $\BIINDoc\B$ and $\BIINDbc\B$ of $\BIIND\B$ of vertical transformations $\tau$ such that $\tau_*$, or $\tau^*$, or both of them are natural isomorphisms (one can indeed check that these conditions are preserved by vertical and horizontal composition and identities). The indices $c$, $oc$ and $bc$ denote respectively the words Cartesian, opCartesian and biCartesian. The following result is then clear.

\begin{Cor}
The 2-equivalences of the theorem restricts to 2-equivalences
\begin{equation}\label{eq:biindc}
\xymatrix{\BIINDc\B  \ar[r]^{\sim} & \INDc\B \ar[r]^\sim & \BIFIB_c(\B)}
\end{equation}
\begin{equation}\label{eq:biindoc}
\xymatrix{\BIINDoc\B  \ar[r]^\sim & \OPINDoc\B \ar[r]^\sim & \BIFIB_{oc}(\B)}
\end{equation}
\cqfd
\end{Cor}

Consider then the case of $\BIINDbc\B$. Let $\tau\colon\Omega\Rightarrow\Omega'$ be a morphism in this 2-XL-category. Its image $\G(\tau^*)$ under the composite \hyperref[eq:biindc]{(\ref*{eq:biindc})} is thus a Cartesian functor. Again, by naturality of the isomorphism of Theorem \ref{thm:Bifib}, the following diagram commutes:
$$\xymatrix@=1.3cm{
\E_{\Omega^*}\ar[r]^{\G\tau^*}\ar[d]_\cong & \E_{\Omega'^*}\ar[d]^\cong\\
\E_{\Omega_*}\ar[r]_{\G_{op}\tau_*} & \E_{\Omega'_*}.
}$$
Thus, as $\G\tau_*$ is an opCartesian functor, so is $\G\tau^*$. One similarly shows that any biCartesian functor $F\colon\E_{\Omega^*}\to\E_{\Omega'^*}$ over \B comes from a vertical transformation $\tau\colon\Omega\Rightarrow\Omega'$ whose both components $\tau_*$ and $\tau^*$ are natural isomorphisms. In conclusion, the following result holds.

\begin{Cor}
The composites (\ref*{eq:biindc}) and (\ref*{eq:biindoc}) restrict to 2-isomorphic 2-equivalences
$$\xymatrix@1@C=1.5cm{
**[l]{\BIINDbc\B\quad}\rtwocell<5>^{\sim}_{\sim}{\ \cong} & **[r]{\qquad\BIFIB_{bc}(\B).}
}$$
\cqfd
\end{Cor}

\begin{Exs}[We refer to the examples of fibrations of Exemples \ref{ex:fibration}.]\label{ex:bifibration}
\item
Any category \C can be viewed as a bifibration $\C\to\one$ over the terminal object $\one$ of \CAT, or over itself $\C\ra{Id_C}\C$.
\item
Any product $\prod \A_i\to\A_i$ of categories onto one of the factors is a bifibration.
\item
The canonical fibration $\C^\two\to\C$ is a bifibration, with a canonical opcleavage. If $f\colon A\to B$ is in \C and $\xi=(E\ra{p}A)$ is over $A$, then an opCartesian lift of $f$ at $\xi$ is given by post-composing by $f$:
$$
(E\ra{p}A)\ra{(1_E,f)}(E\ra{f\circ p}B).
$$
One thus obtains an adjunction between direct and inverse image functors relative to $f$:
$$\xymatrix{
	f_*:\C/A \ar@<.9ex>[r]\ar@{}[r]|{\perp} & \C/B:f^*.\ar@<.9ex>[l]
}$$
\item
The functor $Ob\colon\Cat\to\Set$ is a bifibration. Given a function $f\colon X\to Y$ and a category $\C$ with set of objects $X$, one construct a category $f_*\C$ this way. Define a graph $\C_Y$ with set of objects $Y$. For each pair $(x,x')$ of elements of $X$, $\C_Y(f(x),f(x'))=\C(x,x')$. For elements $y,y'$ of $Y$ not in the image of $f$, when $y=y'$, $\C_Y(y,y)=\{1_y\}$ and $\C_Y(y,y')=\emptyset$ otherwise. The category $f_*\C$ is the category generated by this graph and the following relations. Let $x\ra a x'\ra b x''$ be a composable pair in \C. The following pairs of paths in the graph are relations.
$$\xymatrix@C=1cm{
f(x) \ar[r]^a 		& f(x') \ar[r]^-b 	& f(x'')	& = 	&f(x)\ar[r]^-{b\circ a} & f(x'')\\
f(x) \ar[r]^-{1_x}	& f(x)		&		&=	& f(x)&
}$$
The functor $\underline f_{\C}\colon\C\to f_*\C$ is defined in the obvious way.
\item
The functor of modules over monoids $\Mod(\V)\to\Mon(\V)$ for a monoidal category $\V$ satisfying some mild conditions is a bifibration, as explained in part \ref{ssec:modmon}.
\item\label{ex:bifibration:presheaves}
When \A is a cocomplete category, the indexed category of presheaves with values in \A (\thref{ex:presheafA}) is actually a bi-indexed category \cite[3.7]{BorI94}. In other words, for every functor $F\colon\B\to\C$, the pre-composition with $F$ functor, $F^*$, has a left adjoint:
$$\xymatrix{
	F_*:\PSh(\B;\A) \ar@<.9ex>[r]\ar@{}[r]|{\perp} & \PSh(\C;\A):F^*.\ar@<.9ex>[l]
}$$
There is no canonical choice of the left adjoint in general. In fact, it is constructed ``by hand'', more precisely by a choice of a universal arrow $P\ra{\eta_P} F^*F_*P$ for each $P\in\PSh(\B;\A)$, the action of $F_*$ on morphisms being then defined by universality. Of course, there are in some particular cases, a canonical choice. For instance, the adjunction can always be chosen to be the identity adjunction when $F=Id$ is the identity of a site.

The Grothendieck construction $PSh(\A)\to\Cat$ is thus a bifibration. This result applies also to the pullback along $\T\colon\Top^{op}\to\Cat$ (see \thref{ex:presheafTop}). One has a bi-indexed category of presheaves on topological spaces with values in \A. The adjunction induced by a map $f\colon X\to Y$ is written
$$\xymatrix{
	f_*:\PSh(X;\A) \ar@<-.9ex>[r]\ar@{}[r]|{\perp} & \PSh(Y;\A):f^*.\ar@<-.9ex>[l]
}$$
In the notation of the bi-indexed category of presheaves (on categories), one has $f_*=(f^{-1})^*$ and $f^*=(f^{-1})_*$. There should not be confusion though, because $f$ is a map, not a functor. The bi-indexed category of presheaves on topological spaces is obtained by applying the opposite 2-functor to this adjunction, exchanging right and left adjoints (but not the direction of functors!). One then obtains a bifibration $$PSh_{\Top}(\A)\to\Top^{op}.$$

\end{Exs}

\subsection{Fibrations over a fibration}\label{ssec:FiboverFib}

In this part, we investigate the situation of a composable pair $\E\ra Q\D\ra P\B$ of two fibrations $P$ and $Q$ from different angles. The material is new, unless explicitly stated otherwise.
\subsubsection{Composable pairs of fibrations}
Fibrations can be viewed as morphisms of a category. Indeed, they are particular functors and it is therefore natural to ask whether they are closed under the composition and identity of the category \CAT. The fact that they are is a classical result of fibration theory (e.g., \cite{BorII94}).

\begin{Lem}\label{lem:compfib}
	Let $\E \ra{Q}\D \ra{P}\B$ a composable pair of fibrations. Let $f\colon A\to B$ in \B and $E\in\E_B$.\\
	A Cartesian lift of $f$ at $E$ can be obtained by the following procedure.
	Take the image $Q(E)$ of $E$ and choose a Cartesian lift $D \xrightarrow{h}QE$ of $f$ at $QE$. Then, choose a Cartesian lift $E' \xrightarrow{k}E$ of $h$ at $E$ in \E. This is Cartesian over $f$.
	$$\xymatrix{
	E'\ar[r]^k & E & \E\ar[d]^Q\\
	D\ar[r]^h & QE & \D\ar[d]^P\\
	A\ar[r]^f & B & \B
	}$$
	\cqfd
\end{Lem}

\begin{Prop}
	Objects of \CAT together with fibrations form an sub-XL-category of \CAT.\cqfd
\end{Prop}

\begin{Rem}
	Given a composable pair of functors $\E \ra{Q}\D \ra{P}\B$, we have just seen that the composite $\E \ra{P\circ Q}\B$ inherits the fibration property of $Q$ and $P$. Moreover, cleavages on $P$ and $Q$ determine a cleavage on the composite, using the preceding lemma. We call this cleavage the \emph{composite cleavage}.
\end{Rem}

\subsubsection{Fibrations internal to \texorpdfstring{$\FIB_c$}{FIBc}}
One usually sees a fibration not as a particular morphism in \CAT, but as a particular object of $\CAT^{\two}$. Consequently, we don't look at a composable pair of fibrations as some particular pair of morphisms in \CAT, but rather as a ``two floor object'': a fibration standing on another. We make this precise now. 

The preceding lemma implies that a composable pair of fibrations $$\E\ra Q\D\ra P\B$$ is the same thing as a morphism $Q$ in $\FIB(\B)$ that is itself a fibration
$$\xymatrix@C=.5cm{
\E\ar[rr]^Q\ar[dr] && \D\ar[dl]\\
&\B.&
}$$
Moreover, using the same lemma, one easily shows that a morphism of fibrations over \B that is a fibration, is automatically Cartesian (see Remark \ref{rem:CartMor}). Therefore, it is a morphism in $\FIBc(\B)$. 

There is an internal notion of a fibration in a 2-category with appropriate 2-categorical limits (modern references are \cite{Her99,Web07}, but this goes back to Street \cite{Str74a}). Internal fibrations in \CAT are precisely the Grothendieck fibrations. For a category \B, the category $\FIBc(\B)$ admits these limits. A fibration internal to $\FIBc(\B)$ is a morphism in $\FIBc(\B)$, that is a Cartesian functor over \B 
\begin{equation}\label{eq:fiboverfib}\xymatrix@C=.5cm{
\E\ar[rr]^Q\ar[dr]_R && \D\ar[dl]^P\\
&\B,&
}\end{equation}
such that all vertical arrows in \D admits all Cartesian lifts in \E. Note that all such Cartesian lifts must be also vertical. Moreover, it implies that each fibre restriction $\E_A\ra{Q_A}\D_A$ is a fibration (in \CAT). We call these fibrations \emph{fibre fibrations}. 

If $Q$ is a fibration, then it is Cartesian as already remarked, and thus it is clearly a fibration in $\FIBc(\B)$. The converse is also true \cite{Web07,Her99}, as the following lemma expresses it. This is very useful, because it simplifies the task to prove that $Q$ is a fibration.
\begin{Lem}\label{lem:intfib}
	A morphism of fibrations over \B is a fibration internal to $\FIBc(\B)$ if and only if it is a fibration.
\end{Lem}
\begin{Pf}
	Let us adopt the notations of \eqref{eq:fiboverfib}. There is a very conceptual proof in the first reference. A concrete proof is sketched in the second reference. We give it here because it tells how to construct a Cartesian morphism of $Q$ knowing the Cartesian morphisms of $P\circ Q$ and of the fibre fibrations. 
	
	Let $h\colon D'\to D$ be a morphism in \D over $f\colon A\to B$ in \B. Let $E\in\E$ be an object of \E over $D$. We construct a Cartesian lift of $h$ at $E$. $R$ is a fibration over \B and $E$ is over $B$ with respect to this fibration. Therefore, there is a Cartesian lift $\bar f_E\colon f^*E\to E$ of $f$ at $E$. $Q$ is Cartesian by hypothesis, and thus, $Q(\bar f_E)$ is Cartesian over $f$ in \D. Let $\bar h\colon D'\to Q(f^*E)$ be the morphism induced by $h$ and $Q(\bar f_E)$. This morphism is vertical in \D over $A$. By hypothesis, it admits a Cartesian lift $l\colon{\bar h}^*f^*E\to f^*E$ at $f^*E$ in \E. We prove now that this is a Cartesian lift of $h$ at $E$.
	
	It is clearly over $h$. Let us prove it is Cartesian. The following diagram should help visualizing the proof. Let $m\colon E'\to E$ be in \E and $k\colon Q(E')\to D'$ such that $h\circ k=Q(m)$. As $R(m)=f\circ P(k)$, there exists a unique lift $\overline{P(k)}\colon E'\to f^*E$ in the fibration $R$ such that $m=\bar f_E\circ\overline{P(k)}$. Now $Q(\bar f_E)\circ (\bar h \circ k)=Q(\bar f_E)\circ Q(\overline{P(k)})$. Since $Q(\bar f_E)$ is Cartesian, it implies that $\bar h \circ k=Q(\overline{P(k)})$. By Cartesianness of $l$, there is a unique lift $\bar k\colon E' \to \bar h^*f^*E$ of $k$ such that $l\circ \bar k=\overline{P(k)}$. This shows the existence of a lift of $k$ with respect to the functor $Q$. Uniqueness is easy.
	
$$\xymatrix@C=1.5cm@R=1cm{
&{\bar h}^*f^*E\ar[d]^{l}& \\
E' \ar@/_2pc/[rr]_{m}\ar@{-->}[ur]^{\exists!\bar k} \ar@{-->}[r]_{\exists!\overline{P(k)}} & f^*E \ar[r]_{\bar f_E} & E\\
&&\\
Q(E') \ar[r]^k \ar[dr]_{Q(\overline{P(k)})}\ar@/^2pc/[rr]^{Q(m)} & D'\ar[r]^h \ar@{-->}[d]_{\exists!\,\bar h} & D\\
& Q(f^* E) \ar[ur]_{Q(\bar f_E)} &\\
R(E')\ar@/_2pc/[rr]_{R(m)} \ar[r]^{P(k)} & A\ar[r]^f & B
}$$
\cqfd
\end{Pf}

\subsubsection{Fibrations over fibrations}
A composable pair $\E \ra{Q}\D \ra{P}\B$ is thus precisely a fibration
$$
(\E\ra{P\circ Q}\B)\ra{Q}(\D\ra{P}\B)
$$
internal to $\FIBc(\B)$ over the fibration $(\D\ra{P}\B)$. It explains the terminology \emph{fibration over a fibration} (introduced by Hermida in the cited article). In the following, we prefer the latter terminology to ``composable pair of fibrations'', because of this formal reason, but also because it corresponds better to our intuition, which we explain next. This observation also leads us to introduce the following terminology: $\E\ra{P\circ Q}\B$ is the \emph{domain fibration} and $\D\ra{P}\B$ is the \emph{codomain fibration} of $\E \ra{Q}\D \ra{P}\B$.

\begin{Lem}
	Let $\E \ra{Q}\D\ra{P}\B$ be a fibration over a fibration.\Par
	Cartesian morphisms of $Q$ can be decomposed as a composite of a Cartesian morphism of $P\circ Q$ followed by a Cartesian morphism of a fibre fibration. The former are called \emph{horizontal lifts}, the latter \emph{vertical lifts}.
\end{Lem}

\begin{Pf}
This lemma is a direct consequence of the proof of \thref{lem:intfib}. But we prove it a little bit differently, using the fact that $Q$ is a fibration. This gives another construction of a Cartesian arrow of $Q$. 

Let $g\colon D'\to D$ be a morphism in \D, $E\in\E_D$ and $A \ra{f}B=P(g)$. We construct a Cartesian lift $g^*E\ra{\bar g_E}E$ of $g$ at $E$ in the fibration $Q$. In the following diagram, vertical morphisms are vertical with respect to the fibrations over \B. We use strictly the notation introduced for Cartesian lifts so that one should be able to decode the role of each arrow in the diagram (but at the price of awkwardness). We encourage the reader to consider it from downstairs to upstairs.

$$\xymatrix@C=1.5cm@R=1cm{
{\bar g}^*(\bar f_D)^*E\ar[d]_{\overline{\bar g}_{(\bar f_D)^*E}}\ar@{-->}[dr]^{\bar g_E} & \\
(\bar f_D)^*E \ar[r]^-{(\overline{\bar f_D})_E} & E\\
D'\ar[r]^g \ar@{-->}[d]_{\exists!\,\bar g} & D\\
f^*D \ar[ur]_{\bar f_D} &\\
A\ar[r]_f & B
}$$	
A composite of two Cartesian arrows is Cartesian and therefore $\bar g_E$ is Cartesian over $g$. Moreover, by \thref{lem:compfib}, the Cartesian lift $(\overline{\bar f_D})_E$ of $\bar f_D$ at $E$ is a Cartesian lift of $f$ at $E$. Finally, the Cartesian lift of $\bar g$, which is in $\D_A$, is in $\E_A$.
\end{Pf}

Fibrations internal to \CAT admit a representation as pseudo-functors into \CAT. Similarly, fibrations internal to $\FIBc(\B)$ admit a representation as pseudo-functors $${\B^{op}\longrightarrow \FIBc}.$$

\begin{Prop}\label{prop:pseudoover}
	Let $\E \ra{Q}\D \ra{P}\B$ be a fibration over a fibration. Let cleavages on $P$ and $Q$ be fixed and consider the composite cleavage on $P\circ Q$.\Par
	Then, the following correspondence is a pseudo-functor:
	$$\begin{array}{rcl}
	\B^{op}& \xlongrightarrow{\Xi_Q} & \FIBc\\
	A &\longmapsto & \xymatrix{\E_A\ar[d]^{Q_A}\\ \D_A}\\
	A\ra{f}B & \longmapsto & \xymatrix{\E_B\ar[d]_{Q_B}\ar[r]^{f^*} & \E_A\ar[d]^{Q_A}\\
								\D_B \ar[r]_{f^*} & \D_A.}
	\end{array}$$	
\end{Prop}
\begin{Pf}
	We first check that this correspondence is well defined. We have already seen that the restrictions $Q_A$ are fibrations. What about its effect on morphisms of \B? The square is commutative because one has chosen the composite cleavage on $P\circ Q$. Moreover, it is a Cartesian morphism of fibrations. Indeed, let $g\colon D'\to D$ in $\D_B$ and $\bar g_E\colon g^*E\to E$ a Cartesian lift of $g$ at $E$. Its image under $f^*$ is given by:
	$$\xymatrix{
	(\bar f_{D'})^*(g^*E) \ar[r]^-{\overline{\bar f_{D'}}} \ar@{-->}[d]_{f^*(\bar g_E)}& g^*E\ar[d]^{\bar g_E}\\
	(\bar f_{D})^*E \ar[r]_{\overline{\bar f_{D}}} & E\\
	f^*D' \ar[r]^{\bar f_{D'}}\ar@{-->}[d]_{f^*g} & D'\ar[d]^g\\
	f^*D \ar[r]_{\bar f_{D}} & D\\
	A \ar[r]^f & B
	}$$
	The upper square is mapped to the lower square by the functor $Q$. In the upper square, the upper composite is Cartesian in $Q$, because it is a composite of Cartesian arrows in $Q$. Therefore, the lower composite is Cartesian in $Q$. By Proposition \ref{Prop:Cartesian}, since $\overline{\bar f_{D}}$ is Cartesian in $Q$, $f^*(\bar g_E)$ is Cartesian in $Q$, and therefore in $Q_A$.
	
	Now, let $\dom,\cod\colon \CAT^{\two}\to \CAT$ be the domain and codomain 2-functors. The composites $\dom\circ \Xi_Q=\Phi_{P\circ Q}$ and $ \cod\circ \Xi_Q=\Phi_{Q}$ are pseudo-functors by lemma \ref{lem:pseudo}. Moreover, thanks to the choice of the composite cleavage on $P\circ Q$, the structure isomorphisms $\gamma$ and $\delta$ of $\Phi_{P\circ Q}$ are over those of $\Phi_{Q}$. In conclusion, $\Xi_Q$ is a pseudo-functor, because axioms can be checked levelwise.
\end{Pf}
\begin{Def}
	A pseudo-functor $\B^{op}\to\FIBc$ is called an \emph{indexed category over an indexed category} (\emph{over \B}).
\end{Def}

\begin{Rem}
	There are 2-functors $\dom,\cod\colon\CAT^{\two}\to\CAT$ called \emph{domain} and \emph{codomain}. Thus, an indexed category over an indexed category $\Xi\colon\B^{op}\to\FIBc$ determines two indexed categories $\dom\circ\Xi$ and $\cod\circ\Xi$, which we call the \emph{domain \op resp. codomain\fp indexed category} of $\Xi$. We will see below that, in fact, $\Xi$ determines a strict morphism between these two indexed categories.
\end{Rem}

We could now define a 2-category of fibrations over fibrations over \B and prove that there is a 2-equivalence between it and the 2-category $\FIBc^{\B^{op}}$. We don't want to prove such a theorem here. Yet, we want to remark two important facts in this direction. 

First, we notice that there is no ``loss of information up to isomorphism'' in the process of going from a fibration over a fibration to its corresponding pseudo-functor. Indeed, the fibrations $P$ and $P\circ Q$ can be recovered from their indexed categories $\Phi_P$ and $\Phi_{P\circ Q}$ and these are precisely, as we have already noticed, the composite $\cod\circ \Xi_F$ and $\dom\circ \Xi_F$, respectively. What about the functor $Q$? The pseudo-functor contains its restriction $Q_A$, for all $A\in\B$. As $P\circ Q$ is a fibration, any morphism $E\ra{k}E'$ in \E over $A\ra{f}B$ in \B factors as vertical morphism of $P\circ Q$ followed by a Cartesian arrow of the cleavage of $P\circ Q$. The image of the latter is obtained by applying the restrictions $Q_A$. The image of the former is, since the cleavage of $P\circ Q$ is the composite cleavage, the Cartesian lift of $f$ at $Q(E)$ of the cleavage of $P$, which is encoded in the inverse image functor $f^*$.

A slightly more precise and complete way of stating the equivalence between the two descriptions consists of defining a Grothendieck construction $\cal G$ in this context and then showing that $Q\cong\mathcal G(\Xi_Q)$ (after having defined the category of fibrations over fibrations over \B). 

We first discuss what the 2-category of fibrations over fibrations over \B should be. The 2-category of fibrations (in \CAT) is the full sub-2-XL-category $\FIB\subset CAT^{\two}$ of fibrations. In the same manner, the 2-category of fibrations over fibrations over \B is the full sub-2-XL-category $\FIB(\FIBc(\B))\subset \FIBc(\B)^\two$ of fibrations. Therefore, it consists of the following data.
\begin{itemize}
	\item
Objects: Commutative triangles in \CAT
$$\xymatrix@C=.5cm{
\E\ar[rr]^Q\ar[dr]_R && \D\ar[dl]^P\\
&\B.&
}$$
where all functors are fibrations (in \CAT). As before, one also just writes $$\E\ra{Q}\D\ra{P}\B.$$
\item
Morphisms: A morphism from $\E\ra{Q}\D\ra{P}\B$ to $\E'\ra{Q'}\D'\ra{P'}\B$ is a tetrahedron that stands on its summit \B, which is the codomain of the edges related to it, and whose base is a commutative square
$$\xymatrix{
\E \ar[r]^Q \ar[d]_{F} & \D \ar[d]^{G}\\
\E' \ar[r]_{Q'}& \D'.
}$$ 
whose vertical edges are Cartesian functors. Thus, $F$ and $G$ are Cartesian functors over \B. One also writes $(F,G)\colon(\E\ra{Q}\D\ra{P}\B)\to(\E'\ra{Q'}\D'\ra{P'}\B)$.
\item
2-cells: A 2-cell 
$$\xymatrix@=2cm{
**[l](\E\ra{Q}\D\ra{P}\B)\rtwocell^{(F,G)}_{(F',G')}{\quad(\alpha,\beta)} &**[r](\E'\ra{Q'}\D'\ra{P'}\B)
}$$
is a pair of natural transformations \emph{over} \B making the following square commute:
$$\xymatrix@=1.5cm{
{\E \ar[r]^Q} \dtwocell^{F'}_{F}{^\alpha} & \D \dtwocell^{G'}_{G}{^\beta}\\
\E' \ar[r]_{Q'}& \D'
}$$
\end{itemize}
We now describe a Grothendieck construction in this context, which we call \emph{2-level Grothendieck construction}\indexb{Grothendieck!2-level -- construction}. Let $\Xi\colon\B^{op}\to \FIBc$ be an indexed category over an indexed category. One can show that it determines a strict transformation of pseudo-functors:
$$\xymatrix@=1.5cm{
\B^{op} \rtwocell^{\dom\circ \Xi}_{\cod\circ\Xi}{\ \ \alpha_{{}_\Xi}}& **[r]\CAT
}$$
with $(\alpha_{{}_\Xi})_A=\Xi(A)$. The Grothendieck construction 2-functor for indexed categories thus determines a Cartesian functor over \B:
$$\xymatrix@C=.5cm@R=1cm{
\E_{\Xi} \ar[dr]_{R_{\Xi}}\ar[rr]^{Q_{\Xi}} && \D_{\Xi}\ar[dl]^{P_{\Xi}}\\
&\B.&
}$$
One can then show, using the fact that each $\Xi(A)$ is a fibration, $A\in\B$, and each $\Xi(f)$ a Cartesian morphism, $f$ in \B, that $Q_{\Xi}$ is itself a fibration. Note that the definition of the Cartesian functor $Q_{\Xi}$ only requires a pseudo-functor $\Xi\colon\B^{op}\to\CAT^{\two}$. However, it must have precisely $\FIBc$ as codomain for $Q_{\Xi}$ to be a fibration.

It remains to show that, given a fibration over a fibration, the Grothendieck construction of its pseudo-functor is isomorphic to the former. One already knows that a fibration over \B is isomorphic over \B to the Grothendieck construction of its indexed category. Moreover, any isomorphic functor is Cartesian and it is easy to verify that the following diagram commutes.
$$\xymatrix{
\E \ar[r]^Q \ar[d]_{\cong} & \D \ar[d]^{\cong}\\
\E_{\Xi} \ar[r]_{Q_{\Xi}}& \D_{\Xi}
}$$ 
\begin{Prop}
	To each indexed category over an indexed category $\B^{op}\ra{\Xi}\FIBc$, there corresponds a fibration over a fibration ${\E_{\Xi}\ra{Q_{\Xi}}\D_{\Xi}\ra{P_{\Xi}}\B}$ over \B and if $\B^{op}\ra{\Xi_Q}\FIBc$ is the pseudo-functor corresponding to $\E\ra{Q}\D\ra{P}\B$ under the construction \ref{prop:pseudoover}, then its corresponding fibration over a fibration over \B is isomorphic to $\E\ra{Q}\D\ra{P}\B$.\cqfd
\end{Prop}

Remark that, for every $\B^{op}\ra{\Xi}\FIBc$, there is an isomorphism of fibrations between the fibre fibration of ${\E_{\Xi}\ra{Q_{\Xi}}\D_{\Xi}}$ at $A\in\B$ and the fibration $\Xi(A)=\E_A\ra{Q_A}\D_A$:
$$\xymatrix@=1cm{
\E_A \ar[r]^{\cong} \ar[d]_{Q_A} & (\E_{\Xi})_A\ar[d]^{(Q_{\Xi})_A}\\
\D_A \ar[r]_{\cong}& (\D_{\Xi})_A.
}$$

\subsubsection{Opfibred and bifibred situations}
We now briefly consider the opfibred and bifibred cases. 

The opfibred setting is an obvious dualization of the fibred case. Let us just mention that given an opfibration over an opfibration $\E\ra{Q}\D\ra{P}\B$, one obtains a pseudo-functor $\B\to\OPFIB_{oc}$.

The bifibred case requires more care. We explore only the process of associating a pseudo-functor to a bifibration over a bifibration and the opposite process of 2-level Grothendieck construction of a bi-indexed category over a bi-indexed category, concept that we still have to define.

Let $\E\ra{Q}\D\ra{P}\B$ be a bifibration over a bifibration. Let us choose a bicleavage of $P$ and of $Q$ and consider the composite bicleavage of $P\circ Q$. Then each restriction $Q_B\colon\E_B\to\D_B$ over $B\in\B$ is a bifibration. Moreover, an arrow $f\colon A\to B$ in \B induces the following structure:
\begin{equation}\label{eq:adjbifibbifib}\begin{aligned}\xymatrix@=1.5cm{
\E_A\ar[d]_{Q_A}\ar@<-.9ex>[r]_{f_*}\ar@{}[r]|{\top} & \E_B\ar[d]^{Q_B}\ar@<-.9ex>[l]_{f^*}\\
\D_A\ar@<-.9ex>[r]_{f_*}\ar@{}[r]|{\top} & \D_B\ar@<-.9ex>[l]_{f^*}
}\end{aligned}\end{equation}
The morphism of bifibrations $(f_*,f_*)$ is opCartesian, whereas $(f^*,f^*)$ is Cartesian. Moreover, this square is an adjoint square, meaning that the pair $(\iota,\iota)$ of identity natural transformations is a mate pair. This can be readily verified using the characterization of adjoint morphisms under an adjunction $f_*\dashv f^*$ of direct and inverse image functors given in \thref{prop:bifibadjoint}, and the fact that $P\circ Q$ is endowed with the composite bicleavage. The condition that the identity pair be mates has a simple expression. Let $\eta^f_\D$, $\epsilon^f_\D$ and $\eta^f_\E$, $\epsilon^f_\E$ be the respective units and counits of the direct and inverse image along $f$ adjunction in \D and \E. Then, the condition that the identity pair $(\iota,\iota)$ be mates is equivalent to $\eta^f_\D\cdot Q_A=Q_A\cdot\eta^f_\E$, or to the condition $\epsilon^f_\D\cdot Q_B=Q_B\cdot\epsilon^f_\E$. This is can be verified from the axiom of mate pairs given, but it is more easily obtained from the diagrammatic form in \cite{KS74}. Note that this is precisely the condition that the square \hyperref[eq:adjbifibbifib]{(\ref*{eq:adjbifibbifib})} be an adjunction in $\CAT^{\two}$. 

There is a 2-XL-category whose objects are bifibrations and morphisms adjoint commutative squares as just described. Two-cells are pairs of mate pairs, one mate pair $(\alpha_*,\alpha^*)$ for the upper adjunction and one mate pair $(\alpha'_*,\alpha'^*)$ for the lower adjunction, such that $\alpha_*$ is over $\alpha'^*$ and similarly for their mates. Let us write it $\BIFIBADJ$\index[not]{BIFIBADJ@\BIFIBADJ}. Bifibrations over bifbrations over \B determine thus pseudo-functors $\B\to\BIFIBADJ$.

Conversely, let us start with such a pseudo-functor $\B\to\BIFIBADJ$. By forgetting the right adjoints, one obtains an opindexed category over an opindexed category and the 2-level Grothendieck op-construction provides us with an opfibration over an opfibration $\E\ra{Q}\D\ra{P}\B$. We verify that this indeed is  a bifibration over a bifibration. 

By Proposition \ref{prop:bifibadjoint}, the opfibrations $P$ and $P\circ Q$ are bifibrations (because their direct image functors admit a right adjoint). We moreover know that $Q$ is an opfibration that is fibrewise a bifibration. We still have to prove it is a bifibration. By Lemma \ref{lem:intfib}, one can decompose the proof in two steps. First, show that $Q$ is a Cartesian functor over \B from $P\circ Q$ to $P$. Then, show that each vertical morphism in \D admits all Cartesian lifts in \E.

Firstly, $Q$ is Cartesian, because a Cartesian arrow in $\E$ over $f\colon A\to B$ is given by $$(f,\epsilon^f_\E)\colon(A,f^*E)\to (B,E)$$ (Proposition \ref{prop:bifibadjoint}) and because the squares \hyperref[eq:adjbifibbifib]{(\ref*{eq:adjbifibbifib})} are supposed to be adjoint squares. In order to prove the second statement, one uses of course that the fibres $\E_A\to\D_A$ are bifibrations. Note there is something to prove though, because one must show that the Cartesian lifts in $\E_A$ are also Cartesian in the whole category $\E$.
\paragraph{}
We therefore have the following definition.

\begin{Def}\label{def:binoverbind}
	\emph{Bi-indexed categories over bi-indexed categories} over \B are pseudo-functors $\B\to\BIFIBADJ$.
\end{Def}

\section{Grothendieck topologies}

In this part, we follow Johnstone \cite{Joh202} and Vistoli \cite{Vis08}. See also \cite{MM92}. We have a more minimal approach to axioms though (we use a weakened form and try to use the least possible), and we give some more details about the relationship between the sifted and unsifted structures. In particular, we divide the generation of a Grothendieck topology from a pretopology into two steps, one consisting in generating the sifted version of the pretopology as described in \cite{Joh202}, the other in saturating the pretopoloy as described in \cite{Vis08}. We study these questions in a more general context of \emph{covering functions}. If a result appears without a proof, then it either means that the proof is easy or that the proof can been found in these references. 

\subsection{Coverings and covering functions}
One first defines the very general notions of covering and refinement of a covering. The example one might have in mind is that of the category $Top$ of topological spaces. A notion of covering for a space $X$ is given by a family $\U=\{U_i\hookrightarrow X\mid i\in I\}$ of inclusions of open subsets of $X$ whose union is $X$. One can compare coverings of this type of a given topological space $X$ by the notion of refinement. A covering $\V$ of $X$ refines a covering $\U$ of $X$ if each open set $V\in\V$ is included in an open set $U\in\U$. This is equivalent to saying that each inclusion $V\hookrightarrow X$ factors through an inclusion $U\hookrightarrow X$ in $Top$. \begin{Defs}
[Let \C be a category]
\item
A \emph{covering}\indexb{Covering} of an object $C\in\C$ is a class of arrows in \C with codomain $C$. It is called a \emph{set-covering}\indexb{Covering!Set-{}--} if $R$ is a set%
\footnote{Authors often define a covering in a category \C not as a mere family but as an \emph{indexed} family $$(C_i\to C)_{i\in I},$$ i.e., a \emph{function} $I\to\Mor\C$. We do not follow this convention, because we think it quite complicates the matter to differentiate coverings that have the same class of arrows but different indexing. Defining coverings as classes of arrows also allows us to consider sieves as particular coverings. 

Now, it is sometimes practical, if only for notational reasons, to consider indexed families, and we do so from time to time. In fact, a covering $R$ has always an canonical indexing: the identity function $1_R$ and thus, one can assume without loss of generality that a covering comes with an index class. Yet, the reader should remember that, ultimately, the covering will always be the range of the index function, that is the class $$\{C_i\to C\mid i\in I\}.$$ In this work, this notation always means the class corresponding to the indexed class $(C_i\to C)_{i\in I}$.}.
\item
A covering $R$  of an object $C$ \emph{refines}\indexb{Refinement|see{Covering}} a covering $S$ of $C$ if each arrow of $R$ factors through an arrow of $S$, i.e., for each $f\colon D\to C$ of $R$, there exists an arrow $g\colon E\to C$ of $S$ and a factorization:
$$\xymatrix{
& E\ar[d]^{g}\\
D\ar[r]_{f}\ar@{-->}[ur]^{\exists} & C.
}$$
One says that $R$ is a \emph{refinement}\indexb{Covering!Refinement of a --} of $S$.
\end{Defs}
The refinement relation determines a preorder structure on the conglomerate of all coverings of an object $C\in\C$. 

On the whole category $Top$, one has a function $K$ that assigns to each topological space the set of all its open subset coverings. One can consider another type of covering of a space $X$: sets $\{Y_i\xhookrightarrow{f_i}X\mid i\in I\}$ of open embeddings into $X$ that are \emph{collectively surjective} in the sense that the induced map $\coprod_{i\in I} Y_i\to X$ is surjective. This defines another function $K'$ on the class of objects of $\Top$. Of course, one has $K\subset K'$. Thus, any covering $R$ of $K$ admits a refinement in $K'$, $R$ itself. One says that $K$ is \emph{subordinated} to $K'$. On the other hand, $K'$ is subordinated to $K$. Indeed, if $R'=\{Y_i\xhookrightarrow{f_i}X\mid i\in I\}\in K'(X)$, then consider the factorization
$$\xymatrix{
&f_i(Y_i)\ar@{^{(}->}[d]\\
Y_i\ar@{^{(}->}[r]_{f_i}\ar@{-->}[ru]_{\cong}^{\bar f_i} & X
}$$
The family $R=\{f_i(Y_i)\hookrightarrow X\mid i\in I\}$ belongs to $K(X)$. It refines $R'$ because, for each $i\in I$, one has the factorization
$$\xymatrix{
&Y_i\ar@{^{(}->}[d]^{f_i}\\
f_i(Y_i)\ar@{^{(}->}[r]\ar@{-->}[ru]_{\cong}^{{\bar f_i}^{-1}} & X
}$$
In this situation, one says that $K$ and $K'$ are equivalent. We define now these concepts in a general context.
\begin{Defs}\label{def:site}
\item\label{def:site:covfunc}
A \emph{covering function on \C}\indexb{Covering function}\index[not]{K@$K$} is a function $$K\colon\Ob\C\to\power(\power(Mor\C))$$ that assigns to each object $C\in\C$ a conglomerate $K(C)$ of coverings of $C$%
\footnote{$Mor\C$ is a class, and in $\NBG$, the power class $\mathcal{P}(A)$\index{Class!Power --} of a class $A$ is well-defined. It is the class of all \emph{subsets} of $A$. What we denote $\power Mor\C$ here is not the power class, it is the \emph{power conglomerate}\index{Conglomerate!Power --} of $Mor\C$, i.e., the conglomerate of all subconglomerates of $Mor\C$ (so it is the usual ``power class'', but in the ambient \NBI). Note that any subconglomerate of a class is a class, since classes are given by all subconglomerates of the universe of sets.}.
Coverings of $K$ are sometimes called \emph{$K$-coverings}. 
\item
A category equipped with a covering function on it is called a \emph{site}.\indexb{Site}

\item\index{Subordinated|see{Covering function}}
A covering function $K'$ on \C is \emph{subordinated}\indexb{Covering function!Subordinated -- } to a covering function $K$, and we write $K'\preceq K$, if each covering of $K'$ admits a refinement belonging to $K$. One says that $K$ and $K'$ are \emph{equivalent}\indexb{Covering function!Equivalent --}, and write $K\equiv K'$\index[not]{Equivalent@$\equiv$}, if $K\preceq K'$\index[not]{Subordinated@$\preceq$} and $K'\preceq K$.
\end{Defs}
Being subordinate gives a preorder structure on the conglomerate of covering functions on \C, and being equivalent is thus an equivalence relation.

\paragraph{}
One can also define the dual notion of a covering.

\begin{Def}
	A \emph{opcovering}\indexb{Opcovering (and related concepts)} $R$ of an object $C\in\C$ is a covering of $C$ in $\C^{op}$.
\end{Def}

One obtains similarly all the dual notions of \emph{refinement of a opcovering}, \emph{opcovering functions} and \emph{subordination of opcovering functions}. In the sequel, we do not formulate the dual definitions and theorems.

\subsection{Axioms on covering functions}
We now consider three different possible axioms on a covering function $K$. The reason for the choice of the letter that names them will become clearer later. Let \C be a category.

\begin{Axns}
\begin{enumerate}
\item[\normalfont(M)]\index[not]{AxiomM@(M)}
For all objects $C\in\C$, the \emph{identity covering} $\{1_C\colon C\to C\}$ refines some covering of $K$.
\item[\normalfont(C)]\index[not]{AxiomC@(C)}
For each $C\in\C$, each $K$-covering $R$ of $C$ and each arrow $g\colon D\to C$ in \C, there exists a $K$-covering $S$ of $D$ such that the \emph{composite covering} $$g\circ S\coloneqq\{g\circ h\mid h\in S\}$$ of $C$ refines $R$.
\item[\normalfont(L)]\index[not]{AxiomL@(L)}
Given a $K$-covering $R$ and, for each $f\in R$, a $K$-covering $R_f$ of $\dom f$, there exists a $K$-covering that refines the \emph{composite covering}\indexb{Covering!Composite --}
$$
\bigcup_{f\in R} f\circ R_f=\{f\circ g\mid f\in R\text{ and }g\in R_f\}.
$$
\end{enumerate}
\end{Axns}
\begin{Rem}
	Axiom (M) is equivalent to axiom
	\begin{itshape}
	\begin{enumerate}
\item[\normalfont(M')]
For each $C\in\C$, there is a $K$-covering $R$ of $C$ that contains a split epi.
\end{enumerate}
\end{itshape}
\end{Rem}

\begin{Defs}
	[Let \C be a category]
	\item After Johnstone, a covering function $K$ on \C that satisfies axiom (C) is called a \emph{coverage}\indexb{Coverage}. This is where the name of axiom (C) comes from. A \emph{coverage-site}\indexb{Site!Coverage{-}--} is a site whose covering function is a coverage
	\item A covering function $K$ on \C that satisfies all axioms (M), (C) and (L) is called a \emph{\op Grothendieck\fp pretopology}\indexb{Pretopology}\index{Covering function!Grothendieck pretopology|see{Pretopology}}. A \emph{Grothendieck site}\indexb{Site!Grothendieck --}\index{Covering function!Grothendieck site} is a site whose covering function is a pretopology.
\end{Defs}

The following result does not appear in the literature, as far as we know. Indeed, Vistoli seems to be the only one to introduce this notion of equivalence and he restricts his attention to Grothendieck pretopologies.
\begin{Prop}\label{prop:equivcovfun}
	Let $K$ and $K'$ be covering functions. Then,
	\begin{enumerate}
\item
If $K$ and $K'$ are equivalent, then $K$ satisfies axiom (C) \ssi $K'$ does.
\item
If $K$ and $K'$ are equivalent coverages, then $K$ satisfies axiom (L) \ssi $K'$ does. 
\end{enumerate}
\end{Prop}
\begin{Pf}
	The first statement follows quite directly from the hypotheses. Let us prove the second assertion. The situation being perfectly symmetric, it is enough to prove only one direction. Suppose $K$ is a coverage that satisfies (L). Consider a $K'$-covering $R'$ and for each $f'\in R'$, a $K'$-covering $R'_{f'}$ of $\dom f'$. 
	
	Since $K'\preceq K$, there exist $K$-coverings $R$ that refines $R'$ and $R_{f'}$ that refines $R'_{f'}$, for all $f'\in R'$. Let $f\colon D\to C$ be an element of $R$. Then, there exists $f'\colon D'\to C\in R'$ and $g\colon D\to D'$ such that $f=f'\circ g$. Now, by axiom (C), there exists a $K$-covering $S_f$ of $D$ such that the covering $f\circ S$ refines $R_{f'}$. One can apply axiom (L) for $K$ to the composite covering $\bigcup_{f\in R} f\circ S_f$, to obtain a $K$-covering $T$ of $C$ that refines it. Since $\bigcup_{f\in R} f\circ S_f$ refines $\bigcup_{f'\in R'} f'\circ R'_{f'}$, $T$ refines $\bigcup_{f'\in R'} f'\circ R'_{f'}$. One uses the fact that $K\preceq K'$ to conclude.
\end{Pf}

Here is a result that shows the importance of axioms (C) and (L).

\begin{Lem}\label{lem:comref}
	Let $K$ be a coverage on a category \C satisfying axiom (L). Then any two $K$-coverings of the same object admit a common refinement by a $K$-covering.
\end{Lem}
\begin{Pf}
	Let $R$ and $S$ be $K$-coverings of an object $C$ of \C. By axiom (C), there exists a $K$-covering $S_f$ for every $f\in S$ such that the composite covering $f\circ S_f$ refines $R$. Now, by axiom (L), there exists a $K$-covering $T$ that refines the composite covering $\bigcup_{f\in S}f\circ S_f$. Yet, the latter refines both $R$ and $S$, and therefore, $T$ refines also $R$ and $S$. 
\end{Pf}

The general notion of a covering function on a category was first introduced for the sake of defining \emph{sheaves} on the category.

\begin{Def}
	Let \C be a category. A presheaf $P$ on \C, i.e., a functor $P\colon \C^{op}\to \Set$\indexb{Sheaf!Presheaf}\index{Presheaf|see{Sheaf}}, satisfies the \emph{sheaf axiom}\indexb{Sheaf!Sheaf axiom} (or \emph{has the sheaf property}\indexb{Sheaf!Sheaf property}) for a covering $R$ if the following property holds.
	\begin{enumerate}
\item[]
If an indexed family $(x_f)_{f\in R}$ of elements $x_f\in P(\dom f)$ is \emph{compatible} in the sense that, 
whenever $h\colon D\to \dom f$ and $k\colon D\to \dom g$ satisfy $f\circ h=g\circ k$, one has $$P(h)(x_f)=P(k)(x_g),$$ 
then there exists a unique $x\in P(C)$ such that $P(f)(x)=x_f$, for all $f\in R$.
\end{enumerate}
If $K$ is a covering function on \C and if $P$ satisfies the sheaf axiom for each $K$-covering, then one says that $P$ is a \emph{$K$-sheaf on \C}\indexb{Sheaf!K-sheaf on \C@$K$-sheaf on \C}\indexb{Sheaf}, or simply a \emph{sheaf} if the covering function and the category it is defined on are clear from the context.
\end{Def}

Covering functions contain too much information with respect to the sheaf axiom. Indeed, there are different covering functions that have the same sheaves. One first considers the sheaf property with respect to coverings and then to covering functions. 

\subsection{Sifted and saturated covering functions}
Given an object $C\in\C$, two coverings of $C$ have the same sheaves if they generate the same \emph{sieves}. To get an intuition for this notion, consider the category $\mathcal{O}(X)$ of open subsets of a topological space $X$ whose arrows are inclusions. Consider on this category the covering function $K$ of open subset coverings. Given a covering $\{U_i\hookrightarrow V\}$ of an open subset $V\subset X$, it generates a particular $K$-covering, the covering $$\{W\hookrightarrow V\text{ in }\mathcal{O}(X)\mid \exists\ i\in I \text{ with } W\subset U_i\}.$$ This covering is closed under pre-composition. Moreover, it resembles a ``material sieve'' because every neighbourhood that is ``smaller'' than one of its holes $U_i$ goes through it, but not otherwise. So, in this particular case, it happens that given a covering $\U$ of an open subset $V$, it changes nothing as far as sheaves are concerned to add subsets $W$ of the elements of \U to the covering \U.

\begin{Def}
Let \C be a category. A \emph{sieve}\indexb{Sieve} on an object $C$ of \C is a covering $R$ of $C$ that is a right ideal in \C, i.e., such that if $f\in R$, then, for any right composable arrow $g$ in \C, $f\circ g\in R$.
\end{Def}
 
 Given a covering $R$ of an object $C$, there exists a smallest sieve (with respect to inclusion) that contains $R$, called the \emph{sieve generated by $R$}\indexb{Sieve!-- generated by a covering} and denoted $\overline R$\index[not]{Rbar@$\overline R$}. It is given by
 $$
 \overline R=\{f\circ g\mid f\in R,\ g\text{ composable with }f\}.
 $$
 Therefore, each covering function $K$ on \C generates a covering function $\overline{K}$ whose coverings are the sieves generated by the coverings of $K$. A covering function all of whose coverings are sieves is called \emph{sifted}\indexb{Covering function!Sifted --}. In particular, one considers sifted coverages.
 
 \begin{Term}
Sieves are particular coverings, and we apply to them the same terminology. In particular, given a covering function $K$ on a category \C, a \emph{$K$-sieve} is a sieve $S$ that belongs to $K$. We warn the reader that Vistoli uses this terminology with a different meaning: that there is a covering of $K$ that refines $S$.
\end{Term}
 
\begin{Lem}
	Let $R$ be a covering of a category \C and $P$ a presheaf on \C.\\
	Then, $P$ has the sheaf property with respect to $R$ \ssi it has it with respect to the sieve $\overline R$ generated by $R$.\\
	In particular, if $K$ is a covering function, then $K$ and $\overline{K}$ have the same sheaves.\cqfd
\end{Lem}

So, the sheaf property does not distinguish a covering and the sieve it generates. Note that when passing from $K$ to $\overline{K}$, one does not add coverings. One only adds arrows to existing coverings. On the contrary, one might have lost coverings, as different coverings can generate the same sieve. Remark that given a sifted covering function $K$, there is a maximal covering function $K_{max}$ such that $\overline{K_{max}}=K$. It is defined by
$$
R\in K_{max}\iff\overline R\in K.
$$ 
\begin{Def}
	The covering function $K_{max}$\index[not]{Kmax@$K_{max}$} defined above is called the \emph{maximal generating covering function}\indexb{Covering function!Maximal generating -- of a sifted one} of the sifted covering function $K$.
\end{Def}
We will see later in examples that coverings that appear naturally are not sieves. Nevertheless, sieves are important because they give a more conceptual definition of a covering by means of subobjects of representable functors (see \cite{Joh202,MM92}). They also simplify the theory, as we will observe throughout the text.

We first take a look at the relationship between a covering $R$ and its corresponding sieve. Consider a category \C and an object $C$ of \C. The sieve generated by the identity covering $\{1_C\colon C\to C\}$ is the \emph{maximal sieve}\indexb{Sieve!Maximal --} $M_C$ on $C$, i.e., the sieve of all arrows in \C with codomain $C$. This explains the name of axiom (M) (see also Lemma \ref{lem:sifted}).

If $R$ and $R'$ are coverings, then $R\subset R'$ implies that $R$ refines $R'$. The converse is wrong in general, but true for sieves. By definition, a covering $R$ refines a covering $R'$ \ssi $R\subset \overline{R'}$. In fact, one has the following.
\begin{Lem}
	A covering $R$ refines a covering $R'$ \ssi $\overline{R}\subset \overline{R'}$. In particular, if $S$ and $S'$ are sieves, then $S$ refines $S'$ \ssi $S\subset S'$.\cqfd
\end{Lem}

We consider now some relationships between $K$ and its associated sifted covering function $\overline{K}$.

\begin{Not}
	Let $S$ be a sieve on an object $C$ in a category \C and $f\colon D\to C$ a morphism in \C. Then, one defines the covering
	$$
	f^*(S)\coloneqq\{g\text{ in }\C\mid \cod g=D \text{ and } f\circ g\in S\},
	$$
	which is clearly a sieve.
\end{Not}
We first give equivalent forms of the axioms (M), (C) and (L) in the sifted setting (the lowercase ``s'' added to an axiom's name is for ``sifted'').
\begin{Lem}\label{lem:sifted}
Let \C be a category and $K$ a \emph{sifted} covering function on \C.	
\begin{enumerate}
\item
Axiom (M) is equivalent to the axiom
\emph{\begin{enumerate}
\item[(M\subscript s)]\index[not]{AxiomMs@(M\subscript s)}
For each $C\in\C$, the maximal sieve $M_C$ belongs to $K(C)$.
\end{enumerate}}
\item
Axiom (C) is equivalent to the axiom
\emph{\begin{enumerate}
\item[(C\subscript s)]\index[not]{AxiomCs@(C\subscript s)}
If $R$ is a $K$-sieve on $C\in\C$ and $g$ a morphism in \C with codomain $C$, then $g^*(R)$ has a refinement in $K$, i.e., contains a $K$-sieve. 
\end{enumerate}}
\item
Axiom (L) is equivalent to the axiom
\emph{\begin{enumerate}
\item[(L\subscript s)]\index[not]{AxiomLs@(L\subscript s)}
Let $C\in \C$. If $S\in K(C)$ and $R$ is any sieve on $C$ such that $f^*(R)$ is a $K$-sieve for all morphisms $f\in S$, then $R$ contains a $K$-sieve.
\end{enumerate}}
\end{enumerate}
\cqfd
\end{Lem}

\begin{Def}
	A covering function (resp.\ a sifted covering function) $K$ on \C is \emph{saturated}\indexb{Covering function!Saturated --} (resp.\ sieve-saturated)\indexb{Covering function!Sieve-saturated sifted --} if any covering having a refinement in $K$ is in $K$ (resp.\ any sieve containing a $K$-sieve is a $K$-sieve).
\end{Def}
For saturated covering functions, the preorder relation of being subordinate coincides with the inclusion relation.

\begin{Lem}
	Let $K$ be a saturated (resp.\ sieve-saturated sifted) covering function. Then, for every covering function (resp.\ sifted covering function) $K'$,
	$$
	K'\preceq K\iff K'\subset K.
	$$
	\cqfd
\end{Lem} 

If $K$ is a sieve-saturated sifted covering function, then it is clear that axioms (C\subscript s) and (L\subscript s) are equivalent to the following axioms:
\emph{\begin{enumerate}
\item[(C\subscript s')]\index[not]{AxiomCs'@(C\subscript s')}
If $R$ is a $K$-sieve of $C\in\C$ and $g$ a morphism in \C with codomain $C$, then $g^*(R)$ is a $K$-sieve. 
\end{enumerate}}

\emph{\begin{enumerate}
\item[(L\subscript s')]\index[not]{AxiomLs'@(L\subscript s')}
Let $C\in \C$. If $S\in K(C)$ and $R$ is any sieve on $C$ such that $f^*(R)$ is a $K$-sieve for all morphisms $f\in S$, then $R$ is a $K$-sieve.
\end{enumerate}}

Axiom (L\subscript s') gives an intuition for the axiom (L)'s name. The letter (L) refers to the \emph{local character condition}. In its sifted saturated version (L\subscript s'), it expresses the fact that if a sieve is ``locally a $K$-sieve'' over a $K$-sieve, then it must be a $K$-sieve itself.

\begin{Lem}\label{lem:topsat}
	 A sifted covering function on a category \C that satisfies (M\subscript s) and (L\subscript s') is sieve-saturated.
	 
	 \cqfd
\end{Lem}

\begin{Def}
A \emph{Grothendieck topology}\indexb{Grothendieck!-- topology|see{Covering function}}\indexb{Covering function!Grothendieck topology} on a category \C is a sifted covering function $J$ on \C satisfying one of the following equivalent conditions:
\begin{enumerate}[(i)]
\item
It satisfies axioms (M\subscript s), (C\subscript s') and (L\subscript s').
\item
It is a coverage that satisfies axioms (M\subscript s) and (L\subscript s').
\item
It is sieve-saturated and satisfies (M\subscript s), (C\subscript s) and (L\subscript s).
\end{enumerate}
\end{Def}

Now we turn to comparison between $K$ and $\overline{K}$.
\begin{Prop}\label{prop:sifted}
Let \C be a category and $K$ a covering function on \C.
\begin{enumerate}
\item
$\overline{K}$ is equivalent to $K$.
\item
$\overline K$ is sifted. Moreover, $K$ is sifted \ssi $K=\overline K$.
\item
$K$ satisfies axiom (M) or (C) \ssi $\overline K$ does, and therefore, \ssi $\overline K$ satisfies (M\subscript s) and (C\subscript s).
\item
If $K$ satisfies (L), then $\overline K$ also, and therefore, $\overline K$ satisfies (L\subscript s). Moreover, $K$ satisfies (C) and (L) \ssi $\overline K$ does, and therefore, \ssi $\overline K$ satisfies (C\subscript s) and (L\subscript s).
\item
If $K$ is saturated, then $\overline{K}$ is sieve-saturated. Conversely, if $\overline{K}$ is sieve-saturated, then $K_{max}$ is saturated.\label{prop:sifted5}
\item
$K\preceq K'\iff\overline{K}\preceq\overline{K'}$. 
\end{enumerate}	
\end{Prop}

\begin{Pf}
These results are mostly easy to obtain or follow directly from Proposition \ref{prop:equivcovfun}. The fact that if $K$ satisfies (L), then so does $\overline K$ is, on the other hand, not a consequence of the latter proposition. 

Let $\overline R$\ be the $\overline K$-sieve generated by a $K$-covering $R$ of $C\in\C$. For each $f\in\overline R$, let $\overline R_f$ be a $\overline K$-covering of $\dom f$, generated by a $K$-covering $R_f$. One considers the composite sieve $\bigcup_{f\in\overline R}f\circ\overline R_f$ and show that it admits a refinement in $\overline K$. Let us write $R=\{C_i\ra{f_i}C\mid i\in I\}$ and $R_{f_i}=\{C_{ij}\ra{g_{ij}}C_i\mid j\in J_i\}$. Since $K$ satisfies (L), the composite covering $\bigcup_{i\in I}f_i\circ R_{f_i}$ admits a refinement $S$ in $K$. Let $\overline{S}$ be its generated sieve. Then, $\overline S$ refines $S$, which refines $\bigcup_{i\in I}f_i\circ R_{f_i}$, which is included in $\bigcup_{f\in\overline R}f\circ\overline R_f$. Therefore $\bigcup_{f\in\overline R}f\circ\overline R_f$ is refined by $\overline S$, which is in $\overline K$.

\cqfd
\end{Pf}
\subsection{Sifted saturation of a covering function}
The sheaf property happens also to be independent under addition of some sort of coverings.

\begin{Lem}
	Let $C$ be an object of a category \C. Then any presheaf on \C has the sheaf property with respect to the identity covering $\{C\ra{1_C}C\}$ of $C$.\cqfd
\end{Lem}

\begin{Lem}\label{lem:sheafrefinement}
	Let $F$ be a sheaf for a coverage $K$ on a category \C. Then $F$ has the sheaf property for every covering in \C that has a refinement in $K$.\\
	In particular, if $K'$ is any covering function on \C subordinated to $K$, then every $K$-sheaf is a $K'$-sheaf.\\
	Finally, equivalent coverages have the same sheaves.\cqfd
\end{Lem}

Thus, given a coverage $K$, as far as sheaves are concerned, it changes nothing to add to $K$ coverings that have refinement in $K$ (and in particular identity coverings). 

Each equivalence class of covering functions has a largest element when seen as a partially ordered conglomerate under inclusion, as we now verify. 

\begin{Def}\label{def:saturation}
	Let \C be a category and $K$ a covering function (resp. a sifted covering function) on \C. The \emph{saturation}\indexb{Covering function!Saturation of a --} (resp. \emph{sifted saturation})\indexb{Covering function!Sifted saturation of a sifted --} of $K$, written $K_{Sat}$\index[not]{KSat@$K_{Sat}$}, is the covering function whose coverings are all coverings having a refinement in $K$ (resp. the sifted covering function whose sieves are all sieves that contain a $K$-sieve). 
\end{Def} 
The following proposition provides the important properties of $K_{Sat}$ (we leave to the reader the task of expressing the corresponding properties in the sifted case).
\begin{Prop}\label{prop:saturation}
	Let $K$ be a covering function, 
	\begin{enumerate}
\item
$K\subset K_{Sat}$ and is equivalent to $K_{Sat}$. 
\item
$K_{Sat}$ satisfies axiom (M) if $K\neq\emptyset$. Moreover, $K_{Sat}$ satisfies axiom (C) \ssi $K$ does. If $K_{Sat}$ satisfies (L), then so does $K$. Finally, $K$ is a coverage satisfying (L) \ssi $K_{Sat}$ is such a coverage.
\item
$K_{Sat}$ is saturated. Moreover, $K$ is saturated \ssi $K=K_{Sat}$.
\item
For each covering function $K'$, $K'\preceq K\iff K'\subset K_{Sat}\iff K'_{Sat}\subset K_{Sat}$,
\item
For each covering function $K'$, $K\equiv K'\iff K_{Sat}=K'_{Sat}$.
\end{enumerate}
\end{Prop}
\begin{Pf}
	Most of the results follow directly from definitions or from Proposition \ref{prop:equivcovfun}. Yet, there still remains to prove that if $K_{Sat}$ satisfies (L) then so does $K$ (without the assumption that they are coverages). But this follows readily from the facts that $K\subset K_{Sat}$ and that a covering admits a refinement in $K_{Sat}$ \ssi it admits a refinement in $K$.
\end{Pf}

$K_{Sat}$ is the saturation of $K$ under coverings that have refinement in $K$. $\overline{K}$ is the saturation of the coverings of $K$ under pre-composition. When one combines the two saturation processes, one gets the \emph{sifted saturation}\indexb{Covering function!Sifted saturation of a --} of the covering function $K$. This is coherent with Definition \ref{def:saturation} when applied to a sifted covering function and does not depend on the order of the composition of these two saturation processes.

\begin{Prop}
	Let \C be a category and $K$ any covering function on \C. Then,
	$$
	(\overline K)_{Sat}=\overline{K_{Sat}}.
	$$
	We denote this covering function $J_K$\index[not]{JK@$J_K$} and call it the \emph{sifted saturation} of $K$. It is characterized by the following property: for any sieve $S$,
	$$
	S\in J_K\iff S\text{ contains a }K\text{-covering.}
	$$
	\cqfd
\end{Prop}
One can then combine Propositions \ref{prop:sifted} and \ref{prop:saturation} to study $J_K$. Note that it is sifted and sieve-saturated. 

If $K$ is a pretopology, then $J_K$ is a Grothendieck topology. This is the \emph{Grothendieck topology generated by the pretopology $K$}\indexb{Covering function!Grothendieck topology!-- generated by a pretopology}. 

Conversely, if a covering function $K$ generates a Grothendieck topology $J_K$, then $K$ must satisfy (C) and (L) (but not necessarily (M)). Given a Grothendieck topology $J$, there is a maximal covering function that generates it, and it is a pretopology. It is the maximal generating covering function $J_{max}$ of the sifted covering function $J$ (indeed, $J_{max}$ is automatically saturated, use Lemma \ref{lem:topsat} and \ref{prop:sifted5} of Proposition \ref{prop:sifted}).

\begin{Def}
A \emph{morphism of sites}\indexb{Site!Morphism of --s} from $(\C,K)$ to $(\D,L)$ is a functor $F\colon\C\to\D$ such that the \emph{image} $$F(R)\coloneqq\{F(f)\mid f\in R\}$$ of a covering $R$ of $K$ belongs to $L_{Sat}$, i.e., admits a refinement in $L$. 

For instance, if $K$ and $L$ are two covering functions on the same category \C, then $K$ is subordinated to $L$ if and only if the identity functor is a morphism of sites $$Id\colon(\C,K)\to(\C,L).$$ Sites and morphisms of sites with the usual composition of functors determine an XL-category denoted \SITE\index[not]{SITE@\SITE}. The category of small sites is written \Site\index[not]{Site@\Site}. Their full subcategories consisting of coverage-sites are written \textit{CSITE}\index[not]{CSITE@\textit{CSITE}} and \textit{CSite}\index[not]{CSite@\textit{CSite}} respectively.

We occasionally consider functors between sites that strictly preserve coverings.  Let $(\C,K)$ and $(\D,L)$ be sites. A functor $F\colon\C\to\D$ is a \emph{strict morphism of sites}\indexb{Site!Strict morphism of --s} if, given a $K$-covering $R$, the covering $F(R)=\{F(f)\mid f\in R\}$ belongs to $L$. One obtains the sub-XL-category \SITEs\index[not]{SITEs@\SITEs} of \SITE consisting of sites and strict morphisms of sites, and similarly for the other variants of \SITE.
\end{Def}

\subsection{Sites with pullbacks}\label{ssec:pretoppullback}
\subsubsection{Axioms of pretopologies revisited}
	The axioms we stated for a Grothendieck pretopology are somehow weaker than the usual axioms in the literature. Indeed, the following axioms are more frequent for a pretopology $K$ on a category \C.
	\begin{Axns}\label{axiom:strongpretop}
	\begin{itemize}
\item[\normalfont($\tilde{\text{M}}$)]\index[not]{AxiomM2@($\tilde{\text M}$)}
For every isomorphism $f\colon D\ra{\cong}C$ in \C, the covering $\{f\}$ belongs to $K$.
\item[]\textnormal{One also encounters the following variant}
\item[\normalfont($\tilde{\text{M}}$')]\index[not]{AxiomM3@($\tilde{\text M}$')} For every $C\in \C$, the identity covering $\{1_C\colon C\to C\}$ belongs to $K$.
\item[\normalfont($\tilde{\text{C}}$)]\index[not]{AxiomC2@($\tilde{\text C}$)}
For every $K$-covering $R$, every $C_f\ra{f}C\in R$ and every arrow $g\colon D\to C$, there exists a pullback
\begin{equation}\label{eq:pullcovering}
\begin{aligned}\shorthandoff{;:!?}
\xymatrix@!{
g^*(C_f) \ar[d]_{\bar f} \ar[r] & C_f \ar[d]^{f}\\
D \ar[r]_g & C
}\end{aligned}\end{equation}
such that the covering $\{g^*(C_f)\ra{\bar f} D\mid f\in R\}$ belongs to $K$.
\item[\normalfont($\tilde{\text{L}}$)]\index[not]{AxiomL2@($\tilde{\text L}$)}
Given a $K$-covering $R$ and, for each $f\in R$, a $K$-covering $R_f$ of $\dom f$, the composite covering
$$
\bigcup_{f\in R} f\circ R_f=\{f\circ g\mid f\in R\text{ and }g\in R_f\}.
$$
belongs to $K$.
\end{itemize}
\end{Axns}

\begin{Defs}
	\item
Given any covering $R$ of an object $C$ and an arrow ${g\colon D\to C}$, if pullbacks \hyperref[eq:pullcovering]{(\ref*{eq:pullcovering})} exist for each $f\colon C_f\to C$ of $R$, then we call the covering $\{\bar f\colon g^*(C_f)\to D\mid f\in R\}$ a \emph{pullback-covering of the covering $R$ along $g$}. If pullback-coverings of $R$ exist along every arrow with codomain $C$, then we say that the covering $R$ is a \emph{covering with pullbacks}.
\item
A site $(\C,K)$ \emph{has pullbacks of coverings} or is a \emph{site with pullbacks}\indexb{Site!-- with pullbacks} if all pullback-coverings of $K$-coverings exist (we do not require that they be $K$-coverings though; this is axiom ($\tilde{\text C}$)).
\end{Defs}

\begin{Rems}
	\item First, we remark that for a category \C with pullbacks, if $(\C,K)$ is a Grothendieck site, then $K_{Sat}$ is a Grothendieck site in this stronger sense. The converse is true, except for axiom (M): $K$ will satisfy axioms (C) and (L). We explain this now.
	
	First, for any site $(\C,K)$ with $K(C)\neq\emptyset$ for all $C\in\C$, the saturation $K_{Sat}$ of $K$ satisfies ($\tilde{\text{M}}$), since all coverings of an object $C$ refine a given isomorphism into $C$. Again, a Grothendieck site has coverings for all of its objects because of axiom (M). Now, recall that a covering function $K$ satisfies (C) and (L) \ssi $K_{Sat}$ satisfies (C) and (L) (Proposition \ref{prop:saturation}). Since $K_{Sat}$ is saturated, it satisfies (L) \ssi it satisfies ($\tilde{\text L}$) and when \C has pullbacks, it satisfies (C) \ssi it satisfies ($\tilde{\text{C}}$).
	
	Note that axiom (M) is not a consequence of $K_{Sat}$ being a Grothendieck topology in the stronger sense. However, this axiom is more a question of convention. In effect, given a covering function $K$, one can always decide to add to it all identities. This has no effect on the other axioms (C) and (L), and (M) is then satisfied. Thus, for categories with pullbacks, our definition of a pretopology is equivalent to the usual one. Nevertheless, the former gives more freedom for defining the equivalence class representative of covering functions. For example, there is a very natural pretopology on the category of topological spaces that consists of the open subset coverings. It does not satisfy ($\tilde{\text{M}}$), but (M).
	\item Next, we note that in case $K$ satisfies ($\tilde{\text{M}}$) and ($\tilde{\text{L}}$), then coverings are closed under pre-composition with isomorphisms. In particular, the condition ($\tilde{\text{C}}$) does not depend on the choice of the pullback in this situation.
\item
Finally, we remark that given a sieve $R$ on $C\in\C$ and an arrow $g\colon D\to C$ such that the pullback-covering $\{g^*(C_f)\ra{\bar f} D\mid f\in R\}$ of $R$ along $g$ exists, the sieve $g^*(R)$ is nothing but the sieve generated by the pullback-covering. This explains the notation $g^*(R)$.
\end{Rems}

\subsubsection{Sheaves with values in a category}
In sites with set-coverings and pullbacks of coverings, there is a generalization of the notion of sheaf. Let \C be any category, $R=\{C_i\ra{f_i}C\mid i\in I\}$ a set-covering in \C and suppose that pullbacks
\begin{equation}\label{eq:pullcov}\begin{aligned}
\xymatrix@=1.3cm{
**[l]{C_i\times_C C_j}\ar[r]^{\pi_{ij}^2}\ar[d]_{\pi_{ij}^1} & C_j\ar[d]^{f_j}\\
C_i\ar[r]_{f_i} & C
}\end{aligned}\end{equation}
exist for all $i,j\in I$. Let us make a choice of such a pullback for each pair $(i,j)\in I\times I$. Then, it is well known and easy to check that a presheaf $P$ on $\C$ is a sheaf for the covering $R$ \ssi, in the following diagram, $e$ is an equalizer of $p_1$ and $p_2$ in \Set.
\begin{equation}\label{eq:sheafpull}\begin{aligned}\xymatrixnocompile@C=1.5cm{
										& P(C_i)\ar[r]^{P(\pi_{ij}^1)}	& P(C_i\times_C C_j)\\
P(C)\ar[ru]^{P(f_i)}\ar[rd]_{P(f_j)}\ar@{-->}[r]^e	& \prod_{i\in I} P(C_i)\ar[u]_{pr_i}\ar[d]^{pr_j}\ar@<+.7ex>@{-->}[r]^-{p_1}\ar@<-.7ex>@{-->}[r]_-{p_2}		& \prod_{(i,j)\in I\times I}P(C_i\times_C C_j)\ar[u]_{pr_{ij}}\ar[d]^{pr_{ij}}\\
										& P(C_j)\ar[r]_{P(\pi_{ij}^2)}		& P(C_i\times_C C_j)
}\end{aligned}\end{equation}
\begin{Defs}
	[Let \C be a category and \A a category with products.]
	\item
A \emph{presheaf with values in \A}\indexb{Sheaf!Presheaf with values in \A} is a functor $P\colon\C^{op}\to \A$.
\item
Let $R=\{C_i\ra{f_i}C\mid i\in I\}$ be a set-covering in \C with pullbacks.\Par
A presheaf $P$ on \C with values in \A \emph{has the $\A$-sheaf property for the covering $R$}\indexb{Sheaf!A-sheaf property@\A-sheaf property} if $e$ is the equalizer of $p_1$ and $p_2$ in the diagram (\ref{eq:sheafpull}).

\item Given a site $(\C,K)$ with pullbacks and set-coverings, an \emph{$\A$-sheaf on \C in $K$}\indexb{Sheaf!A-sheaf@\A-sheaf} is a presheaf on \C that has the $\A$-sheaf property for every $K$-covering. We denote $Sh(\C,K;\A)$ the full subcategory of $\PSh(\C;\A)$ whose objects are $\A$-sheaves in $K$.
\end{Defs}
\begin{Rem}
	If the site is small, then there is a way to define sheaves with values in \A without pullbacks. Indeed, the sheaf property for a covering $R$ can be expressed in the same form as diagram (\ref{eq:sheafpull}), where the second product is not indexed by the pullbacks $C_i\times_C C_j$ but by commutative diagrams of the form
	$$\xymatrix{
	&C_i\ar[rd]^{f_i} &\\
	D\ar[ru]^g\ar[dr]_h && C\\
	& C_j\ar[ru]_{f_j} &
	}$$
	See also \cite{MM92} for a sifted version. When the site is large, the problem with these versions of the axiom is that they include a product indexed by a possibly proper class (take for example the singleton covering $\{id_C\}$). Now, by definition, given a set $A$ and sets $x_{\alpha}$ indexed by $A$, $$\prod_{\alpha\in A}x_{\alpha}=(\bigcup_{\alpha\in A}x_{\alpha})^A.$$
	A function with domain a proper class is a proper class. So, in \NBG, when $A$ is a proper class, the second member of this equality is empty. In \NBI, unless $\bigcup_{\alpha\in A}x_{\alpha}$ is empty, it is a proper conglomerate. This shows that the usual form of the product in \Set can at least not be a product of a proper class of (not all empty) sets. In fact, one can prove in \NBI that \Set does not have all products indexed by classes \cite[p. 114]{McL97}.
\end{Rem}

\begin{Def}
	A \emph{morphism of sites with pullbacks}\indexb{Site!Morphism of --s with pullbacks} is a morphism of sites that preserves pullbacks of coverings.\Par
	Sites with pullbacks and set-coverings (resp. small sites with pullbacks), and their morphisms, form an $XL$-category $\PSITE$\index[not]{PSITE@\PSITE} (resp. a category $\PSite$\index[not]{PSite@\PSite}). There are their coverage versions, \PCSITE\index[not]{PCSITE@\PCSITE} and \PCSite\index[not]{PCSite@\PCSite}.
\end{Def}

The following results are classical, but we haven't seen them in this generality. The next lemma is a generalization of lemma \ref{lem:sheafrefinement} for sheaves with values in a category \A (recall though that \A-sheaves are defined only for \emph{set}-coverings \emph{with pullbacks}). Indeed, the proof of lemma \ref{lem:sheafrefinement} in \cite{Joh202} can be written diagrammatically.

\begin{Lem}\label{lem:Asheafrefinement}
	Let \A be a category with products. Let $F$ be an \A-sheaf on a site $(\C,K)\in\PCSITE$.\Par
	Then $F$ has the $\A$-sheaf property for every set-covering in \C with pullbacks that has a refinement in $K$.\Par 
	In particular, if $(\C,K')$ is a site in $\PSITE$ such that $K'$ is subordinated to $K$ , then every $\A$-sheaf in $K$ is a $\A$-sheaf in $K'$. Finally, equivalent sites in \PCSITE have the same $\A$-sheaves.\cqfd
\end{Lem}

The following result is a direct corollary.
\begin{Lem}\label{lem:precompsheaves}
	Let $F\colon(\C,K)\to(\D,L)$ a morphism in \PCSITE and \A a category with products.\Par
	If $P$ is a $\A$-sheaf on $(\D,L)$, then $P\circ F$ is a $\A$-sheaf on $(\C,K)$.
\end{Lem}
\begin{Pf}
	Let $R$ be a covering of $(\C,K)$. Then $F(R)$ admits a refinement in $L$. Since $P$ is a sheaf in $L$ and $L$ is a coverage, by \thref{lem:Asheafrefinement}, $P$ is a sheaf for the covering $F(R)$. Now, since $F$ preserves pullbacks of coverings, $P$ is a sheaf for $F(R)$ \ssi $P\circ F$ is a sheaf for $R$.
\end{Pf}

\subsubsection{The bi-indexed category of sheaves}
Recall that we have defined the indexed category $\PSh(-;\A)\colon\Cat^{op}\to\CAT$ of presheaves with values in \A in \thref{ex:presheafA}. Given a morphism $F\colon(\C,K)\to(\D,L)$ in \PCSITE, the functor $$F^*=\PSh(F;\A)\colon\PSh(\D;\A)\to\PSh(\C;\A)$$ restricts, by the preceding lemma, to a functor $\Sh(\D,L;\A)\to\Sh(\C,K;\A)$, which we also denote $F^*$. The following proposition is now a direct consequence.

\begin{Prop}
	The indexed category $\PSh(-;\A)\colon\Cat^{op}\to\CAT$ of presheaves with values in \A induces the \emph{indexed category of sheaves with values in \A} over the category of small coverage-sites with pullbacks.
$$\begin{array}{rcl}
\Sh(-;\A)\colon\PCSite^{op}	& \longrightarrow 	& \CAT\\
(\C,K) 					& \longmapsto		& \Sh(\C,K;\A)\\
(\C,K)\ra{F}(\D,L)			& \longmapsto		& F^*\colon\Sh(\D,L;\A)\to\Sh(\C,K;\A)
\end{array}$$
	\cqfd
\end{Prop}
Its Grothendieck construction is written $Sh(\A)\to\PCSite$. When $\A$ is the category of rings (or commutative rings), this is called the fibration of \emph{ringed sites}.

The functor $\T\colon\Top^{op}\to\Cat$ induces a functor into $PCSite$. Observe that pullbacks in the category $\T(X)$ are given by intersections. One therefore obtains the \emph{opindexed category of sheaves on topological spaces with values in \A} by the composition:
$$
\Sh_{\Top}(-;\A)^{op}\colon\Top\ra{\T}\PCSite^{op}\ra{Sh(-;\A)}\CAT\ra{(-)^{op}}\CAT
$$
(see \thref{ex:presheafTop}). By abuse of terminology, given a topological space $X$, one calls a \emph{sheaf on $X$} a sheaf on the small site of $X$.
The total category $Sh_{\Top}(\A)^{op}$\footnote{Recall $Sh_{\Top}(\A)^{op}$ is equal to the dual of the pullback of $\Sh(\A)\to\PCSite$ over $\Top^{op}\to\PCSite$, whence the notation.} of the Grothendieck op-construction of the opindexed category $\Sh_{\Top}(-;\A)^{op}$ has the following form (we adopt the usual algebraic geometry notations):

\begin{itemize}
\item
\textbf{Objects}: Pairs $(X,\Oc X)$ where $X$ is a topological space and $\Oc X\in\Sh(X;\A)$ is a sheaf on $X$ with values in \A.
\item
\textbf{Morphisms}: Pairs $(X,\Oc X)\ra{(f,f^\sharp)}(Y,\Oc Y)$, where $f\colon X\to Y$ is a continuous map, $f^\sharp\colon \Oc Y\to f_*\Oc X$ is a morphism in $\Sh(Y;\A)$ and $f_*\Oc X=\Oc X\circ (f^{-1})^{op}$.

This indexed category is in fact a 2-functor and therefore, composition and identities of its Grothendieck construction have a simple form.
\item
\textbf{Composition}: The composite of the pair $$(X,\Oc X)\ra{(f,f^\sharp)}(Y,\Oc Y) \ra{(g,g^\sharp)}(Z,\Oc{Z})$$
has first component the composite $g\circ f$ of the continuous maps and second component given by
$$
\Oc Z \ra{g^\sharp}g_*\Oc Y \ra{g_*(f^\sharp)} (g\circ f)_*\Oc X.
$$
\item
\textbf{Identities}: $(X,\Oc X) \ra{(1_X,1_{\Oc X})}(X,\Oc X)$.
\end{itemize}

We call it the \emph{category of \A-spaces}\indexb{A-space@\A-space!Category of --s}. When $\A$ is the category of rings (it is often supposed that it is the category of commutative rings), this is called the category of \emph{ringed spaces}\indexb{Ringed space}. We denote it \Ringed\index[not]{Ringed@\Ringed}.

\paragraph{Bi-indexed structure}
We would like now to study the opindexedness of the indexed category $Sh(-;\A)$. In other words, we look for a left adjoint to the functor \begin{equation}\label{eq:invimsheaf}F^*=\Sh(F;\A)\end{equation} associated to a morphism $F\colon(\B,K)\to(\C,L)$ of coverage-sites with pullbacks. In fact, given such a morphism, the restriction to sheaves of a left adjoint $F_*$ to the functor $$F^*\colon\PSh(\C;\A)\to\PSh(\B;\A),$$ when it exists  (see Examples \hyperref[ex:bifibration:presheaves]{\ref*{ex:bifibration}(\ref*{ex:bifibration:presheaves})}), does not in general take values in sheaves. One therefore needs a \emph{sheafification functor}, that is, a left adjoint to the inclusion of the category of sheaves in the category of presheaves:
$$\xymatrix{
	a:\PSh(\C;\A) \ar@<.9ex>[r]\ar@{}[r]|-{\perp} & \Sh(\C,K;\A):i\ar@<.9ex>[l]
}$$
Once this is guaranteed, one obtains a composable pair of adjunctions
$$\xymatrix{
	\PSh(\B;\A) \ar@<.9ex>[r]^{F_*}\ar@{}[r]|{\perp} & \PSh(\C;\A)\ar@<.9ex>[l]^{F^*} \ar@<.9ex>[r]^-a\ar@{}[r]|-{\perp} & \Sh(\C,K;\A)\ar@<.9ex>[l]^-i
}$$
Using the composite adjunction and the preceding proposition, it is now not difficult to prove that it induces an adjunction between the respective categories of sheaves, which we denote:
$$\xymatrix{
	F_\sharp\colon\Sh(\B,K;\A) \ar@<.9ex>[r]\ar@{}[r]|-{\perp} & \Sh(\C,K;\A):F^*.\ar@<.9ex>[l]
}$$
Schapira and Kashiwara give a set of conditions on the category \A that guarantees the existence of a sheafification functor when $K$ is a Grothendieck topology on a small category \cite{KS06}\footnote{This therefore also applies when $K$ is a pretopology, since a pretopology is equivalent to the Grothendieck topology it generates, and therefore has the same $\A$-sheaves by Lemma \ref{lem:Asheafrefinement}. We haven't checked if the result of \cite{KS06} relies on the supplementary axioms of a pretopology or if it would already work for coverages. Note, by the way, that for each coverage, there exists a Grothendieck topology (in general not equivalent to the coverage, by Proposition \ref{prop:equivcovfun}) that has the same sheaves of sets \cite{Joh202}. This might also be true for \A-sheaves.}. See also the discussion in the article ``Sheafification'' in \cite{nLab} and references therein. Important examples of such categories $\A$ are the category \Set of sets, $\Grp$ of groups, $K$-\textit{Alg} of (commutative) algebras over a commutative ring $K$, $Mod_R$ of modules over a ring $R$. In particular, one obtains also the categories \Ab of abelian groups, \Ring of rings and \Comm of commutative rings.

When the category \A is such that the sheafification functor exists for small Grothendieck sites, one has in particular a bi-indexed category of $\A$-sheaves on topological spaces $\Sh_{\Top}(-;\A)^{op}\colon\Top\to\CAT$. A map $f\colon X\to Y$ gives rise to the following adjunction (we use the notation  of the algebraic geometry literature) with right adjoint on the left side
\begin{equation}\label{eq:adjsheavestop}\xymatrix{
	f_{*}\colon\Sh(X;\A) \ar@<-.9ex>[r]\ar@{}[r]|-{\perp} & \Sh(Y;\A):f^{-1}.\ar@<-.9ex>[l]
}\end{equation}
Let us write the left adjoint $f^{-1}$ by $\tilde f$ just for the next sentence. In the notation of the bi-indexed category of sheaves (on sites) and in terms of the morphism of sites ${f^{-1}\colon\T(Y)\to\T(X)}$, the functors of this adjunction are defined by $f_*\coloneqq(f^{-1})^*$ and $\tilde f\coloneqq(f^{-1})_\sharp$. The direct image functor seems to be a right adjoint. Nevertheless, these adjunctions determine the \emph{bi-indexed category of sheaves on topological spaces with values in \A}\label{page:bincatsheavestop} after having applied the opposite 2-functor, which reverses the roles of the adjoint functors (but not the directions of functors!). 

Consequently, the Grothendieck op-construction $\Sh_{\Top}(\A)^{op}\to\Top$ of the op-indexed category $\Sh_{\Top}(-;\A)^{op}\colon\Top\to\CAT$ is a bifibration, the \emph{bifibration of \A-spaces}\label{page:bifibAspace}\indexb{A-space@\A-space!Bifibration of --s}. An opCartesian lift of a map $f\colon X\to Y$ at an \A-space $(X,\Oc X)$ is given by the canonical opcleavage: $$(f,1_{f_*\Oc X})\colon(X,\Oc X)\to (Y,f_*\Oc X).$$ Once an adjunction $f^{-1}\dashv f_*$ of unit $\eta^f$ is chosen, a Cartesian lift of $f$ at $(Y,\Oc Y)$ is given by \begin{equation}\label{eq:invimAspace}(f,\eta^f_{\Oc Y})\colon(X,f^{-1}\Oc Y)\to(Y,\Oc Y).\end{equation} This is an application of \thref{prop:bifibadjoint}, but one should not be misled by the fact that the opfibred category is in fact built from the \emph{dual} of the adjunction $f^{-1}\dashv f_*$. Note that $(\eta^f)^{op}\colon (f_*)^{op}(f^{-1})^{op}\Rightarrow Id$ is the \emph{counit} of this dual adjunction.

This bifibration is isomorphic to the Grothendieck construction of the associated indexed category (whose inverse image functors are the dual of functor $f^{-1}$) (see Theorem \ref{thm:Bifib}). But one prefers the op-construction, because the op-indexed category is a 2-functor, whereas the indexed category is a mere pseudo-functor. Indeed, for each direct image functor $F_*$ associated to a morphism of sites, there is a non-canonical choice of a left adjoint $F_{\sharp}$ (when it exists). One can construct an indexed category from these left adjoints using an axiom of choice (see Remark \ref{rem:choiceadjpseudo}). 

We observe now that one can, in particular cases, make use of the freedom in the choice of the left adjoint functor to have particularly simple inverse image functors, instead of taking the adhoc construction via colimits. We state first the following toy result. Its proof is quite straightforward. We just remark that point (i) comes from the fact that the coverings of the empty open subset $\emptyset\in\T(X)$ of a space are the identity covering and the empty covering. The latter covering induces the condition.

\begin{Lem}\label{lem:sheafonpoint}
	Let \A be a category with products. Let $X$ be a space and $*$ the one-point space.
	
	\begin{enumerate}[(i)]
\item\label{lem:sheafonpoint:terminal}
For any \A-sheaf $F$ on $X$, $F(\emptyset)$ is a terminal object in \A.
\item
$\PSh(*;\A)\cong\A^{\two}$ and there is an strictly surjective equivalence of categories
$$\Sh(*;\A)\simeq\A,$$
given by taking the image $F(*)$ of a sheaf $F$ at the point $*$. Under these identifications, the inclusion functor of sheaves into presheaves amounts, given a choice $*_{\A}$ of a terminal object of \A, to associate the arrow $A\ra!*_{\A}$ to the object $A\in\A$.
\item
The sheafification functor $\PSh(*;\A)\to\Sh(*;\A)$ exists and is given, under these identifications, by the functor domain $\dom\colon\A^{\two}\to\A$. 
\end{enumerate}
\cqfd
\end{Lem}

\begin{Exs}\label{ex:invimsheaves}
	\item
	When $f=1_X$ is the identity map, the adjunction \hyperref[eq:adjsheavestop]{(\ref*{eq:adjsheavestop})} can be chosen to be the identity adjunction, because $(1_X)_*=Id_{\Sh(X;\A)}$.
	\item
	\label{ex:invimsheaves:stalk} Suppose \A has products and colimits. 
	
	The adjunction \hyperref[eq:adjsheavestop]{(\ref*{eq:adjsheavestop})}, when the map $f$ is the inclusion $i_x\colon\{x\}\hookrightarrow X$ of a point $x$ of a space $X$, is very important in algebraic geometry. The functor $$(i_x)_*\colon\A\to\Sh(X;\A)$$ is called the \emph{skyscraper sheaf functor over $x$} and is also written $Sky_x$. The sheaf $Sky_x(A)$ has value $A$ on any open neighbourhood $U$ of $x$ and value $*$ (a fixed terminal object of \A) on open subsets of $X$ not containing $x$. Its left adjoint $(i_x)^{-1}$, which exists with no further assumptions on \A, is called the \emph{stalk functor} and is also written $Stalk_x$. To describe it, observe that there is a subcategory $\T(X)_x\xhookrightarrow{i}\T(X)$ of open neighbourhoods of $x$, which is a sub-coverage-site with pullbacks. The inclusion functor, by \thref{lem:precompsheaves}, induces a functor between the categories of sheaves
	$$
	\Sh(\T(X);\A)\ra{i^{*}}\Sh(\T(X)_x;\A).
	$$
	The stalk functor is the following composite of functors
	$$
	\Sh(\T(X);\A)\ra{i^{*}}\Sh(\T(X)_x;\A)\ra{\text{colim}|_{\Sh}}\A,
	$$
	where $\text{colim}|_{\Sh}$ is the restriction to sheaves of the usual colimit functor.
	Finally, one has the adjunction
	$$\xymatrix{
	Stalk_x\colon\Sh(X;\A) \ar@<.9ex>[r]\ar@{}[r]|-{\perp} & \A:Sky_x.\ar@<.9ex>[l]
	}$$
	We introduce standard notation: Given sheaves $F$ or $\Oc X$ on $X$, we write $$F_x\coloneqq Stalk_x(F)\text{ and }\Oc{X,x}\coloneqq Stalk_x(\Oc X)$$ and given a morphism $\phi$ of sheaves on $X$, we denote $\phi_x\coloneqq Stalk_x(\phi)$.
	\item \label{ex:invimsheaves:toterminal}
	Consider now another extreme case, where the map $f$ is the unique map $!\colon X\to*$ from a space $X$ to a singleton space. Suppose the sheafification functor $$\PSh(X;\A)\to\Sh(X;\A)$$ exists. The functor $!_*\colon\Sh(X;\A)\to\A$ takes a sheaf $F$ to its value $F(X)$ at $X$. It is called the \emph{global section functor} and denoted $\Gamma(X;-)$. We now calculate its left adjoint. Observe that the global section functor is in fact the restriction to sheaves of the \emph{limit functor} $\lim\colon\A^{\T(X)^{op}}\to\A$, which has a left adjoint, the constant functor $\Delta$. One obtains a composite adjunction:
	$$\xymatrix{
	\A \ar@<.9ex>[r]^-{\Delta}\ar@{}[r]|-{\perp} & \PSh(X;\A)\ar@<.9ex>[l]^-{\lim} \ar@<.9ex>[r]^{a}\ar@{}[r]|-{\perp} & \Sh(X;\A)\ar@<.9ex>@{_(->}[l]^{I}.
	}$$
	Thus, a left adjoint of the global section functor acts as follows on objects $A\in\A$: take the sheafification of the constant presheaf $\Delta(A)$. This functor is called the \emph{constant sheaf functor}\indexb{Sheaf!Constant -- functor} and we write it $\tilde\Delta_X$\index[not]{Deltatilde@$\tilde\Delta$}. Note that $\Delta(A)$ is never a sheaf, unless $A$ is a terminal object of \A, by \hyperref[lem:sheafonpoint:terminal]{Lemma \ref*{lem:sheafonpoint}(\ref*{lem:sheafonpoint:terminal})}. In fact, when $*_{\A}$ is a terminal object, $\Delta(*_{\A})$ is a sheaf, and it is terminal in $\Sh(X;\A)$. In order to avoid confusion, remark that in the algebraic geometry literature, the constant presheaf $\Delta(A)$ is not the image of $A$ under the constant functor. It has the same values as our $\Delta(A)$, except at the empty open set $\emptyset$, where it is defined to have value a terminal object $*_{\A}$ of \A (and so $A\ra!*_\A$ on $\emptyset\hookrightarrow U$ for an open $U$). Then, $\tilde\Delta_X(A)$ is also the sheafification of this presheaf. See \cite{Har77,Qin02} for further information on constant sheaves.
	\item\label{ex:invimsheaves:openembedding}
Suppose \A has products and consider the case where the map $f$ is an \emph{open embedding}, i.e., an open injective map, $f\colon Y\hookrightarrow X$ into a space $X$. In this situation, the inverse image functor always exists (that is, with no further assumptions on \A) and is defined as follows. The map $f$ induces a morphism of coverage-sites with pullbacks $f\colon\T(Y)\to\T(X)$ and therefore an inverse image functor between the categories of sheaves on these sites $f^*\colon\Sh(X;\A)\to\Sh(Y;\A)$, like in \hyperref[eq:invimsheaf]{(\ref*{eq:invimsheaf})}. This functor is in fact a left adjoint to the functor $f_*\colon\Sh(Y;\A)\to\Sh(X;\A)$, as one may readily verify, and is therefore denoted $f^{-1}$. The unit $\eta^f$ of this adjunction at a sheaf $F\in\Sh(X;\A)$ is determined by the restriction morphisms $$\res_{\scriptscriptstyle V,U}=F(U\hookrightarrow V)\colon F(V)\to F(U).$$ It is the morphism of sheaves $\eta^f_F\colon F\to f_*f^{-1}F$ defined by
$$
(\eta^f_F)_U\coloneqq\res_{\scriptscriptstyle U,ff^{-1}(U)}\colon F(U)\to F(ff^{-1}(U)).
$$
In case of an inclusion of an open subset $i\colon V\hookrightarrow X$, one writes $F|_V\coloneqq i^{-1}F$. The unit in this case is thus given by the restrictions $(\eta^i_F)_U=\res_{\scriptscriptstyle U,U\cap V}\colon F(U)\to F(U\cap V)$.

\end{Exs}

\subsection{Examples of covering functions}\label{ssec:covfunex}
We now turn to examples of covering functions, some of general use and some in particular categories that appears in the next chapters. We also show that some apparently different covering functions are in fact equivalent. We do not give examples of covering functions on $\Cat$, but we would like to mention the article \cite{BBP99} that is devoted to the study of all different kinds of epis in $\Cat$. 

\subsubsection{General categories}\label{sssec:gentop}
\paragraph{Trivial topologies}\label{par:trivtop}
We first consider some extreme cases.
\begin{enumerate}
\item
The \emph{no-covering coverage} that has no coverings at all. 
$$
\begin{array}{rcl}
K\colon\Ob\C & \longrightarrow & \power(\power(Mor\C))\\
C & \longmapsto & \emptyset
\end{array}
$$
It is trivially a (sifted) coverage, but not a (pre-)topology because of axiom (M). It is subordinated to all covering functions.
\item
The \emph{empty-covering coverage} defined by
$$
\begin{array}{rcl}
K\colon\Ob\C & \longrightarrow & \power(\power(Mor\C))\\
C & \longmapsto & \{\emptyset\}.
\end{array}
$$
It has the same properties as the previous one.
\item\label{indtop}
The covering function that has all and only identity coverings is a pretopology, called the \emph{coarsest {\normalfont or} indiscrete pretopology}. Its coverings do not belong to every pretopology in our axiom system, but they belong to their saturation. It belongs to every covering functions satisfying axiom $(\tilde{\text{M}})$ though, and that is where its name comes from (see part \ref{ssec:pretoppullback}). Idem for its sifted version.
\item
The covering function that has all possible coverings (resp. sieves) is also trivially a pretopology (resp. a topology), called the \emph{finest {\normalfont or} indiscrete pretopology} (resp. the \emph{finest topology.}). All covering functions are subordinated to it.
\end{enumerate}

\paragraph{Pretopologies of epis}
Let us consider a general category \C. Recall that axiom (M) is equivalent to the condition that for each object, there is some covering containing a split epi. Important classes of pretopologies come from covering functions whose coverings are all singletons that belong to some class of epis containing the split epis. Recall that there is a hierarchy of epis in a category \C (see \cite{BorI94,PT04} for this part):
$$
Split Epis\subset Regular Epis \subset Strong Epis \subset Extremal Epis \subset Epis.
$$ 
Of course, this hierarchy starts with $Identities\subset Isomorphisms$. In categories with pullbacks, extremal and strong epis coincide. In regular categories, regular, strong and extremal epis coincide and are moreover stable under pullbacks. 

Identities, isomorphisms, and split epis are stable under pullbacks in general categories, but this is not true for the other weaker notions of epis. One says that a morphism is a \emph{stably-extremal epi} or \emph{universally extremal epi} if all its existing pullbacks are extremal epis, and similarly for the other kind of epis. Moreover, using the fact that the composite of two pullback squares is a pullback square, one obtains that these classes of stable types of epi are themselves stable under pullbacks. Together with assumption that \C has pullbacks, this guarantees axiom ($\tilde{\text{C}}$) of the covering functions. 

Axiom ($\tilde{\text{L}}$) requires stability under composition. This in fact requires pullbacks and we therefore do not need to separate the extremal and strong cases. If \C has pullbacks, then the class of stably-regular epis and of stably-strong epis are closed under composition (the regular case is proved in \cite{PT04} and the strong case is not difficult). 

One therefore has the following pretopologies on a category \C \emph{with pullbacks}:
\begin{enumerate}
\item
The \emph{pretopology $SplitEpi(\C)$ of split epis}, which is equivalent to the coarsest topology.
\item
The \emph{pretopology $SRegEpi(\C)$ of stably-regular epis}. In a regular category (in particular in a Barr-exact category and in a topos), all regular epis are stable. One calls it the \emph{regular pretopology}.
\item
The \emph{pretopology $SStrongEpi(\C)$ of stably-strong epis}. It coincides with the regular pretopology in a regular category.
\item
The \emph{pretopology $SEpi(\C)$ of stable epis}. It coincides with the regular topology in a topos, since all epis are regular in there.
\end{enumerate}


\paragraph{Subcanonical covering functions}
 \begin{Def}
	A covering function $K$ on a locally small category \C such that every representable presheaf on \C is a $K$-sheaf is called \emph{subcanonical}.
\end{Def}

Observe that, by Lemma \ref{lem:sheafrefinement}, if $K$ is a subcanonical coverage and $K'$ is subordinated to $K$, then $K'$ is subcanonical. In particular, if $K$ and $K'$ are equivalent coverages, then one is subcanonical \ssi the other is.

Let \C be a locally small category with pullbacks and $K$ a covering function on \C. If the coverings of $K$ are singletons, then $K$ is subcanonical \ssi its coverings are regular epis. If, moreover, $K$ is a coverage, then it is subcanonical \ssi its coverings are stably-regular epis \cite{Joh202}. Therefore, under the hypotheses on \C, the stably regular pretopology is the largest subcanonical coverage whose coverings are singletons. In the following, when we talk about subcanonical covering functions on a category \C, it is implicitly supposed that \C is locally small. For further study of subcanonical coverages, see the next paragraph on jointly epimorphic coverings. See also the latter reference and \cite{DLORV}.

\paragraph{Jointly epimorphic coverings}
\begin{Def}
	Let $\C$ be a category. A covering $R=\{C_f\ra{f}C\}$ of an object $C\in\C$ is \emph{jointly epimorphic} if, given two arrows $g,h\colon C\rightrightarrows D$ in \C such that, for all $f$ in $R$,
	$$g\circ f=h\circ f,$$
then $g=h$.
\end{Def}
If $\C$ has (small or large, depending on the size of the covering) coproducts, this is equivalent to the condition that the induced arrow $\coprod C_f\to C$ to be epimorphic. Important examples come from colimits: any colimiting cone is indeed jointly epimorphic.

Many usual covering functions have jointly epimorphic coverings. Of course, all singleton coverings that are epimorphic are jointly so. So all the pretopologies of epis that we have described give examples. 

There is a many-arrow generalization of the notion of strong epi, called \emph{jointly strongly epimorphic} coverings. For a category with (large coproducts), it is equivalent to the condition that the induced arrow from the coproduct is a strong epi. There is also the corresponding pullback-stable notion, called \emph{stably jointly strongly epimorphic} coverings. Under the preceding hypothesis of existence of coproducts, such a covering induces a stably strong epi from the coproduct, but the converse seems not to be true in general (We have not searched for counter-examples. It is equivalent in $\Top$). We have proved that the stably jointly strongly epimorphic coverings form a pretopology. 

All subcanonical covering functions have jointly epimorphic coverings, as one can readily check. More precisely, any covering for which every representable presheaf is a sheaf is jointly epimorphic. They are called \emph{effective-epimorphic} coverings. They can be seen as many-arrow generalization of regular epis
\footnote{See the paragraph on subcanonical covering functions. In fact, effective-epimorphic coverings $R$ on an object $C$ are also characterized by a the fact that the sieve $\overline{R}$ they generate is a colimit cone on $C$ for the diagram consisting on the subcategory of \C whose objects are the domains of arrows in $\overline R$ and morphisms, arrows over $C$ \cite{Joh202}.}.
The largest subcanonical covering function is of course the covering function all of whose coverings are effective-epimorphic. Now, what about the largest subcanonical coverage ? 

If $K$ is a subcanonical coverage, then its coverings have a stronger property. In a general category, it is more easily stated in the sifted setting. If $P$ is a representable presheaf, then it is a sheaf on any $K$-sieve. Given a $K$-sieve $R$ of an object $C$ and an arrow $g\colon D\to C$, there exists a $K$-sieve $S$ on $D$ contained in $g^*(R)$ and $P$ must be a sheaf on $S$. So, by Lemma \ref{lem:sheafrefinement}, $P$ is a sheaf on $g^*(R)$. Thus $R$ has the additional property that given an arrow $g$ with codomain $C$, any presheaf is a sheaf for $g^*(R)$. In general, a sieve $R$ that has this property is called \emph{universally effective-epimorphic}. It happens that the sifted covering function consisting of all sieves of this type is a subcanonical sifted coverage, and therefore, the largest one \cite{Joh202}. It is called the \emph{canonical topology} on \C and written $J_{can}$. In fact, it is even a topology. We now describe the canonical coverage for a general category and relate it to the characterization of \cite{DLORV} when \C has pullbacks. To my knowledge, this is new, but might be obvious to experts.
\begin{Prop}
	Let \C be a category. There exists a largest subcanonical coverage on \C. It is given by the maximal generating covering function $(J_{can})_{max}$ of the canonical topology $J_{can}$. It is a pretopology. We call it the \emph{canonical pretopology} on \C and write it $K_{can}$.\Par
	If \C has pullbacks, then $K_{can}$ is the covering function consisting of all pullback-stable effective-epimorphic coverings.
\end{Prop}

\begin{Pf}
	$K_{can}$ is a pretopology because it is the maximal generating covering function of a topology. In particular, it is a coverage. It is subcanonical because being a sheaf on a covering $R$ is equivalent to being a sheaf on the sieve $\overline R$ it generates. It is maximal among subcanonical coverages by the discussion above.
	
	Now, let \C have pullbacks. Let $R$ be a covering of an object $C$, $g\colon D\to C$ a morphism in \C, and S a pullback-covering of $R$ along $g$. We have already remarked that $\overline S= g^*R$. Let $P$ be a representable presheaf. Then $P$ is a sheaf on $S$ \ssi it is on $\overline S=g^*(\overline R)$. Therefore, $\overline R$ is universally effective-epimorphic \ssi $R$ is a pullback-stable effective epimorphic covering.
\end{Pf}

\paragraph{Induced covering functions on slice categories}

Let $K$ be a covering function on a category \C and $C\in\C$. It induces a covering function on the slice category $\C/C$, called the \emph{slice covering function} and written $K_C$, in the following manner. A covering $R_C=\{\xi_i\ra{f_i}\xi\}$ of $\xi=(A\ra{p}C)$ (with $\xi_i=(A_i\ra{p_i}C)$) belongs to $K_C$ if and only if $R=\{A_i\ra{f_i}A\}$ belongs to $K$. Thus, there is a bijection between the coverings of $\xi$ in $K_C$ and the coverings of $A$ in $K$. $K_C$ inherits properties of $K$. For instance, it not difficult to show that if $K$ is sifted, respectively satisfies axioms (M), (C) or (L), or is subcanonical, then so is $K_C$.

Covering functions $K$ are often defined to be the collections of all coverings whose arrows satisfy some categorical property in \C. It happens sometimes that the induced covering function $K_C$ is exactly the covering function defined by same property, but in $\C/C$. For instance, consider \C \emph{with all finite limits} and the pretopology $SEpi(\C)$ of stable epis in $\C$. Then, the slice pretopology $SEpi(\C)_B$ over $B$ is equal to the pretopology $SEpi(\C/B)$ of stable epis in $\C/B$. One can check this using the fact that the forgetful functor $U\colon\C/B\to\C$ is a left adjoint (see Examples \ref{ex:fibmonoids}), and therefore preserves epis, and that the pullback in $\C/B$ is given by the pullback in $\C$ (after applying $U$) with the unique possible arrow to $B$. We have also checked this result for the pretopologies of stably-strong (=stably-extremal) epis and of split epis (\C only needs to have pullbacks for the latter). 
\subsubsection{Covering functions on topological spaces}
Firstly, we consider a single topological space and its preorder category $\T(X)$ of open subsets with inclusions. This category admits a pretopology. Its coverings are the open subset coverings in the topological meaning: a covering of an open subset $V\subset X$ is a family of open subsets of $V$ whose union is $V$. $\T(X)$ with this pretopology is called the \emph{small site associated to X}.\indexb{Site!Small - associated to a space}
\paragraph{}
We now turn to covering functions on the category $\Top$ of topological spaces.
\paragraph{Epi coverings}
In the category \Top of topological spaces, epis are the surjective maps. 

Moreover, extremal, strong and regular epis all coincide and are the quotient maps. Yet, there exists non regular epis, for example all bijective maps that are not homeomorphisms. 

The stably-regular(=stably-strong=stably-extremal) epis in \Top are characterized in an article of Day and Kelly \cite{DK70}. See also \cite{Mic68}. They are called \emph{biquotient maps}. Note that all open surjections are biquotient. Moreover, they are stable under composition and pullbacks. Consequently, they form a sub-pretopology of the stably-regular topology.


\paragraph{Jointly epimorphic open coverings}
We provide examples of pretopologies in the category of topological spaces whose coverings are particular set of jointly epimorphic open maps. Note first that in $\Top$, a covering $\{Y_i\ra{f_i}X\}$ is a set of jointly epimorphic open maps \ssi the induced map $\coprod_{i\in I}Y_i\ra{f}X$ is a surjective open map. In the following examples, each pretopology is contained in the following.
\begin{enumerate}
\item
The coverings of the \emph{open subset pretopology} are the open coverings in the topological sense: if $X$ is any topological space, then an open covering of $X$ is a family ${\{U_i\hookrightarrow X\mid i\in I\}}$ of inclusions of open subsets $U_i\subset X$ such that their union $\bigcup_{i\in I}U_i$ is $X$. \Top together with this covering function is a coverage-site with pullbacks. Note however that only particular choices of pullbacks of coverings of the coverage belongs to the coverage.
\item
The \emph{open embedding pretopology} has the following coverings. They are set of arrows $\U=\{Y_i\ra{f_i}X\}$ such that
\begin{enumerate}[(i)]
\item
Each $f_i\colon Y_i\to X$ is an \emph{open embedding}, i.e., an injective open map. 
\item
The family $\U$ is jointly surjective.
\end{enumerate}

\item
The \emph{étale pretopology} has the following coverings. They are set of arrows $$\U=\{Y_i\ra{f_i}X\}$$ such that
\begin{enumerate}[(i)]
\item
Each $f_i\colon Y_i\to X$ is \emph{étale} (also called \emph{local homeomorphism}): for all $y\in Y_i$, there exists an open neighbourhood $U$ of $y$ such that the restriction of $f_i$ to $U$ is a homeomorphism onto an open subset of $X$. 
\item
The family $\U$ is jointly surjective.
\end{enumerate}
\item
The \emph{open pretopology} has coverings that are sets of jointly surjective open maps.
\end{enumerate}

\begin{Lem}
	The open subset, open embedding and étale pretopologies are equivalent. Moreover, the open pretopology is subcanonical, and therefore, all these pretopologies are subcanonical.
\end{Lem}
The first statement, that can be found in \cite{Vis08}, is an exercise of topology. A proof of the second can be found in \cite{MM92}. 

Given a space $X$, the \emph{big site associated to $X$}\indexb{Site!Big - associated to a space} is the slice category $Top/X$ together with the slice covering function induced by one of the first three equivalent pretopologies of this lemma.

\subsubsection{Covering functions on smooth manifolds}
We now consider the category $\Diff$ of smooth manifolds \cite{NS03,Nab97,War83}. We give two pretopologies on this category. One should be aware that \Diff does not have all pullbacks. Nevertheless, particular types of morphisms admit pullbacks. For instance, pullbacks of open embeddings and of surjective submersions exist, and are still open embeddings, resp. surjective submersions.
\paragraph{The pretopology of open subset coverings}
The coverings of a manifold $M$ in this pretopology are those of the pretopology of open subsets when $M$ is seen as a topological space. Each open subset is seen as an open submanifold and the inclusion is then automatically an open embedding. This pretopology is equivalent to the pretopology of collectively surjective open smooth embeddings (topological embeddings that are immersions).
\paragraph{The pretopology of surjective submersions}
Coverings of this pretopology are singletons, whose map is a surjective submersion. The induced map $\coprod U_i\to M$ of an open subset covering $\{U_i\hookrightarrow M\}$ of a manifold $M$ is a surjective submersion. Moreover, this pretopology is subordinated to the previous one.

\subsubsection{Covering functions on commutative rings}\label{sssec:comring}
Let $\Comm$ denotes the category of commutative rings and ring homomorphisms. The \emph{Zariski pre-optopology} on \Comm has coverings $R=\{A \ra{f_i}A_i\mid i\in I\}$ indexed by sets $I$ such that:
\begin{enumerate}[(i)]
\item
Each $A\ra{f_i}A_i$ is a localization of $A$ at an element $a_i$ of $A$.
\item
The ideal $(a_i)_{i\in I}\subset A$ generated by the set $\{a_i\mid i\in I\}$ is the whole ring $A$.
\end{enumerate}
There is a proof in \cite{MM92} that this covering function is a pre-optopology in the stronger sense of \thref{axiom:strongpretop}, but under the assumption that all coverings are finite. One can adapt the proof to the more general coverings we have defined though, by means of the following well-known observation: any covering of this covering function admits a finite subcovering. Indeed, suppose $(a_i)_{i\in I}= A$. Then, $$1=\sum_{k=1}^nr_ka_{i_k}$$ for a natural $n$ and $r_k\in A$, $1\leq k\leq n$. Consequently, $A=(a_{i_1},\dots,a_{i_n})$. 

\section{Fibred sites}\label{sssec:covinfib}

There are different nonequivalent ways of mixing the notions of site and of fibration. 

On one hand, one may consider a type of fibred site that one could call a ``category fibred in sites''. It is a fibration whose fibres are sites and whose inverse image functors are strict morphisms of sites. This notion is well-defined (that is, independent of choices of inverse image functors) if the sites of the fibres satisfy both axioms \tm and (L). Note that the base category does not need to be endowed with a covering function. This kind of concept is considered in SGA 4 \cite[Exposé VI]{SGAIV2}. Informally, one can see it as a ``site internal to the 2-XL-category of fibrations''. 

On the other hand, recall the following definition\footnote{By internal fibrations, we mean the ``strict'' version as defined in \cite{Str74a}. See also \cite{Web07,nLab}}.

\begin{Def}
	Let $\dA$ be a 2-category. A \emph{fibration internal to} $\dA$ is a morphism $p\colon E\to B$ of \dA such that
	\begin{enumerate}[(i)]
\item
for all $X\in \dA$, the functor
$$p_*\colon\dA(X,E)\to \dA(X,B)$$
is a fibration (of categories),
\item
for all morphism $f\colon X\to Y$ in \dA, the following commutative square is a Cartesian morphism of fibrations
$$\xymatrix@=1.5cm{
\dA(X,E)\ar[r]^{f^*}\ar[d]_{p_*} & \dA(Y,E)\ar[d]^{p_*}\\
\dA(X,B)\ar[r]_{f^*} & \dA(Y,B).
}$$
\end{enumerate}
\end{Def}

One might thus consider (formally now) fibrations internal to the 2-XL-category of sites, strict morphisms of sites and natural transformations, which we denote \SITEs\index[not]{SITEs@\SITEs}, as the corresponding XL-category. Let us denote $$(\B,K_\B)^{(\A,K_\A)}$$ the hom-category $\SITEs((\A,K_\A),(\B,K_\B))$. If the base site satisfies axiom \tmp and the total site axioms \tm and (L), then such an object is also a ``category fibred in sites''. 

What we call a fibred site is a third kind of object, related to the previous ones. A sifted version of this notion exists in Jardine's paper \cite{Jar06} (see below after \thref{prop:KandKE}).
\begin{Def}
	A \emph{fibred site}\indexb{Site!Fibred--} consists in a Grothendieck fibration equipped with a covering function on its base. 
\end{Def}
Let us now study the relationship between this notion and the two previously defined notions of fibred sites. Consider a Grothendieck fibration $\E\ra P\B$ and a covering function $K$ on the base \B. It induces a covering function $K_P$\index[not]{KP@$K_P$} on the total category $\E$ whose coverings are all the Cartesian lifts of coverings of $K$. More precisely, the coverings of $K_P$ are the families of Cartesian arrows
$$\{{f_{i}}^*E\ra{\overline{f_{i}}_{E}}E\mid i\in I\},$$
where $R=\{f_{i}\colon B_{i}\to B\}$ is a $K$-covering. We say that such a covering is a $K_P$-covering \emph{over} the $K$-covering $R$. The functor $P$ becomes a strict morphism of sites, by definition of $K_P$.

\begin{Def}
  \label{def:indcov}
The covering function $K_P$ on the total category \E, also denoted $K_\E$ when this does not induce confusion, is called the \emph{covering function induced by the fibration $P$}\indexb{Covering function!-- induced by a fibration}.
\end{Def}

When the base site satisfies \tmp, this structure provides us with a ``category fibred in sites'', but in a rather trivial way. Indeed, the restriction of $K_P$ to a fibre $\E_B$ consists in isomorphism singleton coverings. This is a pretopology equivalent to the indiscrete pretopology, see \hyperref[indtop]{Exemple \ref*{indtop} on page \pageref*{indtop}}. 

A fibred site is more interestingly seen as a internal fibration in \SITEs. Indeed, a fibred site $P\colon\E\to(\B,K)$ gives rise to such an object when \E is equipped with the induced covering function $K_P$. Moreover, under some mild assumptions, the covering function $K_P$ is, up to equivalence, the least element among covering functions on \E that make $P$ an internal fibration in \SITEs.

\begin{Lem}
	Let $P\colon (\E,K_\E)\to (\B,K_\B)$ be an internal fibration in \SITEs and $(\D,K_\D)$ a site. Suppose that $K_\B$, $K_\D$ and $K_\E$ satisfy \tmp. Then all components of a Cartesian morphism in $P_*\colon (\E,K_P)^{(\D,K_\D)}\to (\B,K)^{(\D,K_\D)}$ are Cartesian in $P$.
\end{Lem}

\begin{Pf}
	Consider a site $(\one,id)$ with one object, one arrow and one identity covering. Under the hypotheses, the constant functor $\Delta(D)\colon(\one,\id)\to(\D,K_\D)$ at an object $D\in \D$ is a strict morphism of sites. Moreover, on has the equality of categories
	$$(\B,K_\B)^{(\one,\id)}=\B^\one,$$
	where the right-hand side is the functor category, and similarly for $(\E,\id)$. Consider the following commutative diagram.
	\begin{equation}\begin{aligned}\label{eq:cartcomp}\xymatrix@=1.3cm{
(\E,K_\E)^{(\D,K_\D)}\ar[d]_{P_*}\ar[r]^{\Delta(D)^*} &(\E,K_\E)^{(\one,id)}\ar@{}[r]|(.6)=\ar[d]_{P_*} &\E^\one\ar@{}[r]|{\cong}\ar[d]_{P_*} &\E\ar[d]^{P}\\
(\B,K_\B)^{(\D,K_\D)}\ar[r]^{\Delta(D)^*} &(\B,K_\B)^{(\one,id)}\ar@{}[r]|(.6)= &\B^\one\ar@{}[r]|{\cong} &\B
}\end{aligned}\end{equation}
Since $P$ is an internal fibration in \SITEs and $\Delta(D)$ is in \SITEs, this is a Cartesian morphism of fibrations. Again, applying it to a Cartesian morphism $\beta$ in $$P_*\colon (\E,K_P)^{(\D,K_\D)}\to (\B,K)^{(\D,K_\D)}$$ produces the $D$-component of $\beta$, which has to be Cartesian.
\end{Pf}

\begin{Prop}
	If $P\colon\E\to(\B,K)$ is a fibred site, then $K_P$ makes $P$ an internal fibration in \SITEs.\Par	
	Conversely, let $P\colon (\E,K_\E)\to (\B,K_\B)$ be an internal fibration in \SITEs. If both $K_\B$ and $K_\E$ satisfy \tmp, then $P$ is a fibration (of categories). If moreover $K_\E$ satisfies axioms \tm and \textnormal{(L)}, then $K_P\preceq K_\E$.
\end{Prop}

\begin{Pf}
	Let $P\colon\E\to(\B,K)$ be a fibred site and $(\D,K_\D)$ a site. Let $\alpha\colon F\Rightarrow G$ in $(\B,K)^{(\D,K_\D)}$ and $H\in (\E,K_P)^{(\D,K_\D)}$ over $G$ with respect to the functor \begin{equation}\label{eq:P*}P_*\colon (\E,K_P)^{(\D,K_\D)}\to (\B,K)^{(\D,K_\D)}.\end{equation} One obtains a Cartesian morphism $\bar{\alpha}_H\colon\alpha^*H\Rightarrow H$ in $P_*$ by defining it component $(\bar\alpha_H)_D$ at $D\in\D$ to be a Cartesian lift of the component $\alpha_D$ of $\alpha$ at $H(E)$ in the fibration $P\colon\E\to\B$. One checks the fact that $\alpha^*H$ is not only a functor but a strict morphism of sites by means of \hyperref[Prop:Cartesianii]{Proposition \ref*{Prop:Cartesian}(\ref*{Prop:Cartesianii})}.\Par
	Let $P\colon (\E,K_\E)\to (\B,K_\B)$ be an internal fibration in \SITEs. The fact that $P$ is a fibration if $K_\B$ and $K_\E$ satisfy \tmp results directly from the two last squares in diagram \hyperref[eq:cartcomp]{(\ref*{eq:cartcomp})}.\Par
Now let $R=\{B_i\ra{f_i} B\mid i\in I\}$ be a $K_\B$-covering and $S=\{f_i^*E\ra{\overline{f_i}_E} E\}$ a lift of it in $P$. We build from the index set $I$ a category $\mathcal I$, the category generated by the graph $\{i\to*\mid i\in I\}$. We moreover endow it with the covering function $K_\I$ consisting of all identity coverings and the covering $\{i\to*\mid i\in I\}$. The coverings $R$ and $S$ give then rise to strict morphisms of sites $R\colon(\I,K_\I)\to(\B,K_\B)$ and $S\colon (\I,K_\I)\to(\E,K_\E)$. Let us denote by $\Delta(B)\colon\I\to\B$ the constant functor at $B$ and similarly for $\Delta(E)\colon\I\to\E$. They also determine strict morphisms of the sites just considered. Besides, one obtains a natural transformation $\alpha\colon R\Rightarrow\Delta(B)$ whose components are $\alpha_i=f_i$ and $\alpha_*=1_B$. Taking its lift $\bar\alpha_{\Delta(E)}\colon\alpha^*\Delta(E)\Rightarrow\Delta(E)$ at $\Delta(E)$ in the fibration $P_*$ of \hyperref[eq:P*]{(\ref*{eq:P*})}, one gets a strict morphism of sites $\alpha^*\Delta(E)\colon(\D,K_\D)\to (\E,K_E)$. Since $K_\I$ satisfies \tmp, by the preceding lemma, all components of $\bar\alpha_{\Delta(E)}$ are Cartesian morphisms in $P$. Therefore, the following diagram is commutative. 
$$\shorthandoff{;:!?}\xymatrixnocompile@!C@R=1.3cm{
 & E\ar[rrd]^{1_E}  & &\\
 f_i^*(E)\ar[ru]^{\overline{f_i}_E}\ar[rrd]^{\overline{f_i}_E} & \alpha^*\Delta(E)(*)\ar[rr]_{(\bar\alpha_{\Delta(E)})_{*}}\ar@{}[u]|*[@]{\cong} & & E\\
\alpha^*\Delta(E)(i)\ar[rr]_{(\bar\alpha_{\Delta(E)})_{i}}\ar@{}[u]|*[@]{\cong}\ar[ru] & & E\ar[ur]_{1_E} &\\
 & & &\\
 & B\ar[rr]^{1_B}& &B\\
B_i\ar[ur]^{f_i}\ar[rr]_{f_i} & & B\ar[ur]_{1_B} &
}$$
Axiom (L) leads to the conclusion.
\end{Pf}

\begin{Prop}\label{prop:KandKE}
	Let $\E\ra P(\B,K)$ be a fibred site. Let $K_P$ be the induced covering function on \E.
	\begin{enumerate}[(i)]
\item
If $K$ satisfies axiom \textnormal{(M)}, \tm, \tmp, \textnormal{(C)}, \textnormal{(L)} or \tl, so does $K_P$. The converse is true if $P$ is surjective on objects.
\item
If $(\B,K)$ is a site with pullbacks or, more strongly, satisfies \tc, so does $(\E,K_P)$. Moreover, $P$ is a morphism of sites with pullbacks. Conversely, if $P$ is surjective on objects and if $(\E,K_P)$ satisfies \tc, then so does $(\B,K)$. 
\item
If $K'$ is a covering function on \B such that $K\preceq K'$, then $K_P\preceq K'_P$. In particular, if $K'\equiv K$, then $K_P\equiv K'_P$. The converse of these implications are true if $P$ is surjective on objects.
\end{enumerate}
\end{Prop}
\begin{Pf}
	The proofs of (i) and (iii) are not difficult. The questions of pullbacks are a little trickier, and we therefore give their constructions. 
	
	Let $R_E=\{{f_i}^*E \ra{\overline{f_i}_E}E\mid i\in I\}$ be a $K_P$-covering of $E$ over the $K$-covering $R=\{B_i\ra{f_i}B\mid i\in I\}$ and $h\colon D\to E$ an arrow in \E. Suppose a pullback
$$\shorthandoff{;:!?}\xymatrix@!=1cm{
**[l]{A\times_BB_i}\ar[r]^{l_i}\ar[d]_{g_i} & B_i\ar[d]^{f_i}\\
A\ar[r]_{P(h)} & B
}$$
exists for each $i\in I$. Then, for each $i\in I$, a pullback of $\bar {f_i}_E$ along $h$ is given by the following square, where $\bar g_i$ is a Cartesian lift of $g_i$ at $D$ and $\bar l_i$ is induced by Cartesianness of $\bar f_i$.
$$\shorthandoff{;:!?}\xymatrix@!=1cm{
g_i^*D\ar[d]_{\bar g_i}\ar@{-->}[r]^{\bar l_i} & f_i^*E\ar[d]^{\bar {f_i}_E}\\
D\ar[r]_h & E
}$$
Conversely, given a $K$-covering $R=\{B_i\ra{f_i}B\mid i\in I\}$ and an arrow $f\colon A\to B$, consider a $K_P$-covering $R_E=\{{f_i}^*E \ra{\overline{f_i}_E}E\mid i\in I\}$ of an object $E$ over $B$ over $R$ ($P$ being surjective on objects by hypothesis). Let $\bar f_E\colon f^*E\to E$ a Cartesian lift of $f$ at $E$. Suppose there exists a pullback-covering of $R_E$ along $\bar f_E$ \emph{that is a $K_P$-covering}. Then, its image under $P$ is a pullback-covering of $R$ along $f$ that is a $K$-covering.
\end{Pf}

This is another example that shows the relevance of weakening the definition of a site with pullbacks (we only require the existence of pullback of coverings, not of all pullbacks). Indeed, $(\B,K)$ having pullbacks of coverings implies that $(\E,K_P)$ has such pullbacks. On the contrary, there seems to be no reason that \B having pullbacks implies that \E does (but I haven't looked for a counter-example). 

This construction can also be done in the sifted context. A sifted covering function $K$ on \B induces a sifted covering function $K_P$ on \E this way. For all $E\in\E$ and all coverings $S$ of $E$,
$$
S\in K_P(E)\iff \exists R\in K\text{ such that }S=P^{-1}(R)\cap\{h\text{ in }\E\mid\cod h=E\}.
$$
This definition is coherent with the non-sifted definition in the sense that given a plain covering function $K$ on \B, $\overline{K_P}=\overline{K}_P$.

\thref{prop:KandKE} applies especially in the situation of a fibration over a fibration $$\E\ra Q\D\ra P\B,$$ as described in \autoref{ssec:FiboverFib}. Indeed, given a covering function on \B, it provides the category \D with the structure of a site, making for instance the fibration $\E\ra Q\D$ suitable for the study of locally trivial objects. A general case of interest is the fibration of modules over the fibration of monoids $$\Mod(\E)\to\Mon(\E)\to\B$$ of a monoidal fibred category $\E\to\B$ (see \autoref{fibmodovfibmon}). This is how \emph{quasi-coherent\textnormal{ or }locally free sheaves of modules} are defined for instance (see \autoref{ssec:exloctriv}).


\chapter{Local triviality}\label{cha:LocTriv}
We now turn to the question of locally trivial objects in a fibred site (or dually in an opfibred cosite). Suppose one considers some objects in the base category as ``trivial''. Suppose that we also consider some objects in the fibres over these trivial objects as ``trivial''. For instance, in algebra, these might be (finitely generated) free modules over all rings. In bundle theory, topological or geometrical, with vector space or $G$-action structure for a group $G$ on fibres, the trivial objects are product bundles over all spaces. In algebraic geometry, in the fibration of sheaves of modules over locally ringed spaces, one considers affine schemes as trivial in the base and ``affine modules'' as trivial in the fibres. If the base category admits a terminal object, then trivial objects in fibres are often defined by inverse images of trivial objects in the fibre over the terminal.

Suppose now that the base category is endowed with a structure of a site. Then an object in the base is locally trivial if it can be covered by trivial objects. An object $E$ in the total category is locally trivial if there is a covering of the object in the base that it sits over such that the inverse images of $E$ over each arrow of the covering are trivial. 

We first define locally trivial objects in a site, in order to get some important tools and intuition for the more involved notion of locally trivial in a fibred site. One then study some properties and examples. In particular, we consider conditions that insure that the restriction of the fibration to locally trivial objects is again a fibration, which is an important property for application to $K$-theory.

We discovered the notion of locally trivial objects in Street's articles \cite{Str82,Str04}, where he uses it for the purpose of characterizing stacks by means of torsors. The article \cite{Str04} is a very important inspiration for our work. Street already notes in these articles that vector bundles are particular examples of his notion and also suggests an application to $K$-Theory. We explain below (\autoref{ssec:firstexloctriv}) how Street's notion of locally trivial object is a particular case of ours.

\section{Locally trivial objects in a site}

\subsection{Trivial objects}
\label{ssec:triv-objects-categ}

\begin{Def}
\label{def:subcattriv}
  Let \C be a category. A \emph{subcategory of trivial objects}\indexb{Trivial objects!-- in a category} is a replete subcategory $\Triv\subset\C$.
\end{Def}

In practice, subcategories of trivial objects are often defined via a functor $F\colon\A\to\C$ from an auxiliary category \A. This is for instance the case when trivial objects are determined by a ``free functor'', i.e., a left adjoint to a forgetful functor

\adjoint{\A}{\C}{F}{U.}

Now, a functor $F\colon\A\to\C$ determines several subcategories of \C, corresponding to the different possible factorizations of $F$ through the inclusion of a subcategory. Before studying these possible factorizations, we first recall some properties of functors.

\subsubsection{Conditions on functors}\label{ssec:conditions-functors}
We found some of these definitions in  \cite{Tay99,nLab}. The definition of an ``isofull'' functor is new.

\paragraph{Warning} When $F\colon\A\to\D$ is a functor, the notation
$$
F(A)\ra{F(f)}F(A')
$$
denotes, as usual in the category \D, an arrow that has domain $F(A)$ and codomain $F(A')$, and, specifically, an arrow that is the image of an arrow $f$ of \C. It does \emph{not} indicate that $A$ and $A'$ are the respective domain and codomain of $f$.

\begin{Defs}
  [Let $F\colon\A\to\D$ be a functor.]\label{def:conditions-functors}
	\item
$F$ is \emph{replete}\indexb{Functor!Replete --} if, given an isomorphism $g\colon F(A) \ra{\cong}C$ in \C, there exists an isomorphism $f\colon A\ra\cong A'$ in \A such that $g=F(f)$. A full and replete functor is sometimes called \emph{strictly full}. 
\item
$F$ is \emph{isofull}\indexb{Functor!Isofull --} if for each pair $(A,A')$ of objects of \A, the function
$$F_{A,A'}\colon\A(A,A')\to\C(F(A),F(A'))$$
is surjective on isomorphisms., i.e., for each isomorphism $g\colon F(A)\ra\cong F(A')\text{ in }\C$, there exists an arrow $f\colon A\to A'$ in \A such that $g=F(f)$. It is \emph{pseudomonic} if, moreover, it is faithful.
\item
$F$ \emph{reflects isomorphisms} if, $F(f)$ being an isomorphism implies that $f$ is isomorphism, for all $f$ in \A.
\item
A subcategory $\A\subset\C$ is \resp \emph{replete\textnormal{,} pseudomonic \textnormal{or} essentially wide},\index{Replete subcategory|see{Category}}\index{Category!Replete sub--}\index{Pseudomonic subcategory|see{Category}}\index{Category!Pseudomonic sub--}\index{Essentially wide subcategory|see{Category}}\index{Category!Essentially wide sub--} if its inclusion functor is \resp replete, isofull (and therefore pseudomonic) or essentially surjective. A full and replete subcategory is sometimes said to be \emph{strictly full}. 
\end{Defs}

\begin{Exs}
\item
By proposition \hyperref[Prop:Cartesian:liftiso]{\ref*{Prop:Cartesian}(\ref*{Prop:Cartesian:liftiso})}, fibrations and opfibrations are replete.
\item
The subcategory $\Sh(\C,K;\A)$ of sheaves on the site $(\C,K)$ with values in \A of the category $\PSh(\C;\A)$ of presheaves on \C with values in \A is full and replete.
\item
Recall that a commutative ring is \emph{local} if it has a unique maximal ideal. A homomorphism of local commutative rings $\phi\colon A\to B$ is a \emph{local homomorphism} if its associated map $\phi^{-1}\colon \Spec B\to\Spec A$ (that associates to a prime ideal $\frak p$ of $B$ the prime ideal $\phi^{-1}(\frak p)$ of $A$) preserves the maximal ideal. The subcategory of local commutative rings and local homomorphisms of rings is replete (and not full) in the category of commutative rings.
\end{Exs}

\begin{Rem}
	Pseudomonic functors reflect isomorphisms, as one readily verifies. This implies that a pseudomonic functor $F\colon\C\to\D$ respects the following stronger isofullness property. Given objects $C,D\in\C$, it induces a bijection between the set of isomorphisms from $C$ to $D$ and the set of isomorphisms from $F(C)$ to $F(D)$:
$$
F_{C,D}\colon \mathrm{Iso}_{\C}(C,D)\ra\cong \mathrm{Iso}_\D(F(C),F(D)).
$$
Moreover, recall that a functor is monic in \CAT \ssi it is faithful and injective on objects. Thus, the pseudomonic condition is an ``up to isomorphism'' weakening of the monic condition. 
\end{Rem}

\subsubsection{Notions of image of a functor}\label{ssec:factsyst}
Let $F\colon \A\to\C$ be a functor. There are different notions of image of $F$ available, depending on the \emph{factorization system} \cite{AHS90,nCatLab,BorI94} $(\fL,\fR)$ of the XL-category \CAT one considers. We denote $F(\A)$ the set-theoretic image of the functor $F$, which is not necessarily a subcategory of \C.

\paragraph{Image}
\label{par:image-factorization}\indexb{Functor!Image of a --}\index{Image!-- of a functor@(Replete, full) -- of a functor|see{Functor}}\index[not]{imF@\im F}

\fL is the class of extremal epimorphisms and \fR the class of monomorphisms \cite{BBP}. Recall that monomorphisms in \Cat are those functors that are injective on morphisms, or, equivalently, injective on objects and faithful. Recall moreover that extremal (= strong) epimorphisms in \CAT are the functors $G\colon \D\to\E$ such that every morphism in $\E$ is a composite of finitely many morphisms in $G(\D)$ (in particular, they are surjective on objects). Such a factorization of the functor $F$ can be obtained via the smallest subcategory of \C containing $F(\A)$, called the \emph{image} of $F$ and denoted \im F. Its objects are the objects of $F(\A)$ and its arrows are the finite composites of arrows in $F(\A)$.
 
\paragraph{Replete image}
\label{par:repl-image-fact}\indexb{Functor!Replete image of a --}\index[not]{rimF@\rim F}

\fL is the class of essentially surjective functors $G\colon\D\to\E$ such that every morphism in \E is a finite composite of isomorphisms in \E and morphisms in $G(\D)$ and \fR the class of replete monomorphisms (I do not know any reference for this system). One obtains such a factorization of $F$ via the smallest replete subcategory of $\C$ containing $F(\A)$. Its objects are all objects of \C isomorphic to some object of F(\A) and its morphisms are all finite composite of isomorphisms in \C and morphisms in F(\A). This subcategory is called the \emph{replete image} of $F$ and is denoted $\rim F$.

\paragraph{Full image}
\label{par:full-image-fact}\indexb{Functor!Full image of a --}\index[not]{fimF@\fim F}\index{Category!Associated full sub--}\index[not]{fA@\full\A}

\fL is the class of functors surjective on objects and \fR the class of full monomorphisms \cite{nCatLab}. One obtains such a factorization via the full subcategory of \C whose objects are the objects of $F(\A)$, called the \emph{full image} of $F$ and written \fim F. In case $F$ is the inclusion of a subcategory $\A\subset\C$, then \fim F is called the \emph{associated full subcategory} of \A and denoted $\full\A$.

\paragraph{Replete full image}
\label{par:replete-full-image}\indexb{Functor!Replete full image of a --} \index[not]{rfimF@\rfim F}\indexb{Category!Full repletion of a sub--}\indexb{Repletion@(Full) repletion|see{Category}}

\fL is the class of essentially surjective functors and \fR the class of strictly full monomorphisms \cite{nCatLab}. One obtains such a factorization via the full subcategory of \C whose objects are all objects that are isomorphic to an object of $F(\A)$. This is called the \emph{replete full image} and denoted \rfim F. When $F$ is the inclusion of a subcategory $\A\subset \C$, \rfim F is called the \emph{full repletion} of \A and is denoted $\frep\A$.

\begin{Rems}
[Let $F\colon\A\to\C$ be a functor.]\label{rem:replete-full-image}
\item
If $F$ is isofull, its image and replete image take a much simpler form. Indeed, under this hypothesis, $\im F=F(\A)$ (this is also true if $F$ is monic on objects). Moreover,
the replete image of $F$ is the subcategory of \C consisting in all objects isomorphic to some object of $F(\A)$, and all morphisms $h\colon C\to C'$ that admit a factorization
\begin{equation}\label{eq:repimage}\begin{aligned}\xymatrix{
C\ar[d]_\cong \ar[r]^h & C' \ar[d]^\cong\\
F(A)\ar[r]_{F(f)} & F(A').
}\end{aligned}\end{equation}
If a morphism $h$ in \C admits such a factorization, then, thanks to isofullness, given any choice of isomorphisms $C\cong F(A)$ and $C'\cong F(A')$, there exists a factorization \eqref{eq:repimage} and the morphism $f$ can be chosen such that its domain is $A$ and its codomain is $A'$.

\item
The different notions of image of a functor $F$ are coherent in the sense that:
$${}_r(\fim F)={}_f(\rim F)=\rfim F.$$

\item
If $F$ is replete, then $\rim=\im F$ and $\rfim F=\fim F$. If $F$ is full, then $\fim F=\im F=F(\A)$ and $\rfim F=\rim F$. If $F$ is full and replete, then $\rfim F=\fim F=\rim F=\im F=F(A)$.

\item
The subcategory \rim F (resp.\ $\rfim F$) of \C is replete (resp.\ strictly full). Moreover, the replete (full) image of $F$ is equal to the (full) repletion of $\im F$. Therefore, the notations for functors and subcategories are coherent. Moreover, \rim F is a pseudomonic subcategory of \C and if $F$ is isofull, then $\im F$ also.

\item \label{rem:replete-full-image:equ}Given the factorization $\xymatrix@1{\A\ar@{->>}[r]^-{\bar F} & \,\rim F\, \ar@{^{(}->}[r] & \C}$ of the functor $F$, the functor $\bar F$ is an equivalence \ssi $F$ is pseudomonic \cite{nLabR}.
\end{Rems}

Since we want trivial objects to be defined up to isomorphisms, we have two candidates for defining the subcategory of trivial objects of \C associated to a functor $F\colon\A\to\C$: the replete and the replete full image. The former has be chosen when we care about morphisms of \A, whereas the latter is usualy chosen when we just want to remember the objects. In both cases, one gets a factorization:
\[\xymatrix@1{\A\ar@{->>}[r]^-{\bar F} & \,\Triv\, \ar@{^{(}->}[r] & \C}.\]

\begin{Exs}
[Here are a few examples of functors determining subcategories of trivial objects via the replete full image.]
	\item Free $R$-modules and free $R$-modules of finite type are determined respectively by functors $\Set\to\Mod_R$ and $\N\to\Mod_R$.
	\item Trivial bundles over a space $X$ are determined by the ``cofree bundle'' functor $$\Top\to\Top/X.$$
\end{Exs}

\subsection{Locally trivial objects}
\label{ssec:locally-triv-objects}
In this part, we define the full subcategory of locally trivial objects in a site with trivial objects. In some situations, one could be interested in non-full subcategories of locally trivially objects. Indeed, when the subcategory of trivial objects in not full, one might consider morphisms that are locally trivial in some sense. Yet, we will not study this topic here (see discussion in \autoref{sssec:morph-locally-triv}, page \pageref{sssec:morph-locally-triv} though).

\begin{Defs}
[\indexb{Locally trivial objects!-- in a site}Let $\Triv\subset\C$ be a subcategory of trivial objects.] 
\item An object $C\in\C$ is \emph{locally trivial for a covering $R$} (with respect to \Triv) if $R$ covers $C$ and has trivial domains.
\item Suppose now that $(\C,K)$ is a site that satisfies axiom \tmp. An object $C\in C$ is \emph{locally trivial} (in the site $(C,K)$ with respect to \Triv) if it is locally trivial for a $K$-covering. The \emph{subcategory of locally trivial objects}, denoted $\Loc=\Loc(P,\Triv,K)$ is the full subcategory of \C consisting of locally trivial objects.\\
There is an induced covering function on \Loc, denoted $K_{l}$, whose coverings are $K$-coverings with trivial domains.
\end{Defs}

\begin{Rem}
  The category \Loc is itself replete if $K$ satisfies either \tm and \tl, or \tc.
\end{Rem}

\begin{Exs}
\item Euclidean spaces are spaces homeomorphic to some $\R^{n}$. Let us consider the full subcategory $\mathit{Euc}$ of Euclidean spaces in the category \Top of topological spaces. When \Top is endowed with the pretopology of open subset coverings, locally trivial objects are called \emph{locally Euclidean spaces}.
\item For this example, we anticipate some material that is detailed further down. There is a contravariant full functor $(\Spec,\Oc{})\colon\Comm^{op}\to\LRinged$ from commutative rings to locally ringed spaces (see \autoref{ssec:Firstexmorfib}). The replete (full) image of this morphism is the category of affine schemes. When one put on \LRinged the pretopology of open subringed spaces coverings, then the full subcategory of locally trivial objects is the category of \emph{schemes}.
\item There also is subcategory of $\LRinged$ consisting of locally ringed spaces isomorphic to some $(\R^{n},\C^{\infty}_{n})$, where $\C^{\infty}_{n}$ is the sheaf of smooth real functions on $\R^{n}$. When considering the same covering function as in the preceding example, the subcategory of locally trivial objects is the category of \emph{smooth manifolds}.
\end{Exs}
Properties of the covering function $K_{l}$ happen to be crucial in the following developments. There are situations where most properties of $K$ pass to $K_l$. This happens when every $K$-covering of a locally trivial object admits a refinement in $K_l$. Indeed, in this situation, one readily verifies that if $K$ satisfies either (C) or (L), then so does $K_{l}$. The next lemma asserts that under some mild assumptions on $K$, it is actually enough to check this condition on coverings of trivial objects of the base.

\begin{Lem}\label{lem:reftrivbase}
	Suppose $K$ is a pretopology satisfying axiom ($\tilde{L}$).\Par
	If every $K$-covering of trivial objects of \C admits a refinement in $K_l$, then this is also true for every $K$-covering of locally trivial objects.
\end{Lem}
\begin{Pf}	
	Now, let $C\in\Loc$ and $R$ a $K$-covering of $B$. Then, $B$ admits a $K_l$-covering $S$. Let $g\colon \tilde B\to B$ in $S$. By axiom (C), there exists a $K$-covering $S'_g$ of $\tilde B$ such that the composite covering $g\circ S'_g$ refines $R$. Moreover, by definition of a $K_l$-covering, ${\tilde B\in\Triv_b}$. Therefore, by hypothesis, there exists a $K_l$-covering $S_g$ of $\tilde B$ that refines $S'_g$. Thus, $\bigcup_{g\in S}(g\circ S_g)$ refines $\bigcup_{g\in S}(g\circ S'_g)$, which refines $R$. Consequently, $\bigcup_{g\in S}(g\circ S_g)$ refines $R$. Furthermore, since $K$ satisfies ($\tilde{\mathrm L}$), $\bigcup_{g\in S}(g\circ S_g)$ is a $K$-covering, and thus a $K_l$-covering.
	
	\cqfd
\end{Pf}

\begin{Ex}\label{ex:spec}
We consider locally trivial objects in the category \Top in the pretopology of open subset coverings, with respect to the subcategory of trivial objects induced by the \emph{prime spectrum functor} $$\Spec\colon\Comm\to\Top.$$ Recall that for a commutative ring $A$, the space $\Spec A$ is the set of all prime ideals of the commutative ring $A$ topologized by defining the closed subsets to be the subsets of the form $V(\frak a)=\{\frak p\in\Spec A\mid\frak a\subset\frak p\}$ for any ideal $\frak a$ of $A$. A basis for this topology is given by the open subsets $D(f)=\Spec A-V((f))$ for each $f\in A$. Given a homomorphism of commutative rings $\phi\colon A\to B$, the continuous function $\Spec \phi\colon\Spec B\to\Spec A$ is defined by $\Spec\phi(\frak p)=\phi^{-1}(\frak p)$. Locally trivial objects in this situation are thus topological spaces that can be covered by open subsets homeomorphic to the prime spectrum of a commutative ring.

Not every space can be covered by open subsets homeomorphic to prime spectra of commutative rings. Indeed, spectra of commutative rings have the nice property of being compact (and not every space is locally compact). On the other hand, prime spectra have bad separability properties: the prime spectrum of a commutative ring $A$ is almost never Hausdorff, not even T1, since the closed points are precisely the maximal ideals of $A$.

Every $K$-covering of a space homeomorphic to the prime spectrum $\Spec A$ of a ring $A$ admits a $K_l$-refinement. This is due to the fact that the topology on $\Spec A$ is generated by the basis elements $D(f)$, $f\in A$, and that, for each $f\in A$, $D(f)\cong \Spec A_f$ \cite{Qin02}. Note that the pretopology $K$ satisfies the hypotheses of \thref{lem:reftrivbase}. Therefore, $K_{l}$ safisfies axioms (C) and (L). 
\end{Ex}

\section{Trivial objects in a fibration}
Let $P\colon\E\to\B$ be a fibration. We want to define an analog of a subcategory of trivial objects in that context. As usual in mathematics, the notion of triviality should be isomorphism-independent\footnote{In a homotopy theoretic framework, one could work with a notion that is defined up to weak equivalence. We do not consider this framework in our thesis.}. Since \E is a fibred category, there are two notions of isomorphisms available: isomorphisms and vertical isomorphisms. When we need to distinguish isomorphisms in \E from vertical ones, we call the former \emph{global isomorphisms}. 

What kind of isomorphisms should the notion of triviality be independent of? This will be decided by what \E is designed for. If it is mainly a tool for studying \B, then one should consider vertical isomorphisms. On the other hand, if one wants to study \E, or its subcategory of (locally) trivial objects, for its own sake, then one should opt for global isomorphisms. For example, affine schemes are trivial objects defined up to global isomorphism, because one is interested in (affine) schemes for their own sake (or for the study of commutative rings), not for studying the topological spaces they sit over. On the other hand, one studies vector bundles up to vertical isomorphism because one is mainly interested in them as a tool for probing the space they live over.

\begin{Def}
  A subfunctor $Q\subset P$
\[\xymatrix@=1.3cm{
\D\ar@{^{(}->}[r]\ar[d]_{Q} & \E\ar[d]^{P}\\
\A \ar@{^{(}->}[r] & \B,
}\]
is \emph{replete (resp.\ globally replete)} if each restriction $\D_{A}\subset\E_{A}, A\in\A$, is a replete subcategory (resp.\ if both $\D\subset\E$ and $\A\subset \B$ are replete subcategories).
\end{Def}
Note that global repleteness implies repleteness.

\begin{Def}
  A \emph{subfunctor of trivial objects}\index{Trivial objects!Subfunctor of --} in a fibration $P$ is a replete subfunctor
\[\xymatrix@=1.3cm{
\Triv_{t}\ar@{^{(}->}[r]\ar[d] & \E\ar[d]^{P}\\
\Triv_{b} \ar@{^{(}->}[r] & \B.
}\]
Objects of $\Triv_{b}$, resp.\ $\Triv_{t}$ are called \emph{trivial objects of the base}, resp.\ \emph{of the total category}, or just \emph{trivial objects} if the level where they are is clear from the context.
\end{Def}
Note that when $\Triv$ is a \emph{globally} replete subfunctor, then both $\Triv_{b}\subset\B$ and $\Triv_{t}\subset\E$ are subcategories of trivial objects as defined in \thref{def:subcattriv}. When \Triv is only replete, then it is a fibrewise subcategory of trivial objects, in the sense that for each $B\in\Triv_{b}$, $(\Triv_{t})_{B}\subset\E_{B}$ is a subcategory of trivial objects.

Recall that trivial objects in a category \C are often defined via a functor from an auxiliary category \A. Similarly, trivial objects in a fibration are often naturally defined ``externally'', by means of a morphism of fibrations:
\begin{equation*}
 \xymatrix@=1.3cm{
\D \ar[r]^F \ar[d]_{Q}& \E\ar[d]^P\\
\A \ar[r]_G & \B.
} 
\end{equation*}

The XL-category \CAT, and thus also its arrow XL-category $\CAT^{\two}$, have good factorization system features and therefore, it makes in there no difference to define trivial objects internally or externally. On the other hand, the XL-category \FIB does not have such good factorization properties. For that reason, it is sometimes necessary to define trivial objects externally, if one desires having a fibration of trivial objects, not a mere functor.

We first extend the conditions on functors studied in \autoref{ssec:conditions-functors} to morphisms in $\CAT^{\two}$. To each condition on functors, there corresponds a ``global'' and a ``fibred'' notion for commutative squares of functors.  We then study two important examples, coming from algebraic geometry. Next, we study factorization systems in $\CAT^{\two}$ and \FIB. This will lead us to categories and fibrations of trivial objects induced by a morphism of fibrations. Finally, we describe classes of examples that will have applications later.

\subsection{Conditions on morphisms of fibrations}

We now consider the conditions on morphisms in $\CAT^{\two}$ corresponding to conditions  on functors in \thref{def:conditions-functors}. But now one has two variants for each notion, the ``global'' notion, that is defined componentwise, and the ``fibred'', that is defined fibrewise.

\begin{Defs}
	[{Consider a morphism $(F,G)\colon Q\to P$ in $\CAT^{\two}$
	\begin{equation*}\begin{aligned}\xymatrix@=1.3cm{
\D \ar[r]^F \ar[d]_{Q}& \E\ar[d]^P\\
\A \ar[r]_G & \B,
}\end{aligned}\end{equation*}}]\label{def:fibredversion}
\item
$(F,G)$ is \emph{replete} if $F$ is fibrewise replete, i.e., each restriction $F_A\colon\D_A\to\E_{G(A)}$, $A\in\A$, is replete. In detail, this means that if $k\colon F(D)\ra\cong E$ is a \emph{vertical} isomorphism in \E, then there exists a \emph{vertical} isomorphism $h\colon D\ra\cong D'$ in \D such that $k=F(h)$. This is called the \emph{fibred notion} of repleteness. The morphism $(F,G)$ is \emph{globally replete} if $F$ and $G$ are replete (as functors). This is called the \emph{global notion} of repleteness.
\item
Similarly, $(F,G)$ is \emph{full\textnormal{,} faithful\textnormal{,} isofull\textnormal{,} pseudomonic\textnormal{,} monic on objects\textnormal{,} monic on morphisms \textnormal{or} essentially surjective} if $F$ is respectively fibrewise so. It is \resp \emph{globally} so if $F$ and $G$ are respectively full, faithful, isofull, pseudomonic, monic on objects, monic on morphisms or essentially surjective (as functors). For each of these notions, one thus has the \emph{fibred\textnormal{ and the }global notion}.
\end{Defs}
When the morphism is not only in $\CAT^{\two}$ but in \FIB, the global and fibred notions coincide under some hypotheses.
\begin{Lem}\label{lem:repfullglob}
Let $(F,G)\colon(\D\ra Q\A)\to(\E\ra P\B)$ be a morphism of fibrations.
	\begin{enumerate}[(i)]
        \item Suppose $(F,G)$ is globally faithful, then it is faithful.
	\item\label{lem:repfullglob:globfib}
Suppose that $G$ is faithful. Then, $(F,G)$ is full, isofull, or pseudomonic, if it is globally so.
\item Suppose that $G$ is monic and $(F,G)$ globally replete. Then, $(F,G)$ is replete. 
\item\label{lem:repfullglob:ess} Suppose that $G$ is isomorphism reflecting isofull. Then $(F,G)$ is essentially surjective if it is globally so.
\item\label{lem:repfullglob:i} If $(F,G)$ and $G$ are \resp replete, monic on objects or essentially surjective, then $(F,G)$ is globally so. If $(F,G)$ is isofull and $G$ is isomorphism reflecting isofull, then $(F,G)$ is globally isofull.
\item\label{lem:repfullglob:cart} Suppose that $(F,G)$ is Cartesian. If $(F,G)$ and $G$ are \resp isofull, full, faithful, pseudomonic or monic on morphisms, then $(F,G)$ is globally so.
\item Suppose that $Q$ is surjective on objects. If $F$ is \resp replete or replete and isofull, then $(F,G)$ is globally so.
\item Suppose that $P$ is surjective on objects. If $F$ is essentially surjective, then $(F,G)$ is globally so.
\end{enumerate}
\end{Lem}
\begin{Pf}
Recall first the important property of Cartesian arrows (\thref{Prop:Cartesian}) that every isomorphism is Cartesian and every Cartesian lift of an isomorphism is also an isomorphism. This is why the Cartesianness is not required for \ref{lem:repfullglob:i}, but it is for \ref{lem:repfullglob:cart}. We prove \ref{lem:repfullglob:i} for repleteness, the other points being similar. 

Let $k\colon E\ra\cong F(D)$ be an isomorphism in \E. Define $$A\coloneqq Q(D)\quad\text{and}\quad (g\colon B\ra\cong G(A))\coloneqq P(k).$$
	One can follow the proof on the diagram below. Since $G$ is replete, there is an isomorphism $f\colon A'\ra\cong A$ in \A such that $g=G(f)$. Let $\bar f_D\colon f^*D\ra{\cong} D$ be a Cartesian lift of $f$ at $D$ in \D. It is an isomorphism, since it lifts an isomorphism. Now, $F(\bar f_D)\colon F(f^*D)\to F(D)$ and $k\colon E\ra\cong F(D)$ are two Cartesian lifts of the same arrow $g=G(f)$, because they are both isomorphisms. In consequence, there exists a vertical isomorphism $\bar k$ making commute the triangle in the following diagram.
	$$
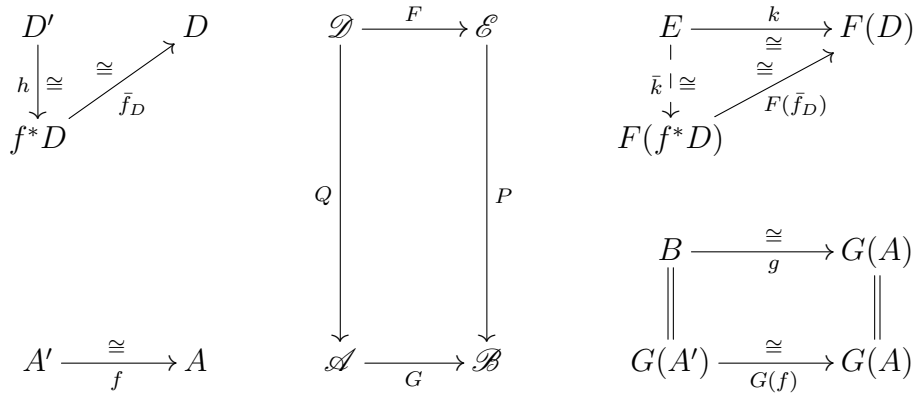
\xymatrix@C=1.3cm{
	D'\ar[d]_h^\cong			& D	& \D\ar[ddd]_Q\ar[r]^F	& \E\ar[ddd]^P	& E\ar[r]_\cong^k\ar@{-->}[d]_{\bar k}^\cong	& F(D)\\
	f^*D\ar[ru]_{\bar f_D}^\cong 	&		&					& 			& F(f^*D) \ar[ur]_{F(\bar f_D)}^\cong 			&\\
				 			&  		& 					& 			& B\ar@{=}[d]\ar[r]^\cong_g 			& G(A)\ar@{=}[d]\\
	A'\ar[r]^\cong_f				& A		& \A \ar[r]_G			& \B			&G(A')\ar[r]^\cong_{G(f)} 				& G(A)
	}$$
	Now, $(F,G)$ is supposed to be replete, and therefore, there exists a vertical isomorphism $h\colon D'\to f^*D$ in \D such that $F(h)=\bar k$. Finally, $k=F(\bar f_D\circ h)$.
\end{Pf}

\begin{Rems}
	\item
	One should not be misled by the word ``global'': there seem to be no reason in general that a global notion implies its fibred notion (yet, one should provide a counterexample).
	\item
	In the situation where $G=Id$, i.e., for a functor over \B, \thref{lem:repfullglob} shows first that the global notions implies the fibred. Consequently, it also shows that there is no distinction between being a replete as a functor over \B and replete as a mere functor. The two very close terminologies do therefore not conflict. This is also true for the isofull and essentially surjective conditions. Moreover, for Cartesian functors over \B, there is no distinction between the global and the fibred notion for all the notions introduced.
	
	For non Cartesian functors over \B, the notions of isofull functor over \B and isofull functor might differ. In order to avoid confusion, we will insist that a functor over \B is isofull as a functor over \B, and not as a mere functor, by saying it is \emph{fibrewise isofull}. The same applies to the full, faithful and pseudomonic conditions.
\end{Rems}

We can now define more or less strong notions of \emph{subfibration} of a fibration $P\colon\E\to\B$. The strength of the notion partly depends on the ambient category of $P$ one considers: \FIB, \FIB(\B), \FIBc or $\FIBc(\B)$. 

\begin{Defs}
[Let $P\colon\E\to\B$ be a fibration.]
	\item The weakest notion we consider, called \emph{subfibration}, consists of a morphism of fibration $(I,J)$ with codomain $P$ whose functors $I$ and $J$ are inclusions of subcategories:
$$\xymatrix{
\D \ar@{^(->}[r]^I \ar[d]_{\bar P}& \E\ar[d]^P\\
\A \ar@{^(->}[r]_J & \B,
}$$
The functor $I$ is called the \emph{total inclusion functor}, whereas $J$ is called the \emph{base inclusion functor}. We do not impose that the morphism $(I,J)$ be Cartesian. In other words, the Cartesian morphisms of $\bar P$ are not necessarily Cartesian morphisms of \E that belong to \D. 
	\item A \emph{subfibration over \B} is a subfibration whose base inclusion is the identity $$Id_{\B}\colon\B\to\B.$$ 
	\item
A \emph{strong subfibration \op over \B \fp}is a subfibration (over \B) that is a Cartesian morphism of fibrations (over \B).
\item
A subfibration is \emph{\op globally\fp replete\textnormal{,} \op globally\fp isofull\textnormal{,} \op globally\fp pseudomonic\textnormal{,} \op globally\fp full}   if its morphism of fibrations $(I,J)$ is respectively so. It is \emph{\op globally\fp essentially wide} if $(I,J)$ is (globally) essentially surjective. A (globally) replete full subfibration is called (globally) \emph{strictly full}.
\end{Defs}

\begin{Cor}\label{cor:subfib}
Let $(I,J)\colon(\D\to\A)\to(\E\ra P\B)$ be a subfibration.
\begin{enumerate}[(i)]
\item
The isofullness and pseudomonicity properties are equivalent. Idem for their global versions.
\item
$(I,J)$ is replete, resp.\ full, faithful, isofull, or pseudomonic, if it is globally so.
\item If $\A$ is a pseudomonic subcategory of \B, then $(I,J)$ is essentially wide if it is globally so.
\item
If $(I,J)$ is replete, then it is pseudomonic. Idem for the global notions.
\item\label{cor:subfib:strong}
If $(I,J)$ is strong and (globally) replete, then any Cartesian lift in \E of an arrow in $\A$ at an object $D$ of \D is in \D. In particular, it is a Cartesian lift in the subfibration.
\end{enumerate}
\end{Cor}
\begin{Pf}
  (i) to (iii) are obvious.
  (iv) The hypothesis implies that $(I,J)$ is (globally) replete and monic, and thus (globally) isofull.
	
  (v) Let $A'\ra f A$ in \A and $D$ over $A$. Let $f^*D\ra{\bar f_D}D$ be a Cartesian lift of $f$ at $D$ in the subfibration. It is also Cartesian in the ambient fibration $P$, since $(I,J)$ is strong. Consider now a Cartesian lift $h\colon E\to D$ of $f$ at $A$ in the fibration $P$. Then, there exists a vertical isomorphism $k\colon E\ra\cong f^*D$ in \E such that $\bar f_d\circ k=h$. Since $\D$ is replete, the isomorphism $k$ is in  \D. Finally, $h$ is in \D as a composite of arrows in \D, and is Cartesian in the subfibration, since it is the composite of a Cartesian arrow in the subfibration and an isomorphism in \D.
\end{Pf}

\subsection{First examples}\label{ssec:Firstexmorfib}

\subsubsection{Locally ringed spaces}\label{sssec:locringspaces} Recall that a ringed space $(X,\Oc X)$ is a \emph{locally ringed space}\indexb{Ringed space!Locally --} if each stalk $\Oc {X,x}$, $x\in X$, is a local ring. Recall that the stalk functor at $x$ is defined only up to isomorphism (being defined as a left adjoint, or more concretely, as a colimit functor). Yet, the definition of a locally ringed space is meaningful because the subcategory of local commutative rings is replete in the category of commutative rings. 
	
	Consider now a morphism of ringed spaces $(f,f^\sharp)\colon(X,\Oc X)\to (Y,\Oc Y)$. Let 
	\begin{equation}\label{eq:adjlocring}
	f_\sharp\colon f^{-1}\Oc Y\to\Oc X
	\end{equation}
	be the transpose map of $f^\sharp\colon\Oc Y\to f_*\Oc X$ under some adjunction $f^{-1}\dashv f_*$. This is the only morphism of sheaves on $X$ making the following diagram in \Ringed commutes
	$$\xymatrix{
(X,\Oc X)\ar@{-->}[d]_{(1_X,f_\sharp)}\ar[r]^{(f,f^\sharp)} & (Y,\Oc Y)\\
(X,f^{-1}\Oc Y)\ar[ur]_{(f,\eta^f_{\Oc Y})},&
}$$
where $\eta^f$ is the unit of the adjunction $f^{-1}\dashv f_*$ (by \thref{prop:bifibadjoint}). The stalk at $x\in X$ of $f^{-1}\Oc Y$ has a simple form, thanks to Example \hyperref[ex:invimsheaves:stalk]{\ref*{ex:invimsheaves}(\ref*{ex:invimsheaves:stalk})}. Indeed, recall that taking the stalk at $x$ is, up to the equivalence $\Sh(*;\A)\simeq\A$, the functor $(i_x)^{-1}$ of the map $i_x\colon*\hookrightarrow X$ such that $i_x(*)=x$. Thus ${i_x}^{-1}\circ f^{-1}$ is left adjoint to $f_*\circ (i_x)_*=(f\circ i_x)_*=(i_{f(x)})_*$. The functor ${i_x}^{-1}\circ f^{-1}$ is therefore just the stalk functor $Stalk_{f(x)}=(i_{f(x)})^{-1}$ at $f(x)$. Now, given two locally ringed spaces $(X,\Oc X)$ and $(Y,\Oc Y)$, a morphism of ringed spaces $(f,f^\sharp)$ between them is a \emph{morphism of locally ringed spaces}\indexb{Ringed space!Morphisms of locally --s} if the homomorphism of rings $$(f_\sharp)_x\colon\Oc{Y,f(x)}\to\Oc{X,x}$$ is local%
\footnote{
In algebraic geometry books, this homomorphism is defined as follows: apply the functor $Stalk_{f(x)}$ (defined via colimits) to the morphism $f^\sharp\colon\Oc Y\to f_*\Oc X$ to obtain an homomorphism of rings $$\Oc{Y,f(x)}\ra{(f^\sharp)_{f(x)}}(f_*\Oc X)_{f(x)}.$$ Then it is said that there is a homomorphism $(f_*\Oc X)_{f(x)}\to\Oc{X,x}$, and one consider the composition of this homomorphism with the previous one. That this map exists already at the level of sheaves is never mentioned to my knowledge. Indeed, it seems not to make sense at the level of sheaves, since $f_*\Oc X$ and $\Oc X$ do not live in the same category. Yet, as explained above $(f_*\Oc X)_{f(x)}=(f^{-1}f_*\Oc X)_x$. So, at the level of sheaves, one is really looking at the composite $$f^{-1}\Oc Y\ra{f^{-1}(f^\sharp)}f^{-1}f_*\Oc X\ra{\epsilon_{\Oc X}}\Oc X,$$ where $\epsilon$ is the counit of the adjunction $f^{-1}\dashv f_*$. In other words, this composite is the transpose morphism of the morphism $f^\sharp$ under the latter adjunction. This fact is maybe obvious to the authors, but we think it is worth mentioning, being useful both conceptually and for calculations (as we shall see). 
}.
This definition does also not depend on the choice of the stalk because of the repleteness of the subcategory of local commutative rings. It does not depend either on the choice of the left adjoint $f^{-1}$ of $f_*$, because the transpose morphisms \hyperref[eq:adjlocring]{(\ref*{eq:adjlocring})} under two different adjunctions are isomorphic over $\Oc X$.

Now locally ringed spaces and their morphisms form a strong subfibration over \Top of the fibration of ringed spaces.  First of all, they form a subcategory of \Ringed. It is clearly closed under identities and we briefly explain the closure under composition. Consider a composable pair
$$
(X,\Oc X)\ra{(f,f^\sharp)}(Y,\Oc Y)\ra{(g,g^\sharp)}(Z,\Oc Z).
$$
Its composite has first component $g\circ f$ and second component $$\Oc Z\ra{g^\sharp}g_*\Oc Y\ra{g_*(f^\sharp)}g_*f_*\Oc X=(g\circ f)_*\Oc X.$$
A left adjoint to $(g\circ f)_*$ is given by $f^{-1}\circ g^{-1}$ where $f^{-1}$ and $g^{-1}$ are left adjoints of $f_*$ and $g_*$ respectively. The transpose morphism of $g_*(f^\sharp)\circ g^\sharp$ under this adjunction is $$f^{-1}g^{-1}\Oc Z\ra{f^{-1}(g_\sharp)}f^{-1}\Oc Y\ra{f_\sharp}\Oc X,$$ where $f_\sharp$ and $g_\sharp$ are the transpose morphisms of $f^{\sharp}$ and $g^\sharp$. The stalk at $x$ of this morphism is thus the following homomorphism 
\begin{equation}\label{eq:compositestalk}
\Oc{Z,g\circ f(x)}\ra{(g_\sharp)_{f(x)}}\Oc{Y,f(x)}\ra{(f_\sharp)_x}\Oc{X,x}.
\end{equation}

We denote \LRinged\index[not]{LRinged@\LRinged} this category. The functor $\Ringed\to\Top$ restricts to $$\LRinged\to\Top.$$ 

This is a strong subfibration over \Top. Indeed, let $f\colon X\to Y$ be a continuous map and $(Y,\Oc Y)$ a locally ringed space over $Y$. Recall that once we have chosen an adjunction $f^{-1}\dashv f_*$, we obtain a Cartesian lift of $f$ at $(Y,\Oc Y)$ in \Ringed (see \hyperref[eq:invimAspace]{(\ref*{eq:invimAspace})}). It is defined by 
\begin{equation}\label{eq:Cartliftringed}
(X,f^{-1}\Oc Y)\ra{(f,\eta^f_{\Oc Y})}(Y,\Oc Y),
\end{equation}
where $\eta^f\colon Id_{\Sh(Y;\Comm)}\Rightarrow f_*f^{-1}$ is the unit of the adjunction $$f^{-1}\dashv f_*\colon\Sh(Y;\Comm)\to\Sh(X;\Comm).$$ Now, if $(Y,\Oc Y)$ is locally ringed, then its inverse image along the Cartesian arrow \hyperref[eq:Cartliftringed]{(\ref*{eq:Cartliftringed})} is also, because the stalk at $x$ of $f^{-1}\Oc Y$ can be chosen to be $\Oc{Y,f(x)}$. Moreover, the Cartesian morphism \hyperref[eq:Cartliftringed]{(\ref*{eq:Cartliftringed})} is in \LRinged because the adjoints of components of the unit of an adjunction are identities. We still have to prove this is a Cartesian arrow in the subfibration. Consider the following situation, where $(h,h^\sharp)$ is a morphism of locally ringed spaces.
$$\xymatrix@=1.3cm{
(Z,\Oc Z)\ar@{-->}[r]_-{(g,h_{\sharp f})}\ar@/^3pc/[rr]^{(h,h^\sharp)} & (Y,f^{-1}\Oc X)\ar[r]_{(f,\eta^f_{\Oc X})} & (X,\Oc X)\\
Z\ar@/^-3pc/_h[rr]\ar[r]_g & Y\ar[r]_f & X
}$$
The dotted arrow is the induced arrow in the indexed category \Ringed. The morphism $h_{\sharp f}\colon f^{-1}\Oc Y\to g_{*}\Oc Z$ is the transpose under $f^{-1}\dashv f_*$ of the morphism $h^\sharp$. The transpose of $h_{\sharp f}$ under $g^{-1}\dashv g_*$ is the transpose $h_\sharp$ of $h^\sharp$ under the composite adjunction $g^{-1}\circ f^{-1}\dashv f_*\circ g_*=h_*$. Since $(h,h^\sharp)$ is a morphism of locally ringed space, the stalk at $z\in Z$ of $h_\sharp$ is local. In conclusion, $(g,h_{\sharp f})$ is a morphism of locally ringed spaces.

This subfibration is replete. It is therefore globally so and (globally) pseudomonic, by \thref{cor:subfib}. In order to prove it is replete, let $$(1_X,(1_X)^\sharp)\colon(X,\Oc X)\ra{\cong}(X,\Oc X')$$ be a vertical isomorphism in \Ringed, with $(X,\Oc X)$ locally ringed. Since $(1_X,(1_X)^\sharp)$ is an isomorphism in \Ringed, $(1_X)^\sharp\colon\Oc X'\to\Oc X$ is an isomorphism in $\Sh(X;\Comm)$. Therefore, its stalk at a point $x\in X$ is an isomorphism of rings. By repleteness of the category of local commutative rings, both the ring $\Oc{X,x}'$ and the isomorphism ${(1_X)^\sharp}_x$ are local. Moreover, the adjunction over the identity map can be chosen to be the identity adjunction $(1_X)^{-1}=Id\dashv Id=(1_X)_*$. Thus, we have also shown that $(1_X,(1_X)^\sharp)$ is a morphism of locally ringed spaces. Note that the result that \LRinged is a strong subfibration can now be slightly strengthened, thanks to repleteness of \LRinged. For, any Cartesian lift of a continuous map at a locally ringed space in \Ringed is a Cartesian lift in \LRinged.

This subfibration is not full though. For an example of a morphism of ringed spaces between locally ringed spaces that is not in \LRinged, see for instance \cite[Example 2.3.2]{Har77}. Moreover, $\LRinged\to\Top$ does not seem to be an opfibration. In any case, it is not a sub-opfibration of $\Ringed\to\Top$. Indeed, let $A\in\Comm$ be a \emph{non local} commutative ring. Then, the associated affine scheme $(\Spec A,\Oc A)$ is a locally ringed space, as it is for any ring. But its direct image over $\Spec A\to *$ is the global section ring $\Oc A(\Spec A)\cong A$. This is not a locally ringed space over $*$, by hypothesis.

\subsubsection{Affine schemes}\label{sssec:affsch}\index{Scheme! Affine --}
There is a morphism of fibrations
\begin{equation}\label{eq:affschemes}\begin{aligned}\xymatrix@=1.3cm{
\Comm^{op} \ar[r]^{(\mathrm{Spec},\Oc{})} \ar[d]_{Id_{\Comm^{op}}}& \LRinged\ar[d]\\
\Comm^{op} \ar[r]_{\Spec} & \Top.
}\end{aligned}\end{equation}
The functor $(\mathrm{Spec},\Oc{})$ is defined in any book of algebraic geometry and we briefly recall its construction. The base functor $\Spec$ has already been defined in \thref{ex:spec}. In order to define $(\mathrm{Spec},\Oc{})$, recall that it is enough to define a sheaf on a space $X$ by its values on a basis of the topology of $X$ \cite{MM92,Qin02}. The value of the sheaf $\Oc A$ on the basis element $D(f)$ is the localization $A_f$ of the ring $A$ at the element $f$ (for a direct definition see \cite{Har77}). It is a locally ringed space for every ring \cite{Har77} and its global section ring is $$\Oc A(\Spec A)=\Oc A(D(1))=A_1=A.$$ 
Let now $\phi\colon A\to B$ be a homomorphism of commutative rings. It induces a continuous function $\Spec \phi\colon\Spec B\to\Spec A$ (\thref{ex:spec}) and a morphism of sheaves $\Oc \phi\colon\Oc A\to (\Spec \phi)_*\Oc B$ by the localization $(\Oc\phi)_{D(f)}\coloneqq\phi_f\colon A_f\to B_{\phi(f)}$. Then, 
$$(\Spec\phi,\Oc\phi)\colon(\Spec B,\Oc B)\to(\Spec A,\Oc A)$$ 
is a morphism of locally ringed spaces \cite{Har77}. 

It is a basic classical result of algebraic geometry that the functor $(\mathrm{Spec},\Oc{})$ is fully faithful \cite{Har77,Qin02}\footnote{It is not full into the category $\Ringed$ of ringed spaces, though (see \cite{Har77} for a counterexample).}. On the contrary, the base functor $Spec\colon\Comm^{op}\to\Top$ is not full. Here is a counterexample\footnote{We found this counterexample on a blog called ``blogarithm''. We reproduce it here, because internet references are not stable enough.}. Since $\Q$ is a field, $\Spec\Q=\{0\}$. Moreover, since $\Z$ is a PID, $\Spec\Z=\{0\}\cup\{(p)\mid p\text{ is prime}\}$. A prime number $p\in\Z$ defines a continuous map $\Spec\Q\to\Spec\Z$ by $(0)\mapsto(p)$. This continuous map is not the image of the unique homomorphism of rings $\phi\colon\Z\to\Q$ (\Z is an initial ring) under the functor Spec. Indeed, suppose there is such a homomorphism. Then, $\phi^{-1}((0))=(p)$, and thus $\phi(p)=0$. Again, $0=\phi(p)=p\phi(1)$. Thus, $\phi(1)=0$, since \Q is a field. In consequence, $\phi^{-1}((0))=\Z$, a contradiction.

The morphism of fibrations \hyperref[eq:affschemes]{(\ref*{eq:affschemes})} is not Cartesian. Recall first that a morphism of ringed spaces into a locally ringed space is Cartesian in \LRinged \ssi it is Cartesian in \Ringed. Let now $\phi\colon A\to B$ be a morphism of commutative rings and $(\Spec \phi,\Oc\phi)\colon(\Spec B,\Oc B)\to (\Spec A,\Oc A)$ the induced morphism of locally ringed spaces. As explained above, a Cartesian morphism over $\Spec\phi\colon\Spec B\to\Spec A$ at $(\Spec A,\Oc A)$ is given by $(\Spec B,(\Spec\phi)^{-1}\Oc A)\to(\Spec A,\Oc A)$. Now if $(\Spec \phi,\Oc\phi)$ were Cartesian, then there would be a vertical isomorphism $(\Spec B,\Oc B)\ra{\cong}(\Spec B,(\Spec\phi)^{-1}\Oc A)$, and therefore, for each point $\frak p\in\Spec B$, an isomorphism of the stalks $A_{\phi^{-1}(\frak p)}\ra{\cong}B_{\frak p}$. Consider, for instance, the morphism $\phi$ given by the inclusion $i\colon\Z\hookrightarrow\Z[x]$ and the prime ideal $(0)\in\Z[x]$ ($\Z[x]$ is a domain because $\Z$ is). Then, $i^{-1}((0))=(0)$, $\Z[x]_{(0)}=\Q(x)$ and $\Z_{(0)}=\Q$. However, the field of rational functions $\Q(x)$ over \Q is not isomorphic to $\Q$.

The \emph{category of affine schemes} is the replete (full) image of the functor $(\mathrm{Spec},\Oc{})$ \hyperref[eq:affschemes]{(\ref*{eq:affschemes})}. In other words, its objects are all locally ringed spaces $(X,\Oc X)\cong (\Spec A,\Oc A)$ isomorphic in \LRinged to the image of a ring $A$. Its morphisms are all morphisms of locally ringed spaces between those. 

Note that it is also the replete (non full) image of the composite functor $$\Comm^{op}\xlongrightarrow{(\mathrm{Spec},\Oc{})}\LRinged\hookrightarrow\Ringed.$$ Indeed, by repleteness of \LRinged in \Ringed and fullness of $(\mathrm{Spec},\Oc{})$, this functor, although not full, is isofull. Moreover, a ringed space isomorphic to a $(\Spec A,\Oc A)$ is automatically locally ringed and isomorphic by an isomorphism of locally ringed spaces. Finally, a morphism is in the replete image of the composite functor \ssi it is a morphism of locally ringed space (use the fullness of $(\mathrm{Spec},\Oc{})$ into \LRinged).

\subsection{Notions of image of a morphism of fibrations}\label{ssec:imagemorphfib}

Factorization systems have the nice property (compare to the weak ones) that they can always be made functorial. More precisely, this means that given a factorization system $(\fL,\fR)$ of a category \C, there exists a triple $(E,\epsilon,\mu)$ where:
\begin{itemize}
\item $E\colon\C^{\two}\to\C$ is a functor,
\item $\epsilon\colon\dom\Rightarrow E$ and $\mu\colon E\Rightarrow\cod$ are natural transformations,
\end{itemize}
such that, for all $f$ in $\C$,
\begin{enumerate}[(i)]
\item $\rho_{f}\circ\epsilon_{f}=f$,
\item $\epsilon_{f}\in\fL$, 
\item $\mu_{f}\in\fR$.
\end{enumerate}
Moreover, the triple $(E,\epsilon,\mu)$ is ``uniquely determined up to isomorphism'' by $(\fL,\fR)$ in the obvious sense. One obtains such a functorial factorization system via a choice of an $(\fL,\fR)$-factorization for each arrow in \C and by using the \emph{orthogonality property} of the factorization system. 

Now, given a small category \D, any functorial factorization system $(\fL,\fR,E,\epsilon,\mu)$ on a category \C determines a componentwise factorization system $(\fL^{\D},\fR^{\D})$ on the functor category $\C^{\D}$: the elements of $\fL^{\D}$ are the natural transformations in $\C^{\D}$ whose components are in $\fL$, and similarly for $\fR^{\D}$. Indeed, a natural transformation ${\alpha\colon F\Rightarrow G\colon\D\to\C}$ canonically determines a functor $A\colon\D\to\C^{\two}$ and the factorization of $\alpha$ is given by the composition 

\[\xymatrix{
F \ar@{=>}[rr]^{\alpha}\ar@{=>}[rd]_{\epsilon\cdot A}& & G \\
& E\circ A\ar@{=>}[ru]_{\mu\cdot A} &}
\]

In particular, any factorization system $(\fL,\fR)$ on a category \C determines (via a choice of a factorization functor) a factorization system $(\fL^{\two},\fR^{\two})$ on the arrow category $\C^{\two}$ of \C. Therefore, the factorization systems studied on \CAT pass to $\CAT^{\two}$. But what about the XL-category \FIB?

In this part, we explore the different notions of image of a morphism of fibrations. As for properties of the morphisms, there corresponds to each concept a global version, and these will fortunately coincide under some natural conditions on the morphism. As we shall see, images of morphisms of fibrations are not always subfibrations. We give conditions under which they are subfibrations, strong or not. All these concepts are new. 

\paragraph{}

Let
\begin{equation}\label{eq:morfibtriv}\begin{aligned}\xymatrix@=1.3cm{
\D \ar[r]^F \ar[d]_{Q}& \E\ar[d]^P\\
\A \ar[r]_G & \B,
}\end{aligned}\end{equation} 
be a morphisms in $\CAT^{\two}$.

We first consider factorization systems of $\CAT^{\two}$ determined by factorization systems in \CAT via the procedure explained above. We refer to the factorization systems on \CAT defined in \autoref{ssec:factsyst} with the choice of factorization given by the corresponding image notion.

\paragraph{Image}

The \emph{image of $(F,G)$} is the restriction  of $P$ to the respective images of $F$ and $G$:
$$\im F\to\im G.$$ 
Note that if $(F,G)$ is (not necessarily globally) isofull, then $\im F=F(\A)$.

\paragraph{Global replete image}
The \emph{global replete image of $(F,G)$} is the restriction of $P$ to the replete images of $F$ and $G$:
$$
\rim F\to\rim G.
$$

\paragraph{Full image} 
The \emph{full image} of $(F,G)$ is the restriction of $P$ to the respective full images of $F$ and $G$:
$$\fim F\to\fim G.$$

\paragraph{Global replete full image}
The \emph{global replete full image of $(F,G)$} is the restriction of $P$ to the replete full images of $F$ and $G$:
$$
\rfim F\to\rfim G.
$$

We now consider factorization system in $\CAT^{\two}$ that are not determined componentwise by factorization systems in \CAT. Vertical isomorphisms are the isomorphisms relevant here, not general isomorphisms.

\paragraph{Replete image}

The left class of this factorization system consists in essentially surjective morphisms $(F,G)$ such that every morphism in \E is the composite of finitely many vertical isomorphisms and morphisms in $F(\D)$ and such that $G$ is an extremal epimorphism. It right class consists in morphisms that are replete and globally monic on morphisms.

The corresponding image notion is called the \emph{replete image of $(F,G)$} and is denoted
\begin{equation}\label{eq:repim}
\rim{(F,G)}\to\im G,
\end{equation}
where the subscript $v$ is for ``vertical''. The total category $\rim{(F,G)}$ is defined as follows. Its objects are all $E\in\E$ that admit a \emph{vertical} isomorphism $E\cong F(D)$. In other words, they are objects of the replete images $\rim{F_A}$ of the fibre restrictions $F_A\colon\D_A\to\E_{GA}$, $A\in\A$. Its morphisms are all morphisms in \E that are finite composites of \emph{vertical} isomorphisms in \E and arrows in $F(\D)$. The base category \im G is the usual image of the functor $G$ and the functor \hyperref[eq:repim]{(\ref*{eq:repim})} is the restriction of $P$.

In case $(F,G)$ is isofull, the morphisms of the replete image have the following description. They are all morphisms in \E that admit a factorization
$$\xymatrix{
E\ar[d]_\cong \ar[r]^k & E' \ar[d]^\cong\\
F(D)\ar[r]_{F(h)} & F(D'),
}$$
whose legs are vertical isomorphisms and where $h$ is a morphism in \D. This condition does not depend on the choice of the vertical isomorphisms and $h$ can always be supposed to have $D$ and $D'$ as domain and codomain.

For a Cartesian morphism of fibration $(F,G)$ with $G$ isofull, the morphisms of the category $\rim{(F,G)}$ can be equivalently defined by conditions completely expressed in the fibres. For, note that in this case, the objects of $\rim{(F,G)}$ are closed under inverse images along arrows of \im G (use the fact that $G$ is isofull). Now, a morphism $k\colon E'\to E$ in $\E$ between objects $E, E'$ equipped with vertical isomorphisms $E\cong F(D),E'\cong F(D')$ is in $\rim{(F,G)}$ \ssi it sits over a morphism $g$ in \im G and the induced morphism $\bar k\colon E'\to g^*E$ is in \rim{F_{Q(D')}}. The fact that this condition is sufficient is easy. One shows it is necessary by a recurrence on the number of arrows in $F(\D)$ that occur in the arrow $k$ (which is a composite of arrows in $F(\D)$ and of vertical isomorphisms).

The replete image $\rim{(F,G)}\subset\rim F$ is a subcategory of the global replete image, which is full when $(F,G)$ is isofull. Moreover, when $(F,G)$ is a morphism of fibrations with $G$ replete, then they coincide.

\paragraph{Replete full image}

The left class of this factorization system consists in essentially surjective morphisms $(F,G)$ with $G$ surjective on objects. Its right class consists in morphisms that are replete, globally full and globally monic on morphisms.

The  corresponding image notion is called the \emph{replete full image} of $(F,G)$ and is denoted 
\begin{equation}\label{eq:refulim}
\rfim{(F,G)}\to\fim G.
\end{equation}
The category $\rfim{(F,G)}$ is the full subcategory of $\E$ of objects $E\in \E$ admit a \emph{vertical} isomorphism $E\cong F(D)$. The category $\fim G$ is the usual full image of $G$, and the functor \hyperref[eq:refulim]{(\ref*{eq:refulim})} is the restriction of the functor $P$.

The replete full image $\rfim{(F,G)}\subset\rfim F$ is a full subcategory of the global replete full image. Moreover, when $(F,G)$ is a morphism of fibrations with $G$ replete, then they coincide.

\subsubsection{Fibred properties of the image}
We now give the main result of this part, which describes the fibred properties of these different notions of image when the morphism $(F,G)$ is a morphism of fibrations.

\begin{Prop}\label{prop:fibpropim}
Let $(F,G)\colon (\D\ra Q\A)\to (\E\ra P\B)$ be a morphism of fibrations.
\begin{enumerate}[(i)]
\item
If $(F,G)$ is globally pseudomonic, then the image of $(F,G)$ is a subfibration of $P$ and the replete, \resp global replete, image is a replete, \resp a globally replete, subfibration.
\item
If $(F,G)$ is globally isofull Cartesian and $G$ faithful, then the image of $(F,G)$ is a strong subfibration of $P$ and the replete, \resp global replete, image is a replete, resp.\ globally replete, strong subfibration.
\item
If $(F,G)$ is Cartesian and $G$ full, then the full, replete full and global replete full images of $(F,G)$ are full ($\iff$globally full) strong subfibrations. The replete full, \resp global replete full, image subfibration is in addition replete, \resp globally replete.
\end{enumerate}
\end{Prop}
\begin{Pf}\NoEndMark
(i) We prove the most difficult case, which is the global replete image. The diagrams \hyperref[eq:sitQ]{(\ref*{eq:sitQ})} and \hyperref[eq:sitP]{(\ref*{eq:sitP})} summarize the development and may help in following the proof. Let ${E\in\rim F}$ be over $B\in\B$ and $f\colon B'\to B$ in $\rim G$. We construct a Cartesian lift of $f$ at $E$ in $\rim F$ with respect to the functor $\rim Q$. 

By definition of the replete image of a functor, there exists an isomorphism $h\colon F(D)\ra\cong E$, and a factorization
\begin{equation}\label{eq:fact1}\begin{aligned}\xymatrix@=1.3cm{
G(A)\ar[r]^{G(a)}\ar[d]_{g^{-1}}^\cong & GQ(D)\ar[d]^{P(h)}_\cong\\
B'\ar[r]_f & B
}\end{aligned}\end{equation}
with $a\colon A\to Q(D)$ a morphism in \A (recall that one can chose freely the isomorphisms). Let $\bar a_D\colon a^*D\to D$ be a Cartesian lift of $a$ at $D$ in the fibration $Q$. Let $\bar g_{F(a^*D)}\colon g^*F(a^*D)\ra\cong F(a^*D)$ be a Cartesian lift of $g$ at $F(a^*D)$ in the fibration $P$. The composite
\begin{equation}\label{eq:Cartmorreplete}\xymatrix@1@=1.3cm{
g^*F(a^*D)\ar[r]^{\bar g_{F(a^*D)}}_\cong & F(a^*D)\ar[r]^{F(\bar a_D)} & F(D)\ar[r]^-{h}_-{\cong} & E
}\end{equation}
is Cartesian over $f$ with respect to \rim Q. 

In order to verify this, let $l\colon E'\to E$ be a morphism in \rim F and $b\colon P(E')\to B'$ in \rim G such that $f\circ b=P(l)$. Thus, there are factorizations
\begin{equation}\label{eq:fact2}\begin{aligned}\xymatrix@=1.3cm{
F(D')\ar[d]_k^\cong\ar[r]^{F(d)} & F(D)\ar[d]^{h}_\cong\\
E'\ar[r]_{l} & E
}\quad\text{and}\quad\xymatrix@=1.3cm{
GQ(D')\ar[d]_{P(k)}^\cong\ar[r]^{G(a')} & G(A)\ar[d]^{g^{-1}}_\cong\\
P(E')\ar[r]_{b} & B'
}\end{aligned}\end{equation}
with $d\colon D'\to D$ in \D and $a'\colon Q(D')\to A$. Now, by combining diagrams \hyperref[eq:fact1]{(\ref*{eq:fact1})} and \hyperref[eq:fact2]{(\ref*{eq:fact2})}, one obtains both equalities
$$
P(l)\circ P(k)=P(h)\circ GQ(d)\quad\text{and}\quad P(l)\circ P(k)=P(h)\circ G(a)\circ G(a').
$$
Since $P(h)$ is an isomorphism and $G$ is faithful, one concludes $Q(d)=a\circ a'$. Thus, there is the following situation in the fibration $Q$:
\begin{equation}\label{eq:sitQ}\begin{aligned}\xymatrix@C=1.5cm@R=1.3cm{
D'\ar@/^2pc/[rr]^d\ar@{-->}[r]_{\exists!\bar a} & a^{*}D\ar[r]_{\bar a_D} & D\\
Q(D')\ar@/_2pc/[rr]_{Q(d)}\ar[r]^{a'} & A \ar[r]^a & QD
}\end{aligned}\end{equation}
The situation in $P$ is now the following.
\begin{equation}\label{eq:sitP}\begin{aligned}\xymatrixnocompile@C=1cm@R=3cm{
F(D')\ar@/_2pc/[rrr]^{F(\bar a)}\ar[r]_\cong^k\ar@/_3pc/[rrrr]_{F(d)}& E' \ar@/^3pc/[rrrr]^l \ar@{-->}[r]^-{\bar b}	& g^*F(a^*D) \ar[r]_\cong^-{\bar g_{F(a^*D)}} & F(a^*D) \ar[r]^{F(\bar a_D)} & F(D)\ar[r]^h_\cong & E\\
GQ(D')\ar@/^2.5pc/[rrr]^{G(a')}\ar[r]^{P(k)}_\cong & P(E')\ar@/_3pc/[rrrr]_{P(l)}\ar[r]^b & B'\ar@/_2pc/[rrr]^f\ar[r]^g_\cong & G(A)\ar[r]^{G(a)} & GQ(D)\ar[r]^{P(h)}_\cong & B
}\end{aligned}\end{equation}
In the bottom row, everything commutes by definition. Consider the top row. It includes, on the upper part, a pentagon containing $l$, a quadrilateral containing $F(\bar a)$ and a triangle containing $F(d)$. The lower triangle commutes by definition of $\bar a$. The morphism $\bar b$ is defined to be the only one to make the quadrilateral commute, and is therefore a morphism in \rim F. The morphism $l$ is also part of an external quadrilateral, which is precisely the first square of \hyperref[eq:fact2]{(\ref*{eq:fact2})}. This implies the commutativity of the pentagon, since $k$ is an isomorphism. Thus $\bar b$ has the commutativity property required in the proof of the Cartesianness. Finally, using the fact that $Q(\bar a)=a'$, one concludes that $P(\bar b)=b$. We thus have proved the existence part of the Cartesian property.

We now turn to the uniqueness of $\bar b$. Suppose $e\colon E'\to g^*F(a^*D)$ makes the pentagon commute, sits over $b$ and belongs to \rim F. By the last condition, there exists a morphism $\tilde a \colon D'\to a^*D$ making the following square commute:
\begin{equation*}\begin{aligned}\xymatrix@=1.3cm{
F(D')\ar[d]_k^\cong\ar[r]^{F(\tilde a)} & F(a^*D)\ar[d]^{\bar g_{F(a^*D)}^{-1}}_\cong\\
E'\ar[r]_{e} & g^*F(a^*D)
}\end{aligned}\end{equation*}
We can now replace $\bar b$ by $e$ and $\bar a$ by $\tilde a$ in diagram \hyperref[eq:sitP]{(\ref*{eq:sitP})}. The lower triangle commutes, because its postcomposition with $h$ does. We have thus $F(\bar a_D)\circ F(\tilde a)=F(d)$. By faithfulness of $F$, \begin{equation}\label{eq:atildea}\bar a_D\circ\tilde a=d.\end{equation} Moreover, $GQ(\tilde a)=PF(\tilde a)=G(a')$, because $P(e)=b$ by hypothesis. By faithfulness of $G$, this implies that $Q(\tilde a)=a'$. Combining the latter equality with \hyperref[eq:atildea]{(\ref*{eq:atildea})}, by Cartesianness of $\bar a_D$, $\tilde a=\bar a$. Thus $e=\bar b$.

(ii) We prove the case of the global replete image. The other cases work the same way. First note that $F$ is isofull and $G$ pseudomonic. Now, the first part of the proof of (i), that is the construction of the morphism \hyperref[eq:Cartmorreplete]{(\ref*{eq:Cartmorreplete})}, and the existence part of the Cartesian property of this morphism only requires the latter condition: $F$ isofull and $G$ pseudomonic. Yet, the morphism \hyperref[eq:Cartmorreplete]{(\ref*{eq:Cartmorreplete})} is Cartesian with respect to $P$, because we have supposed that $(F,G)$ is Cartesian. Since we have shown that the unique morphism $\bar b$ lies in \rim f in this situation, \hyperref[eq:Cartmorreplete]{(\ref*{eq:Cartmorreplete})} is Cartesian with respect to \rim Q, the restriction of $P$ to \rim F and \rim G.

(iii) Again, the global replete full image is more complicated, and we show this case. Note that $G$ full implies that $G$ is in particular isofull. The other cases work with the same ideas. Let $E\in\rfim F$ over $B\in\rfim G$, $h\colon E\ra\cong F(D)$ an isomorphism and $f\colon B'\to B$ a morphism in \rfim G. Then, there is a factorization as in \hyperref[eq:fact1]{(\ref*{eq:fact1})} because $G$ is full, and therefore a morphism is defined as in \hyperref[eq:Cartmorreplete]{(\ref*{eq:Cartmorreplete})}. The latter morphism has its domain in \rim F by construction and, since $F$ is Cartesian, it is Cartesian with respect to $P$. Again, \rfim F is a full subcategory, and thus this a Cartesian morphism with respect to \rim Q.

\cqfd
\end{Pf}
Let us denote by $\FIBgp\subset\FIB$, $\FIBcgp\subset\FIBc$ the wide sub-XL-categories whose morphisms are globally pseudomonic and $\FIBcf\subset\FIBc$ the wide sub-XL-category whose morphisms have a full base component functor.
\begin{Cor}
  \begin{enumerate}[(i)]
  \item   The image, replete image et global replete image factorization systems on $\CAT^{\two}$ restrict to the sub-XL-categories \FIBgp and \FIBcgp.
  \item   The full, replete full and global replete full image factorization systems on $\CAT^{\two}$ restrict to the sub-XL-category \FIBcf. 
  \end{enumerate}
\cqfd
\end{Cor}

We now consider some particular cases.

\subsubsection{Subfibrations}
Let 
$$\xymatrix@=1.3cm{
\D \ar@{^(->}[r]^I \ar[d]_{\bar P}& \E\ar[d]^P\\
\A \ar@{^(->}[r]_J & \B,
}$$
be a subfibration. The notion of image always exists and is trivial: it is the subfibration itself. The replete image of $(I,J)$ is called the \emph{repletion} of the subfibration. The global replete image of $(I,J)$ is called the \emph{global repletion} and is thus the restriction of $P$ to the repletions of the respective subcategories:
$$
\rep{\bar P}\colon\rep\D\to\rep \A.
$$
The full image of the subfibration is just the restriction of $P$ to the respective full subcategories of \E and \B generated by the objects of $\D$ and $\A$, $$\full{\bar P}\colon\full\D\to\full\A.$$ It is called the \emph{associated full subfibration}. The (global) replete full image of $(I,J)$ is called the \emph{\op global\fp full repletion}. The global full repletion of $(I,J)$ is thus the restriction of $P$ to the respective full repletions of $\D$ and \A:
$$
\frep{\bar P}\colon\frep\D\to\frep \A.
$$
\begin{Cor}
Let $(I,J)\colon(\D\to\A)\hookrightarrow(\E\ra P\B)$ be a subfibration.
\begin{enumerate}[(i)]
\item
If it is a globally pseudomonic subfibration, then its repletion, \resp global repletion, is a replete, \resp globally replete, subfibration.
\item
If it is a globally pseudomonic strong subfibration, then its repletion, \resp global repletion, is a replete, \resp globally replete, strong subfibration.
\item
If it is a strong subfibration with \A full in \B, then its associated full subfibration, its full repletion and its global full repletion are strong subfibrations. The full repletion, \resp global full repletion, is in addition replete, \resp globally replete.
\end{enumerate}
\end{Cor}

\subsubsection{Functors over \texorpdfstring{\B}{B}}
Let 
$$\xymatrix@C=.5cm{
\D\ar[rr]^F \ar[dr]_Q && \E\ar[dl]^P\\
&\B&
}$$
be a functor over \B. Then, all the global and fibred notions of images coincide and have \B as codomain. Thus, there is no conflict of notations and terminology between some sort of image of $F$ as a functor over \B and as a plain functor. For instance, the full image of $F$ as a functor over \B is the restriction $\fim F\to\B$ of $P$, where $\fim F$ is the image of $F$ as a mere functor.
\begin{Cor}\label{cor:imfunovB}
Let $F\colon(\D\ra Q\B)\to(\E\ra P\B)$ be a functor over \B. Then,
\begin{enumerate}[(i)]
\item
If $F$ is a pseudomonic functor, the restrictions of $P$ to the image and replete image of $F$ are respectively subfibrations and replete ($\iff$globally replete) subfibrations of $P$ over \B.
\item
If $F$ is isofull and Cartesian, the restrictions of $P$ to its image and replete image are strong subfibrations of $P$ over \B.
\item\label{cor:imfunovB:fullreplete}
If $F$ is Cartesian, the restriction of $P$ to its full image is a full ($\iff$globally full) strong subfibration of $P$ over \B. The restriction of $P$ to its replete full image is a strictly full ($\iff$globally strictly full) strong subfibration of $P$ over \B.
\end{enumerate}
\end{Cor}

\subsection{Categories and fibrations of trivial objects}
We explain now how a morphism of fibrations $$(F,G)\colon(\D\ra Q\A)\to(\E\ra P\B)$$ determines a subfunctor of trivial objects $\Triv\subset P$. There are different choices available. As explained above, we want to have, depending on the situation, either a replete or a globally replete subfunctor. 

In the former case, there are two different subfunctors of interest, the replete image and the replete full image 
\begin{equation}
  \label{eq:6}
\rim(F,G)\to\im G,\qquad\rfim(F,G)\to\fim G.  
\end{equation}
In the latter case, the two choices available are the global replete image and the global replete full image
\begin{equation}
  \label{eq:7}
\rim F\to\rim G,\qquad\rfim F\to\rfim G.  
\end{equation}

Recall that if the base functor $G$ is replete, then there is no difference between \eqref{eq:6} and \eqref{eq:7}: an object $E\in\E$ is isomorphic to some object of $F(\D)$ \ssi it is vertically isomorphic to some object of $F(\D)$. In the examples we have in mind, when the choice does matter, then it is the global notion that should be chosen. As for the question of fullness of the image, in both cases \eqref{eq:6} and \eqref{eq:7}, when $(F,G)$ is globally full, the full and non-full versions coincide, but otherwise, they might differ. The non-full versions remember something of the morphisms of in $Q$, whereas the full versions do not. In fact, if $F$ is pseudomonic functor, then \rim F is even equivalent to \D (see~\hyperref[rem:replete-full-image:equ]{Remark~\ref*{rem:replete-full-image}(\ref*{rem:replete-full-image:equ})}). Recall finally that \thref{prop:fibpropim} helps to determine whether these different subfunctors of trivial objects are (strong) subfibrations of $P$.

\begin{Exs}[We refer to examples of \autoref{ssec:Firstexmorfib}.]
	\item (Locally ringed spaces) Consider the subfibration over \Top of locally ringed spaces in the fibration of ringed spaces. Since it is over \Top, the global and fibrewise notions coincide. Moreover, the subfibration is replete (thus pseudomonic), and therefore is equal to its repletion $\LRinged=\rep\LRinged$. On the other hand, it is different than its full repletion $\frep\LRinged$.
	\item (Affine schemes) One has a morphism of fibrations $$((\mathrm{Spec},\Oc{}),\mathrm{Spec})\colon(\Comm^{op}\ra{Id}\Comm^{op})\to(\LRinged\to\Top).$$ Here, one should consider the global notion of isomorphism in order to obtain affine schemes, and therefore global versions of images of this morphism of fibrations. The functor $(\mathrm{Spec},\Oc{})$ being full, its replete full and replete images coincide. The subcategory of trivial objects $\Triv_{t}$ in this context is the category of affine schemes. As explained in the last paragraph of \autoref{sssec:affsch}, it is also equal to the replete image of $(\mathrm{Spec},\Oc{})$ in \Ringed. 
	\end{Exs}

\subsection{More examples}\label{ssec:moreextrivobjects}
We now describe classes of examples that will be of great use for locally trivial objects that one meets in nature.

\subsubsection{Product fibrations}\label{ssec:prodfib}
This is a particular situation that gathers a whole class of examples.
\begin{equation}\label{eq:funovBprodfib}\begin{aligned}\shorthandoff{;:!?}\xymatrix@!=.7cm@R=1.5cm{
\B\times\C\ar[dr]_{Pr_{\B}}\ar[rr]^F & & \E\ar[dl]^P\\
&\B&
}\end{aligned}\end{equation}
Here are some examples.
\paragraph{Free modules}
Consider, for each ring $R$, a fixed \emph{free $R$-module adjunction}
$$\xymatrix{
F_R\colon\Set\ar@<.4pc>[r] \ar@{}[r]|{\perp}& Mod_R:\!U\ar@<.4pc>[l]
}$$
with unit $\eta^R$.
One defines an opCartesian functor over \Ring
$$\shorthandoff{;:!?}\xymatrix@!=.7cm@R=1.5cm{
**[l]{\Ring\times\Set}\ar[dr]_{Pr_{\Ring}}\ar[rr]^-F & & **[r]{\Mod}\ar[dl]\\
&\Ring&
}$$
this way. 
$$ \begin{array}{rcl}
F(R,X)&\coloneqq&(R,F_R(X))\\
F((R,X)\ra{(\phi,f)}(S,Y))&\coloneqq&(R,F_R(X))\ra{(\phi,\phi_f)}(S,F_S(Y)),
\end{array}$$
where $\phi_f\colon F_R(X)\to F_S(Y)$ is the following morphism of $R$-modules, seeing $F_S(Y)$ as an $R$-module via $\phi$.
$$\xymatrix{
X\ar[d]_f\ar[r]^-{\eta^R_X} & UF_R(X)\ar@{-->}[d]^{U(\phi_f)}\\
Y\ar[r]_-{\eta^S_X} & UF_S(Y).
}$$ 
The opCartesianness can be directly verified. Alternatively, one may use the following composite adjunction.
$$\xymatrixnocompile@C=2cm@R=.1cm{
&\scriptstyle\cong&\\
\Set\ar@/^2.5pc/[rr]^{F_S}\ar@<.4pc>[r]^{F_R} \ar@{}[r]|{\perp} & \Mod_R\ar@<.4pc>[l]^{U_R} \ar@<.4pc>[r]^{-\otimes_R S} \ar@{}[r]|{\perp}& \Mod_S\ar@<.4pc>[l]^{\phi^\sharp}\ar@/^2.5pc/[ll]^{U_S}\\
&\scriptstyle\stackrel{}{=}&
}$$
One might restrict the category \Set to subcategories like the full subcategory \textit{FinSet} of finite sets, or even to the discrete subcategory \N of natural numbers. By the dual of Corollary \hyperref[cor:imfunovB:fullreplete]{\ref*{cor:imfunovB}(\ref*{cor:imfunovB:fullreplete})}, the replete full image of $F$ yields a strong sub-opfibration. It is the full subcategory of \Mod of free modules. One would obtain the full subcategories of finitely generated free modules with both other choices \textit{FinSet} and \N.

\subsubsection{Subcategory of the fibre over the terminal object}
We found the inspiration for this class of examples in Street's notion of locally trivial objects, but these trivial objects are in some sense both more general and more particular than Street's ones (see discussion in \autoref{sssec:Street}). The examples we explain now are particular cases of the preceding examples where the fibration is a product fibration. In fact, when \B has a terminal object, all Cartesian functors over \B of the form \hyperref[eq:funovBprodfib]{(\ref*{eq:funovBprodfib})} can be seen as of the type described now. Moreover, important categories of locally trivial objects, such as those involved in $K$-theory of schemes and of spaces, are actually defined from trivial objects of this type.

Let $P\colon\E\to\B$ be a cloven fibration, and suppose $\B$ has a terminal object $*$. Let $\C\subset\E_*$ be a subcategory of the fibre $\E_*$ of $P$ over the terminal object $*$. There is a Cartesian functor
$$\shorthandoff{;:!?}\xymatrix@!=.7cm@R=1.5cm{
\B\times\C\ar[dr]_{Pr_{\B}}\ar[rr]^{F_\C} & & \E\ar[dl]^P\\
&\B&
}$$
defined as follows. We denote by $!_B\colon B\to *$ the unique morphism to the terminal object from $B\in\B$ and, as usual, $$\overline{(!_B)}_C\colon\oc_B^*C\to C$$ the Cartesian lift of $!_B$ at $C\in\C$ belonging to the cleavage of $P$.
\begin{equation}\label{eq:funsubcatterm}\begin{array}{rcl}
(B,C)				& \longmapsto	& F_\C(B,C)=\oc_B^*C\\
&&\\
(B,C)\ra{(f,h)}(B',C')	& \longmapsto 	& \xymatrixnocompile{!_B^{*}C\ar[rr]^{\overline{(!_B)}_C}\ar@{-->}[dr]_{F_\C(f,h)}&	& C\ar[d]^h\\
							& !_{B'}^{*}C'\ar[r]_{\overline{(!_B')}_C'}	& C'\\
B\ar@/_1pc/[rr]_{!_B}\ar[r]^f			& B' \ar[r]^{!_{B'}}  	& {*}}
\end{array}\end{equation}
Thus, the restriction of $F_\C$ to the fibre $Pr_\B^{-1}(B)$ is just the restriction to \C of the inverse image functor $!_B^*\colon\E_*\to\E_B$ induced by the cleavage.
In order to show that this functor is Cartesian, use Proposition \hyperref[Prop:Cartesianii]{\ref*{Prop:Cartesian}(\ref*{Prop:Cartesianii})}.

By Corollary \ref{cor:imfunovB}, its replete full image is a strictly full strong subfibration of $P$. Its total category $\rfim F$ has the following characterization in \E. It is the full subcategory of \E of objects $E$ such that there exists an object $C\in \C$ and a Cartesian arrow as in the following diagram.
\begin{equation}\begin{aligned}\label{eq:trivterminal}\xymatrix{
E\ar[r] & C\\
P(E)\ar[r]^{!} & {*}
}\end{aligned}\end{equation}
Let us return to the general situation of \hyperref[eq:funovBprodfib]{(\ref*{eq:funovBprodfib})}. Suppose the functor $$F\colon\B\times\C\to\E$$ is \emph{Cartesian over \B} and consider a cleavage of $P$. Then, one can compare the functor $F$ with the composite 
$$\B\times \C\ra{Id_{\B}\times F(*,-)}\B\times \E_*\xlongrightarrow{F_{\E_*}}\E,$$ 
where the second functor is defined as in \hyperref[eq:funsubcatterm]{(\ref*{eq:funsubcatterm})}. One readily verifies that $F$ is isomorphic to this composite functor. Observe that in general $F(*\times \C)$ might not be a subcategory, though. These observations show us that the replete full image of the functor $F$ has a simple shape. It is the full sub-category of \E of objects $E$ that admit a Cartesian arrow $E\to F(*,C)$, for some $C\in\C$. 

\section{Locally trivial objects in a fibred site}\label{sec:loctriv}

\subsection{Definitions}
We first study the notion of local triviality at the level of objects, and then turn to morphisms.

\subsubsection{Objects} 
\begin{Defs}[\indexb{Locally trivial objects!-- in a fibred site}]
\item 
  Let $P\colon\E\to\B$ be a fibration and $\Triv\subset P$ be a subfunctor of trivial objects.
\begin{itemize}
  \item An object $B\in\B$ is \emph{locally trivial in the base for a covering $R$} if $R$ covers $B$ and has trivial domains.
  \item An object $E\in\E$ is \emph{locally trivial in the total category for a covering $R$} if $R$ covers $P(E)$ and for each $f\in R$, the inverse image of $E$ along $f$ is trivial. In other words, for each $f\in R$, there is a trivial object $T_{f}\in\Triv_{t}$ and a Cartesian arrow over $f$:
\[\xymatrix{
T_{f}\ar[r] &E\\
A\ar[r]^{f} & B
}\] 
  \end{itemize}
\item  Let $P\colon\E\to(\B,K)$ be a fibred site such that $K$ satisfies axiom \tm and $\Triv\subset P$ a subfunctor of trivial objects. An object of the base or of the total category is \emph{locally trivial} if there exists a $K$-covering for which it is locally trivial. 
\end{Defs}

\begin{Rems}[\label{rem:proploctriv}]
\item Note that the functor $P$ restricts to locally trivial objects in the total and base categories.
  \item Locally trivial objects in both base and total categories have a common characterization as objects than can be covered by trivial ones. Indeed, an object $E\in\E$ is locally trivial for a covering $R$ if and only if it is can be $K_{\E}$-covered by trivial objects, where $K_{\E}$ is the covering function on \E induced by $P$ and $K$ (\thref{def:indcov}). Moreover, if $E\in\E$ is locally trivial for a covering $R$, then $P(E)$ also is.
  
  \item Let us now, and for the following remarks, consider a subfunctor of trivial objects induced by a morphism of fibrations
\[\xymatrix{
\D\ar[r]^{F}\ar[d]_{Q} & \E\ar[d]^{P}\\
\A\ar[r]_{G} & \B.
}\]	
Consider trivial objects determined by the replete (full) image. An object $E$ is locally trivial for a covering $R$ \ssi there exists, for each $f\in R$, a Cartesian lift of $f$ at $E$ with domain in $F(\D)$:
	$$\xymatrix{
	F(D)\ar[r]^-{\mathrm{cart.}} & E\\
	GQ(D)\ar[r]^f & P(E).
	}$$
	\item Consider trivial objects determined by the globally replete (full) image of $(F,G)$. An object $E$ is locally trivial for a covering $R$ if and only if there exists, for each $f\in R$, a Cartesian lift of $f$ at $E$ such that:
$$\xymatrix{
	F(D)\ar[r]^h_\cong & f^*E\ar[r]^{\bar f_E} & E\\
	GQ(D)\ar[r]^{P(h)}_\cong & B\ar[r]^f & P(E).
	}$$
	\item Whether one has a replete or a globally replete subfunctor \Triv, an object $E$ is trivial \ssi it is locally trivial for the identity covering of $P(E)$, because a morphism is Cartesian over an identity if and only if it is a vertical isomorphism. Moreover, in the global case, trivial objects have the following characterization: an object $E$ is trivial \ssi it is locally trivial for an isomorphism singleton covering $\{B\ra\cong P(E)\}$. From this point of view, one sees that the question of deciding between fibrewise or global images amounts to the following. Should we consider that having a property locally over an isomorphism is equivalent to having the property itself, or is this true only over an identity?

\item
Choosing the global replete image instead of the replete image of $(F,G)$  is in some sense just a means to gain freedom in the choice of covering functions $K$ on \B that define the locally trivial objects. Indeed, if a covering function $K$ satisfies both axioms ($\tilde{\mathrm{M}}$) and ($\tilde{\mathrm{L}}$), then locally trivial objects with respect to both sorts of images coincide. For, if $E$ is locally trivial with respect to the global replete image of $(F,G)$ for a $K$-covering $R=\{f_i\colon B_i\to B\mid i\in I\}$, then, for each $i\in I$, we have
	$$\xymatrix{
	F(D_i)\ar[r]^{h_i}_\cong & f_i^*E\ar[r]^{\overline{(f_i)}_E} & E\\
	GQ(D_i)\ar[r]^{P(h_i)}_\cong & B_i\ar[r]^{f_i} & P(E).
	}$$
Therefore, $E$ is locally trivial with respect to the replete image for the composite covering $$\{GQ(D_i)\ra{P(h_i)}B_i\ra{f_i}B\mid i\in I\},$$ which is in $K$ under these hypotheses (see in \autoref*{par:schemes} the \hyperref[par:schemes]{example of schemes}).
\end{Rems}

\subsubsection{Morphisms of locally trivial objects}
\label{sssec:morph-locally-triv}
We now turn to the question of morphisms between locally trivial objects. If the subfunctor of trivial objects is globally full, then it is natural just to consider the globally full subfunctor of locally trivial objects. On the other hand, when the subfunctor is not globally full, one might consider morphisms that themselves are ``locally trivial'' in some sense. Let us consider an example: morphisms of locally ringed spaces. One can indeed see locally ringed spaces as locally trivial objects in the following manner.

Let $\Comm_l$ denote the category of local commutative rings and local homomorphisms and let $*$ be a chosen singleton. Recall that $\Sh(*;\Comm)\cong\Comm$ and thus, the fibre $\Ringed_{*}$ of the bifibration $\Ringed\to\Top$ is isomorphic to $\{*\}\times\Comm^{op}$. There is a morphism of fibrations
$$\xymatrix{
\Comm_l^{op}\ar[r]^{({*},Id)}\ar[d] & \Ringed\ar[d]\\
\one\ar[r]_{*} & \Top
}$$
The whole image of this morphism in thus concentrated in the fibre over the singleton *. This morphism of fibrations is biCartesian and (globally) isofull and faithful (for isofullness, see ``Locally ringed spaces'' in \autoref{ssec:Firstexmorfib}). Its global replete image is thus a strong subfibration (\thref{prop:fibpropim}). Its total category has as objects all ringed spaces whose space is a point and sheaf a local commutative ring (by repleteness of $\Comm_l$ in \Comm). Its morphisms are morphisms of ringed spaces $(\{x\},A)\ra{(f,f^\sharp)}(\{y\},B)$ between to such objects such that $f^\sharp\colon B\to A$ is a morphism of local commutative rings.

Now, consider on \Top the coverage whose coverings are $\{\{x\}\hookrightarrow X\mid x\in X\}$. Recall that the inverse image functor over an inclusion of a point $x\in X$ is the stalk functor at $x$. Thus, with the previously defined subfibration of trivial objects, the locally trivial objects in \Ringed in this fibred site are nothing but the locally ringed spaces. The category $\LRinged$ is not full in \Ringed however. Its morphisms are precisely those that are locally trivial, in the following sense. A morphism $(X,\Oc X)\ra{(f,f^\sharp)}(Y,\Oc Y)$ between two locally ringed spaces is locally trivial if the inverse image of the induced vertical morphism $(1_X,f_\sharp)\colon (X,\Oc X)\to (X,f^{-1}\Oc Y)$ along each arrow $\{x\}\hookrightarrow X$, $x\in X$, of the covering of $X$ is trivial (see the following commutative diagram).
$$\xymatrix{
(\{x\},\Oc{X,x})\ar[d]_{(f|_{\{x\}},f_{\sharp,x})}\ar[r] & (X,\Oc X)\ar[r]^{(f,f^\sharp)}\ar@{-->}[d]_{(1_X,f_\sharp)} & (Y,\Oc Y)\\
(\{f(x)\},\Oc{Y,f(x)})\ar[r] & (X,f^{-1}\Oc Y)\ar[ur]&
}$$
Yet, this is precisely the definition of a morphism of locally ringed spaces (see the part ``Locally ringed spaces'' in \autoref{ssec:Firstexmorfib}).

Here is another example, close to the latter. Consider the replete strong subfibration over \B of locally ringed spaces and their morphisms in the fibration of ringed spaces
$$\xymatrix@C=.1cm@R=1.5cm{
\LRinged\ \ar@{^(->}[rr]\ar[dr] && \Ringed\ar[dl]\\
&\Top&
}$$
This subfibration being replete, consider it as the subfibration of trivial objects. We put on \Top the pretopology of open subset coverings. In this situation, one has the peculiar property that all locally trivial objects are trivial. Indeed, the stalk of the restriction $\Oc X|_U$ on an open subset $U\subset X$ at a point $x\in U$ can be chosen to be $\Oc{X,x}$. If one would choose the globally full subfunctor of locally trivial objects, one would not obtain the subfibration $\LRinged\to\Top$ though. One need to consider again the morphisms that are locally trivial. Indeed, in the same way as for objects, a morphism in \Ringed between locally ringed spaces that is locally in \LRinged is so globally.
\paragraph{}
These examples show the utility of having a condition of local triviality also on morphisms. Nevertheless, in the examples that we study, this is rarely needed, or if there is a need to restrict the class of morphisms, this can already be done at the level of the fibred site $P\colon\E\to(\B,K)$. For instance, schemes, as we shall see, are locally trivial ringed spaces with respect to affine schemes. Their morphisms are morphisms of locally ringed spaces, not mere morphisms of ringed spaces (and so one should restrict to locally trivial morphisms in some sense). However, schemes are also locally trivial objects in the fibration of locally ringed spaces (still with respect to affine schemes), and form, in the category of locally ringed spaces, a full subcategory. This is a technical trick though. We really think morphisms of schemes as locally trivial morphisms. 

For these reasons and for the sake of simplicity, since local triviality of morphisms seems to be a quite difficult condition to handle, we restrict up to now our attention to full subfunctor of locally trivial objects. Locally trivial morphisms should be interesting structures to study though, and Proposition \ref{prop:fibpropim}, as well as the preceding examples would be a starting point.

\begin{Def}\indexb{Locally trivial objects!-- in a fibred site}
  Let $P\colon\E\to(\B,K)$ be a fibred site such that $K$ satisfies axiom \tm and let $\Triv\subset P$ be a subfunctor of trivial objects.\\
An object of the base or of the total category is \emph{locally trivial} if there exists a $K$-covering for which it is locally trivial.\\
The restriction of $P$ to the respective full subcategories of locally trivial objects of the base and the total categories determines the \emph{subfunctor of locally trivial objects} $\Loc\subset P$:
\[\xymatrix@=1.3cm{
\Triv_{t}\ar[d]\ar@{^{(}->}[r] &\Loc_{t}\ar@{^{(}->}[r]\ar[d] & \E\ar[d]\\
\Triv_{b}\ar@{^{(}->}[r] & \Loc_{b}\ar@{^{(}->}[r] & \B.
}\]
When we speak about locally trivial objects in a fibred site $(P,K)$, we \emph{always} suppose that $K$ satisfies axiom \tm.
\end{Def}

\begin{Rem}
\label{rem:proploctriv:replete} The subfunctor of locally trivial objects is obviously replete. When is it globally replete? If an object $E'$ is isomorphic, by an isomorphism $h\colon E\ra\cong E'$, to an object $E$ that is locally trivial for a covering $R$, then $E'$ is locally trivial for the covering $P(h)\circ R$. Thus, if $K$ satisfies axioms ($\tilde{\mathrm{M}}$) and ($\tilde{\mathrm{L}}$), then $E'$ is also locally trivial.  
\end{Rem}

\subsection{First examples}\label{ssec:firstexloctriv}
\subsubsection{Extreme cases}
Consider a fibred site $\E\ra P(\B,K)$ equipped with a subfunctor \Triv of trivial objects.
\begin{enumerate}
\item
The no-covering coverage: no object is locally trivial.
\item
The empty-covering coverage: all objects are locally trivial.
\item
The finest pretopology: all objects are locally trivial, because of the preceding example.
\item
Coarsest pretopology: the locally trivial objects are precisely the trivial objects.
\end{enumerate}
Thus, when one makes vary the covering function $K$ from the coarsest to the finest pretopology, one makes vary locally trivial objects from the trivial objects to all objects:
$$
\shorthandoff{;:!?}\xymatrix@!0@C=2cm@R=1.5cm{
\Triv_t\ar[d] \ar@{^{(}->}[r]&\Loc_t \ar@{^{(}->}[r] \ar[d] & \E\ar[d]^{P}\\
\Triv_b \ar@{^{(}->}[r]&\Loc_b \ar@{^{(}->}[r] & \B\\
*[l]{\text{Coarsest pretopology}} \ar@{}[r]|{\preceq} & K \ar@{}[r]|{\preceq} & *[r]{\text{Finest pretopology}.}
}$$
Moreover, we will see further down that the construction of the subfunctor of locally trivial objects preserves, under mild assumptions, the subordination of covering functions (see \thref{prop:ltrivsub}).
\subsubsection{Street's notion of locally trivial objects}
\label{sssec:Street}
We now explain the notion of locally trivial objects introduced by Street in \cite{Str04}. We adapt his notation to ours. The author considers an indexed category $\Phi\colon\B^{op}\to\CAT$ over a finitely complete category \B. Given a singleton covering $\{f\colon A\to B\}$ in \B and an object $T$ (the trivial object) over some object $B_0$, he defines locally trivial objects over $B$ as objects $E$ over $B$ for which there exists an arrow $g\colon A\to B_0$ such that the inverse image of $E$ along $f$ is isomorphic to the inverse image of $T$ along $g$.

These locally trivial objects can be seen in our framework this way. Consider a fibration $\E\to\B$. The subfunctor of trivial objects is defined this way. There is no condition on morphisms in Street's paper and he is so considering a full subfunctor of trivial objects. Trivial objects in the total category are all objects $E\in\E$ such that there exists a morphism $g\colon P(E)\to B_0$ and a Cartesian arrow $E\to T$ over $g$. They thus determine a replete, full, strong subfibration over \B:
$$\xymatrix@C=.5cm@R=1.3cm{
*[l]{Triv\ }\ar[dr]_{\bar P}\ar@{^(->}[rr]^I && \E\ar[dl]^P\\
&\B&
}$$
Observe that, being over \B, this subfibration is automatically globally replete and globally full. 

There is also an external description. Given a cleavage of $P$, one has a Cartesian functor over \B
$$\xymatrix@C=.5cm@R=1.3cm{
*[l]{\B/B_0}\ar[dr]_{\dom}\ar[rr]^F &&\E\ar[dl]^P\\
&\B&
}$$
defined by $F(B\ra{g}B_0)=g^*T$ and by universality on morphisms. Its (global) replete full image is the subfibration previously defined. Notice that in case $B_0=*$ is a terminal object of \B, then trivial objects are characterized by \hyperref[eq:trivterminal]{(\ref*{eq:trivterminal})} for $\C=\{T\}$.

Now, Street's notion of locally trivial objects is precisely our notion of locally trivial object for the covering $\{f\}$ with respect to the previously defined trivial objects. Nevertheless, his notion is different in the sense that he uses a full sub-category $Loc_A(\{f\})$ of locally trivial objects for the fixed covering $\{f\}$, whereas we tend to consider locally trivial objects in a covering function $K$ (and therefore only suppose the existence of a $K$-covering). Moreover, Street seems to only consider singleton coverings. Yet, the notion of local triviality is not necessarily equivalent on a covering and on the singleton it generates via the coproduct of its domains. For instance, the category $\Triv_b$ might not be closed under coproducts. In fact, even when $G=Id_\B$, the two notions might differ, as we shall see later.

\subsubsection{Locally trivial \texorpdfstring{\A}{A}-spaces}

The following examples take place in the fibrations of \A-spaces $\Sh_{\Top}(\A)^{op}\to\Top$ or more specifically of ringed spaces $\Ringed\to\Top$ and of locally ringed spaces $\LRinged\to\Top$. We consider on \Top the pretopologies of open subset coverings and of open embeddings. Example \hyperref[ex:invimsheaves:openembedding]{\ref*{ex:invimsheaves}(\ref*{ex:invimsheaves:openembedding})} describes an adjunction $f^{-1}\dashv f_*\colon\Sh(X;\A)\to\Sh(Y;\A)$ for an open embedding $f$. Combining this example with \hyperref[eq:invimAspace]{(\ref*{eq:invimAspace})}, one obtains that a Cartesian lift of $f$ at an \A-space $(X,\Oc X)$ is given by $(f,\res)\colon(Y,f^{-1}\Oc X)\to (X,\Oc X)$, where the morphism $\res\colon \Oc X\to f_*f^{-1}\Oc X$ is defined by the restriction morphisms of the sheaf $\Oc X$. When $(X,\Oc X)$ is a locally ringed space, then this is also a Cartesian lift in \LRinged, since $\LRinged\to\Top$ is a replete strong subfibration of $\Ringed\to\Top$ (\thref{cor:subfib} and \hyperref[sssec:locringspaces]{``Locally ringed spaces''} in \autoref*{sssec:locringspaces}).

In case of an inclusion of an open subset $U\xhookrightarrow i X$, a Cartesian lift of $i$ at $(X,\Oc X)$ is thus given by $(i,\res)\colon(U,\Oc X|_U)\to(X,\Oc X)$, where $\res\colon\Oc X\to i_*(\Oc X|_U)$ is defined by the restriction morphisms $$\res_{\scriptscriptstyle V,V\cap U}=\Oc X(U\cap V\hookrightarrow V)\colon\Oc X(V)\to\Oc X(U\cap V)$$ of the sheaf $\Oc X$.

\paragraph{Schemes}\label{par:schemes}
Consider the fibration of locally ringed spaces $\LRinged\to\Top$ and the morphism of fibrations
\begin{equation}\label{eq:morfibscheme}\begin{aligned}\xymatrix@=1.3cm{
\Comm^{op} \ar[r]^{(\mathrm{Spec},\Oc{})} \ar[d]_{Id_{\Comm^{op}}}& \LRinged\ar[d]\\
\Comm^{op} \ar[r]_{\Spec} & \Top.
}\end{aligned}\end{equation}
(see \hyperref[sssec:affsch]{``Affine schemes''} in \autoref*{sssec:affsch}). The category of affine schemes is the total category of the global replete full image of this morphism. Consider on \Top the pretopology of open subset coverings. Then, the \emph{category of schemes}\indexb{Scheme}\indexb{Scheme!Category of --s}, written $\Sch$\index[not]{Sch@\Sch}, is precisely the total category of the subfunctor of locally trivial objects in this fibred site.

Note that a locally ringed space is a scheme \ssi it is a locally trivial object in the pretopology of open embeddings. In order to see this, use the fact that an open embedding $f\colon Y\to X$ may be factored as an homeomorphism followed by an open inclusion
$$\xymatrix{
Y\ar[rr]^f\ar[dr]_{\bar f}^\cong &&X\\
&f(Y).\ar@{^(->}[ru]_i&
}$$
Moreover, in the pretopology of open embeddings, a locally ringed space is a scheme \ssi it is locally trivial with respect to the (non-global) replete full image of the morphism \eqref{eq:morfibscheme}, because this pretopology satisfies axioms $(\tilde{\mathrm{M}})$ and $(\tilde{\mathrm{L}})$ (see \thref{rem:proploctriv:replete}). This is an example of a more general fact that we have already pointed out. Defining trivial objects as objects globally isomorphic to the image of the functor $(\mathrm{Spec},\Oc{})$ gives the freedom to deal with the more natural pretopology of open subset coverings. One does in effect not even need to know about Grothendieck topologies in order to understand the definition of a scheme.

A scheme is also a locally trivial object in the fibration of ringed spaces. Indeed, in this situation, trivial objects in \Ringed and in \LRinged coincide. Moreover, affine schemes are locally ringed, and a ringed space that is, in the pretopology of open subset coverings, locally a locally ringed space is so globally, as already remarked. The category of locally trivial objects in this setting does differ at the level of morphisms from the category \Sch though, since \LRinged in not full in \Ringed. 

\paragraph{Locally constant \texorpdfstring{\A}{A}-spaces}
This is an example where trivial objects of the total category are determined by a subcategory of the fibre over a terminal object of the base (the base category of trivial objects is the whole base category of the fibred site). We are now in the fibration of \A-spaces $\Sh_{\Top}(\A)^{op}\to\Top$, endowed with the pretopology of open subset coverings. Trivial objects are determined by the whole fibre of the fibration over a terminal object of \Top, that is, a singleton $*$. It is thus the category $\Sh_{\Top}(\A)^{op}_*\cong\Sh(*;\A)^{op}\cong\A^{op}$. Recall \hyperref[eq:trivterminal]{(\ref*{eq:trivterminal})} that a trivial object is defined in this situation to be an object $(X,\Oc X)$ that is the domain of a Cartesian arrow $(X,\Oc X)\to(*,A)$ over $!\colon X\to *$. In other words, it is an \A-space $(X,\Oc X)$ with $\Oc X\cong\tilde\Delta_X(A)$ in \Sh(X;\A), where $\tilde\Delta\colon\A\to\Sh(X;\A)$ is (a choice of) the constant sheaf functor (see Example \hyperref[ex:invimsheaves:toterminal]{\ref*{ex:invimsheaves}(\ref*{ex:invimsheaves:toterminal})}). A locally trivial object in this context is thus an \A-space $(X,\Oc X)$ satisfying the following condition. There exists an open subset covering $\{U_\alpha\}_{\alpha\in I}$ of $X$ and for each $\alpha\in I$, an object $A_\alpha\in\A$, such that $\Oc X|_{U_\alpha}\cong\tilde\Delta_{U_\alpha}(A_\alpha)$ for each $\alpha\in I$. We call these objects \emph{locally constant \A-spaces}. In the literature, they are called \emph{locally constant sheaves}, because they are not defined in the opfibration of \A-spaces, but in the corresponding opindexed category of \A-sheaves. Yet, beside the change of framework, it is exactly the same concept.

\subsection{First properties}
We turn to the study of the basic properties of the subfunctor of locally trivial objects. Let $(P\colon\E\to(\B,K)$ be a fibred site and $\Triv\subset P$ a subfunctor of trivial objects.

We first study the stability properties of locally trivial objects under inverse images in $P$. These, of course, depends on the stability of trivial objects. When $\Triv$ is a strong subfibration of $P$, since it is replete, any inverse image in $P$ of a trivial object over a morphism in $\Triv_b$ is again trivial (see Corollary \hyperref[cor:subfib:strong]{\ref*{cor:subfib}(\ref*{cor:subfib:strong})}). This is a key property for the development of the theory. The properties of locally trivial objects depend, as we shall see, directly on $K_l$, rather than on $K$.

In the following lemma, the coverings are not supposed to belong to $K$.
\begin{Lem}\label{lem:refineloctriv}
	Suppose $\Triv\subset P$ is a replete globally full strong subfibration of $P$ and let $E$ be a locally trivial object for a covering $R$. Then, it is locally trivial for any refinement of $R$ having trivial domains.
\end{Lem}
\begin{Pf}
	This is a direct consequence of the definition of a refinement and of the fact that, under the hypothesis, trivial objects are stable under inverse images in $P$ over a morphism in $\Triv_b$. We nonetheless give the details, since this lemma is at the heart of the next results. 
	
	Let $S$ be a refinement of $R$ with domains in $\Triv_b$ and $f\colon B'\to B$ in $S$. Then, there exists a $g\colon \tilde B\to B$ in $R$ and a factorization of $f$ through $g$ like in the following diagram.
	$$\xymatrix{
	& \tilde B\ar[d]^{g}\\
B'\ar[r]_{f}\ar@{-->}[ru]^{\exists h} & B
}$$
Therefore, $h$ is in $\Triv_b$, since $\Triv_b$ is full in \B. In the following diagram of Cartesian arrows, $g^*E$ is trivial, because $E$ is locally trivial for $R$ and thus $h^*g^*E$ is also trivial. Furthermore, the top row is Cartesian over $f$. 
$$\xymatrix@C=1.5cm{
h^*g^*E\ar[r]^{\bar h_{g^*E}} & g^*E \ar[r]^{\bar g_E} & E\\
B'\ar[r]^h\ar@/_1pc/[rr]_f & \tilde B \ar[r]^g & B
}$$
\cqfd
\end{Pf}
We point out some interesting consequences of this lemma (we thus suppose its hypothesis is satisfied). 

Firstly, suppose $E$ and $E'$ are locally trivial objects over $B$ of \E for coverings $R$ and $R'$. Suppose, moreover, that $S$ has trivial domains and refines both $R$ and $R'$. Then both $E$ and $E'$ are locally trivial for the covering $S$. In particular, if $K_l$ is a coverage satisfying (L), then, by \thref{lem:comref}, any two locally trivial objects in $K$ over the same object can be made locally trivial for the same $K$-covering.

We now suppose that $\Triv_b=\B$, so that the lemma applies to any refinement. Since a covering $R$ refines the sieve $\overline R$ it generates and that, conversely, $\overline R$ refines $R$, they determine the same locally trivial objects.

Still with the hypothesis that $\Triv_b=\B$, consider a covering $\{B_i\ra{f_i}B\mid i\in I\}$ in \B and suppose that the coproduct $\coprod_{i\in I}B_i$ exists. Let $\coprod B_i\ra{f}B$ denote the induced morphism. Then, if an object $E\in\E$ is locally trivial with respect to the singleton covering $\{\coprod B_i\ra{f}B\}$, it is so with respect to $\{B_i\ra{f_i}B\mid i\in I\}$. The converse is not true in general, but let us first state a criterion that guarantees it.

\begin{Prop}
  Let $P\colon\E\to\B$ be a fibration and $\Triv\subset P$ a subfunctor of trivial objects.
\begin{enumerate}[(i)]
\item
P has global coproducts and \Triv is closed under them.
\item
Cartesian morphisms are preserved by coproducts, in the sense that, if $$\{\bar f_j\colon f_j^*E\to E\mid j\in J\}$$ is a set of Cartesian lifts at $E\in\E$, then the induced morphism ${\coprod_{j\in J}f_j^*E\to E}$ is Cartesian;
\end{enumerate}
Then, if $E\in\E$ is locally trivial for a covering $\{B_i\to B\mid i\in I\}$, it also is for the coproduct covering $\{\coprod_{i\in I} B_i\to B\}$.
\end{Prop}

\begin{Pf}
  Let $E$ be locally trivial for $R=\{B_i\ra{f_i} B\mid i\in I\}$ and $\coprod_{i\in I}B_i\ra{f}B$ the induced morphism. By hypothesis, it has a trivial domain. Let $T_i\ra{h_{i}}E$ a Cartesian lift of $f_i$ at $E$ (with $T_i$ trivial) and $\coprod_{i\in I}T_i\ra{h}E$ the induced morphism. By hypothesis, this morphism is Cartesian and has a trivial domain. Therefore, $E$ is locally trivial for the singleton covering $\{P(\coprod T_i)\ra{P(h)}B\}$. Again, one readily verifies that $\{\coprod_{i\in I}B_i\ra{f}B\}$ refines the latter covering. By \thref{lem:refineloctriv}, $E$ is locally trivial over $\{\coprod_{i\in I}B_i\ra{f}B\}$.
\end{Pf}
	
\begin{Cor}\label{prop:coprodCart}
Let $P\colon\E\to\B$ be a fibration and consider the trivial objects determined by a subcategory \C of the fibre over a terminal object of \B.  If
\begin{enumerate}[(i)]
\item
$P$ has global coproducts;
\item
$\C=\{t\}$ is a singleton;
\item\label{prop:coprodCart:ii}
Cartesian morphisms are preserved by coproducts, in the sense that, if $$\{\bar f_j\colon f_j^*E\to E\mid j\in J\}$$ is a set of Cartesian lifts at $E\in\E$, then the induced morphism ${\coprod_{j\in J}f_j^*E\to E}$ is Cartesian;
\end{enumerate}
then, locally trivial objects for a covering $\{B_i\to B\mid i\in I\}$ are locally trivial for its coproduct covering $\{\coprod_{i\in I} B_i\to B\}$.
\end{Cor}
\begin{Pf}
	We just have to check that $\Triv$ is closed under coproducts in $P$. Yet, $\Triv_{b}=\B$ here. Moreover, by definition of trivial objects in this situation and by hypothesis (\ref{prop:coprodCart:ii}), $\Triv_{t}$ is closed under coproducts in \E.
\end{Pf}

\begin{Ex}
	Let us consider the category $\Top$ and the canonical fibration $${\Top^{\two}\ra{cod}Top}.$$ Conditions (i) and (iii) are satisfied. Thus, the proposition applies to locally trivial bundles with constant fibres. If we let the fibre vary, then it is no longer true. Indeed, consider $\C=\{F_1,F_2\}$ and suppose $F_1$ and $F_2$ have different cardinalities. Let $R=\{X_1\ra{f_1}X,X_2\ra{f_2}X\}$ be a covering of a space $X$ by two non empty open subsets and $\xi=(E\ra{p}X)$ a bundle over $X$ that is locally trivial over $R$. The pullback of $\xi$ along $X_1\coprod X_2\to X$ is the bundle
	$$
	(X_1\times F_1)\textstyle\coprod(X_2\times F_2)\ra{pr\coprod pr}X_1\textstyle\coprod X_2
	$$
	This is not trivial over $X_1\coprod X_2$, because it has different cardinality from $(X_1\coprod X_2)\times F_1$ and $(X_1\coprod X_2)\times F_2$. 
\end{Ex}

We now turn to one of the main questions concerning locally trivial objects: when do they form a strong subfibration of $P$?
\begin{Prop}\label{prop:ltsubfib}
	Suppose $\Triv\subset P$ is a replete globally full strong subfibration of $P$. Suppose $K_l$ is a coverage.
\Par
Then, $\Loc\subset P$ is a replete globally full strong subfibration of $P$.
\end{Prop}
\begin{Pf}
	Let $E$ be locally trivial for a covering $R$, $B=P(E)$, and $g\colon B'\to B$ a morphism in $\Loc_b$. Since $K_l$ satisfies (C), then there exists a $K_l$-covering $R'$ of $B'$ such that the composite covering $g\circ R'$ refines $R$. $E$ is locally trivial for $R$, and is therefore locally trivial for $g\circ R'$ by \thref{lem:refineloctriv}. Let $\bar g_{E}\colon g^*E\to E$ be a Cartesian lift of $g$ at $E$ in the fibration $P$ and $f'\colon\tilde B'\to B'$ an arrow of $R'$. Since the composite of Cartesian arrows $f'^*g^*E\to g^*E\to E$ is Cartesian over $g\circ f'$, $f'^*g^*E$ is trivial. Thus $g^*E$ is a locally trivial object. The morphism $\bar g_E\colon g^*E\to E$, being Cartesian in $P$, is Cartesian in $\Loc$, since $\Loc_t\subset\E$ is full.
\end{Pf}

We now state a result that relates the subfunctors of locally trivial objects with respect to two covering functions when one is subordinated to the other. It is a direct consequence of \thref{lem:refineloctriv} and \thref{prop:ltsubfib}.
\begin{Prop}\label{prop:ltrivsub}
	Suppose \Triv is a replete globally full strong subfibration of $P\colon\E\to\B$. Let $K'$ be a covering function on \B such that $K_l\preceq K'_l$.
	\begin{enumerate}[(i)]
\item
Then, every locally trivial object in $K$ is locally trivial in $K'$. In particular, when ${K_l\equiv K'_l}$, then $K$ and $K'$ have the same locally trivial objects.
\item
If, moreover, $K_l$ and $K'_l$ are coverages , then $\Loc(K)$ is a replete, globally full, strong subfibration of $\Loc(K')$.
\end{enumerate}
	\cqfd
\end{Prop}

\begin{Def}\label{def:fibsitetriv}
  A \emph{fibred site with trivial objects}\indexb{Site!Fibred -- with trivial objects}\indexb{Trivial objects!Fibred site with --} is a triple $(P,K,\Triv)$ where
  \begin{itemize}
  \item $(P,K)$ is a fibred site whose covering function satisfies \tm,
  \item \Triv is a subfunctor of trivial objects that is a globally full, strong subfibration and such that $K_{l}$ is a coverage.
  \end{itemize}
\end{Def}

Suppose there is over the fibration $P$ a fibration $P'\colon\E'\to\E$ and that trivial objects in $P'$ are defined such that $\Triv'_{b}=\Triv_{t}$:
\begin{equation}\label{eq:twolevel}\begin{aligned}\xymatrix@=1.3cm{
\Triv'_{t}\ar[d]_{Q'}\ar[r]^{F'} & \E'\ar[d]^{P'}\\
\Triv_{t}\ar[r]^F\ar[d]_Q & \E\ar[d]^P\\
\Triv_{b}\ar[r]_G & \B
}\end{aligned}\end{equation}
Recall (\autoref{sssec:covinfib}) that a covering function $K$ on \B induces a covering function $K_\E$ on \E. Now, an object $E\in\E$ is locally trivial in $K$ \ssi it can be $K_{\E}$-covered by objects in $\Triv_{t}$. Consequently, $\Loc(P,K,\Triv)_t=\Loc(P',K_\E,\Triv')_b$. One thus also obtains a composable pair of functors of locally trivial objects $$\Loc(P',K_\E,\Triv')_t\to\Loc(P',K_\E,\Triv')_b=\Loc(P,K,\Triv)_t\to\Loc(P,K,\Triv)_b.$$

\begin{Ex}
	Consider the fibration of sheaves of modules over the fibration of ringed spaces $\Oc{}$-$\Mod\to\Ringed\to\Top$ (see Example \hyperref[ex:fibmodules:sheaves]{\ref*{ex:fibmodules}(\ref*{ex:fibmodules:sheaves})}). It restricts to the fibration of sheaves of modules over the subfibration of locally ringed spaces $$\Oc{}\text{-}\Mod_l\to\LRinged\to\Top.$$ Now, there is a 2-level situation
	\begin{equation}\label{eq:omod}\begin{aligned}\xymatrix@=1.3cm{
\Mod_c^{op}\ar[d]\ar[r]^{(\mathrm{Spec},\Oc{},\tilde{\ })} & \Oc{}\text{-}\Mod_l\ar[d]\\
\Comm^{op}\ar[r]^{(\mathrm{Spec},\Oc{})}\ar[d]_{Id} & \LRinged\ar[d]\\
\Comm^{op}\ar[r]_{\mathrm{Spec}} & \Top
}\end{aligned}\end{equation}
Given a commutative ring $A$ and an $A$-module $M$, the $\Oc A$-module $\tilde{M}$ in $\Sh(\Spec A;\Ab)$ is defined on a basic open subset by $\tilde M(D(a))=M_a$, the localization of $M$ at $a\in A$ \cite{Qin02,Har77}. The restriction morphisms are defined via the universality of localization and of extension of scalars, using the fact that if $D(b)\subset D(a)$, then $a$ is inversible in $D(g)$ (see \cite[p.\ 42]{Qin02}).

Given a morphism $(\phi,\alpha)\colon(A,M)\to (B,N)$ in $\Mod_c$, one defines 
$$(\Spec\phi,\Oc\phi,\tilde\alpha)\colon(\Spec B,\Oc B,\tilde N)\to(\Spec A,\Oc A,\tilde M).$$
The morphism $(\Oc\phi,\tilde\alpha)\colon(\Oc A,\tilde M)\to(\Spec\phi_{*}\Oc B,\Spec\phi_{*}\tilde N)$ in $\Mod_{c}(\Sh(X,\Ab))$ is defined locally on each basic open set $D(a),a\in A$, this way. Note that $(\Spec\phi_{*}\tilde N)(D(f))=\tilde N(\Spec\phi^{-1}(D(f)))=\tilde N(D(\phi(f)))=N_{\phi(f)}$. The morphism is defined by universality of localizations of commutative rings and opCartesianness of the morphism of extension of scalar \eqref{eq:opCartmodules}
\[\xymatrix{
(A,M)\ar[rrd]_{(\phi,\alpha)}\ar[r]^-{(\lambda_{a},\eta_{a})} &(A_{a},M_{a}) \ar@{-->}[rr]^{(\phi_{a},\alpha_{a})} && (B_{\phi(a)},N_{\phi(a)})\\
&& (B,N)\ar[ur]_{(\lambda_{\phi(a)},\eta)} &\\
A\ar[drr]_{\phi}\ar[r]^{\lambda_{a}} & A_{a}\ar@{-->}[rr]^{\phi_{a}}&& B_{\phi(a)}\\
&& B\ar[ru]_{\lambda_{\phi(a)}}
}\]

The pretopology $K$ of open subset coverings on \Top induces a covering function $K_{\LRinged}$ on $\LRinged$ that we denote by $K_{L}$. Sheaves of modules that are locally trivial with respect to global replete full image of the morphism of fibrations \hyperref[eq:omod]{(\ref*{eq:omod})} and the topology $K_{L}$ are called \emph{quasi-coherent}. They are \emph{coherent} in one restricts to the subcategory of $\Mod_{c}$ of finitely presented modules.

The upper square in \eqref{eq:omod} induces, via the global replete full image of the morphism, a fibred site with trivial objects when equipped with the covering function $K_{L}$. Indeed, the covering function $K_{L}$ is a pretopology satisfying the stronger forms of axioms ($\tilde{\mathrm L}$) and ($\tilde{\mathrm M}$'), since so does $K$ and by \thref{prop:KandKE}. Moreover, every $K_{L}$-covering of an affine scheme admits a $(K_{L})_{l}$-refinement. Thus, \thref{lem:reftrivbase} applies to show that $(K_{L})_{l}$ satisfies (C) and (L). The global replete full image of the upper square in \eqref{eq:omod} is a strong subfibration by \thref{prop:fibpropim} and the following result.

\end{Ex}
\begin{Lem}
  The upper square in \eqref{eq:omod} is Cartesian.
\end{Lem}
\begin{Pf}
  Let $\phi\colon A\to B$ be a morphism of commutative rings, $M$ an $A$-module and $(\phi,\eta)\colon (A,M)\to (B,M\otimes_{A}B)$ the morphism of extension of scalar \eqref{eq:extscal}. Applying functor $(\mathrm{Spec},\Oc{},\tilde{\ })$, one obtains a morphism 
\[(\Spec\phi,\Oc \phi,\tilde\eta)\colon(\Spec B,\Oc B,\widetilde{M\otimes_{A} B})\to(\Spec A,\Oc A,\tilde M).\]
Consider now the following situation in the fibration $\Oc{}\text{-}Mod_{l}\to\LRinged$:
\begin{equation}
  \label{eq:4}
  \begin{aligned}
     \xymatrix@=1.8cm{
(X,\Oc X,\F)\ar@/^3pc/[rr]^{(f,f^{\sharp},\alpha)}\ar@{-->}[r]_-{(g,g^{\sharp},\beta)} &(\Spec B,\Oc B,\widetilde{M\otimes_{A} B})\ar[r]_-{(\Spec\phi,\Oc \phi,\tilde\eta)} & (\Spec A,\Oc A,\tilde M)\\
(X,\Oc X)\ar[r]^{(g,g^{\sharp})}\ar@/_3pc/[rr]_{(f,f^{\sharp})} & (\Spec B,\Oc B)\ar[r]^{(\Spec\phi,\Oc\phi)} & (\Spec A,\Oc A)
}
  \end{aligned}
\end{equation}
We show that, given an $\Oc X$-module $\F$, the data of a morphism of ringed spaces $$(f,f^{\sharp})\colon(X,\Oc X)\to(\Spec A,\Oc A)$$ into an affine scheme and of a morphism of modules $$(f^{\sharp}_{\Spec A},\alpha_{\Spec A})\colon(A,M)\to(\Oc X(X),\F(X))$$ determines a unique morphism of sheaves of modules $(X,\Oc X,\F)\to(\Spec A,\Oc A,\tilde M)$. Indeed, let $a\in A$, then, in $\Mod_{c}\to\Comm$,
\[\shorthandoff{;:!?}\xymatrix@!0@R=2.5cm@C=6cm{
&(A_{a},M_{a})\ar@{-->}[r]^-{\exists\oc(f^{\sharp}_{D(a)},\alpha_{D(a)})} & (\Oc X(f^{-1}D(a)),\F(f^{-1}D(a)))\\
(A,M)\ar[r]_{(f^{\sharp}_{\Spec A},\alpha_{\Spec A})}\ar[ru]^{(\lambda_{a},\eta_{a})} & (\Oc X(X),\F(X))\ar[ru]_{\quad(\text{res},\text{res})} &\\
&A_{a}\ar[r]^{f^{\sharp}_{D(a)}} & \Oc X(f^{-1}D(a))\\
A\ar[r]_{f^{\sharp}_{\Spec A}}\ar[ru]^{\lambda_{a}} & \Oc X(X)\ar[ru]_{\qquad\text{res}_{D(a)\subset\Spec A}} &\\
}\]
and one checks by universality that these morphisms respect restrictions over ${D(b)\subset D(a)}$.

Now let us come back to \eqref{eq:4}. We first define a morphism of modules $(g^{\sharp}_{\Spec B},\beta_{\Spec B})$ by:
\begin{equation}
  \label{eq:5}
  \begin{aligned}
    \xymatrix@=2cm{
(A,M)\ar[r]^{(\phi,\eta_{\phi}} & (B,M\otimes_{A}B)\ar@{-->}[r]^{(g^{\sharp}_{\Spec B},\alpha^{\sharp}_{\Spec B})} & (\Oc X(X),\F(X))\\
A\ar[r]^{\phi} & B\ar[r]^{g^{\sharp}_{\Spec B}} & X
}
  \end{aligned}
\end{equation}
By the preceding remark, this determines a morphism of sheaves of modules $(g,g^{\sharp},\beta)$ and one can check that it makes the diagram \eqref{eq:4} commute. Moreover, any such morphism of sheaves of modules make diagram \eqref{eq:5} commute and therefore must be the morphism we have defined.
\end{Pf}

\begin{Rem}
  The morphism $(\Spec,\Oc{},\tilde{\ })$ is also Cartesian into the fibration $${\Oc{}\text{-}\Mod\to\Ringed},$$ since $\Oc{}\text{-}\Mod_{l}\to\LRinged$ is a strong subfibration of the latter, by \thref{lem:pullfib}. 
\end{Rem}

\subsection{More examples}\label{ssec:exloctriv}
The following examples all are fibrations $P\colon\E\to\B$ with trivial objects determined by a subcategory $\C\subset\E_*$ of the fibre over a terminal object $*$ of \B (or dually, opfibrations with trivial objects determined by a subcategory of the fibre over the initial object), as explained in \autoref{ssec:moreextrivobjects}. We have proved in the same subsection that trivial objects form then a globally replete and globally full strong subfibration of $P$ over \B.

Recall that trivial objects can be characterized in this case as all objects $E\in\E$ that admits a Cartesian arrow $E\to C$ into an object $C\in\C$. This is equivalent to the following. Given, for each $B\in\B$ and each $C\in\C$, a choice $$\oc_B^*C\xlongrightarrow{\bar\oc_C} C$$ of a Cartesian lift of $\oc\colon B\to*$ at $C$, trivial objects consist in all objects $E\in\E$ isomorphic to a $\oc_B^*C$, and one obtains all of them by restricting to vertical isomorphisms. Depending to the situation, both characterizations are useful.

\subsubsection{Principal \textit{G}-bundles and \textit{G}-torsors} 

Let $\C$ be a category with finite limits and fix a Cartesian monoidal structure $(\C,\times,*)$ on it. For an an internal group $G$ in $\C$, consider the fibration $G$-$Bun(\C)\to\C$ of $G$-bundles \C (Example \hyperref[ex:fibration:gbundles]{\ref*{ex:fibration}(\ref*{ex:fibration:gbundles})}). The trivial objects are determined by the singleton subcategory $\{G\to*\}$ of the fibre over a terminal object $*\in\C$. Note that the latter fibre $G$-$Bun(\C)_*$ is isomorphic to the category $Mod_G$ of modules over the monoid $G$ in \C, which we have called $G$-objects (e.g., $G$-sets, $G$-spaces, \ldots). 

Trivial objects over $C$ are therefore all $G$-bundles isomorphic over $C$ to the product bundle $C\times G\to C$, whose $G$-action is induced by the multiplication of $G$. In other words, they are all $G$-bundles whose total space is a product in \C of $C$ and $G$ (not necessarily the one of the Cartesian monoidal structure) and whose $G$-action is induced by the multiplication of $G$.

Since we are in a situation where the base category has pullbacks, there is no harm in considering the stronger version of a pretopology given in \thref{axiom:strongpretop}. So in this example, pretopology implicitly means this stronger notion. Now, locally trivial objects of this fibration in a subcanonical%
\footnote{One requires the pretopology to be subcanonical in order to have the following fundamental proposition. This proposition might be true under milder conditions. The proof seems to rely mainly on the fact that the coverings of a subcanonical pretopology \emph{collectively reflect isomorphisms}, in the sense that if the pullbacks of a morphism $f$ along each arrow of the covering is an isomorphism, then $f$ is an isomorphism. For now, I haven't yet dug into that.}
pretopology $K$ on \C are called \emph{$G$-torsors \op in $K$\fp} \cite{Vis08}. Vistoli gives an important characterization of $G$-torsor in the latter reference:
\begin{Prop}\label{prop:VisTor}
	Let \C be a category with pullbacks and $K$ a subcanonical pretopology on \C. Then, a $G$-bundle $\xi=(E\ra p B)$ is a $G$-torsor in $K$ \ssi the following conditions are satisfied:
	\begin{enumerate}
\item
There is a $K$-covering $R$ of $B$ that refines the singleton covering $\{E \ra{p}B\}$.
\item
The induced morphism $(pr_1,\kappa)\colon E\times G\to E\times_B E$ in \C in the following diagram
$$\xymatrix{
E\times G\ar@{-->}[dr]^{(pr_1,\kappa)}\ar@/^1pc/[drr]^\kappa\ar@/_1pc/[ddr]_{pr_1}&&\\
&E\times_B\ar[d]\ar[r] E & E\ar[d]^p\\
&E\ar[r]_p & B
}$$
is an isomorphism.
\end{enumerate}
\end{Prop}
The idea for proving these conditions are sufficient is that the morphism $(pr_1,\kappa)$ is in fact a morphism of $G$-bundles over $E$ (where $E\times G$ has the $G$-action induced by multiplication of $G$), and therefore an isomorphism in $G\text{-}Bun(\C)_E$. Thus, $\xi$ is locally trivial over $p$, and then also over $R$, since it refines $p$. The other way is more difficult. See \cite{Vis08}.

In case of the $\C=Top$ with the open subset pretopology, the $G$-torsors are precisely the principal $G$-bundles (as so called usually in the literature). Note that Husemoller, in his classical book \cite{Hus94}, call them locally trivial principal bundles. He defines principal bundles as $G$-bundles satisfying conditions (4) in the next proposition. We now state a new result, that characterizes the $G$-torsors for the open pretopology.
\begin{Prop}
	Let $\C=\Top$. Let $G$ be a topological group and $\xi=(E\ra p B)$ a $G$-bundle. The following conditions are equivalent:
	\begin{enumerate}[(1)]
	
\item
$\xi$ is a $G$-torsor in the pretopology of collectively surjective open maps.

\item
$\xi$ is a $G$-torsor in the pretopology of surjective open maps.

\item The two following conditions are satisfied:
\begin{itemize}
\item
$p\colon E\to B$ is a quotient map,
\item
The induced map $(pr_1,\kappa)\colon E\times G\to E\times_B E$ in
$$\xymatrix{
E\times G\ar@{-->}[dr]^{(pr_1,\kappa)}\ar@/^1pc/[drr]^\kappa\ar@/_1pc/[ddr]_{pr_1}&&\\
&E\times_B\ar[d]\ar[r] E & E\ar[d]^p\\
&E\ar[r]_p & B
}$$
is a homeomorphism.
\end{itemize}

\item
The following three conditions are satisfied:
\begin{itemize}
\item
The action of $G$ on $E$ is free,
\item
The induced map $\bar p$ in
$$\xymatrix{
E\ar[d]_q \ar[dr]^p &\\
E/G \ar@{-->}[r]_{\bar p} & B,
}$$
where $q$ is the quotient with respect to $G$-action, is a homeomorphism,
\item
The function $\delta\colon E\times_B E\to G$ defined by $x\delta(x,y)=y$ is continuous.
\end{itemize}
\end{enumerate}
\end{Prop}
\begin{Pf}
	(2) $\Rightarrow$ (1) is clear.
	
	(1) $\Rightarrow$ (2) because hypotheses of \thref{prop:coprodCart} are verified. First, we verify that coproducts exist in $G$-Bun(\Top). Recall that the product functors $\_\times X\colon\Top\to\Top$ for a space $X$ preserve coproducts. Given $G$-bundles $\xi_i=(E_i\to B_i)$ with $G$-actions $\kappa_i$, this allows us to define a $G$-action $$\textstyle(\coprod E_i)\times G\cong\coprod(E_i\times G)\ra{\coprod \kappa_i}\coprod E_i$$ and provides a coproduct $\coprod\kappa_i\colon\coprod E_i\to\coprod B_i$ of the $G$-bundles $\xi_i$. Now, given any covering $\{B_i\ra{f_i}B\}$ in \Top and a $G$-bundle $E\ra p B$, consider the pullback $G$-bundles 
	$$\xymatrix{
	f_i^*E\ar[r]^{\overline{f_i}}\ar[d]_{\overline{p_i}} & E\ar[d]^p\\
	B_i \ar[r]_{f_i}& B
	}$$
The following square, with horizontal arrows induced by the coproduct, is a pullback in \Top:
	$$\xymatrixnocompile{
	\coprod f_i^* E_i\ar[d]_{\coprod \overline{p_i}}\ar[r]^{\bar f} & E\ar[d]^p\\
	\coprod B_i\ar[r]_{f} & B.
	}$$
One readily verifies that the pullback $G$-action and the coproduct $G$-action defined above coincide.

	(2) $\Rightarrow$ (3): Apply \thref{prop:VisTor}. We get the homeomorphism. Let us verify that $p$ is quotient. We know that there is an open surjective map $f\colon A\to B$ that factors through $p$ via an arrow $g\colon A\to E$. This implies first that $p$ is surjective. Now, let $U\subset B$ be a subset of $B$ such that $p^{-1}(U)$ is open in $E$. Then $g^{-1}(p^{-1}(U))=f^{-1}(U)$ is open. Therefore, $f(f^{-1}(U))=U$ is open, since $f$ is an open map.
	
	(3)$\iff$(4): $(pr_1,\kappa)$ is injective (resp.\ surjective) \ssi the action is free (resp.\ fibrewise transitive). Moreover, $\bar p$ is a quotient map (resp.\ $\bar p$ is injective) \ssi $p$ is a quotient map (resp.\ the $G$-action is fibrewise transitive). Using the fact that an injective quotient map is a homeomorphism, one obtains the first part. We leave to the reader the rest of this equivalence.
		
	(3) and (4) $\Rightarrow$ (2): We leave to the reader the task to prove that $(pr_1,\kappa)$ is a morphism of $G$-bundle. Plus, it is a well-known topological fact that the quotient map $E\ra{q}E/G$ of a $G$-space is open \cite{Hus94}. Consequently, since $\bar p$ is a homeomorphism, $p$ is open. (3) and (4) together show then that $\xi$ is locally trivial for the cover $\{E\ra{p}B\}$, which is a surjective open map.
\end{Pf}

Torsors appear in different areas of mathematics. In case the group $G$ is abelian, they give, under conditions, a ``geometrical'' description of the first cohomology group with coefficients in $G$. When $G$ is not abelian, they provide a way of defining the first \emph{non-abelian cohomology} group. We won't go further on this topic here. We end this part on torsors by noting a relationship between torsors in $\C$ and in the slices $\C/C$ and then giving a few examples of categories where torsors are applied.

The trivial $G$-bundle $\theta_C^G=(C\times G\ra{pr}C)$ has a natural structure of group in $\C/C$. In fact, $G$-bundles over $C$ are precisely the $\theta_C^G$-bundles of the category $\C/C$ over the terminal object $1_C$ of $\C/C$ (see Example \ref{ex:fibmodules}(\ref{ex:fibmodules:gbundles})). Moreover, it is not difficult to see that a $G$-bundle in \C over $C$ is a $G$-torsor in a covering function $K$ \ssi the corresponding $\theta_C^G$-bundle in $\C/C$ over $1_C$ is a $\theta_C^G$-torsor in the slice covering function $K_C$. Note that a torsor over the terminal object in $\C/C$ does not necessarily come from a torsor in \C, because one can consider non trivial groups in $\C/C$.

In the category of differentiable manifolds%
\footnote{As already mentioned, the category \Diff does not have all pullbacks, but it has finite products. One can define the functor $G\text{-}Bun(\Diff)\to\Diff$, but it is not a fibration. In order to define torsors, we only need Cartesian arrows above the arrows of the coverings and over the arrows to the terminal objects. And indeed, these exist.} 
equipped with the open subset pretopology, torsors correspond to the classical notion of smooth principal bundles for a Lie group $G$. On the other hand, Moerdijk in \cite{Moe03} considers torsors in the slice category $\Diff/M$ for a manifold $M$ over the terminal object $1_M$ and in the slice pretopology of the surjective submersion pretopology on $\Diff$\footnote{Use the fact that if $g\circ f$ is a surjective submersion, then so is $g$.}. Torsors also appear in algebraic geometry. One considers torsors in a slice category $\Sch/S$ of schemes over a scheme $S$ over the terminal object $1_S$ in various pretopologies \cite{Mil80,Vis08}. Finally, in topos theory, one considers torsors over the terminal object in the pretopology of epimorphisms\cite{Joh77}. 

\subsubsection{Locally trivial categories}
Consider the fibration $Ob\colon\Cat\to\Set$ of Example \hyperref[ex:fibration:cat]{\ref*{ex:fibration}(\ref*{ex:fibration:cat})}. Consider on \Set the coverage of inclusions of elements, i.e., each set has only one covering, defined by the family ${\{\{x\}\xhookrightarrow{i_x}X\mid x\in X\}}$. Finally, define trivial objects to be induced by the singleton subcategory $\{\one\}$ of the fibre over *. Then, trivial objects are \emph{codiscrete} categories, i.e., categories \C that have exactly one morphism for each pair $(C,C')$ of its objects. The locally trivial objects are therefore categories $\C$ that have only trivial sets of endomorphisms of their objects, i.e., such that $\C(C,C)=\{1_C\}$ for all $C$ in \C. These categories are called indeed \emph{locally trivial categories} in the literature \cite{Til87,PST88} and seem to be of interest mainly in computer science.

\subsubsection{Projective modules over commutative rings}
Consider the category $\Comm$ of commutative rings. In this part, when talking about rings, we therefore mean \emph{commutative} rings. Applying results of \autoref{sec:bifibmodmon} to the symmetric monoidal category of abelian groups and tensor product of abelian groups, one obtains the opfibration of modules over commutative rings, that we denote $$\Mod\to\Comm.$$ Recall that a direct image of a $K$-module $M$ along a ring homomorphism $f\colon K\to L$ is given by \emph{extension of scalars} $M\to M\otimes_KL$. Put on \Comm the Zariski pre-optopology and let trivial objects be induced by the full sub-category of free abelian groups of finite rank. One might as well consider the discrete subcategory $\{\Z^n\mid n\in \N\}$ of $\Mod_\Z$. Then, trivial modules over a ring $K$ are free $K$-modules of finite rank\footnote{%
Recall that non-zero commutative rings have invariant basis number (IBN)\cite{Rot09}. Zero rings $\{0\}$ are terminal objects in the category of (commutative) rings. Thus, they do not have IBN, because $\{0\}\cong\{0\}^n$ as $\{0\}$-modules, for every $n\in\N$. In fact, any $\{0\}$-module $M$ is isomorphic to the (null) module $\{0\}$, because, for every $m\in M$, one has $1\cdot m=m=0\cdot m=0$. In conclusion, all $\{0\}$-modules are trivial over $\{0\}$ in our sense.
}. Locally trivial modules over a ring $K$ are precisely finitely generated projective $K$-modules. This indeed follows from a classical theorem of commutative algebra \cite{Bou98}.
\begin{Thm}
	Let $K$ be a commutative ring and $M$ a $K$-module.\Par
	Then, $M$ is finitely generated projective \ssi there exists a finite set $\{a_i\mid i\in I\}$ of elements of $K$ generating $K$ as an ideal and such that, for each $i\in I$, the $K[a_i^{-1}]$-module $M[a_i^{-1}]$ is free of finite rank.
	
	\cqfd
\end{Thm}
Note that the coverings of the Zariski pre-optopology are not supposed to be finite. Yet, we have remarked in \autoref{sssec:comring} that any of these coverings admits a finite subcovering belonging to the Zariski pre-optopology.

\subsubsection{Locally free modules over ringed spaces}
Consider the restriction of the fibration of sheaves of modules over locally ringed spaces to the category of schemes: $\Oc{}$-$\Mod\to\Sch$. Consider it together with the pretopology of open subscheme coverings. $(\Spec\Z,\Oc\Z)$ is terminal in \Sch.  Let the trivial objects be determined by the modules $(Spec\Z,\Oc\Z,\Oc\Z^n)$. The the locally trivial objects are called \emph{locally free sheaves of modules}.

\subsubsection{Vector bundles}
Recall the fibration of bundles of vector spaces $\mathrm{VBun}\to\Top$ defined in Example \hyperref[ex:fibration:bunvect]{\ref*{ex:fibration}(\ref*{ex:fibration:bunvect})}. Let the trivial objects be determined by the bundles of vector spaces over * given by the vector spaces $\R^n$. Put on \Top the pretopology of open subset coverings. Then, locally trivial objects are \emph{vector bundles}.


\chapter{Modules in a monoidal fibred category}\label{cha:FibMod}
\section{The bifibration of modules over monoids}\label{sec:bifibmodmon}
After recalling some basic material, we explain a classical result : under some mild assumptions, monoids and modules in a monoidal category \V organize into a bifibration 
$${\Mod(\V)\to \Mon(\V)}.$$
We end this section by exploring 2-functoriality of this correspondence between monoidal categories and bifibrations%
\footnote{There are other ways to organize together monoids and modules. In fact, monoids and modules, under the same assumptions but for both right and left tensors $\_\otimes A$ and $A \otimes\_$, also form a bicategory $\mathds{M}od(\V)$ whose 0-cells are monoids, 1-cells bimodules and 2-cells morphisms of bimodules \cite{Lei04,Lac07}. But the monoid morphisms are not present (at least explicitly). One can have these in sight if one consider a \emph{double category} $\frak{Mod}(\V)$ whose horizontal bicategory is $\mathds{M}od(\V)$ \cite{Lei04,PS09,Shu08a}. The latter references give also some discussion on the relation between these three different points of view.}. This is a new insight of the matter as far as we know. 

\subsection{Basic notions}
We recall briefly basic notions and results in order to fix the notations and the terminology and to have them ready to use for the sequel (see e.g. \cite{McL97,Lei04} for more details).

\paragraph{The 2-XL-category of monoidal categories}

A \emph{monoidal category}\indexb{Monoidal category} is sextuplet \linebreak $(\V,\otimes,I,\alpha,\lambda,\rho)$, where \V is a category, $\otimes\colon {\V_\times\V\to\V}$ is a functor (called the \emph{tensor product}), $I$ is an object of \V (called the \emph{unit} of \V, also denoted as a functor $u\colon \one\to\V$), $\alpha_{A,B,C}\colon {(A \otimes B) \otimes C}\ra{\cong}A \otimes(B \otimes C)$ is a natural isomorphism (called the \emph{associator}), and $\lambda_A\colon I \otimes A\ra{\cong}A$ and $\rho_A\colon A \otimes I\ra{\cong}A$ are natural isomorphisms (called \emph{left} and \emph{right unitors}). This data is required to satisfy coherence axioms. When talking about a monoidal category, we freely specify only the triple $(\V,\otimes,I)$ or even just the underlying category \V, if it induces no confusion. A \emph{symmetric monoidal category}\indexb{Monoidal category!Symmetric --} consists of a pair $(\V,\sigma)$ where \V is a monoidal category (whose tensor product we denote $\otimes$) and $\sigma_{A,B}\colon A\otimes B\ra{\cong}B\otimes A$ is an isomorphism natural in $A$ and $B$. This data is subject to coherence axioms.

A \emph{monoidal functor}\indexb{Monoidal functor} from a monoidal category $(\V,\otimes,I)$ to a monoidal category \linebreak $(\V',\otimes',I')$ is a triple $(F, \phi,\psi)$ where $F\colon \V\to \V'$ is a functor, $\phi\colon \otimes'\circ (F \times F)\Rightarrow F\circ \otimes$ is a natural transformation, $\psi\colon I'\to F(I)$ is a morphism in $\V'$ (also written as a natural transformation $\psi\colon u'\Rightarrow F\circ u$). This data is subject to coherence axioms. A \emph{symmetric monoidal functor}\indexb{Monoidal functor!Symmetric --} from $(\V,\sigma)$ to $(\V',\sigma')$ is a monoidal functor $(F,\phi,\psi)\colon\V\to\V'$ such that the following diagram commutes.
$$\xymatrix@=1.5cm{
F(A)\otimes F(B)\ar[d]_{\phi_{A,B}}\ar[r]^{\sigma'_{FA,FB}} & F(B) \otimes F(A)\ar[d]^{\phi_{B,A}}\\
F(A\otimes B)\ar[r]_{F(\sigma_{A,B})} & F(B\otimes A)
}$$

A (symmetric) monoidal functor $(F, \phi,\psi)$ is \emph{strong}\indexb{Monoidal functor!Strong (symmetric) --} (resp. \emph{strict}\indexb{Monoidal functor!Strict (symmetric)}) when its structure morphisms \(\phi\) and \(\psi\) are isomorphisms (resp. identities).

A \emph{monoidal natural transformation}\indexb{Monoidal natural transformation} from a monoidal functor $(F,\phi,\psi)$ to a monoidal functor $(F',\phi',\psi')$ is a natural transformation of the underlying functors $\alpha\colon F\Rightarrow F'$ subject to coherence axioms with respect to monoidal structures. We state them in an object-free manner, in order to have them ready for the more general context of monoidal objects in a monoidal 2-category.
\begin{equation}\label{eq:monnattransf:phi}\begin{aligned}
\xymatrixnocompile@C=1.5cm@R=1cm{
\V \times \V\ar@/_1pc/[dr]_{\otimes}\rtwocell<5>^{F \times F}_{F' \times F'}{\omit}\ar@{}[r]|(.5){\text{\small$\Downarrow$} \alpha \times \alpha} & \V' \times \V'\ar[r]^(.6){\otimes'} & \V'\\
 \rrtwocell<\omit>{<-4>\ \ \phi'}& \V\ar@/_1pc/[ur]_{F'} & }
=
\xymatrixnocompile@C=1.5cm@R=1cm{
\V \times \V\ar[r]^{F \times F}\ar@/_1pc/[dr]_{\otimes} & \V' \times \V' \ar[r]^(.6){\otimes'} & \V'\\
\rrtwocell<\omit>{<-4>\ \ \phi} & \D\urtwocell^F_{F'}{\alpha} &}
\end{aligned}\end{equation}
and
\begin{equation}\label{eq:monnattransf:psi}\begin{aligned}
\xymatrix@=1.5cm{
\one\ar[r]^u\ar[dr]_{u'} & \V\dtwocell^{F'}_F{^\alpha}\\
& \V'\ultwocell<\omit>{<3>\psi}}
=
\xymatrix@=1.5cm{
\one\ar[r]^u\ar[dr]_{u'} & \V\ar[d]^{F'}\\
& \V'\ultwocell<\omit>{<3>\psi'}}
\end{aligned}\end{equation}

Monoidal categories, monoidal functors and monoidal natural transformations form a 2-XL-category, which we denote $\MONCAT$\index[not]{MONCAT@\MONCAT}. One also consider its 2-cell-full sub-2-XL-category $\MONCAT_s$\index[not]{MONCATs@$\MONCAT_s$} of monoidal categories and strong monoidal functors. This remains true for the corresponding symmetric structures: symmetric monoidal categories, (strong) symmetric monoidal functors and monoidal natural transformations form a 2-XL-category, which we denote $\SYMMON$\index[not]{SYMMON@\SYMMON} (resp. $\SYMMON_s$\index[not]{SYMMONs@$\SYMMON_s$}).
\paragraph{The category of monoids in \V}

A \emph{monoid}\indexb{Monoid in a monoidal category} in \V is a triple $(R,\mu,\eta)$ where $R$ is an object of \V, and $\mu\colon R \otimes R\to R$ and $\eta\colon I\to R$ are morphisms in \V (called the \emph{product} and the \emph{unit} of the monoid). This data must satisfy the associativity and the unit axioms, i.e., the following diagrams must commute.\bigskip
\begin{center}
$\xymatrix@C=0pc@R=4pc{
&**[l]{(R \otimes R)\otimes R}\ar[rr]^{\cong}_{\alpha_{RRR}}\ar[ld]_{\mu\otimes 1_{R}} &&
**[r]{R \otimes(R\otimes R)}\ar[rd]^{1_{R}\otimes\mu} & \\
	R\otimes R\ar[rrd]_{\mu} &&&&	R\otimes R\ar[lld]^{\mu}\\
&&R&&
	}\quad
\xymatrix@=1.5cm{
I\otimes R\ar[r]^{\eta\otimes 1_{R}}\ar[dr]_{\cong}^{\lambda_R} &R\otimes R\ar[d]^{\mu} & R\otimes I\ar[l]_{1_{R}\otimes \eta}\ar[dl]^{\cong}_{\rho_R}\\
&R&
}$
\end{center}
A monoid $(M,\mu,\eta)$ in a symmetric monoidal category $((\V,\otimes,I),\sigma)$ is \emph{commutative}\indexb{Monoid in a monoidal category!Commutative --} if the following diagram commutes.
$$\xymatrix{
M\otimes M \ar[dr]_{\mu}\ar[r]^{\sigma_{M,M}}_{\cong} & M\otimes M\ar[d]^{\mu}\\
&M
}$$

\emph{Morphisms of monoids}\indexb{Monoid in a monoidal category!Morphism of --} are morphisms of the underlying objects of \V commuting with the product and the unit. Monoids in \V together with their morphisms and composition and identities of $\V$ form a category that we denote $\Mon(\V)$\index[not]{MonV@$\Mon(\V)$}. When \V is symmetric, we denote $Comm(\V)$\index[not]{CommV@$\Comm(\V)$} its full subcategory of commutative monoids.
\subparagraph{}
It is an important feature of monoidal functors that they preserve monoids and their morphisms. 

Let $(F,\phi,\psi)\colon \V\to \V'$ be a monoidal functor and $(R,\mu,\eta)$ be a monoid in \V. Then $F$ induces a monoid structure on $F(R)$ defined by 
$$F(R)\otimes F(R)\ra{\phi_{R,R}}F(R \otimes R)\ra{F(\mu)}F(R),\quad\quad I'\ra{\psi}F(I)\ra{F(\eta)}F(R).$$
Moreover, given a morphism of monoids $F\colon R\to S$ in \V, the morphism $F(f)$ in $\V'$ is a morphism of monoids for the induced structures. Therefore, one obtains a functor 
$$\Mon(F)\colon \Mon(\V)\to \Mon(\V').$$

Similarly, symmetric monoidal functors preserve commutative monoids. If $F$ is such a functor, then we denote $Comm(F)$ the restriction of $\Mon(F)$ to $Comm(\V)$:
$$
Comm(F)=Mon(F)\mid_{Comm(\V)}\colon Comm(\V)\to Comm(\V).
$$

Monoidal natural transformations also behave well with respect to monoids. Their components are indeed automatically monoid morphisms for the monoid structures induced by the functors. In consequence, a monoidal natural transformation $\tau\colon F\Rightarrow F'$ determines a natural transformation 
$$\Mon(\tau)\colon \Mon(F)\Rightarrow \Mon(F'),$$ 
with $\Mon(\tau)_R=\tau_R$. One defines in the same manner $Comm(\tau)$ for a natural transformation $\tau$ between symmetric monoidal functors.
\begin{Prop}\label{prop:Mon}
	The correspondences described above determine 2-functors
$$Mon\colon \MONCAT\longrightarrow \CAT$$
and
$$
Comm\colon \SYMMON\longrightarrow \CAT
$$
\cqfd
\end{Prop}
This is what happens on the ground. Now let us gain some altitude.

\paragraph{The category of (right) modules over a monoid in \V}

A \emph{right $R$-module}\indexb{Module in a monoidal category} is a pair $(M,\kappa)$ where $M$ is an object of \V and $\kappa\colon M\otimes R\to M$ is morphism in \V satisfying the following commutative diagram.
$$\xymatrix@=1.5cm{
(M\otimes R) \otimes R\ar[r]^{\cong}_{\alpha_{MRR}}\ar[d]_{\kappa\otimes 1_R} & M \otimes(R \otimes R)\ar[r]^{1_M \otimes\mu} & M \otimes R\ar[d]^{\kappa} & M\otimes I\ar[l]_{1_M \otimes \eta}\ar[dl]^{\rho_M}\\
M \otimes R\ar[rr]_{\kappa} &&M
}$$
A \emph{morphism of right $R$-modules}\indexb{Module in a monoidal category!Morphism of --} is a morphism in \V of the underlying \V-objects that commutes with the actions of $R$. Again, with composition and identities of \V, they altogether form the category of right $R$-modules that we denote $\Mod_R$\index[not]{ModR@$\Mod_R$}. In a similar manner, one defines the category $\leftidx{_R}{\Mod}{}$\index[not]{ModR2@$\leftidx{_R}{\Mod}{}$} of left $R$-modules.

\paragraph{The category of bimodules over a pair of monoids (R,S) in \V}

Given monoids $R$ and $S$, a \emph{$(R,S)$-bimodule}\indexb{Bimodule in a monoidal category} is a triple $(M,\kappa, \sigma)$ where $(M,\kappa)$ is a left $R$-module, $(M,\sigma)$ is a right $S$-module satisfying the following coherence axiom.
$$\xymatrix@C=.3pc@R=5pc{
&**[l]{(R \otimes M)\otimes S}\ar[rr]^{\cong}_{\alpha_{RMS}}\ar[ld]_{\kappa\otimes 1_{S}} &&
**[r]{R \otimes(M\otimes S)}\ar[rd]^{1_{R}\otimes\sigma} & \\
	M\otimes S\ar[rrd]_{\sigma} &&&&	R\otimes M\ar[lld]^{\kappa}\\
&&R&&
}$$
A \emph{morphism of bimodules}\indexb{Bimodule in a monoidal category!Morphism of --} is a morphism of the underlying \V-objects that is both a left module and a right module morphism. Together with bimodules, and composition and identities of \V, they constitute the category of $(R,S)$-bimodules, denoted $\leftidx{_R}{\Mod}{_S}$.
\subparagraph{}
Monoidal functors handle modules well also. Let $(F,\phi,\psi)\colon \V\to \V'$ be a monoidal functor, $(R,\mu,\eta)$ a monoid in \V and $(M,\kappa)$ a right $R$-module. $F$ induces a right $\Mon(F)(R)$-module structure on $F(M)$ by
\begin{equation}\label{eq:ModF}
F(M)\otimes F(R)\ra{\phi_{M,R}}F(M \otimes R)\ra{F(\kappa)}F(M).
\end{equation}
Moreover, if $f\colon M\to N$ is an $R$-module morphism, then $F(f)$ is a $\Mon(F)(R)$-module morphism for the induced structures on $F(M)$ and $F(N)$. 

The relationship between monoidal functors, monoidal natural transformations and modules is understood more deeply and more easily in the bifibred setting.
\subsection{The global view}\label{ssec:modmon}

All (right) modules over any monoid in \V organize into a bifibration over the category of monoids $\Mon(\V)$, as we now explain.

In this part, $(R,\mu,\eta)$ and $(S,\nu,\epsilon)$ are monoids in \V and $f\colon R\to S$ is a morphism of monoids.
\paragraph{Restriction of scalars}\label{par:resscal}

The morphism $f\colon R\to S$ of monoids  induces a functor, called \emph{restriction of scalars}\indexb{Restriction of scalars}\index[not]{fsharp1@$f^\sharp$}\footnote{This terminology is motivated by the case of rings and modules over them. If $R$ is a field, a morphism of rings $f\colon R\to S$ has to be injective, unless $S$ is a trivial ring.},
$$\begin{array}{rcl}
f^{\sharp}\colon \Mod_S & \longrightarrow & \Mod_R\\
(M,\kappa) & \mapsto & (M,M \otimes R\ra{1_M\otimes f}M \otimes S\ra{\kappa}M),
\end{array}$$
which is identity on morphisms, and similarly for left modules.

Moreover, this functor preserves bimodule structures and so induces, for a monoid $T$, a functor 
$$f^{\sharp}\colon \leftidx{_T}{\Mod}{_S}\to\leftidx{_T}{\Mod}{_R},$$
and similarly on the left.

\paragraph{Extension of scalars}

We describe a functor induced by the morphism of monoids $f$ going the direction opposite to the restriction of scalars functor, called \emph{extension of scalars}.\indexb{Extension of scalars}

\emph{Suppose that \V has all reflexive coequalizers}. Then one defines a tensor product $M\otimes_R N$ of a right module $(M,\kappa)$ and a left module $(N,\sigma)$ over a monoid $R$ by the (reflexive) coequalizer:
$$\xymatrix@C=4pc{
(M\otimes R) \otimes N\ar@<.7ex>[r]^-{\kappa\otimes 1_N}\ar@<-.7ex>[r]_-{(1_M \otimes\sigma)\circ\alpha_{MRN}} & M \otimes N\ar[r]^{coeq} & M \otimes_R N.
}$$
A choice of such a coequalizer for each pair of modules determines a functor 
$$\Mod_R\times \leftidx{_R}{\Mod}{}\to\V.$$
A quite long calculation shows that if, in addition, the monoidal endofunctors $-\otimes A$ of \V preserve reflexive coequalizers (which is true, e.g., if \V is right closed), then it induces a functor
\begin{equation}\label{modtens}\otimes_R\colon \Mod_R\times\leftidx{_R}{\Mod}{_S}\to \Mod_S.\end{equation}
Given a right $R$-module $(M,\kappa)$ and a $(R,S)$-bimodule $(N,\sigma,\tau)$, the right $S$-action $\tilde\tau$ on $M\otimes_R N$ is defined by the following diagram
\begin{equation}\begin{aligned}\label{eq:inducedactiontensor}\xymatrix@C=1.5cm{
(M\otimes N)\otimes S\ar[r]^{coeq\otimes 1_S}\ar[d]_{\alpha_{MNS}} & (M\otimes_R N)\otimes S\ar@{-->}[dd]^{\tilde\tau}\\
M\otimes(N\otimes S)\ar[d]_{1_M\otimes\tau} &\\
M\otimes N\ar[r]_{coeq} & M\otimes_R N,
}\end{aligned}\end{equation}
thanks to the fact that the functor $-\otimes S$ preserves reflexive coequalizers.

Now, a monoid $S$ is canonically a $(S,S)$-bimodule. Applying to it the functor 
$$f^{\sharp}\colon \leftidx{_S}{\Mod}{_S}\to\leftidx{_R}{\Mod}{_S},$$
one gets a $(R,S)$-bimodule $f^*(S)$, which, by abuse of notation, we simply denote $S$ again.

The bifunctor \eqref{modtens} restricted to this $(R,S)$-bimodule gives finally the extension of scalars functor:\index[not]{fsharp2@$f_\sharp$}
$$f_{\sharp}:=-\otimes_R S\colon \Mod_R\to \Mod_S.$$

\paragraph{An adjunction}

\begin{Prop}\label{prop:adjointmodules}
	Suppose that \V is a monoidal category with reflexive coequalizers, and that the latter are preserved by the functors $-\otimes A$, for all $A\in\V$.\\
Let $f\colon (R,\mu,\eta)\to (S,\nu,\epsilon)$ be a morphism of monoids in \V.\medskip\\
Then the extension of scalars functor relative to $f$ is left adjoint to the restriction of scalars functor:\medskip
\begin{equation}\label{adjointmodules}\xymatrix@C=4pc@1{
f_{\sharp}:\Mod_R\ar@<-.9ex>[r]\ar@{}[r]|{\top}& \Mod_S : f^{\sharp}\ar@<-.9ex>[l]
}\medskip\end{equation}
Its unit $\eta^f$ is given by ``tensoring by the unit'':
\begin{equation}
  \label{eq:extscal}
  M\ra[\rho_M]{\cong}M \otimes I\ra{1_M \otimes\epsilon}M \otimes S\ra{coeq}M \otimes_R S.
\end{equation}
\end{Prop}
A proof of this proposition can be found in Bruno Vallette's thesis \cite{Val03}.\cqfd

\paragraph{The bifibration}

\emph{Under the conditions of \thref{prop:adjointmodules}} and by means of its adjunctions, we can now construct the \emph{bifibration of right modules over monoids in \V}. It is not difficult to prove that restrictions of scalars give rise to a 2-functor:
$$\begin{array}{rcl}
Mod\colon \Mon(\V)^{op} & \longrightarrow & \CAT\\
R & \longmapsto & \Mod_R\\
R\ra{f}S &\longmapsto & \Mod_S\ra{f^{\sharp}}\Mod_R
\end{array}$$
Via the Grothendieck construction, one obtains a (split) fibration 
\begin{equation}\label{eq:bifibmodmon}
\Mod(\V)\to \Mon(\V).
\end{equation} 
Due to the splitness of the fibration and the particularity of the restriction of scalars functors, the category $\Mod(\V)$ takes a simple shape:
\begin{itemize}
\item
\textbf{Objects of $\Mod(\V)$:} Pairs $(R,M)$ whose first member is a monoid and second is a right module over this monoid.
\item
\textbf{Morphisms of $\Mod(\V)$:} Pairs $(f,\phi)\colon (R,M)\to (S,N)$, where $f\colon R\to S$ is a morphism of monoids in \V and $\phi\colon M\to N$ is a morphism in $\V$ such that the following diagram commutes, where $\kappa$ and $\sigma$ denote the corresponding actions of $R$ and $S$ on $M$ and $N$.
$$\xymatrix@=1.2cm{
M \otimes R\ar[d]_{\kappa}\ar[r]^{\phi \otimes f} & N \otimes S\ar[d]^{\sigma}\\
M\ar[r]_{\phi} & N
}$$
\item
\textbf{Composition and identities:} Those of $\V \times \V$.
\end{itemize}
The Grothendieck construction determines a canonical cleavage. The Cartesian lift of a morphism of monoids $f\colon R\to S$ at an $S$-module $(S,N)$ is given by
$$
(R,f^{\sharp}N)\ra{(f,id)}(S,N).
$$
Now, by \thref{prop:bifibadjoint} and under its hypotheses, the adjunctions \eqref{adjointmodules} define on the fibration 
$$
\Mod(\V)\to \Mon(\V)
$$
a structure of (non split) opfibration. An opCartesian morphism with source $(R,M)$ over a morphism of monoids $f\colon R\to S$ is given by:
\begin{equation}\label{eq:opCartmodules}(R,M)\ra{(f,\eta^f_M)}(S,M \otimes_R S).\end{equation}
For a symmetric monoidal category \V, we are sometimes interested in the restriction 
\begin{equation}\label{eq:bifibmodcomm}
\Mod(\V)|_{\Comm(\V)}\to\Comm(\V)
\end{equation}
of this bifibration to the full subcategory $\Comm(\V)$ of commutative monoids and the full subcategory of modules over them. It is also a bifibration, because it is a pullback of \eqref{eq:bifibmodmon} along the inclusion $\Comm(\V)\hookrightarrow\Mon(\V)$ (see Lemmas \ref{lem:pullfib} and \ref{lem:pullopfib}).

\begin{Not}
	$\Mod_c(\V):=\Mod(\V)|_{\Comm(\V)}$
\end{Not}

\paragraph{The 2-functor ``Modules over monoids''}
We are now in the right context for talking about the relationship between monoidal functors between monoidal categories \V and $\V'$, natural transformations between them and modules in the respective monoidal categories. 

Let $F\colon\V\to\V'$ be a monoidal functor. It determines a functor
$$
\Mod(F)\colon \Mod(\V) \longrightarrow \Mod(\V').
$$
The image of an $R$-module $(M,\kappa)$ is the $\Mon(F)(R)$-module structure on $F(M)$ defined by \hyperref[eq:ModF]{(\ref*{eq:ModF})}. The image of a morphism of modules $(f,\phi)\colon(R,M)\to(S,N)$ in \V is the pair $(F(f),F(\phi))$.

Let now $$\xymatrix{\V\rtwocell^{F}_{F'}{\alpha}& \V'}$$ be a monoidal natural transformation. It determines a natural transformation 
$$\xymatrix@=1.5cm{\Mod(\V)\rtwocell^{\Mod(F)}_{\Mod(F')}{\omit \Downarrow \Mod(\alpha)} &\Mod(\V')}$$
by $$\Mod(\alpha)_{(R,M)}\colon (F(R),F(M))\ra{(\alpha_R,\alpha_M)}(F'(R),F'(M)).$$

\begin{Not}\index[not]{RCMONCAT@$\rcMONCAT$}\index[not]{RCSYMMON@$\rcSYMMON$}
	The 2-cell full sub-XL-2-category of $MONCAT_s$ consisting of monoidal categories having reflexive coequalizers and such that the right tensors $-\otimes A$, $A\in\V$ preserve them, and strong monoidal functors that preserve reflexive coequalizers, is denoted $$\rcMONCAT.$$
	We similarly write $\rcSYMMON$.
\end{Not}

We can now complete the picture of \thref{prop:Mon}. The following result is new, to the best of our knowledge.

\begin{Prop}\label{prop:bifib}
	The correspondence described above yields a 2-functor between the 2-XL-categories of monoidal categories and monoidal functors and of fibrations and Cartesian functors:
	$$\begin{array}{rcl}
	\mathcal{M}\colon\MONCAT& \longrightarrow & \FIB_c\\
	\V &\longmapsto & \xymatrix{\Mod(\V)\ar[d]\\ \Mon(\V)}\\
	\xymatrix{\V\rtwocell^{F}_{F'}{\alpha} & \V'} & \longmapsto & \xymatrix@=1.5cm{\Mod(\V)\ar[d]\rtwocell^{\Mod(F)}_{\Mod(F')}{\omit \Downarrow \Mod(\alpha)} &\Mod(\V')\ar[d]\\
								\Mon(\V) \rtwocell^{\Mon(F)}_{\Mon(F')}{\omit \Downarrow \Mon(\alpha)} & \Mon(\V').}
	\end{array}$$
	It restricts to a 2-functor
	$$
	\leftidx{_{\mathit{rc}}}{\mathcal M}{}\colon\rcMONCAT \longrightarrow \BIFIB_{bc}.
	$$
	Moreover, this induces, by restriction to commutative monoids, 2-functors $$\M_c\colon\SYMMON\to\FIB_c\quad\text{and}\quad\leftidx{_{\mathit{rc}}}{\mathcal M}{_c}\colon\rcSYMMON\to\BIFIB_{bc}.$$
\end{Prop}
\begin{Pf}
	The corresponding statements for commutative monoids being a direct consequence of the general case, it remains to show the existence of the 2-functor $$\leftidx{_{\mathit{rc}}}{\mathcal M}{}\colon\rcMONCAT \longrightarrow \BIFIB_{bc}.$$ We already know that the fibration $\Mod(\V)\to\Mon(\V)$ is a bifibration when \V is in \rcMONCAT. Thus, we have to prove that $(\Mod(F),\Mon(F))$ is an opCartesian morphism of opfibrations when $(F,\phi,\psi)\colon\V\to\V'$ is a strong monoidal functor that preserves reflexive coequalizers. 
	
	Let $f\colon R\to S$ be a morphism of monoids in \V and $(R,M)$ an $R$-module. An opCartesian lift of $f$ at $(R,(M,\kappa))$ is given by $(f,\eta^f_M)\colon(R,M)\to(S,M\otimes_RS)$ (see \hyperref[eq:opCartmodules]{(\ref*{eq:opCartmodules})}). Let us denote by $\sigma$ the left $R$-action on $S$ induced by $f$. Consider the following diagram in $\V'$.
\begin{equation}\label{eq:morphismk}\begin{aligned}\xymatrix@C=1.8cm{
(F(M)\otimes F(R))\otimes F(S)\ar[d]_{\phi\otimes 1_{F(S)}}\ar@<.8ex>[r]^-{(F(\kappa)\circ\phi)\otimes 1_S}\ar@<-.8ex>[r]_-{[1_M \otimes (F(\sigma)\circ\phi)]\circ\alpha'} & F(M) \otimes F(S) \ar[dd]_{\phi}\ar[r]^{coeq} & F(M) \otimes_{F(R)} F(S)\ar@{-->}[dd]^{k}\\
F(M\otimes R)\otimes F(S)\ar[d]_{\phi} &&\\
F((M\otimes R)\otimes S)\ar@<.8ex>[r]^-{F(\kappa\otimes 1_S)}\ar@<-.8ex>[r]_-{F(1_M \otimes\sigma)\circ F(\alpha)} & F(M \otimes S)\ar@{}[luu]|*+[o][F]{1}\ar[r]^{F(coeq)} & F(M \otimes_R S)
}\end{aligned}\end{equation}
The diagram $\xymatrix{*+[o][F]{1}}$ is composed of two squares, one formed of the upper arrows and the other of the lower arrows. The former commutes by naturality of $\phi$. The latter, by naturality of $\phi$ and by the associativity coherence axiom of $F$. Consequently, the morphism $F(coeq)\circ\phi$ induces a unique morphism $k$ making the right square of the diagram commutes. 

Again, the morphism $k$ is a morphism of $F(S)$-modules. We won't draw the diagram, but one checks this fact by using definition \hyperref[eq:inducedactiontensor]{(\ref*{eq:inducedactiontensor})} of the induced action on a tensor product over a ring of a module and a bimodule. One also uses the preservation of reflexive coequalizers by right tensors, naturality of $\phi$ and the associativity axiom for $F$.

Now, one readily verifies that the morphism $(1_{F(S)},k)$ in $\Mod(\V')$ is the unique vertical morphism that makes the following diagram commutes:
$$\xymatrix@C=2cm@R=1.3cm{
(F(R),F(M))\ar[r]^-{(F(f),\eta^{F(f)}_{F(M)})}\ar[dr]_{(F(f),F(\eta^f_M))\ } & (F(S),F(M)\otimes_{F(R)}F(S))\ar@{-->}[d]^{(1_{F(S)},k)}\\
& (F(S),F(M\otimes_RS))\\
F(R)\ar[r]_{F(f)} & F(S)
}$$
The proof uses again the naturality of $\phi$, as well as the right unit coherence axiom for $F$.

Finally, since $F$ is strong monoidal, the vertical arrows in \hyperref[eq:morphismk]{(\ref*{eq:morphismk})} are isomorphisms. Moreover, as $F$ preserves reflexive coequalizers, the lower part of \hyperref[eq:morphismk]{(\ref*{eq:morphismk})} is a coequalizer. In conclusion, $k$ is an isomorphism in \V and thus, $(1_{F(S)},k)$ is an isomorphism in $\Mod(\V)$.
\end{Pf}

\begin{Rem}\label{rem:modmonadjfunctor}
	2-functors preserve adjunctions and therefore, given an adjunction 
	\begin{equation*}\xymatrix@1@C=1.3cm{
	F:\V \ar@<-.9ex>[r]\ar@{}[r]|{\top} & \V':G\ar@<-.9ex>[l]
}\end{equation*}
in \MONCAT, one obtains an adjunction in $\FIBc$:
\begin{equation}\label{eq:modmonadj}\begin{aligned}\xymatrix@=1.5cm{
\Mod(\V)\ar[d]_P \ar@<-.9ex>[r]_{\Mod(F)}\ar@{}[r]|{\top} & \Mod(\V')\ar[d]^{P'}\ar@<-.9ex>[l]_{\Mod(G)}\\
\Mon(\V) \ar@<-.9ex>[r]_{\Mon(F)}\ar@{}[r]|{\top} & \Mon(\V'),\ar@<-.9ex>[l]_{\Mon(G)}
}\end{aligned}\end{equation}
where we have named by $P$ and $P'$ the usual fibrations of modules over monoids in \V and $\V'$ respectively. Let us unpack a little this notion. It means that the squares determined respectively by $F$ and $G$ are morphisms in $\FIB_c$, i.e., Cartesian morphisms of fibrations, that $F\dashv G$ induces two adjunctions, one at the level of monoids, the other at the level of modules, and finally that these adjunctions are coherent in the following sense. Let $\eta$ and $\epsilon$ be the unit and counit of $F\dashv G$. Then, the adjunction at the level of monoids has unit $\Mon(\eta)$ and counit $\Mon(\epsilon)$, and similarly for modules. Now, these adjunctions are coherent in the sense that $$P\cdot \Mod(\eta)=\Mon(\eta)\cdot P\quad\text{and}\quad P'\cdot \Mod(\epsilon)=\Mon(\epsilon)\cdot P'.$$

In fact, a 2-functor between 2-(XL-)categories $\Phi\colon\dA\to\dB$ induces a 2-functor between the corresponding 2-(XL-)categories of adjunctions $\ADJ(\Phi)\colon\ADJ(\dA)\to\ADJ(\dB)$\index[not]{ADJ(Phi)@$\ADJ(\Phi)$}. In particular, one obtains here a 2-functor
$$
\ADJ(\M)\colon\ADJ(\MONCAT)\to\ADJ(\FIB_c).
$$

Suppose now that $\V$, $\V'$ and $F$ are in $\rcMONCAT$\footnote{As a left adjoint, $F$ automatically preserves reflexive coequalizers, and thus, the only requirement on $F$ is that it be strong.} (and $G$ any monoidal functor). Then, in \hyperref[eq:modmonadj]{(\ref*{eq:modmonadj})}, $P$ and $P'$ are bifibrations, the morphism of bifibrations $$(\Mod(F),\Mon(F))\colon P\to P'$$ is opCartesian (and thus biCartesian) and the morphism of bifibrations $$(\Mod(G),\Mon(G))\colon P'\to P$$ is Cartesian. Thus, the restriction of the 2-functor $\ADJ(\M)$ to the sub-2-XL-category of $\ADJ(\MONCAT)$ whose monoidal categories are in $\rcMONCAT$ and whose morphisms have a left adjoint that is strong monoidal takes values in \BIFIBADJ.
\end{Rem}

\section{Internal rings and modules}

In this thesis, we are mainly interested in monoids and modules in a monoidal category. Indeed, as we have seen in the preceding section, they provide, under some mild assumptions, a bifibration. On the other hand, internal modules over internal rings do in general yield only a fibration. Recall that the opfibration property of the bifibration of modules over monoids makes use of the monoidal product. This monoidal product is not available in general in the internal setting. The fibration is enough in some cases, like vector bundles. But the opfibred structure is needed for sheaves of modules on the contrary. 

Yet, many important examples of monoidal categories are in fact categories of internal abelian groups and their monoids and modules have a much more concrete definition as internal rings and modules. Moreover, categories of internal abelian groups have properties, as being additive, and sometimes even abelian. Thus, internal rings and modules are of interest for our theory. The main goal of this section is to study the following questions: when do internal abelian groups admit a tensor product in the internal sense (see \thref{def:tensprod})? When it exists, when does it provide a monoidal structure on the category of internal abelian groups? This section has not been finished though, by lack of time.

\subsection{Context and basic notions}
In this part, we will stay in the realm (that is, in the 2-XL-category) of Cartesian monoidal categories, symmetric strong monoidal functors and monoidal natural transformations. It is essentially the 2-XL-category of categories with finite products, finite product preserving functors and natural transformations, but with the advantage of having tools of symmetric monoidal categories for doing algebra. Let's make this more precise.

Consider a category \C with finite products. Then, making a choice of a product functor $\_\times\_\colon\C\times \C\to\C$ and of a terminal object $*$, one automatically gets a symmetric monoidal structure on \C with unit, associativity and symmetry isomorphisms given by universality of products. This is called a \emph{Cartesian monoidal category} (which has to be distinguished from a \emph{Cartesian category} which is a category with finite limits).

In the same manner, if one starts from a finite product preserving functor $F\colon\C\to~\D$ between categories with finite products, one automatically gets a symmetric strong monoidal functor between the corresponding Cartesian monoidal categories (whatever choice of such structures was made on the categories), with structure isomorphisms given by the universality of products. There are actually no other type of strong monoidal functors between Cartesian monoidal categories: they preserve products and their structure isomorphisms have to be the universal ones. Moreover, they automatically are symmetric in a unique way.

Finally, if \C and \D are categories with finite products, then any natural transformation $\tau\colon F\Rightarrow G$ between product preserving functors is monoidal once a choice is made of a Cartesian monoidal structure on \C and \D.

Now the global picture of the land. There is a 2-XL-category $\SYMMON$ of symmetric monoidal categories, symmetric monoidal functors and monoidal natural transformations. We denote $\CARTMON$ its 2-cell-full sub-2-XL-category whose objects are Cartesian monoidal categories, morphisms (symmetric) strong monoidal functors and 2-cells (monoidal) natural transformations. On the other hand, one has the 2-cell-full sub-2-XL-category \textit{FPCAT} of \CAT of categories with finite products, finite product preserving functors. The following proposition states in short the discussion above.
\begin{Prop}
	The forgetful 2-functor $U\colon \SYMMON\to \CAT$ restricts to a 2-equivalence
	$$U\colon \CARTMON\ra{\simeq} \mathit{FPCAT}.$$
\cqfd
\end{Prop}

\begin{Exs}[We recall here two examples of categories with finite products that will be extensively used in this work.]
	\item Let \C be a category with pullbacks. Then, for every objects $C\in \C$, the slice category has $\C/C$ has products given by pullbacks in \C over $C$ and terminal object given by the identity morphism $1_C$ of $C$.
	\item Given a category $\A$ with (finite) limits, it is well-known that for any small category $\C$, the functor category $\A^\C$ has objectwise (finite) limits. In case of $\A=\Set$, the category of presheaves $Set^{\C^{op}}$ on a small category $\C$ has therefore all limits. Now, consider a pretopology $K$ on \C. The full sub-category $\Sh(\C,K)\subset \Set^{\C^{op}}$ of sheaves on the site $(\C,K)$ is a \emph{localization} of the category $Set^{\C^{op}}$ of presheaves on \C, meaning that it is \emph{reflective in $Set^{\C^{op}}$} (the inclusion functor admits a left adjoint) and that the reflector (the left adjoint of the inclusion) preserves finite limits\cite{BorIII94,MM92}:
	\begin{equation*}
	\xymatrix@C=4pc@1{
a:\Set^{\C^{op}\ }\ar@<-.9ex>[r]\ar@{}[r]|{\top}&\ \Sh(\C,K) : I.\ar@<-.9ex>@{_(->}[l]
}\medskip\end{equation*}
From the fact that  $\Sh(\C,K)$ is a full replete reflective subcategory of the category of presheaves, one deduces that $\Sh(\C,K)$ has all limits, given by the limits in $\Set^{\C^{op}}$, thus by objectwise limits \cite{BorI94}.
\end{Exs}

\begin{Rem}\label{FunctorCat}
	Let $\D$ be a small category. The 2-functor $(-)^{\D}\colon \CAT\to \CAT$ ``lifts'' to a 2-functor $(-)^{\D}:\SYMMON\to \SYMMON$ using objectwise monoidal products. As products in a functor category can be chosen to be the objectwise products, this 2-functor restricts to $\CARTMON$ as domain and codomain.
\end{Rem}

A Cartesian monoidal category has the special feature of being symmetric and of having a diagonal map $C\ra{\triangle} C\times C$ as well as an augmentation map $C\ra{\epsilon}*$ for any object $C$. This is the necessary structure for stating the axioms of internal abelian groups, rings and modules
\footnote{There is a more general concept of a \emph{Hopf monoid} in a symmetric monoidal category \V : these are both monoids and comonoids in \V with the two structure being coherent, i.e. \emph{bimonoids}, and with the additional data of an \emph{antipode} map (playing the role of an inverse).

In a Cartesian monoidal category $(\C,\times,*)$, every object is, with its diagonal and augmentation maps, naturally a comonoid in \C, and it is the unique comonoid structure possible on it. Moreover, this natural comonoid structure on an object of \C is automatically coherent with any monoid structure given on it and an inverse of a group object verifies the axiom of an antipode. This gives rise to an isomorphism of categories $Gr(\C)\cong Hopf(\C,\times,*)$. Thus groups objects are particular case of Hopf monoids. 

Conversely, cocommutative Hopf monoids are particular cases of group objects. Given a symmetric monoidal category \V, the category of cocommutative comonoids $CocComon(\V)$ inherits the symmetric monoidal structure of \V and this monoidal structure on $CocComon(\V)$ happens to be Cartesian. Then one has the the equality $CocHopf(\V)=Gr(CocComon(\V))$\cite{Por08,Por08a}.}.

\begin{Def}\label{IntAb}
Let $(\cal C,\times,*)$ be a Cartesian monoidal category. An \emph{internal group} or \emph{group object in \C} is a quadruple $(G,\mu,\eta,\zeta)$ where 
\begin{itemize}
\item $(G,\mu,\eta)$ is a monoid in $(\cal C,\times,*)$ with binary operation $\mu$ and unit $\eta$,
\item $\zeta\colon G\to G$ is the inverse morphism, i.e., it makes the following diagrams commute.
\end{itemize}
$$\xymatrix@=1.5cm{
	G \ar[r]^-{(1_{G}\times\zeta)\circ\triangle} \ar[d]_{\epsilon}	& G\times G \ar[d]^{\mu}	& G \ar[l]_-{(\zeta\times 1_{G})\circ\triangle} \ar[d]^{\epsilon}\\
	\ast \ar[r]_{\eta} 				& G					& \ast \ar[l]^{\eta}
}$$
An \emph{internal abelian group} is an internal group that is commutative as a monoid (we will then denote additively its binary operation).\\
A \emph{morphism of (abelian) internal groups} is a morphism of the underlying monoids. One also calls \emph{additive} the morphisms of internal abelian groups.\\
With the composition and identities of \C, internal groups (resp. internal abelian groups) and their morphisms form a category $Gr(\C)$ (resp. $Ab(\C)$).
\end{Def}

\begin{Rem}\label{InvPres}
	If a morphism $G\to H$ in \C between two internal groups preserves their multiplications, it will automatically preserve their units and their inverses, as in $Set$. In particular, there is no need of verifying the unit axiom for morphisms of internal groups.
\end{Rem}

Finite product preserving functors preserve categories of internal (abelian) groups. In a similar way as for (commutative) monoids in a (symmetric) monoidal category (see Proposition \ref{prop:Mon}), this gives rise to a 2-functor.

\begin{Prop}
	Let \C be a Cartesian monoidal category.\Par
	There are 2-functors $Ab\subset Gr:\CARTMON\to \CAT.$
\qed
\end{Prop}

\begin{Def}
	Let $A,B,C\in\C$ be internal abelian groups in \C.
	A morphism 
	$$f\colon A\times B\to C$$ 
	in \C is \emph{biadditive} if the two following diagrams commute.
	$$\xymatrix{
		A\times (B\times B) \ar[r]^-{1_A\times +} \ar[d]_{\triangle\times 1_{B\times B}}	& A\times B \ar[ddd]^{f}	& (A\times A)\times B \ar[l]_-{+\times 1_{B}} \ar[d]^{1_{A\times A}\times\triangle}\\
		(A\times A)\times (B\times B) \ar[d]_{\cong}							&					& (A\times A)\times (B\times B) \ar[d]^{\cong}\\
		(A\times B)\times (A\times B) \ar[d]_{f\times f}							&					& (A\times B)\times (A\times B) \ar[d]^{f\times f}\\
		C\times C \ar[r]_{+}												& C					& (C\times C) \ar[l]^{+}
	}$$
\end{Def}

\begin{Def}
	An \emph{internal ring} in a Cartesian monoidal category \C is a sextuplet $(R,+,0,-,\cdot,1)$ where:
	\begin{itemize}
		\item $(R,+,0,-)$ is an internal abelian group in \C,
		\item $(R,\cdot,1)$ is a monoid in \C,
		\item $R\times R\ra{\cdot} R$ is a biadditive map.
	\end{itemize}
	An internal ring is \emph{commutative} if it is so as a monoid. A \emph{morphism of internal rings} is a morphism $f\colon R\to S$ in \C that is both a morphism of the underlying internal abelian groups and of the underlying monoids.\\
	With the composition and identities of \C, they all together organize into the category $Ring(\C)$\index[not]{Ring(C)@$Ring(\C)$}. It has a full subcategory \cRing(\C)\index[not]{cRing(C)@$\cRing(\C)$} of commutative rings.
\end{Def}

\begin{Def}
	Let $R$ be an internal ring in $\C$. An \emph{internal (right) $R$-module} is a quintuple $(M,+,0,-,\kappa)$ where
	\begin{itemize}
		\item $(M,+,0,-)$ is an internal abelian group in \C,
		\item $(M,\kappa)$ is an $R$-module over the monoid $R$ in \C,
		\item $\kappa\colon M\times R\to M$ is biadditive.
	\end{itemize}
	A \emph{morphism of internal $R$-modules} is a morphism $f\colon M\to N$ in \C that is both a morphism of the underlying internal abelian groups and of the underlying $R$-modules.\\
	Given an internal ring, one thus has the category $M_{R}(\C)$ of internal (right) $R-$modules.
\end{Def}

\begin{Rem}
	The definition we have given for a biadditive morphism was meant to emphasize the fact that it is defined in the more general context of Hopf monoids in a symmetric monoidal category. But we will soon exploit the fact that our category is in fact Cartesian monoidal by using the universality of products and terminal object. In this context, a more handy (and equivalent) version is available. \Par
	Let $(A\times B\times B,\pi_{1},\pi_{2},\pi_{3})$ be any product of $A$, $B$ and $B$. Then $f\colon A\times B\to C$ is biadditive if and only if the following diagram commutes.
		$$\xymatrix@R=1cm@C=2cm{
		A\times B\times B \ar[r]^-{(\pi_{1},+\circ(\pi_{2},\pi_{3}))} \ar[d]_{(f\circ(\pi_{1},\pi_{2}),f\circ(\pi_{1},\pi_{3}))}	& A\times B \ar[d]^{f} & A\times A\times B\ar[l]_-{(+\circ(\pi_{1},\pi_{2}),\pi_3)}
		\ar[d]^{(f\circ(\pi_{1},\pi_{3}),f\circ(\pi_{2},\pi_{3}))}\\
		C\times C \ar[r]_{+}												& C & C\times C\ar[l]^{+}
	}$$
\end{Rem}

\subsection{Properties of the category of internal abelian groups}

We first notice some important properties of internal (abelian) groups that are shared with (abelian) groups (in $Set$). A useful result is the fact that a group $G$ in $\C$ induces a ``pointwise'' group structure on the hom-sets $\C(C,G)$ \cite{Freyd64}.

\begin{Lem}
	Let \C be a Cartesian monoidal category and $G$ a internal group in $\C$.\Par
	Then the representable functor $\C(-,G)\colon\C^{op}\to Set$ factors through the forgetful functor of the category of groups $Gr=Gr(Set)$ to the category $Set$:
	$$\xymatrix{
	\C^{op}\ar[rr]^{\C(-,G)}\ar@{-->}[dr] && Set\\
	&Gr\ar[ur]_U&
	}$$
	Moreover, if $G$ is abelian, then it factorizes through the category $Ab=Ab(Set)$ of abelian groups.
\end{Lem}

If a functor $F\colon \D\to Gr(\C)$ has a limit in \C given by a cone $u\colon\lim UF\Rightarrow UF$, then there is a unique structure of internal group on $\lim UF$ such that each $${u_D\colon\lim UF\to F(D)}$$ is a morphism of internal groups. Moreover this is a limit in $Gr(\C)$. The same is true for abelian group objects. This is folklore, which we prove now.
\begin{Prop}\label{LimitCreation}
	Let \C be a Cartesian monoidal category.\Par
	Then the forgetful functors $U\colon Gr(\C)\to \C$ and $U\colon Ab(\C)\to \C$ create limits.\Par
	In particular, $Gr(\C)$ and $Ab(\C)$ have all finite products and their forgetful functors preserve them. Moreover, if $\C$ is complete, then $Gr(\C)$ and $Ab(\C)$ are complete and the forgetful functors preserve limits.
\end{Prop}

\begin{Pf}
	Group categories behave well with functor categories: there is a 2-natural isomorphism $\alpha$
	
	\begin{equation}\label{eq:grfuncat}\begin{aligned}\xymatrix@=1.5cm{
		\CARTMON \ar[r]^{Gr}\ar[d]_{(-)^{\D}}\drtwocell<\omit>{<0>\alpha} & \CAT\ar[d]^{(-)^{\D}}\\
		\CARTMON \ar[r]^{Gr} & \CAT
	}\end{aligned}\end{equation}
	where the left vertical functor is the functor described in \thref{FunctorCat}. The isomorphism $\alpha$ makes the following diagram commute
	$$\xymatrix@=1.5cm{
		Gr(\C)^{\D}\ar[r]^{U_{*}}\ar[dr]_{\alpha_{\C}}^{\cong}	& \C^{\D}\\
											& Gr(\C^{\D})\ar[u]_{\cal U}
	}$$
where $\cal U$ is the forgetful functor and $U_{*}$ is post-composition with $U$. Let $\bar F\colon \D\to Gr(\C)$ be a functor. Denote $F:=U_{*}(\bar F)=U\circ \bar F$, $\bar F(D):=(F(D),m_{D},e_{D},z_{D})$. Then $\alpha_{\C}(\bar F)=(F,m,e,z)$ where the $D$-component of the natural transformations $m$, $e$ and $z$ are $m_{D}$, $e_{D}$ and $z_{D}$ respectively.\\ 
Suppose now that a limit $u\colon\lim F\Rightarrow F$ exists in \C. 

\emph{Uniqueness}: Let $\bar{A}=(A,\mu,\eta,\zeta)\in Gr(\C)$ and $\tau\colon \bar{A}\Rightarrow \bar F$ be a cone in $Gr(\C)$ above $u$, i.e. $U\cdot\tau=u$. We want to show that this cone is uniquely determined. There actually is no choice for $\tau$ since then $A=\lim F$ and $\tau_{D}=u_{D}$, for all $D\in\D$. It remains to show that the group structure of $\bar A$ is uniquely determined.\\ 
Since each $\tau_{D}$ is a morphism of internal groups, the following diagram commutes
\begin{equation}\begin{aligned}\xymatrix{
	\lim F\times \lim F\ar[r]^-{\mu}\ar[d]_{u_{D}\times u_{D}}	& \lim F\ar[d]^{u_{D}}\\
	F(D)\times F(D)\ar[r]_-{\mu_{D}}	& F(D)
}\end{aligned}\end{equation}
One gets the following commutative diagram in $\C^{\D}$.
\begin{equation}\begin{aligned}\xymatrix{
	\lim F\times \lim F\ar[r]^-{\mu}\ar@{=>}[d]_{u\times u}	& \lim F\ar@{=>}[d]^{u}\\
	F\times F\ar@{=>}[r]_-{m}	& F
}\end{aligned}\end{equation}
Thus $\mu$ is uniquely determined. The uniqueness of the unit and inverse morphisms is similar (see \thref{InvPres}): the two following diagram have to commute in $\C^{\D}$.

\begin{equation}\begin{aligned}\xymatrix{
	{*}\ar[r]^-{\eta}\ar@{=>}[dr]_{e}	& \lim F\ar@{=>}[d]^{u} 	& & \lim F\ar@{=>}[d]_{u}\ar[r]^{\zeta}	& \lim F\ar@{=>}[d]^{u}\\
							& F					& & F\ar@{=>}[r]_{z}	& F
}\end{aligned}\end{equation} 

\emph{Existence}: Consider the quadruple $\bar L=(\lim F,\mu,\eta,\zeta)$, where $\mu$, $\eta$ and $\zeta$ are defined by the three preceding diagrams. This is an internal group in $\C$ as one can check directly. Moreover, the cone $u\colon\lim F\Rightarrow F$ determines a cone $\bar u\colon \bar L\Rightarrow \bar F$ in $Gr(\C)$ with $U\cdot \bar u=u$ by definition. This is a limit in $Gr(\C)$ of $\bar F$. One can verify these facts objectwise, but one is really moving back and forth using the isomorphism $Gr(\C)^{\D}\cong Gr(\C^{\D})$.

For internal abelian groups, one first remarks that the 2-natural isomorphism $\alpha$ in (\ref{eq:grfuncat}) restricts to the 2-functor $Ab$. One can then show that if $\bar F$ takes values in internal abelian groups, the limit in internal groups $\bar L$ we have constructed is abelian.
\end{Pf}

\begin{Ex}
By (\ref{eq:grfuncat}), there is a natural isomorphism $Ab(\Set^{\C^{op}})\cong Ab^{\C^{op}}$.
\end{Ex}

Given a Cartesian monoidal category \C, the previous proposition determines a canonical Cartesian monoidal structure on $Gr(\C)$ and $Ab(\C)$: the products and terminal object created by the forgetful functor from the Cartesian monoidal structure of \C. In particular, the group structure on the terminal $*$ is the unique possible such structure and is automatically abelian. The group structure on the product $G\times H$ of two groups in $\C$ is the only group structure such that the projections $G\leftarrow G\times H\to H$ are group object morphisms, that is, componentwise multiplication, unit and inverse. It is abelian if both $G$ and $H$ are.

We turn now more specifically to the category of internal abelian groups $Ab(\C)$. It is a classical fact that it is an additive category \cite{Freyd64,BG04}, but we haven't seen a proof of it.

\begin{Prop}
	Let \C be a Cartesian monoidal category. Then $Ab(\C)$ is an additive category. 
	Given internal abelian groups $A$ and $(B,+,0,-)$, the abelian group structure on the hom-set $Ab(\C)(A,B)$ is given by:
	$$\begin{array}{rllll}
	Ab(\C)(A,B)\times Ab(\C)(A,B) & \xlongrightarrow{\cong}  & Ab(\C)(A,B\times B) & \xlongrightarrow{+_*} & Ab(\C)(A,B)\\
	*  & \xlongrightarrow{\cong} & Ab(\C)(A,*) & \xlongrightarrow{0_*} & Ab(\C)(A,B)\\
	Ab(\C)(A,B) & \xlongrightarrow{-_*} &  Ab(\C)(A,B) &&
	\end{array}$$
	Therefore, the sum of two arrows $f,g\colon A\to B$ is the composite $A\ra{(f,g)}B\times B\ra{+}B$, the unit is the composite $A\to *\ra{0}B$ and the opposite of an arrow $f\colon A\to B$ is the composite $A\ra{f}B\ra{-}B$.
\end{Prop}
We have defined the abelian group structure on the hom-sets in a way that suggests a proof. Indeed, consider a representable functor $Ab(\C)(A,-)\colon Ab(\C)\to Set$. This functor preserves products and therefore group objects. Given an abelian group object $B$ in $Ab(\C)$, the abelian group object structure it induces on $Ab(\C)(A,B)$ is precisely the one given in the Proposition. But abelian group objects in $Ab(\C)$ are precisely abelian group objects in \C. Indeed, as for groups, if an object admits two compatible internal group structures, these structure coincide and are commutative.

\begin{Lem}
	Let \C be a Cartesian monoidal category, and consider $Gr(\C)$ with its induced Cartesian monoidal structure. Then the forgetful functor $U\colon Gr(\C)\to\C$ induces an isomorphism of categories
	$$
	Gr(Gr(\C))\xlongrightarrow[\cong]{Gr(U)}Ab(\C).
	$$
	\cqfd
\end{Lem}

Now things get more complicated when one study the abelianness of $Ab(\C)$. Indeed, if $Top$ denotes the category of topological spaces, then $Ab(Top)$ is not abelian \cite{BG04}. Recall that one way of characterizing an abelian category is that it is both an additive and a Barr-exact category. If \C is Barr-exact, then $Ab(\C)$ inherits its Barr-exactness and thus is abelian \cite{BG04}.

\begin{Lem}
	If \C is a Barr-exact Cartesian monoidal category, then $Ab(\C)$ is abelian.
\end{Lem}

\begin{Exs}
	\item
	The category $Ab(HComp)$ of internal abelian groups in the category of compact Hausdorff spaces is abelian because the latter is Barr-exact \cite{BorII94}.
	\item
	If $\C$ is a topos, for instance a category of sheaves on a space or on a site, $Ab(\C)$ is abelian, because any topos is Barr-exact \cite{BorIII94}.
\end{Exs}

\subsection{Tensor product of internal abelian groups}
\begin{Def}\label{def:tensprod}\indexb{Tensor product of internal abelian groups}
	Let $A$ and $B$ be internal abelian groups in $\C$. A \emph{tensor product} of $A$ and $B$ is a biadditive morphism $A\times B\ra{h} A\otimes B$ in \C that is universal among such, i.e. every biadditive morphism $A\times B\ra f C$ factors uniquely through it:
	$$\xymatrix@R=1.7cm@C=1.7cm{
	A\times B\ar[r]^h\ar[dr]_f	& A\otimes B\ar@{-->}[d]^{\exists!\hat f}\\
					& C
	}$$
\end{Def}

There are two natural questions that arise. Do tensor products exist in general? If so, do they determine a monoidal structure on $Ab(\C)$? We don't know much about these questions in general, and we are not aware of much research on them. Here are nevertheless some indications.

Let us consider at first the existence question. The results below follow quite directly from the definition.

\begin{Lem}
	If $f\colon A\times B\to C$ is biadditive, then 
\begin{enumerate}[(i)]
	\item if $g\colon C\to D$ is additive, then $g\circ f$ is biadditive,
	\item if $g\colon A'\to A$ and $h\colon B'\to B$ are additive, then $f\circ (g\times h)$ is biadditive.
\end{enumerate}

\qed
\end{Lem}

\begin{Prop}
	Let \C be a locally small Cartesian monoidal category.\\
	Then for every $A$, $B\in Ab(\C)$, biadditive morphisms from $A\times B$ determine a functor
	$$Biadd(A,B;-)\colon Ab(\C)\to Set.$$
	and a tensor product of $A$ and $B$ exists if and only if this functor is representable.\\
	Moreover, if \C is small complete, then the functor $Biadd(A,B;-)$ preserves all small limits.
\qed
\end{Prop}

This proposition can help in answering the question of the existence of a tensor product in particular cases. For, if \C is a locally small and small complete Cartesian monoidal category, then $Ab(\C)$ is also locally small and small complete (the latter by \thref{LimitCreation}). Thus, the representability of $Biadd(A,B;-)$ can be proved by means of a \emph{solution set} for the \emph{Representation Theorem} \cite[p. 122]{McL97}. If, in addition, $Ab(\C)$ is \emph{well-powered} and has a \emph{cogenerating set}, then the Special Representation Theorem gives the representability directly [ibid. p. 130]. Let us recall the definition of these concepts (e.g., \cite{McL97,BorI94}).

\begin{Defs}
	\item
A \emph{subobject} of an object $C$ in a category \C is an isomorphism class  in $\C/C$ of monos $A\rightarrowtail C$ in \C. A category is \emph{well-powered} if the class of subobjects of each of its objects is a set.
\item
A set $Q$ of objects of \C is a \emph{cogenerating set} if, given a pair of distinct parallel arrows $f\neq g\colon A\rightrightarrows B$, there is an object $C\in Q$ and an arrow $h\colon B\to C$ such that $h\circ f\neq h\circ g$.
\end{Defs}

\begin{Lem}
	If \C is well-powered and admits a free abelian group object functor, then $Ab(\C)$ is also well-powered.
\end{Lem}
\begin{Pf}
	Let $f\colon A\to B$ be a morphism of internal abelian groups in \C. Observe first two easy facts.
	\begin{enumerate}
\item
$f$ is an isomorphism in \C \ssi it is an isomorphism in $Ab(\C)$.
\item
If $f$ is mono is \C, then it is mono in $Ab(\C)$.

\end{enumerate}

Let us prove now the converse of fact 2 under the hypothesis of existence of a free abelian group functor $F\colon \C\to Ab(\C)$, left adjoint of the forgetful functor $U$. Suppose $f$ is mono in $Ab(\C)$ and consider a diagram
$$\xymatrix{
C\ar@<.8ex>[r]^h\ar@<-.8ex>[r]_k & A\ar[r]^f & B
}$$
in \C with $f\circ h=f\circ k$. Then the transpose $\hat k$ and $\hat h$
$$\xymatrix{
F(C)\ar@<.8ex>[r]^{\hat h}\ar@<-.8ex>[r]_{\hat k} & A\ar[r]^f & B
}$$
 of $h$ and $k$ verify $f\circ \hat h=f\circ \hat k$ (use universality of the unit $C\to UF(C)$). Since $f$ is mono in $Ab(\C)$, $\hat h=\hat k$, and therefore $h=k$.
 
 Now we prove the result. Let us write $Mono_{\C}(A)$ for the class of monos in \C with codomain $A$ and $Sub_{\C}(A)$ the class of subobjects of $A$. The dotted arrow in the following diagram exists and is injective because of the preceding results.
 $$\xymatrix{
 Mono_{Ab(\C)}(A)\ar@{^{(}->}[r]\ar[d] & Mono_{\C}(A)\ar[d]\\
 Sub_{Ab(\C)}(A)\ar@{^{(}-->}[r] & Sub_{\C}(A)
 }$$
 \cqfd
\end{Pf}

Here is a result that can be combined with the preceding \cite{BorIII94}.
\begin{Lem}
	A topos is well-powered.
\end{Lem}

\begin{Exs}
	\item The category $Ab=Ab(Set)$ has a cogenerator, i.e., a cogenerating set of one element, given by $\Q/\Z$ (see \cite[p. 169]{BorI94}).
	\item The category $Ab(Top)$ of topological abelian groups has a cogenerator given by $\Q/\Z$ with the coarsest topology. The category $Ab(k\text{-}Top)$ of compactly generated topological abelian groups has a cogenerator $k(\Q/\Z)$ where $\Q/\Z$ has the coarsest topology. 
	\item The category $Sh(X;Ab)\cong Ab(Sh(X))$ of abelian sheaves over a space $X$ has a cogenerator (see, e.g., \cite{Mit65}). In fact, this remains true for abelian sheaves on a site: if $\E$ is a Grothendieck topos, then the category $Ab(\E)$ has a cogenerator (see \cite[p. 261]{Joh77}).
\end{Exs}

\section{The fibred context}
After having defined the monoidal objects in the world of (op-,bi-)fibrations, we explore the notions of monoids and modules in this context.
\subsection{Monoidal fibred categories}
In this part, we will introduce the fibred version of a monoidal category. We actually define two notions, monoidal fibred categories and strong monoidal fibred categories. In the latter case, the monoidal structure preserves Cartesian morphisms, while in the former, it is required only to preserve the fibre. We explain the theory in the languages of both fibrations and indexed categories. 

The notion of monoidal fibred category seems to appear, expressed in the language of indexed categories, in the 90's in the computer science literature (see, e.g., \cite{CA93}, but there might be older references I am not aware of). Maltsiniotis describes in some detail these kinds of structure in the fibration language (in French, \cite{Mal95}), but the author considered \emph{strict} monoidal objects. Our (strict) monoidal fibred categories are a particular case of his \emph{monoidal categories over \B} where the functor $\E\to\B$ is a fibration and what he called \emph{fibred monoidal categories} are our (strict) strong monoidal fibred categories (see \autoref{fnote:Wordorder} on \autopageref*{fnote:Wordorder} for a remark on the word order). More recent references, both written in terms of indexed categories, are \cite{HM06} and \cite{PS09}\footnote{The first authors make the same distinction that we do between monoidal and strong monoidal objects, whereas the second authors consider only strong (symmetric) monoidal objects without the adjective ``strong''. Shulman in another article \cite{Shu08a}, talks about \emph{monoidal fibrations} which are in general completely different from the structure we consider because both total and base categories are supposed to be monoidal. See also \cite{Gom09} for these kinds of structures.}. Our new contribution in this part is a systematic treatment of the subject, in both frameworks of fibrations and of indexed categories. In particular, we prove the 2-equivalence between the two notions. We finish with a new insight about monoidal bifibred categories.

\subsubsection{The fibration setting}

Recall that $\mathit{OPFIB}(\B)$, $\BIFIB(\B)$ and $\FIB(\B)$ are the respective full sub-2-XL-ca\-te\-go\-ries of the strict slice 2-XL-category $\CAT/\B$ consisting of (op-,bi-)fibrations. We also consider stronger versions of them, the 2-cell-full sub 2-XL-cat\-e\-gories of $\CAT/\B$ consisting of respectively (op-,bi-)fibrations and (op-,bi-)Cartesian functors over \B, written respectively $\mathit{OPFIB}_{oc}(\B)$, $\BIFIB_{bc}(\B)$ and $\FIBc(\B)$. In case of bifibrations, one also encounters functors that are Cartesian, but not opCartesian, or the contrary, and thus we should also consider $\BIFIB_c(\B)$ and $\BIFIB_{oc}(\B)$.

These all admit Cartesian monoidal 2-structures whose product $\times_{\B}$ is a fibre product in \CAT and whose unit is $Id_{\B}$, the identity functor of the category $\B$ (see \thref{prop:2prodfib} and \thref{lem:pullopfib}). For convenience, we will choose the usual fibre product of \CAT. Recall that for fibred products of fibrations for instance, a choice of Cartesian morphisms is then given by pairwise Cartesian morphisms.
\paragraph{}
We now formulate our theory for fibrations, as the corresponding statements can be easily recovered for op- and bifibrations. Recall that a monoidal category is a \emph{monoidal object} (or a \emph{pseudo-monoid}) in the Cartesian monoidal 2-XL-category \CAT \cite{DS03}.  

\begin{Def}\index{Monoidal fibred category|see{Monoidal fibration}}\indexb{Monoidal fibration}\indexb{Monoidal fibration!Strong --}
The monoidal objects in $\FIB(\B)$ (resp. $\FIBc(\B)$) are called \emph{monoidal} (resp. \emph{strong monoidal}) \emph{fibred categories} (or \emph{fibrations}) \emph{over \B}\footnote{I prefer this word order (than the more usual ``fibred monoidal categories'') because these objects are not monoidal categories. Moreover, one should note that the adjective \emph{strong} has nothing to do with the ``laxness'' of the monoidal objects: both are pseudo-monoids (and not lax or strict monoids). Neither has it to do with the ``splitness'' of the fibration: both are general. It indicates in what monoidal 2-category the monoidal object lives and this change of environment has an effect on the strength of the relation between the monoidal product and the Cartesian morphisms. In addition, we will see later that this terminology is indicative of the associated indexed categories.\label{fnote:Wordorder}}.

Thus, a monoidal (resp. strong monoidal) fibred category over \B is a sextuple $(P, \otimes,u,\alpha,\lambda,\rho)$ where: 
\begin{enumerate}[(i)]
\item
$P:\E\to\B$ is a fibration,
\item
$\otimes \colon \E \times_{\B}\E \to \E$ is a functor over \B (resp. a Cartesian functor),
\item
$u\colon\B\to\E$ is a functor over \B, i.e., a section of $P$ (resp. a Cartesian functor),
\item
\(\alpha\) is a natural isomorphism over \B, called \emph{associator}, and filling the following diagram
$$\shorthandoff{;:!?}\xymatrixnocompile@!0@C=4pc@R=6pc{
&(\E \times_{\B}\E)\times_{\B}\E\ar[rr]^{\cong}\ar[ld]_{\otimes\times_{\B}Id_{\E}} &&
\E \times_{\B}(\E\times_{\B}\E)\ar[rd]^{Id_{\E}\times_{\B}\otimes} & \\
	\E\times_{\B}\E\ar[rrd]_{\otimes} &\rrtwocell<\omit>{!\mbox{$\stackrel{\alpha}{\cong}$}}&&&	\E\times_{\B}\E\ar[lld]^{\otimes}\\
&&\E&&
	}$$
\item
\(\lambda\) and \(\rho\) are natural isomorphisms over \B, called respectively \emph{left\textnormal{ and }right unitors}, and filling the following diagram.
$$\xymatrix@=1.5cm{
\B \times_{\B}\E\ar[r]^{u\times_{\B}Id_{\E}}\ar[dr]_{\cong}\drtwocell<\omit>{<-3>\lambda}&\E\times_{\B}\E\ar[d]^{\otimes} & \E\times_{\B}\B\ar[l]_{Id_{\E}\times_{\B}u}\ar[dl]^{\cong}\\
&\E\urtwocell<\omit>{<-3>\rho}&
}$$
\end{enumerate}
These data are subject to the usual \emph{coherence axioms} in each fibre of $P$. 
\end{Def}

Because the whole monoidal fibred structure is \emph{over} \B, it restricts to each fibre $\E_B$, $B\in\B$, and therefore determines a monoidal category structure on each of them. We call these monoidal categories the \emph{fibre monoidal categories}\indexb{Monoidal fibration!Fibre monoidal category of a --} of the monoidal fibred category.

Conversely, it is not sufficient to define a monoidal structure on each fibre, because these structures should be related to each other by the horizontal morphisms. Even in the non-strong setting one has to check that the tensor and the unit are functors and that the associativity and unit are natural, not only on vertical morphisms, but on any type of morphisms. Nevertheless, it is very useful to understand the monoidal structure on the fibres. For example, they might be Cartesian monoidal, allowing us to use the tools of universality. Moreover, all the coherence diagrams take place in the fibres.

\paragraph{}
Since both $\FIB(\B)$ and $\FIB_c(\B)$ are symmetric 2-monoidal (whose symmetry we denote \emph{sym}), one can define \emph{symmetric monoidal objects} in them.

\begin{Def}\indexb{Monoidal fibration!Symmetric (strong) --}
	A \emph{symmetric monoidal \op resp. symmetric strong monoidal\fp fibred category} over a category \B is given by a monoidal (resp. strong monoidal) fibred category $(P\colon\E\to\B,\otimes)$ and a natural isomorphism $\sigma$ over \B as in the following diagram.
	$$\xymatrix@C=.5cm@R=1.5cm{
	\E\times_\B\E\ar[rr]^{\text{sym}}_\cong\ar[rd]_\otimes\drrtwocell<\omit>{!\stackrel{\displaystyle\sigma}{\cong}} & & \E\times_\B\E\ar[dl]^\otimes\\
	&\E&
	}$$
	subject to the usual axioms in each fibre of $P$.
\end{Def}

\begin{Exs}\label{ex:monfib}
	\item
Any monoidal category \V is strong monoidal fibred over the terminal object \one of \CAT.
	\item\label{ex:monfib:can}
If \C has pullbacks, then the canonical fibration $\C^{\two}\ra{\cod}\C$ over \C determines, via a choice of a pullback for each pair $A\to B\leftarrow C$ in $\C$, a \emph{Cartesian monoidal fibred category}. This means that it is a symmetric strong monoidal fibred category such that the induced monoidal structures on its fibres are Cartesian.
\end{Exs}

\paragraph{}
We turn now to morphisms of monoidal objects in $\FIB(\B)$ and $\FIBc(\B)$.
\begin{Def}\indexb{Monoidal functor over \B}\indexb{Monoidal Cartesian functor over \B}
	A \emph{monoidal functor \textnormal{(\resp}monoidal Cartesian functor\fp over \B} from a monoidal fibration $(P\colon \D\to\B,\otimes,u)$ to a monoidal fibration $(Q\colon\E\to\B,\otimes',u')$ is a triple $(F,\phi,\psi)$ where
\begin{enumerate}[(i)]
\item
$F\colon\D\to\E$ is a functor (resp. a Cartesian functor) over \B,
\item
\(\phi\) is a natural transformation over \B filling the following diagram,
$$\xymatrix@R=1.5cm{
\D\times_{\B}\D\ar[r]^{F\times_{\B}F}\ar[d]_{\otimes}\drtwocell<\omit>{\phi} & \E\times_{\B}\E\ar[d]^{\otimes'}\\
\D\ar[r]_{F} & \E
}$$
\item
\(\psi\) is a natural transformation over \B filling the following diagram
$$\xymatrix@=1.5cm{
\B\ar[d]_u\ar[dr]^{u'}\drtwocell<\omit>{<3>\psi} & \\
\D\ar[r]_F & \E
}$$
\end{enumerate}
and this data is subject to the usual \emph{coherence axioms} in each fibre of $Q$. 

Of course, one has, as in the non-fibred context, the special cases of \emph{strong} and \emph{strict monoidal} morphisms.\indexb{Monoidal functor over \B!Strong --}\indexb{Monoidal Cartesian functor over \B!Strong --}\indexb{Monoidal functor over \B!Strict --}\indexb{Monoidal Cartesian functor over \B!Strict --}
\end{Def}
Such a functor $F$ restricts to a monoidal functor $(F_B,\phi_B,\psi_B)$ on the fibres over $B$ for each $B\in\B$. Again, these restrictions are very informative (in particular, all the coherence diagrams stay entirely in the fibres) but not sufficient for defining a monoidal functor over \B. In particular, the possible Cartesianness of $F$ is not captured in the fibre-restrictions.

\begin{Def}\indexb{Monoidal functor over \B!Symmetric (strong, strict) --}
	A \emph{symmetric monoidal \op resp. monoidal Cartesian\fp functor over \B} between symmetric monoidal fibred categories $(\E,\sigma)$ and $(\E',\sigma')$ is a monoidal (resp. monoidal Cartesian) functor $(F,\phi,\psi)$ over \B between the underlying monoidal fibred categories $\E$ and $\E'$ that preserves the symmetry isomorphisms $\sigma$ and $\sigma'$. That is, it is subject to the usual coherence axiom in the fibres of $\E'$. It terms of natural transformations, the axiom is given by
	$$
	(F\cdot\sigma)\bullet\phi=(\phi\cdot\text{sym})\bullet(\sigma\cdot (F\times_\B F)).
	$$
	Symmetric monoidal functors are \emph{strong \op\resp strict\fp} when they are so as monoidal functors.
\end{Def}

\paragraph{}
Finally, we define the 2-cells in $\FIB(\B)$ and $\FIBc(\B)$ that preserve monoidal structures.
\begin{Def}\indexb{Monoidal natural transformation over \B}
	A \emph{monoidal natural transformation \(\alpha\) over \B} from $(F,\phi^F,\psi^F)$ to $(G,\phi^G,\psi^G)$, $F$ and $G$ being monoidal functors over \B (Cartesian or not) from \D to \E, (strong) monoidal fibred categories, is a natural transformation $\alpha\colon F\Rightarrow G$ satisfying the usual coherence axioms, which we state:
$$
\xymatrixnocompile@C=1.5cm@R=1cm{
\D \times_{\B} \D\ar@/_1pc/[dr]_{\otimes}\rtwocell<5>^{F \times_{\B}F}_{G \times_{\B}G}{\omit}\ar@{}[r]|(.5){\text{\small$\Downarrow$} \alpha \times_{\B} \alpha} & \E \times_{\B} \E\ar[r]^(.6){\otimes} & \E\\
 \rrtwocell<\omit>{<-4>\ \ \phi^G}& \D\ar@/_1pc/[ur]_G & }
=
\xymatrixnocompile@C=1.5cm@R=1cm{
\D \times_{\B} \D\ar[r]^{F \times_{\B}F}\ar@/_1pc/[dr]_{\otimes} & \E \times_{\B} \E \ar[r]^(.6){\otimes} & \E\\
\rrtwocell<\omit>{<-4>\ \ \phi^F} & \D\urtwocell^F_G{\alpha} &}
$$
and
$$
\xymatrix@=1.5cm{
\B\ar[r]^u\ar[dr]_u & \D\dtwocell^G_F{^\alpha}\\
& \E\ultwocell<\omit>{<3>\psi^F}}
=
\xymatrix@=1.5cm{
\B\ar[r]^u\ar[dr]_u & \D\ar[d]^G\\
& \E\ultwocell<\omit>{<3>\psi^G}}
$$
\end{Def}
Such a natural transformation determines a natural transformation $\alpha_{B}$ of the fibre-restricted functors, and the coherence axioms hold \ssi they hold fibrewise.
\paragraph{}
Monoidal fibred categories over \B, monoidal functors over \B and monoidal natural transformation over \B determine a 2-XL-category, which we denote $\MONFIB(\B)$\index[not]{MONFIBB@$\MONFIB(\B)$}. Another 2-XL-category of interest is $\MONFIB_s(\B)$\index[not]{MONFIBBs@$\MONFIB_s(\B)$}, the 2-cell-full sub-2-XL-category of strong monoidal fibred categories over \B and strong monoidal functors over \B. Finally, one can also consider sub-2-XL-categories of those by imposing Cartesianness of functors.

Symmetric monoidal fibred categories over \B, symmetric monoidal functors over \B, and monoidal natural transformations over \B also form a 2-XL-category, denoted $\SMONFIB(\B)$\index[not]{SMONFIBB@$\SMONFIB(\B)$}, and one has the same variants, $\SMONFIB_s(\B)$\index[not]{SMONFIBBs@$\SMONFIB_s(\B)$} by imposing strength of objects and morphisms, and by imposing Cartesianness of morphisms.

\subsubsection{The indexed category setting}

We first consider a monoidal fibred category $P\colon\E\to\B$. Recall that a choice of a cleavage gives rise to an indexed category:
$$\begin{array}{rcl}
\Phi_P\colon\B^{op}	&\longrightarrow &\CAT\\
B	& \longmapsto&\E_{B}\\
A\to B & \longmapsto& f^{*}\colon\E_{B}\to\E_A.
\end{array}$$
Now, as noticed above, each fibre comes with a monoidal structure $(\E_B,\otimes_B,I_B,\alpha_B,\lambda_B,\rho_B)$ given by restriction of the global structure. Moreover, as one readily verifies using universality of certain Cartesian morphisms, for each $f\colon A \to B$ in \B, the associated inverse image functor $f^*$, admits a monoidal structure $(f^*,\phi^f,\psi^f)$ with respect to the fibre monoidal categories. It is defined by the following diagrams. Let $D$ and $E$ be objects of $\E_B$.
$$
\xymatrix{
f^*(D \otimes E) \ar[r]^-{\bar f_{D \otimes E}} & D \otimes E\\
f^*D \otimes f^*E \ar[ur]_{\bar f_D \otimes\bar f_E}\ar@{-->}[u]^{\exists!\phi^f_{D,E}} &\\
A \ar[r]_{f}& B}
\qquad\text{and}\qquad
\xymatrix{
f^*(I_B) \ar[r]^-{\bar f_{I_B}} & I_B\\
I_A\ar[ur]_{u(f)}\ar@{-->}[u]^-{\exists!\psi^f}&\\
A \ar[r]_{f} & B}
$$
Finally, the structure isomorphisms of $\Phi_P$ are also monoidal. Therefore, one obtains a pseudo-functor into the 2-XL-category of monoidal categories, monoidal functors and monoidal natural transformations:
$$
\Phi_P\colon\B^{op}\longrightarrow \MONCAT.
$$
Moreover, if the fibred category \E is strong monoidal over \B, then each $f^*$ is strong monoidal and therefore one obtains a pseudo-functor into the 2-XL-category $\MONCAT_s$ of monoidal categories, strong monoidal functors and monoidal natural transformations: 
$$
\Phi_P\colon\B^{op}\longrightarrow \MONCAT_s.
$$
Conversely, a tedious but straightforward verification shows that the Grothendieck construction $\E_{\Phi}\ra{P_{\Phi}}\B$ of a pseudo-functor $\Phi\colon\B^{op}\to \MONCAT$ is a monoidal fibred category, strong monoidal if $\Phi$ takes values in $\MONCAT_s$.

Notice that in case of a monoidal \emph{opfibred} category, the direct image functors $f_*$ are \emph{comonoidal}\footnote{These functors are also called opmonoidal, colax monoidal or oplax monoidal.}. Consequently, they don't preserve monoids and modules unless they are strong.

\paragraph{}

What about monoidal functors over a category? Let $F$ be a monoidal functor over \B making the following triangle commute
$$\xymatrix{
\D\ar[rr]^F\ar[dr]_P & & \E\ar[dl]^Q\\
& \B & \\
}$$
where $P$ and $Q$ are monoidal fibrations. We already know that it induces an oplax natural transformations $\xymatrix@C=1.5cm{\B^{op}\rtwocell^{\Phi_P}_{\Phi_Q}{\omit\Downarrow(\tau_F,\xi_F)} & **[r]{\CAT}}$. Recall that its component $(\tau_F)_A$ at $A\in\B$ is just the restriction $F_A\colon\D_A\to\E_A$ of $F$. It admits therefore a monoidal structure $(F_A,\phi_A,\psi_A)$. Moreover its structure natural transformation $\xi_F^f\colon F_A\circ f^*\Rightarrow f^*\circ F_B\colon \D_B\to \E_A$ at $f\colon A\to B$ in \B is given by the following commutative diagram. Let $D$ be an object of $\D_B$.
$$\xymatrix@=1.5cm{
f^*(F_B(D)) \ar[r]^-{\bar f_{F_B(D)}} & F_B(D) & \E \ar[dd]^P\\
F_A(f^*(D)) \ar[ur]_{F(\bar f_D)}\ar@{-->}[u]^{\exists!(\xi^f_F)_D} & & \\
A \ar[r]_{f}& B & \B
}$$

\begin{Lem}
	The natural transformation $\xi_F^f$ is monoidal.
\end{Lem}

\begin{Pf}
We first prove the coherence with the `` \(\phi\)'s''. Let $D$ and $D'$ be objects of $\D_B$. Consider the following diagram.
$$\xymatrix@=1.3cm{
F_Af^*D\otimes F_Af^*D' \ar[rr]^{\xi^f_F \otimes \xi^f_{F}} \ar[dd]_{\phi_A} \ar[dr]_{F\bar f_D \otimes F\bar f_{D'}\quad} \ar@{}[dddr]|*+[o][F]{1}&  & f^*F_BD \otimes f^*F_BD'\ar[dd]^{\phi^f}\ar[dl]^{\quad\ \bar f_{F_BD} \otimes \bar f_{F_BD'}}\\
& F_BD \otimes F_BD'\ar[dd]^{\phi_B} \ar@{}[dddr]|*+[o][F]{2} & \\
F_A(f^*D \otimes f^* D') \ar[dd]_{F_A(\phi^f)} \ar[rd]^{\ F(\bar f_D \otimes \bar f_{D'})} & & f^*(F_BD \otimes F_BD') \ar[dd]^{f^*(\phi_B)} \ar[lu]_{\ \bar f_{F_BD \otimes F_BD'}}\\
& F_B(D \otimes D') & \\
F_Af^*(D \otimes D') \ar[rr]_{\xi^f_{F}} \ar[ru]^{F(\bar f_{D \otimes D'})\ }& & f^*F_B(D \otimes D') \ar[ul]_{\bar f_{F_B(D \otimes D')}}
}$$
The triangles commute by definition of the natural transformation $\xi_F^f$ and $\phi^f$. The square $\xymatrix@1{*+[o][F]{1}}$ commutes because of the naturality of \(\phi\). The square $\begin{xy}*+[o][F]{2} \end{xy}$ commutes by definition of the inverse image functors $f^*$. This implies that the outer square commutes, by the Cartesianness of the morphism $\bar f_{F_B(D \otimes D')}$ with domain the lower right hand vertex of the latter square.

We turn to the coherence with the ``\(\psi\)'s''. We distinguish the units of \D and \E by a superscript like $I_A^{\D}$ and $I_A^{\E}$. Consider the following diagram.
$$\xymatrix@C=1em{
&&&I_A^{\E} \ar[dl]_{\psi_A} \ar[d]^{u(f)} \ar[dr]^{\psi^f}&&&\\
&& F_A(I_A^{\D}) \ar[dr]_{F(u(f))} \ar[ddll]_{F_A(\psi^f)}& I_B^{\E} \ar[d]^{\psi_B}& f^*(I_B^{\E})\ar[l]^{\bar f_{I_B^{\E}}} \ar[ddrr]^{f^*(\psi_B)}&\\
&&&F_B(I_B^{\D}) &&&\\
F_Af^*(I_B^{\D}) \ar[urrr]^{F(\bar f_{I_B^{\D}})} \ar[rrrrrr]_{\xi_F^f} &&&&&& f^*F_B(I_B^{\D}) \ar[ulll]_{\bar f_{F_B(I_B^{\D})}}
}$$
As in the previous proof, inner triangles commute by definition of the respective natural transformations involved, and the inner squares by respectively naturality of $\psi$ and definition of the inverse image functors. One conclude the proof by Cartesianness of the morphism with source the right-down end of the outer triangle.
\end{Pf}
The monoidal functor $F$ over \B thus determines an oplax natural transformation
$$\shorthandoff{;:!?}\xymatrix@!{\B^{op}\rtwocell^{\Phi_P}_{\Phi_Q}{\ \tau_F} & *\txt{\ \qquad\MONCAT}}.$$

As in the non monoidal case, it is a pseudo-natural transformation \ssi $F$ is Cartesian. 

Moreover, $F$ is strong (resp. strict) if and only if each restriction $F_B$, $B\in\B$, is. If $F$ is strong monoidal between strong monoidal fibred categories, then it determines an oplax natural transformation
$$\shorthandoff{;:!?}\xymatrix@1@!{\B^{op}\rtwocell^{\Phi_P}_{\Phi_Q}{\ \tau_F} & *\txt{\ \qquad$\MONCAT_s$}}.$$

\paragraph{}
Finally, consider a monoidal natural transformation
$$\xymatrix@R=1cm{
\D\rrtwocell^F_G{\alpha}\ar[dr]_P & & \E\ar[dl]^Q\\
& \B. & 
}$$
It induces a modification between the corresponding oplax natural transformations, whose components are monoidal and therefore a modification filling the following diagram
$$\shorthandoff{;:!?}\xymatrix@!{\B^{op}\rtwocell<5>^{\Phi_P}_{\Phi_Q}{\omit\tau_F\big\Downarrow\stackrel{\Omega_{\alpha}}{\Rrightarrow}\big\Downarrow\tau_G} & *\txt{\quad \qquad\MONCAT}}.$$

\paragraph{}
There is a 2-XL-category $\MONCAT^{\B^{op}}$ of pseudo-functors, oplax natural transformations and modifications. The following result is unsurprising, but also desirable for a well-founded theory. Its proof is straightforward from the axioms, but (very) long.

\begin{Thm}
	The correspondence we have described is a 2-equivalence
$$
\MONFIB(\B)\xrightarrow{\sim}\MONCAT^{\B^{op}},
$$
lifting the equivalence between fibrations and indexed categories, i.e., the following diagram commutes
$$\xymatrix@=1.5cm{
\MONFIB(\B)\ar[r]^{\sim}\ar[d]_{U} & \MONCAT^{\B^{op}}\ar[d]^{V_*}\\
\FIB(\B)\ar[r]_{\sim} & \CAT^{\B^{op}}
}$$
where $U$ and $V$ are forgetful functors.\\ 
Moreover, this 2-equivalence restricts to strong objects and morphisms:
$$
	\MONFIB_s(\B)\xrightarrow{\sim}{\MONCAT_s}^{\B^{op}}.
$$
Both 2-equivalences remain true when morphisms of domains are restricted to Cartesian functors and morphisms of codomains to pseudo-natural transformations.\Par
Finally, these 2-equivalences induce 2-equivalences on the respective symmetric variants of these 2-XL-categories.

\cqfd
\end{Thm}

\begin{Def}\indexb{Monoidal indexed category}\indexb{Monoidal indexed category!Strong --}\indexb{Monoidal indexed category!Symmetric (strong) --}
	A pseudo-functor $\B^{op}\to\MONCAT$ is called a \emph{monoidal indexed category}. If it takes its values in $\MONCAT_s$, then it is called a \emph{strong monoidal indexed category}. Similarly, there is the notion of \emph{\op strong\fp symmetric monoidal indexed category}. We denote by \MONIND\B the 2-XL-category $\MONCAT^{\B^{op}}$ and by \MONINDs\B the 2-XL-category $\MONCAT_s^{\B^{op}}$.
\end{Def}

\subsubsection{Monoidal opfibred and bifibred categories}\label{sssec:monopbifib}
We first indicate some important differences between the opfibred and fibred situations. We then explain some facts that arise in the bifibred context and that simplify the matter of showing that something is a monoidal bifibred category. We do not develop the full theory of monoidal bifibrations and monoidal bi-indexed categories as we did for bifibrations and bi-indexed categories. We have a quite pragmatic point of view. Yet, we make use of the ideas developed for bi-indexed categories and try to shed light on the fundamental ideas that would lead to a complete description.

\emph{Monoidal \op\resp opstrong monoidal\fp opfibred categories over} $\B$ are monoidal objects in $\OPFIB(\B)$ (\resp $\OPFIB_{oc}(\B)$). One unpacks easily this definition inspired by the definition of (strong) monoidal fibred categories. 

By duality with the fibred case, a monoidal opfibred category $\E\to\B$ gives rise to a pseudo-functor $\B\to\MONCAT_{oplax}$, where $\MONCAT_{oplax}$ is the 2-XL-category of monoidal categories, oplax monoidal functors and monoidal natural transformations. Such a pseudo-functor is called a \emph{monoidal opindexed category}. If $\E\to\B$ is opstrong monoidal, then its associated monoidal opindexed category takes values in $\MONCAT_{os}$, the 2-cell and object full sub-2-XL-category of $\MONCAT_{oplax}$ consisting of opstrong oplax monoidal functors. Such a monoidal opindexed category is called \emph{opstrong}. Note that an \emph{opstrong oplax monoidal functor}, i.e., an oplax monoidal functor whose structure morphisms are isomorphisms, gives rise to a strong monoidal functor by taking the inverses of the structure morphisms. This correspondence yields a 2-isomorphism $\MONCAT_{os}\cong\MONCAT_s$. Therefore, an opstrong monoidal opindexed category can be considered as a pseudo-functor into $\MONCAT_s$ (just post-compose it with the previous 2-isomorphism).

Monoidal opindexed and opfibred categories are a priori designed for the study of comonoids and comodules\footnote{Furthermore, observe that a lax natural transformation between two pseudo-functors $$\B\to\MONCAT_{oplax}$$ gives rise, via Grothendieck op-construction, not to a monoidal, but indeed an oplax monoidal functor between the corresponding monoidal opfibrations.}. Since we are in this work mainly concerned with their dual objects, monoids and modules, we will mostly consider opstrong opindexed or opfibred categories.

Let us now turn to the bifibred notions. A \emph{monoidal bifibred category over} \B is, equivalently, a monoidal fibration that is opfibred or a monoidal opfibration that is fibred. It can be strong (as a monoidal fibration), opstrong (as a monoidal opfibration), or both, in which case it is called \emph{bistrong}. 

Let us come back to \thref{ex:monfib}.
\begin{Exs}
	\item
Any monoidal category \V is bistrong monoidal bifibred over the terminal object \one of \CAT.
	\item
Consider a category \C with pullbacks. As it is a bifibration, the canonical fibration $\C^{\two}\ra{\cod}\C$ over \C is automatically a monoidal bifibred category, but it happens not to be opstrong monoidal opfibred in general: it is so \ssi the base category is a groupoid. 

For, if $\cod$ is opstrong monoidal, then the unit functor is opCartesian. This implies that for all $f\colon A\to B$ in \C, the morphism $f_*u(A)\ra f u(B)$ in $\C^\two$ is an isomorphism, i.e., \C is a groupoid. Conversely, if \C is a groupoid, then the unit is opCartesian by the previous remark. The product is also opCartesian because then, for all morphisms $f\colon A\to B$ in the base \C, and every object $D\xrightarrow{q}A\xleftarrow{p}E$ of $\C^\two\times_{\C}\C^\two$ over $A$, the pullback $D\times_A E$  is also a pullback of the pair $D\xrightarrow{f\circ q}B\xleftarrow{f\circ p}E$.
\end{Exs}

Let $P\colon\E\to\B$ be a monoidal bifibration and choose a bicleavage of $P$. Then, by \thref{lem:biindexedbifib}, there is a bi-indexed category
$$\begin{array}{rcl}
\Omega_P\colon\B& \longrightarrow	& \ADJ\\
A 		& \longmapsto		& \E_A\\
A \ra{f} B 	& \longmapsto		& \xymatrix{f_*\colon\E_A \ar@<.9ex>[r]\ar@{}[r]|{\perp} & \E_B:f^*,\ar@<.9ex>[l]}
\end{array}
$$
Moreover, we know that the associated opindexed and indexed categories of this bi-indexed category are monoidal. Therefore, the functor $f^*$ is monoidal and the functor $f_*$ is oplax monoidal. One also concludes that the structure isomorphisms $({\gamma_{f,g}}_*,{\gamma_{f,g}}^*)$ and $({\delta_A}_*,{\delta_A}^*)$ of this bi-indexed category are pairs of monoidal natural transformations. This is not all, though. The respective oplax monoidal and monoidal structures of $f_*$ and $f^*$ are closely related, as expressed in the following lemma. This result is new, as far as we know.

\begin{Lem}\label{lem:monbiindcat}
	Let $P\colon\E\to\B$ be a monoidal bicloven bifibration and $f\colon A\to B$ an arrow in \B. Let ${f_*\dashv f^*}$ be the corresponding adjunction of direct and inverse images, whose oplax monoidal and monoidal structure morphisms are respectively $(\phi^{f_*},\psi^{f_*})$ and $(\phi^{f^*},\psi^{f^*})$. Then, the following are mate pairs, i.e., cells in the double category $\frak{ADJ}$.
$$\xymatrix@=2cm{
\E_A\times \E_A\ar[d]_{\otimes}\ar@<-.9ex>[r]_{f_*\times f_*}\ar@{}[r]|{\top} & \E_B\times \E_B\ar@<-.9ex>[l]_{f^*\times f^*}\ar[d]^{\otimes}\\
\E_A \ar@<-.9ex>[r]_{f_*}\ar@{}[r]|{\top} & \E_B\ar@<-.9ex>[l]_{f^*}\ar@{}[lu]|{(\phi^{f_*},\phi^{f^*})}
}\qquad\xymatrix@=2cm{
\one\ar[d]_{I_A}\ar@<-.9ex>[r]_{Id}\ar@{}[r]|{\top} & \one\ar[d]^{I_B}\ar@<-.9ex>[l]_{Id}\\
\E_A \ar@<-.9ex>[r]_{f_*}\ar@{}[r]|{\top} & \E_B\ar@<-.9ex>[l]_{f^*}\ar@{}[lu]|{(\psi^{f_*},\psi^{f^*})}
}$$
\end{Lem}

\begin{Pf}[Sketch]
	One proves that the latter pairs are mate by means of \thref{axiom:matepair} and by the characterization of the transpose morphisms for an adjunction $f_*\dashv f^*$ given in \thref{prop:bifibadjoint}.
\end{Pf}

This lemma leads to the notion of \emph{monoidal bi-indexed category}. This is a bi-indexed category with the following ingredients. Its values on objects are monoidal categories. Its direct and inverse image functors are respectively oplax monoidal and monoidal whose structure morphisms are mates as in the previous lemma. Finally, its structure isomorphisms are conjugate pairs of monoidal natural transformations. This is of course a pseudo-functor from \B to some monoidal version of \ADJ, or a pseudo double functor from $\frak H\B$ to some monoidal version of $\frak{ADJ}$, but as previously said, we will remain pragmatic in this part. A monoidal bi-indexed category is respectively \emph{strong\textnormal{,} opstrong \textnormal{or} bistrong} if its associated monoidal indexed category is strong, its associated monoidal opindexed category is opstrong, or if both conditions are satisfied.

Recall that the notion of bi-indexed category is equivalent to that, for instance, of an indexed category whose inverse image functors admit a left adjoint. One can then, by a choice of these left adjoint for each inverse image functor, construct a bi-indexed category. What about the monoidal context? 

Let $\Phi\colon\B^{op}\to\MONCAT$ be a monoidal indexed category. Suppose moreover that the inverse image functor $f^*$ admits a left adjoint $f_*$ for each $f$ in \B. Then, there exists a unique structure of oplax monoidal functor on each $f_*$ such that the properties of \thref{lem:monbiindcat} hold. This is due to the following classical fact.
\begin{Prop}\label{prop:indmonlefadj}
	Let 
	\begin{equation}\label{eq:adjmonopmon}\xymatrix@1@C=1.3cm{
	F:\V \ar@<-.9ex>[r]\ar@{}[r]|{\top} & \V':G.\ar@<-.9ex>[l]
}\end{equation}
be an adjunction of unit $\eta$, with $(\V,\otimes, I)$ and $(\V',\otimes',I')$ monoidal categories. Let $(G,\phi_G,\psi_G)$ be a monoidal functor.\Par
There exists a structure of oplax monoidal functor on $F$. Its structure morphisms $(\phi_F)_{A,B}$ at objects $A,B\in\V$ and $\psi_F$ are the transpose morphisms under the adjunction $F\dashv G$ of the composite morphisms
\begin{gather}
A\otimes B \xlongrightarrow{\eta_A\otimes \eta_B} GF(A)\otimes GF(B)\xlongrightarrow{(\phi_G)_{FA,FB}}G(F(A)\otimes F(B))\label{eq:oplaxmon1}\\[1em]
I\xlongrightarrow{\quad\psi_G\quad} G(I').\label{eq:oplaxmon2}
\end{gather}
Moreover, when $F$ with this structure is opstrong and is given a structure of strong monoidal functor by inverting $\phi_F$ and $\psi_G$, this becomes an adjunction in $\MONCAT$. Conversely, if \hyperref[eq:adjmonopmon]{(\ref*{eq:adjmonopmon})} is an adjunction in $\MONCAT$ with $F$ strong, then the monoidal structure of $F$ is given by the inverses of the transposes of \hyperref[eq:oplaxmon1]{(\ref*{eq:oplaxmon1})} and \hyperref[eq:oplaxmon2]{(\ref*{eq:oplaxmon2})}.\Par
This proposition also holds in the symmetric setting.
\end{Prop}

\begin{Pf}[Sketch]
Our attention was drawn to this fact in \cite{LM09}, which mentions it without proof. It goes back, as far as we know, to the article of Kelly \cite{Kel74}, where it appears as a particular case of a more general theory. We sketch a direct proof. 

The transpose morphisms of \hyperref[eq:oplaxmon1]{(\ref*{eq:oplaxmon1})} and \hyperref[eq:oplaxmon2]{(\ref*{eq:oplaxmon2})} are precisely the components of the mate natural transformations $\phi_F$ and $\psi_F$ of $\phi_G$ and $\psi_G$ in 
$$\xymatrix@=2cm{
\V\times \V\ar[d]_{\otimes}\ar@<-.9ex>[r]_{F\times F}\ar@{}[r]|{\top} & \V'\times \V'\ar@<-.9ex>[l]_{G\times G}\ar[d]^{\otimes}\\
\V \ar@<-.9ex>[r]_{F}\ar@{}[r]|{\top} & \V'\ar@<-.9ex>[l]_{G}\ar@{}[lu]|{(\phi^{F},\phi^{G})}
}\qquad\xymatrix@=2cm{
\one\ar[d]_{I_A}\ar@<-.9ex>[r]_{Id}\ar@{}[r]|{\top} & \one\ar[d]^{I_B}\ar@<-.9ex>[l]_{Id}\\
\V \ar@<-.9ex>[r]_{F}\ar@{}[r]|{\top} & \V'\ar@<-.9ex>[l]_{G}\ar@{}[lu]|{(\psi^{F},\psi^{G})}
}$$
Now one must check that the coherence axioms of monoidal functors are preserved by the bijection that hold between mates (more precisely, that if $\phi_G$ and $\psi_G$ satisfy the axioms for monoidal functors, then their mates automatically satisfy the axioms for oplax monoidal functors). Recall that it not true that any property of a natural transformation is transmitted to its mate. For instance, the property of being a natural isomorphism is not. In fact, the properties that are transmitted between mates are the properties that can be expressed in terms of the double category $\frak{ADJ}$, because the bijection between mates preserve all the structure of this double XL-category. One verifies that the axioms of a monoidal functor can be expressed in  terms of $\frak{ADJ}$ by writing them ``globally'', that is, as diagrams of the natural transformations $\phi$ and $\psi$, not objectwise as they are usually given. See below \hyperref[eq:conjpairmon]{(\ref*{eq:conjpairmon})} for an example of this kind of argument.

For the last statement of this proposition, one has first to prove that the unit and counit of the adjunction \hyperref[eq:adjmonopmon]{(\ref*{eq:adjmonopmon})} are monoidal natural transformations with respect to the induced strong monoidal structure on $F$. This is a quite direct application of the definitions of $\phi_F$ and $\psi_F$ as transposes of \hyperref[eq:oplaxmon1]{(\ref*{eq:oplaxmon1})} and \hyperref[eq:oplaxmon2]{(\ref*{eq:oplaxmon2})} and, equivalently, as mates of $\phi_G$ and $\psi_G$ (recall also the formulae of transpose morphisms in terms of the unit or the counit). For instance, the axioms of monoidal transformation applied on the unit of the adjunction express precisely that $\phi_F$ and $\psi_F$ are defined by the transposes of \hyperref[eq:oplaxmon1]{(\ref*{eq:oplaxmon1})} and \hyperref[eq:oplaxmon2]{(\ref*{eq:oplaxmon2})}. This proves the last affirmation.
\end{Pf}
\begin{Rem}
	This Proposition shows that an opstrong monoidal bi-indexed category over \B gives rise to a pseudo-functor $\B\to\ADJ(MONCAT)$, the 2-XL-category of adjunctions in \MONCAT. Bistrong monoidal bi-indexed categories over \B are precisely pseudo-functors $\B\to\ADJ(\MONCAT_s)$.
\end{Rem}

Let us come back to our monoidal indexed category $\Phi\colon\B^{op}\to\MONCAT$ whose inverse image functors $f^*$ admit a left adjoint. If we choose a left adjoint $f_*$ for each inverse image functor, we obtain, by \thref{thm:Bifib}, a bi-indexed category $\B\to\ADJ$. In order to prove that this in fact is a monoidal bi-indexed category, the only fact that remains to be proved is the fact that the structure isomorphisms  ${\gamma_{f,g}}_*$ and ${\delta_A}_*$ of the associated opindexed category are monoidal with respect to the oplax monoidal structure just defined on the direct image functors. The latter structure isomorphisms are defined as conjugates of the structure isomorphisms ${\gamma_{f,g}}^*$ and ${\delta_A}^*$ of the original monoidal indexed category $\Phi$, which are thus monoidal. In order to prove that ${\gamma_{f,g}}_*$ and ${\delta_A}_*$ are monoidal, one must verify that the condition of being monoidal for a natural transformation is expressible in $\frak{ADJ}$. We recommend to use the intrinsic, i.e., objectfree, version of the axioms, as given in \hyperref[eq:monnattransf:phi]{(\ref*{eq:monnattransf:phi})} and \hyperref[eq:monnattransf:psi]{(\ref*{eq:monnattransf:psi})}. Given adjunctions $F\dashv G\colon\V\to\V'$, $F'\dashv G'\colon\V\to\V'$ and a conjugate pair 
\begin{equation}\label{eq:conjpairmon}\xymatrix@1@C=3cm{\V\rtwocell^{F\dashv G}_{F'\dashv G'}{\qquad\left(\alpha_*,\alpha^*\right)} & \V'},\end{equation}
the fact that $\alpha_*$ and $\alpha^*$ are monoidal can be expressed entirely in diagrams in $\frak{ADJ}$. The first coherence condition is given by the following equality in $\frak{ADJ}$.
$$\xymatrix@=1.5cm{
\V\times \V\ar@{}[rd]|{(\alpha_*\times\alpha_*,\alpha^*\times\alpha^*)}\ar@{=}[d]\ar@<-.9ex>[r]_{F\times F}\ar@{}[r]|{\top} & \V'\times \V'\ar@<-.9ex>[l]_{G\times G}\ar@{=}[d]\\
\V\times \V\ar[d]_{\otimes}\ar@<-.9ex>[r]_{F'\times F'}\ar@{}[r]|{\top} & \V'\times \V'\ar@<-.9ex>[l]_{G'\times G'}\ar[d]^{\otimes}\\
\V \ar@<-.9ex>[r]_{F'}\ar@{}[r]|{\top} & \V'\ar@<-.9ex>[l]_{G'}\ar@{}[lu]|{(\phi'_{*},\phi'^{*})}
}\quad =\quad
\xymatrix@=1.5cm{
\V\times \V\ar[d]_{\otimes}\ar@<-.9ex>[r]_{F\times F}\ar@{}[r]|{\top} & \V'\times \V'\ar@<-.9ex>[l]_{G\times G}\ar[d]^{\otimes}\\
\V\ar@{=}[d] \ar@<-.9ex>[r]_{F}\ar@{}[r]|{\top} & \V'\ar@<-.9ex>[l]_{G}\ar@{}[lu]|{(\phi_{*},\phi^{*})}\ar@{=}[d]\\
\V \ar@<-.9ex>[r]_{F'}\ar@{}[r]|{\top} & \V'\ar@<-.9ex>[l]_{G'}\ar@{}[lu]|{(\alpha_*,\alpha^*)}
}$$
\begin{Ex}[Sheaves of abelian groups]
Let $X$ be a topological space. Since $\Ab$ is symmetric monoidal for the tensor product, the category of presheaves of abelian groups $\PSh(X;\Ab)$ inherits the corresponding objectwise symmetric monoidal structure (here the objects are the open subsets of $X$). This symmetric monoidal structrure on presheaves induces a symmetric monoidal structure on the category of sheaves of abelian groups $\Sh(X;\Ab)$ by postcomposing it with the sheafification functor. This can be shown using the universality of the unit $\eta_P\colon P\Rightarrow aP$, $P\in \PSh(X;\Ab)$, of the sheafification adjunction $a\dashv I$. In particular, its unit is the constant sheaf $\tilde{\Delta}_X(\Z)$. 

Note that with this monoidal structure, $\Sh(X;\Ab)\in\rcSYMMON$. Indeed, it has all colimits, since it is a reflective subcategory of the cocomplete category of presheaves $\PSh(X;\Ab)$. Moreover, its monoidal product functors $-\otimes F$ preserve colimits because the monoidal product of $\PSh(X;\Ab)$ does and the sheafification functor is a left adjoint \cite{KS06}.

The direct image functor $f_*$ along a map $f\colon X\to Y$ can be shown to be strong symmetric monoidal. One proves this by using the fact that the direct image functor of presheaves is strict monoidal with respect to the objectwise monoidal product, and the universality of the unit of the sheafification adjunction $a\dashv I$. There is thus a strong symmetric monoidal indexed category
$$\begin{array}{rcl}
(\Top^{op})^{op}	&\longrightarrow &\SYMMON_s\\
X	& \longmapsto&\Sh(X;\Ab)\\
f^{op}\colon Y\to X & \longmapsto& f_{*}\colon\Sh(X;\Ab)\to\Sh(Y;\Ab).
\end{array}$$
Each functor $f_*$ admits a left adjoint $f^{-1}$. Let us put on $f^{-1}$ the symmetric monoidal structure induced from the monoidal structure of $f_*$ as in Proposition \ref{prop:indmonlefadj}. This symmetric monoidal functor happens to be strong \cite{LM09}. One finally obtains, as explained above, a bistrong monoidal bi-indexed category over $\Top^{op}$:
\begin{equation}\label{eq:bifibsheaves}\begin{array}{rcl}
\Top^{op}	&\longrightarrow &\ADJ(\SYMMON_s)\\
X	& \longmapsto&\Sh(X;\Ab)\\
f^{op}\colon Y\to X & \longmapsto& \xymatrix@1@C=1.3cm{\Sh(Y;\Ab) \ar@<-.9ex>[r]_{f^{-1}}\ar@{}[r]|{\top} & \Sh(X;\Ab).\ar@<-.9ex>[l]_{f_*}}
\end{array}\end{equation}
\end{Ex}

\subsection{Fibred algebra}
\label{ssec:fibred-algebra}
We define in this part the notions of monoids and modules in a monoidal fibred category and show that they give rise to a fibration of modules over a fibration of monoids. This material is new.
\subsubsection{The fibration of monoids}
Given a monoidal fibred category, we are interested in the categories of monoids in the fibre monoidal categories and in the relationship between them induced by the fibred structure.
\paragraph{Indexed category setting}
In the indexed category setting, it is very easy to define the fibration of monoids in a monoidal indexed category. Consider a monoidal indexed category $\Phi\colon \B^{op}\to \MONCAT$. One obtains fibrewise monoid categories by just post-composing $\Phi$ with the monoid 2-functor $$\Mon\colon\MONCAT\to\CAT$$ of \thref{prop:Mon}. This composition is thus an indexed category:
\begin{equation}\label{eq:monindcat}\begin{array}{rcl}
Mon\circ\Phi\colon\B^{op} & \longrightarrow & \CAT\\
A & \longmapsto & \Mon(\Phi(A))\\
f\colon A\to B & \longmapsto & \Mon(f^*)\colon \Mon(\Phi(B))\to Mon (\Phi(A)),
\end{array}
\end{equation}
which we call the \emph{indexed category of monoids} (\emph{in the monoidal indexed category $\Phi$})\indexb{Monoid in a monoidal indexed category}. We denote its Grothendieck construction 
\begin{equation}\label{eq:fibmon}
\Mon(\Phi)\to\B,
\end{equation} 
and call it the \emph{fibration of monoids in $\Phi$} and its total category the \emph{category of monoids in $\Phi$}.

The artillery of bicategory theory even provides us directly with a \emph{monoid 2-functor}
$$Mon_{\B}\colon \MONCAT^{\B^{op}}\longrightarrow \CAT^{\B^{op}},$$
which is the 2-functor $Mon_*$ given by post-composition with the 2-functor $Mon$. It plays exactly the same role as the functor $Mon\colon\MONCAT\to\CAT$, but in the context of monoidal indexed categories over \B.

When $\Phi$ is symmetric, the 2-functor $$\Comm\colon\SYMMON\to\CAT$$ of \thref{prop:Mon} yields similarly a \emph{fibration of commutative monoids in $\Phi$}\indexb{Monoid in a monoidal indexed category!Commutative --}.

\begin{Ex}[Ringed spaces]\label{ex:bifibringedspaces}
We have defined a bistrong monoidal bi-indexed category $\Top^{op}\to\ADJ(\SYMMON_s)$ in \hyperref[eq:bifibsheaves]{(\ref*{eq:bifibsheaves})}. It its post-composition with the 2-functor $\ADJ(\Comm)$ takes values in $\ADJ$. One can then Grothendieck-construct it, obtaining a bifibration over $\Top^{op}$ and take the dual to finally obtain a bifibration over \Top. Alternatively, one can post-compose it with $^{op}\colon\CAT\to\CAT$ and then Grothendieck op-construct it to finally obtain a bifibration over \Top. These two processes give the same result, as explained in \autoref*{ssec:opfibbifib} on page \pageref{page:grothconsopcons}. Let us consider the second process. One thus first obtains a bi-indexed category:
\begin{equation}\label{eq:biindsheavesab}\begin{array}{rcl}
\Top	&\longrightarrow &\ADJ\\
X	& \longmapsto&\ \Comm(\Sh(X;\Ab))^{op}\\
f\colon X\to Y & \longmapsto& \xymatrix@1@C=2cm{\Comm(\Sh(X;\Ab))^{op}\  \ar@<-.9ex>[r]_{\Comm(f_*)^{op}}\ar@{}[r]|{\top} & \Comm(\Sh(Y;\Ab))^{op}.\ar@<-.9ex>[l]_{\Comm(f^{-1})^{op}}}
\end{array}\end{equation}
It is isomorphic to the bi-indexed category of sheaves on topological spaces with values in \Comm defined in page \pageref{page:bincatsheavestop}. One can show that $\Comm(\Sh(X;\Ab))\cong\Sh(X;\Comm)$ this way. One first proves that the category $\Sh(X;\Comm)$ is isomorphic to the category of commutative internal rings $\cRing(\Sh(X;\Set))$ in the Cartesian monoidal category $\Sh(X;\Set)$. The same is true for sheaves of abelian groups and internal abelian groups in $\Sh(X;\Set)$. This works indeed in the same fashion as for presheaves, since the categorical product of presheaves restricts to a categorical product of sheaves, by reflectivity of the subcategory of sheaves. Again, in the case of presheaves, this comes down from the fact that the categorical product is objectwise (see the proof of \thref{LimitCreation} for more details). 

Secondly, in order to prove that $\cRing(\Sh(X;\Set))\cong\Comm(\Sh(X;\Ab))$, one uses the fact that the monoidal product of $\Sh(X;\Ab)\cong\Ab(\Sh(X;\Set))$ has the universal property of a tensor product (see \thref{def:tensprod}) and \cite{Ten75}.

The Grothendieck op-construction of \hyperref[eq:biindsheavesab]{(\ref*{eq:biindsheavesab})} is thus isomorphic to the bifibration of ringed spaces defined in page \pageref{page:bifibAspace}.
\end{Ex}

\paragraph{Fibration setting}
Now we would like to have an intrinsic, choice-independent, definition of a fibration of monoids.

Let $P\colon \E\to \B$ be a monoidal fibred category. The preceding paragraph gives us a manner to define a fibration of monoids in $P$. Consider, at first, the associated pseudo-functor $\B^{op}\to \MONCAT$, then its image under the monoid 2-functor $Mon_{\B}$ and finally its Grothendieck construction. This does the job. 

Nevertheless, this construction is not completely satisfactory. Indeed, in this way we don't obtain ``the'' category of monoids in $P$, but one for each cleavage of $P$. All these fibrations are of course isomorphic in $\FIB(\B)$, but there is no preferred choice between them unless $P$ has a canonical choice of cleavage, which is the case if $P$ is the Grothendieck construction of a pseudo-functor into \CAT. There is also an objection about the complexity of this definition. The Grothendieck construction has a quite complicated composition law, due to the fact that one has a chosen cleavage. Moreover, one needs such tools as global axiom of choice and bi-XL-category theory%
\footnote{Defining just the inverse image functors $f^*$ requires \emph{global axiom of choice}, an axiom of choice for proper classes. Such an axiom is, in von Neumann axiomatization of $NBG$, a consequence of the axiom of \emph{limitation of size}. But, in the axiomatization of $NBG$ proposed by Mendelson \cite{Men97}, the axiom of choice is independent, and the author imposes only its set version. The global axiom of choice is commonly used in category theory, though, e.g., for defining (co-)limit functors.

Now, if one want to state that these inverse image functors together determine a pseudo-functor into \CAT, one needs a much stronger axiom. Indeed, the category \CAT lives on the third floor of a three-level set theory with \emph{sets}, \emph{classes} and \emph{conglomerates} (this can be axiomatized, for example, by $ZFC$+\textit{Existence of an inaccessible cardinal}. See \autoref{cha:Found}). The problem of global choice vanishes here, because one wants anyway an axiom of choice for the third level of conglomerates, which implies the axiom of choice for classes and, \textit{a fortiori}, sets.\label{fn:GAC}}.%
Why should we need this sophistication for the sake of defining such a simple notion as a the fibred category of monoids in a monoidal fibred category? 
 
There is actually no need for these concepts in order to define the fibration of monoids in $P$, but avoiding them makes the job less straightforward since much more mathematics is done functorially than universally. The fibration, if we suitably define it, should be isomorphic to the fibration obtained by going through the pseudo-functor world.

\begin{Def}\indexb{Monoid in a monoidal fibration}\index[not]{MonE@$\Mon(\E)$}\index[not]{MonP@$\Mon(P)$}
	Let $\E\ra{P} \B$ be a monoidal fibred category.\Par
The \emph{category of monoids in \E} (\emph{or in $P$}), denoted $\Mon(\E)$ (or $\Mon(P)$), is defined by:
\begin{itemize}
\item
$\Ob \Mon(\E)$: Monoids in the fibres $\E_B$, for all $B\in\B$.
\item
$\Mor \Mon(\E)$: A morphism $(R,\mu,\eta) \ra{\phi} (S,\nu,\lambda)$ is a morphism $\phi\colon R\to S$ in \E, such that 
$$\xymatrix@=1.3cm{
R \otimes R \ar[r]^{\phi\otimes\phi}\ar[d]_\mu & S\otimes S\ar[d]^{\nu}\\
R\ar[r]_\phi & S
}\quad\quad
\xymatrix@=1.3cm{
I_{P(R)}\ar[d]_\eta \ar[r]^{u(P(\phi))} & I_{P(S)}\ar[d]^{\lambda}\\
R\ar[r]_\phi & S
}$$

\item
Composition: composition of $\E$.
\end{itemize}
When the monoidal fibration is symmetric, then one also considers the \emph{category commutative monoids in \E}\indexb{Monoid in a monoidal fibration!Commutative}, denoted $Comm(\E)$\index[not]{CommE@$\Comm(\E)$}, which is the full subcategory of $\Mon(\E)$ whose objects are commutative monoids in the fibres.
\end{Def}

There is an obvious projection functor $\Mon(\E)\to\B$ and $\Mon(\E)_{B}=\Mon(\E_{B})$. The following lemma provides us with Cartesian morphisms for it. The lemma can be proved without a choice of inverse image functors for $P\colon\E\to\B$. The notation $f^*(E)\to E$ for \textit{a} Cartesian morphism should not mislead the reader on this point.

\begin{Lem}\label{lem:fibremonoid}
	Let $\E\to \B$ be a monoidal fibred category and $f\colon A\to B$ an arrow in \B.
\begin{enumerate}[(i)]
\item
Let $(R,\mu,\eta)$ be a monoid in $\E_B$ and $f^* R \ra{\bar f_R}R$ a Cartesian arrow over $f$. Then, there is a unique monoid structure on $f^*R$ in $\E_A$ such that $\bar f_R$ is a morphism of monoids in $P$. It is defined by the left legs of the following commutative diagrams.
\begin{equation}\label{eq:uniquemon}\begin{aligned}\xymatrix@=1.5cm{
f^*(R) \otimes f^*(R) \ar[r]^-{\bar f_R \otimes\bar f_R} \ar@{-->}[d]& R \otimes R\ar[d]^{\mu}\\
f^*(R) \ar[r]_{\bar f_R} & R}
\qquad
\xymatrix{
I_A \ar[r]^{u(f)} \ar@{-->}[d]& I_B \ar[d]^{\eta}\\
f^*(R) \ar[r]_{\bar f_R} & R
}
\end{aligned}\end{equation}
Moreover, if $R$ is commutative, then so is $f^*R$.
\item
Let $\phi\colon R\to S$ be a morphism of monoids in $\E_B$. Let $f^*(R)\ra{\bar f_R}R$ and $f^*(S)\ra{\bar f_S}S$ be Cartesian morphisms over $f$ and $f^*(\phi)$ the induced morphism in
$$\xymatrix{
f^*(R) \ar[r]^{\bar f_R} \ar@{-->}[d]_{f^*(\phi)}& R \ar[d]^{\phi}\\
f^*(S) \ar[r]_{\bar f_S} & S.
}$$
Then $f^*(\phi)$ is a morphism of the induced monoids in $\E_A$. \\
In particular, any two inverse images of a monoid over the same arrow have isomorphic induced monoid structures.
\end{enumerate}

\qed
\end{Lem}
\begin{Pf}
	(i) The uniqueness of the monoid structure is clear. Indeed, the multiplication and the unit of such a monoid must make the two squares of (\ref{eq:uniquemon}) commute, by definition of a morphism of monoids in $P$. Since $\bar f_R$ is Cartesian, multiplication and unit morphisms are unique, and given by the left legs of (\ref{eq:uniquemon}). 
	
	The preceding observation gives also the existence of the maps. One then checks that this yields a monoid structure by using the monoid structure of the original monoid $R$ and the naturality of the left unitor and the associator of \E. The commutativity of $f^*R$ when $R$ is commutative follows also readily, using the naturality of the symmetry isomorphism $\sigma$.
	
	(ii) The fact that the inverse image $f^*\phi$ of a morphism of monoids $\phi$ is again a morphism of monoids is very easy to check. Applying this result to the identity morphism of a monoid, one obtains the last affirmation by noticing that a morphism of monoids that is an isomorphism in the ambient monoidal category is an isomorphism of monoids.  
\end{Pf}

In the following, given a Cartesian arrow $f^*(R)\ra{\bar f_R}R$, when talking about ``the monoid'' $f^*(R)$, we refer to the monoid structure defined in the lemma.

\begin{Rem}\label{rem:charcmormon}
Let $R$ and $S$ be monoids in \E. Let $\phi\colon R\to S$ be a morphism in $\E$ over $f\colon A\to B$ in \B. The morphism $\phi$ is a morphism of monoids in $P$ \ssi for some (and then for all) Cartesian lift $f^*(S) \ra{\bar f_S} S$ of $f$ at $S$, the morphism $\bar\phi$ defined by
$$\xymatrix{
R \ar@{-->}[d]_{\bar\phi} \ar[rd]^{\phi} & \\
f^*(S) \ar[r]_{\bar f_S} & S\\
A\ar[r]_f & B
}$$	
is a morphism of monoids in $\E_A$.
\end{Rem}

\begin{Prop}
Let $\E\ra{P} \B$ be a monoidal fibred category.
	\begin{enumerate}
\item
The projection $\Mon(\E)\to \B$ is also a fibration, called the \emph{fibration of monoids in \E} (\emph{or in $P$}).\\
Given a monoid $R$ over $B$ and a morphism $f\colon A\to B$, a Cartesian lift of $f$ at $R$ is given by
\begin{equation}\label{eq:moncarmor}
f^*(R) \ra{\bar f_R}R.
\end{equation}
When \E is symmetric, the projection $\Comm(\E)\to\B$ is a also a fibration, with the same Cartesian morphisms. 
\item
Let $\Mon(\Phi_P)\to\B$ be the fibration of monoids of the indexed category $\Phi_P$ determined by a cleavage of $P$. Then there is an isomorphism of fibrations over \B:
\begin{equation}\label{eq:iso2constr}\begin{aligned}\xymatrix@C=.5em@R=1cm{
\Mon(\Phi_P)\ar[rr]^{\cong}\ar[rd] && \Mon(P)\ar[ld]\\
& \B. &
}\end{aligned}\end{equation}
This result also holds in the commutative setting.
\end{enumerate}
\end{Prop}
\begin{Pf}
	1. We already know that the morphism (\ref{eq:moncarmor}) is a morphism in $\Mon(\E)$ over $f$. In order to check it is Cartesian, consider a morphism $\phi\colon S\to R$ in $\Mon(\E)$ over $g\colon C\to B$ and $k\colon C\to A$ in \B such that $f\circ k=g$. Then, one has the following commutative diagram in $P\colon\E\to\B$.
	$$\xymatrix@=1.5cm{
		S\ar@{-->}[d]_{\overline{\bar\phi}}\ar@{-->}[r]_{\bar\phi}\ar@/^2pc/[rr]^{\phi} &f^*R\ar[r]_{\bar f_R} & R & \E\ar[dd]^P\\
		k^*f^*R\ar[ru]_{\bar k_{f^*R}}&&\\
		C\ar[r]^{k}\ar@/_2pc/[rr]_{g} & A\ar[r]^{f} & B & \B
		}$$
Now, one can use characterization of \thref{rem:charcmormon} in order to show that $(k,\bar \phi)$ is a morphism of monoids. Indeed, the composite of Cartesian morphisms $\bar f_R\circ\bar k_{f^*S}$ is Cartesian. By this remark, $\overline{\bar\phi}$ is a morphism of monoids in $\E_C$. Moreover, $\bar k_{f^*R}$ is a morphism of monoids in \E. Therefore $\bar\phi$ is a morphism of monoids, as a composite of such.

2. Let us make a choice of cleavage on $P$. The isomorphism of fibration \eqref{eq:iso2constr} is defined in the following manner. It acts as identity on objects. As for morphisms, let $(f,\phi)\colon(A,R)\to(B,S)$ be a morphism in $\Mon(\Phi_P)$. Then $f\colon A\to B$ is in \B and $\phi\colon R\to \Mon(f^*)(S)$ is in $\Mon(\E_A)$. Notice that the functor $$\Mon(f^*)\colon\Mon(\E_B)\to\Mon(\E_B)$$ sends a monoid $S$ to a monoid whose underlying object is $f^*S$ and whose monoid structure is precisely the one determined in \thref{lem:fibremonoid}. One gets a morphism $\tilde\phi\colon R\to S$ in $\Mon(\E)$ by the composition $\bar f_S\circ\phi$ of $\phi$ with the Cartesian lift $\bar f_S$ of $f$ at $S$ of the chosen cleavage. It is then routine calculation to check that this indeed determines a functor.
\end{Pf}

\begin{Ex}\label{ex:fibmonoids}
As noted before (see \thref{ex:monfib}), the canonical fibration $$\C^\two\ra{\cod}\C$$ of a category \C with pullbacks is a strong monoidal fibred category. It admits therefore a fibration of monoids $\Mon(\C^\two)\to\C$. We study now a particular class of these monoids.

The canonical fibration over \C is actually a bifibration (see \thref{ex:bifibration}). Suppose that \C has in addition a terminal object $*$ (and therefore all finite limits). The morphism $A\to *$ determines an adjunction between direct and inverse image functors, whose composite with the isomorphism of categories $\C/*\cong\C$ is the following adjunction.
$$\xymatrix{
U:\C/A\ar@<1ex>[r]\ar@{}[r]|-{\perp} & \C:\theta_A^{(-)}.\ar@<1ex>[l]
}$$
$U$ is the forgetful functor. For an object $C\in\C$, $\theta_A^C$, is the product bundle $$\theta_A^C=(A\times C\ra{pr_1} A).$$ For a morphism $f\colon C\to D$ in $\C$, $\theta_A^f=1_A\times f$.

Moreover, each object $C\in\C$ determines a functor
$$
\theta^C_{(-)}\colon\C\to\C^{\two}.
$$
It is defined, for a morphism $f\colon A\to B$, by $\theta_A^C\ra{(f,f\times 1_C)}\theta_B^C$. Since $\theta_A$ is a right adjoint, it preserves internal monoids. One deduces that, for an internal monoid $C$ of \C, the functor $\theta_{(-)}^C$ factors through the forgetful functor $U$ from the category of monoids in $\C^\two$ to $\C^\two$.
$$\xymatrix{
& \Mon(\C^\two)\ar[d]^{U}\\
\C\ar@{-->}[ur]^{\bar\theta^C}\ar[r]_{\theta^C} & \C^\two
}$$
The monoids $\theta_A^G$ for $A\in\C$ and a group object $G$ in \C are of particular interest for the theory $G$-bundles (see \thref{ex:fibmodules} below). In this case, $\theta_A^G$ has the structure of an internal group of $\C/A$.
\end{Ex}

\subsubsection{The fibration of modules}\label{fibmodovfibmon}
We examine now the categories of $R$-modules for monoids $R$ in the fibre monoidal categories of a monoidal fibred category. There are two different kind of connections between these categories. One kind is of algebraic nature and takes place in the fibres: it is the adjunction of \thref{prop:adjointmodules} between categories of modules over different monoids related by a morphism of monoids in a given fibre. The other kind of connection has its origin in the fibred category and connect modules over some arrow in \B.
\paragraph{Indexed category setting} Let us start with a monoidal indexed category 
$${\Phi\colon\B^{op}\to\MONCAT}.$$ 
For every fibre monoidal category $\Phi(A)$, $A\in\B$, one has the fibration of modules over monoids $\Mod(\Phi(A))\to \Mon(\Phi(A))$ as in \eqref{eq:bifibmodmon}. The inverse image functors of $\Phi$ connect these fibrations with each other. Indeed, the 2-functor $$\mathcal{M}\colon\MONCAT\to\FIB_c$$ of \thref{prop:bifib} provides us directly with an indexed category over an indexed category:
$$\begin{array}{rcl}
	\mathcal{M}\circ\Phi\colon\B^{op}& \longrightarrow &\FIB_c\\
	A &\longmapsto & \xymatrix{\Mod(\Phi(A))\ar[d]\\ \Mon(\Phi(A))}\\
	A\ra{f}B & \longmapsto & \xymatrix{\Mod(\Phi(B))\ar[d]\ar[r]^{\Mod(f^*)} & \Mod(\Phi(A))\ar[d]\\
								\Mon(\Phi(B)) \ar[r]_{\Mon(f^*)} & \Mon(\Phi(A)).}
	\end{array}$$	
Recall that post-composition with the 2-functors $\dom,\cod\colon\CAT^{\two}\to\CAT$ determines two indexed categories over \B, called the domain and codomain indexed categories of $\mathcal{M}\circ\Phi$. The codomain indexed category is precisely the indexed category of monoids $Mon\circ\Phi$ in $\Phi$ defined in (\ref{eq:monindcat}). We call the domain indexed category the \emph{indexed category of modules in $\Phi$}, and the whole pseudo-functor $\M\circ\Phi$, the \emph{indexed category of modules over the indexed category of monoids in $\Phi$}. 

One might also consider the corresponding fibration over a fibration via the two-level Grothendieck construction of $\M\circ\Phi$. The codomain fibration is the (one-level) Grothendieck construction of $Mon\circ\Phi$ and is therefore the fibration of monoids $$\Mon(\Phi)\to\B$$ defined in (\ref{eq:fibmon}). We denote the whole structure by
$$
\Mod(\Phi)\to \Mon(\Phi)\to\B,
$$
and call it the \emph{fibration of modules over the fibration of monoids in $\Phi$}. Note that $\Mod(\Phi)\to \Mon(\Phi)$ is a fibrewise opfibration as its fibre at $A\in\B$ is isomorphic to $\Mod(\Phi(A))\to \Mon(\Phi(A))$.

\begin{Ex}[Sheaves of modules]\label{ex:sheavesofmodules}
	This example continues \thref{ex:bifibringedspaces}. Recall that we defined a bistrong monoidal bi-indexed category of sheaves of abelian groups $$\Top^{op}\to\ADJ(\SYMMON_s)$$ in \hyperref[eq:bifibsheaves]{(\ref*{eq:bifibsheaves})}. It has the following important property: all its categories and its left adjoint functors are in \rcSYMMON. This implies that, when post-composing it with the 2-functor $\ADJ(\M_c)$, one ends up in \BIFIBADJ. By post-composing the result with (some version of) the opposite 2-functor, one obtains therefore the following bi-indexed category over a bi-indexed category. We abbreviate $\Sh(X;\Ab)$ by $\Sh(X)$.
\begin{equation}\label{eq:Ocmodules}\begin{array}{rcl}
	\Top & \longrightarrow &\BIFIBADJ\\
	X &\longmapsto & \xymatrix{\Mod_c(\Sh(X))^{op}\ar[d]\\ \Comm(\Sh(X))^{op}}\\
	X\ra{f}Y & \longmapsto & \xymatrix@C=2cm@R=1.5cm{\Mod_c(\Sh(X))^{op}\ar[d]\ar@<-.9ex>[r]_{\Mod_c(f_*)^{op}}\ar@{}[r]|{\top} & \Mod_c(\Sh(Y))^{op}\ar@<-.9ex>[l]_{\Mod_c(f^{-1})^{op}}\ar[d]\\
								\Comm(\Sh(X))^{op}\ar@<-.9ex>[r]_{\Comm(f_*)^{op}}\ar@{}[r]|{\top} & \Comm(\Sh(Y))^{op}.\ar@<-.9ex>[l]_{\Comm(f^{-1})^{op}}}
	\end{array}\end{equation}	
Its 2-level Grothendieck op-construction, which we denote by $$\Oc{}\text{-}\Mod\to\Ringed\to\Top,$$ is therefore a bifibration over a bifibration over \Top. The objects of $\Oc{}$-$\Mod$ over the ringed space $(X,\Oc X)$ are called \emph{sheaves of $\Oc X$-modules}\indexb{Sheaf!-- of $\Oc X$-moduldes}\index{OX-module@\Oc X-module|see{Sheaf}}. We study more deeply this situation in the following \hyperref[ex:fibmodules:sheaves]{Example \ref*{ex:fibmodules}(\ref*{ex:fibmodules:sheaves})}.
\end{Ex}

\paragraph{Fibration setting} Let us now start with a fibration $P\colon\E\to\B$. In the same way as for monoids, one can define the fibration of modules over the fibration of monoids by means of the associated indexed category and the tools explained in the preceding paragraph. We give now an intrinsic definition of it.

\begin{Def}
	Let $\E\ra{P} \B$ be a monoidal fibred category.\Par
The \emph{category of modules in \E} (\emph{or in $P$}), denoted $\Mod(\E)$ (or $\Mod(P)$), is defined by:
\begin{itemize}
\item
$\Ob \Mod(\E)$: Pairs $(R,M)$ where $R$ is a monoid in $\E$ and $M$ is an $R$-module in $\E_{P(R)}$.
\item
$\Mor \Mod(\E)$: Pairs $(R,(M,\kappa)) \ra{(\phi,\alpha)} (S,(N,\sigma))$ where $\phi\colon R\to S$ is in $\Mon(\E)$ and $\alpha\colon M\to N$ is a morphism in \E such that: 
\begin{enumerate}
\item $P(\phi)=P(\alpha)$,
\item $$\xymatrix@=1.3cm{
M\otimes R\ar[d]_\kappa\ar[r]^{\alpha\otimes\phi} & N\otimes S\ar[d]^\sigma\\
M \ar[r]_\alpha & N
}$$
\end{enumerate}

\item
Composition and identities: those of $\E\times \E$.
\end{itemize}
\end{Def}

One has an analog of \thref{lem:fibremonoid}.

\begin{Lem}
	Let $\E\to \B$ be a monoidal fibred category and $f\colon A\to B$ an arrow in \B.

\begin{enumerate}[(i)]
\item
Let $(R,(M,\kappa))$ be in $\Mod(\E_B)$. Let $f^* R \ra{\bar f_R}R$ and $f^* M \ra{\bar f_M}M$ be Cartesian arrows over $f$. Then, there is a unique $f^{*}(R)$-module structure on $f^*M$ such that $(\bar f_R,\bar f_M)$ is a morphism of modules in $P$. It is defined by the left leg of the following commutative diagram.
$$\xymatrix@=1.5cm{
f^*(M) \otimes f^*(R) \ar[r]^-{\bar f_M \otimes\bar f_R} \ar@{-->}[d]& M \otimes R\ar[d]^{\kappa}\\
f^*M \ar[r]_{\bar f_M} & M}$$
\item
Let $(\phi,\alpha)\colon (R,M)\to (S,N)$ be a morphism in $\Mod(\E_B)$. Choose Cartesian lifts of $f$ at $R$, $S$, $M$ and $N$ and let $f^*(\phi)\colon f^*R\to f^*S$, $f^*(\alpha)\colon f^*M\to f^*N$ be the morphisms induced by Cartesianness.\Par Then $(f^*(\phi),f^*(\alpha))$ is a morphism in $\Mod(\E_A)$. \Par
In particular, any two inverse images of a module over the same arrow have isomorphic induced module structures.
\end{enumerate}

\qed
\end{Lem}

\begin{Rem}
Let $(R,M)$ and $(S,N)$ be modules in $\E$. A pair $$(\phi,\alpha)\colon(R,M)\to(S,N)$$ in $\E\times\E$ with $\phi$ and $\alpha$ over $f$ is a morphism of modules in \E \ssi for some Cartesian lifts $f^*(S) \ra{\bar f_S} S$ of $f$ at $S$ and $f^*(N) \ra{\bar f_N} N$ of $f$ at $N$, the pair $(\bar\phi,\bar\alpha)$ defined by
$$\xymatrix{
M \ar@{-->}[d]_{\bar\alpha} \ar[rd]^{\alpha} & \\
f^*(N) \ar[r]_{\bar f_N} & N\\
R \ar@{-->}[d]_{\bar\phi} \ar[rd]^{\phi} & \\
f^*(S) \ar[r]_{\bar f_S} & S\\
A\ar[r]_f & B
}$$
is in $\Mod(\E_A)$.	
\end{Rem}

In order to prove that the projection functor $\Mod(\E)\to \Mon(\E)$ is a fibration, by \thref{lem:intfib}, it is sufficient to prove that the projection $\Mod(\E)\to\B$ is a fibration, that $\Mod(\E)\to \Mon(\E)$ is Cartesian and that it has Cartesian lifts of all vertical morphisms of $\Mon(\E)$. Cartesian lifts of general morphisms are then given by the construction of the proof of the latter lemma. 

A Cartesian lift of $f\colon A\to B$ at a module $(S,N)$ is given by 
$$
(f^*S,f^*N)\ra{(\bar f_S,\bar f_N)}(S,N),
$$
where $\bar f_S\colon f^*S\to S$ and $\bar f_N\colon f^*N\to N$ are Cartesian lifts in $P\colon\E\to\B$ with induced monoid and module structures.

A Cartesian lift of a vertical arrow $\phi\colon R\to S$ in $\Mon(\E)$ over $A\in\B$ at a module $(S,N)$ is given by restriction of scalars in $\E_A$ along the morphism $\phi$ (\autoref{par:resscal}):
$$
(\phi,1_N)\colon(R,\phi^{\sharp}N)\to(S,N).
$$

\begin{Prop}\label{prop:fiboverfibmodules}
Let $\E\ra{P} \B$ be a monoidal fibred category.

\begin{enumerate}
\item
The projection $\Mod(\E)\to \Mon(\E)$ is a fibration, and we call \begin{equation}\label{eq:modovermon}\Mod(\E)\to \Mon(\E)\to\B\end{equation} the \emph{fibration of modules over the fibration of monoids in \E (or in $P$)}.

A Cartesian lift of a morphism $\phi\colon R)\to S$ at a module $(S,N)$ is given by the following composite of a vertical lift and an horizontal lift:
$$\xymatrix@=1.5cm{
(R,\phi^{\sharp}f^*N) \ar[d]_{(\bar\phi,1_N)} \ar[rd]^{\overline{\phi}_{(S,N)}} & \\
(f^*S,f^*N)\ar[r]_{(\bar f_S\bar f_N)} & (S,N)\\
R \ar@{-->}[d]_{\bar\phi} \ar[rd]^{\phi} & \\
f^*(S) \ar[r]_{\bar f_S} & S\\
A\ar[r]_f & B
}$$
\item\label{prop:fibmodmon:2}
The fibre fibration $\Mod(\E)_A\to \Mon(\E)_A$ at $A\in\B$ is isomorphic to the fibration $\Mod(\E_A)\to \Mon(\E_A)$ of modules over monoids in $\E_A$.
\item
Let $\Mod(\Phi_P)\to \Mon(\Phi_P)\to\B$ be the fibration of modules over the fibration of monoids of the indexed category $\Phi_P$ determined by a cleavage of $P$. Then there is an isomorphism of fibrations over fibrations over \B:
$$\xymatrix@C=.5em@R=1cm{
\Mod(\Phi_P)\ar[rr]^{\cong}\ar[d] && \Mod(\E)\ar[d]\\
\Mon(\Phi_P)\ar[rr]^{\cong} && \Mon(\E)\\
}$$
\end{enumerate}
\cqfd
\end{Prop}

\begin{Rems}
\item
	When the fibre monoidal categories of \E satisfy conditions of \thref{prop:adjointmodules}, $\Mod(\E)\to\Mon(\E)\to\B$ is a fibrewise opfibration because of point \ref{prop:fibmodmon:2} in the preceding proposition.
\item
	Let \E be a symmetric monoidal fibred category. The fibration of modules over monoids $\Mod(\E)\to\Mon(\E)$ restricts on commutative monoids to a fibration $$\Mod(\E)|_{\Comm(\E)}\to\Comm(\E)$$ and we denote $\Mod_c(\E):=\Mod(\E)|_{\Comm(\E)}$. One has therefore the \emph{fibration of modules over the fibration of commutative monoids}
        \begin{equation}
          \label{eq:3}
          \Mod_c(\E)\to\Comm(\E)\to\B.
        \end{equation}
In the indexed category setting, one can also restricts one's attention to commutative monoids by post-composing a symmetric monoidal indexed category $\Phi\colon\B^{op}\to\SYMMON$ by the 2-functor $\M_c\colon\SYMMON\to\FIB_c$ defined in \thref{prop:bifib}.
\end{Rems}

\begin{Exs}\label{ex:fibmodules}
	\item \textbf{(\textit{G}-bundles)}\label{ex:fibmodules:gbundles} This is the continuation of \thref{ex:fibmonoids}. The example of $G$-bundles in \thref{ex:fibration} are particular cases of modules in a fibration. 
	The canonical fibration is Cartesian monoidal and therefore admits a fibration of modules over its fibration of monoids $\Mod(\C^\two)\to \Mon(\C^\two)\to\C$. For an internal group $G$ of $\C$, we have defined the internal group $\theta_A^G$ of $\C/A$. It is not difficult to see that $G$-bundles over $A$ are precisely $\theta_A^G$-modules in the Cartesian monoidal category $\C/A$ if one chooses suitably the products in $\C/A$. With this choice of products in the slices, given two $G$-bundles $\xi$ and $\zeta$ over $A$ and $B$ resp., a morphism $\phi\colon\xi\to\zeta$ over $f\colon A\to B$ is a morphism of $G$-bundles as defined in \ref{ex:fibration} if and only if the pair $(\theta_f^G,\phi)\colon(\theta_A^G,\xi)\to(\theta_B^G,\zeta)$ is a morphism of modules in the canonical fibration over \C.

There is a subcategory of $\Mod(\C^\two)$ consisting of these objects and morphisms. The restriction of the fibration $\Mod(\C^2)\to\C$ to this subcategory is isomorphic to the fibration of $G$-bundles introduced in \thref{ex:fibration}. It can also be seen as the pull-back fibration of $\Mod(\C^\two)\to \Mon(\C^\two)$ along the embedding $\bar\theta^G$ of $\C$ into $\Mon(\C^\two)$. The following commutative diagram summarizes the situation, where the square is a pullback.
$$\xymatrix{
*[l]{G\text{-}Bun(\C)}\ar@{^{(}->}[r] \ar[d] & \Mod(\C^\two)\ar[d]\\
\C\ar@{=}[dr]\ar@{^{(}->}[r]^-{\bar\theta^G} & \Mon(\C^\two)\ar[d]\\
&\C
}$$
\item \textbf{(Sheaves of modules)}\label{ex:fibmodules:sheaves} This is continuation of \thref{ex:sheavesofmodules}. The fibration of \Oc{}-modules over the fibration of ringed spaces is a dual example to \thref{prop:fiboverfibmodules}, even though it is a fibration over a fibration and not an opfibration over an opfibration. This comes from the fact that we have applied the dual 2-functor at some point. Before that, the functors $f^{-1}$ were the \emph{direct} image functor of a bifibration over $\Top^{op}$. One could equivalently have first constructed the monoidal bifibration of sheaves of abelian groups $\Sh_{\Top}(\Ab)\to\Top^{op}$ over $\Top^{op}$ and then applied the construction \eqref{eq:3}: \begin{equation}\label{eq:modoverringedspaces}\Mod_{c}(\Sh_{\Top}(\Ab))\to\Comm(\Sh_{\Top}(\Ab))\to\Top^{op}.\end{equation} One finally takes the dual of all this $$\Mod_{c}(\Sh_{\Top}(\Ab))^{op}\to\Comm(\Sh_{\Top}(\Ab))^{op}\to\Top.$$ Thus, one is in fact looking at the opfibration over the opfibration over $\Top^{op}$ \hyperref[eq:modoverringedspaces]{(\ref*{eq:modoverringedspaces})}. This is why we will meet dual behaviour. In particular, it not the restriction, but the extension of scalars that plays a role.

The category \Oc{}-\Mod of sheaves of modules is obtained by the Grothendieck op-construction of the domain opindexed category of \hyperref[eq:Ocmodules]{(\ref*{eq:Ocmodules})}. Its objects are thus triples $(X,\Oc X,\F)$, where $(X,\Oc X)$ is a ringed space and $\F$ an $\Oc X$-module in the monoidal category $\Sh(X;\Ab)$. These modules can be defined as internal modules in the Cartesian monoidal category $\Sh(X;\Set)$, for the same reason that monoids in $\Sh(X;\Ab)$ can be seen as internal rings in $\Sh(X;\Set)$. A morphism $$(f,f^\sharp,\alpha)\colon(X,\Oc X,\F)\to(Y,\Oc Y,\G)$$ in \Oc{}-\Mod is given by a continuous map $f\colon X\to Y$ and a morphism $$(f^\sharp,\alpha)\colon (\Oc Y,\G)\to\Mod_c(f_*)(\Oc X,\F)$$ in $\Mod_c(\Sh(Y;\Ab))$. Equivalently, it is given by a morphism $f^{\sharp}\colon\Oc Y\to f_{*}\Oc X$ in $\Sh(Y;\Comm)$ and a morphism $\alpha\colon\G\to f_{*}\F$ in $\Sh(Y;\Set)$ such that, for all open subsets $U\subset X$, the pair $(f^{\sharp}_{U},\alpha_{U})$ is a morphism of modules. In particular, the pair $(f,f^\sharp)$ is a morphism of ringed spaces.

A Cartesian lift of a morphism $(f,f^\sharp)\colon(X,\Oc X)\to (Y,\Oc Y)$ of ringed spaces at a module $(Y,\Oc Y,\F)$ is given by the following composite of horizontal and vertical lifts.
$$\xymatrix{
(X,\Oc X,f^{-1}\F\otimes_{f^{-1}\Oc Y}\Oc X)\ar[d]_{(1_X,f_\sharp,\eta^{f_\sharp})} &\\
(X,f^{-1}\Oc Y,f^{-1}\F)\ar[r]^-{(f,\eta^f,\eta^f)} & (Y,\Oc Y,\F)\\
(X,\Oc X)\ar@{-->}[d]_{\exists!(1_X,f_\sharp)}\ar[r]^{(f,f^\sharp)} & (Y,\Oc Y)\\
(X,f^{-1}\Oc Y)\ar[ur]_{(f,\eta^f)}&\\
X\ar[r]^f & Y
}$$

In this diagram, $\eta^f$ is the unit of the adjunction $f^{-1}\dashv f_*\colon\Sh(Y;\Ab)\to\Sh(X;\Ab)$ and $\eta^{f_\sharp}$ is the unit of the adjunction of restriction and extension of scalars along $f_\sharp$ in $\Sh(X;\Ab)$. We have omitted to write \Comm and $\Mod_c$ around the functor $f^{-1}$ where needed.
\end{Exs}


\chapter{Prospects}\label{cha:Prospects}
In the introduction we raised the questions motivating this thesis. In this chapter, we first combine the results of the previous chapters to explain our progress in answering these questions. We then discuss what is still missing for a complete answer. Finally, we propose some ideas of research in $K$-Theory, using the categorical framework developed in this thesis. Note that since we finished our thesis, we could answer some of these questions. We are going to explain these progresses in an article in preparation.

The starting point of our thesis was the following question. Given an object $C$ of a category $\C$ , is there a “natural” way to associate categories to it that contain significant information about $C$ and to which the $K$-theory functor could be applied, giving rise to good notions of $K$-theory of the object $C$? This question is a good starting point, but it is in some sense naive, because too general. Asked this way, the answer is probably negative: there is not, in complete generality, a “natural” category one can associate to $C$ in order to have good notion of its $K$-theory. Nevertheless, there is a related question that clarifies the discussion: to what kinds of objects (and thus categories) should K-theory be applied?

After exploring examples such as rings (commutative or not), ringed spaces (commutative or not), topological spaces and ring spectra, we came to the following conclusion. $K$-theory is well suited for being applied to the category \emph{monoids in a monoidal (op-)fibred category} or to some category equipped with a functor into a category of such monoids\footnote{\label{fn:2}Recall that the opfibred situation requires some further conditions. In particular, the monoidal opfibred category should be opstrong (see the discussion in \autoref*{sssec:monopbifib}, at page \pageref{sssec:monopbifib}). One also needs some conditions of existence and preservation of colimits in order to have an opfibration of modules over monoids (see \thref{ex:sheavesofmodules} and \thref{prop:bifib}). The article in preparation will treat these questions in more details.}. In particular, when the monoidal fibred category is symmetric, then commutative monoids are also objects of interest.

In the cases of rings and ring spectra, the monoidal fibred category is trivially fibred, that is, it is just a monoidal category. The (symmetric) monoidal categories underlying rings and ring spectra are respectively those of abelian groups with tensor product and of $S$-modules with smash product \cite{EKMM97}. 

In the cases of spaces and ringed spaces, the fibred context becomes essential. The example of (locally) ringed spaces, and in particular of schemes, has been treated in detail in this text. There, the monoidal context is that of a monoidal bifibred category $\Sh_{\Top}(\Ab)\to\Top^{op}$ of sheaves of abelian groups (it is actually the opfibration that will be of interest later on). Moreover, this example could be easily generalized to ringed sites.

We have not examined the case of topological spaces in this thesis, due to a lack of time. The category associated to a space $X$ that is used for the computation of the $K$-theory of $X$ is the category of (real or complex) vector bundles over $X$\footnote{\label{fn:1}Traditionally, the $K$-theory of a space $X$ is defined in a representable way (via the spectrum BO or BU). Nevertheless, when $X$ is compact, this ad hoc construction can be shown to be equivalent to the $K$-theory of the exact, or symmetric monoidal, category of vector bundles, \emph{when taking the topological enrichment into account}. Indeed, by Swan's theorem, the category of real vector bundles over $X$ is equivalent to the category $\Mod_{C(X)}^{fg,p}$ of finitely generated projective $C(X)$-modules, where $C(X)$ is the ring of continuous functions $X\to\R$. Now this category has a topological enrichment, due to the fact that $C(X)$ is a topological ring. When applying to this category a version of the $K$-theory functor that takes the enrichment into account, one recovers the representable $K$-theory of $X$ \cite{Pal96,Gra94,Gra05,Gil92}. When $X$ is not compact, the category of real vector bundles is not anymore equivalent to that of finitely generated modules over $C(X)$. Yet, this becomes true if one restricts to vector bundles of finite type \cite{Vas86}, but we do not know if one still obtains representable $K$-theory as in the compact case.}. Real vector bundles over $X$ are naturally defined as internal modules in the Cartesian monoidal category $\Top/X$ over the internal ring $\theta_X^\R=(X\times \R\ra{pr_1} X)$ and similarly for the complex case. Are they modules in a monoidal category of internal abelian groups in $\Top/X$?

Let us consider the case of topological vector spaces ($X=*$). The tensor product in the category of topological abelian groups exists, but is not associative \cite{Glo04}. However, it yields a symmetric monoidal structure if one restricts for example to the category $k$-\Top of $k$-spaces, also called compactly generated spaces \cite{Sea05}. One thus has a symmetric monoidal category $\Ab(k\text{-}\Top)$ of compactly generated topological abelian groups and topological vector spaces are $\R$-modules in this monoidal category. This should generalize to the category $\Ab(k\text{-}\Top/X)$ of internal abelian groups in $k\text{-}\Top/X$ for a $k$-space $X$. By the universal property of tensor product, every vector bundle over $X$ can then be seen as a module over the monoid $\theta_X^\R$ in $\Ab(k\text{-}\Top/X)$. The monoidal framework consists thus in the symmetric monoidal fibred category of bundles of topological abelian groups \[\mathit{AbBun}\to k\text{-}\Top.\] The monoids therein can be identified with internal rings in the fibres $k\text{-}\Top/X$. In a similar way as Example \hyperref[ex:fibmodules:gbundles]{\ref*{ex:fibmodules}(\ref*{ex:fibmodules:gbundles})}), the $K$-theory of spaces is then obtained by ``transferring'' the $K$-theory of these bundles of topological rings via the pullback of the fibration $\Mod(\mathit{AbBun})\to\Mon(\mathit{AbBun})$ along the functor $\Top\to\Mon(\mathit{AbBun})$ defined by $X\mapsto \theta_X^\R$\footnote{One has to take the topological enrichment into account though, as already mentioned in \autoref{fn:1} on page \pageref{fn:1}. We come back to this point below.}.

Now that we have identified the kind of objects to which we want to apply $K$-theory, we can return to our first question: to what category over such objects should one apply the $K$-theory functor? A first part of the answer is contained in the identification above of the objects under study. Indeed, when the latter are monoids in a monoidal \linebreak(op-)fibred category, it is natural to consider their associated category of modules. We have seen that a monoidal (op-)fibred category $P\colon\E\to\B$ provides us directly with a fibration $\Mod(\E)\to\Mon(\E)$\footnote{In the opfibred context, the latter functor is an opfibration under some conditions (see \autoref{fn:2} on page \pageref{fn:2}).}. Recall that its fibre over a monoid $R$ (where thus $R$ is a monoid in the monoidal category $\E_{P(R)}$) is isomorphic to the category $\Mod_R$ in the monoidal category $\E_{P(R)}$. When the objects of study are in a category $\C$ together with a functor $F\colon\C\to\Mon(\E)$, then one can take the pullback of $\Mod(\E)\to\Mon(\E)$ along $F$.

In practice, one never applies the $K$-theory functor to the whole category of modules of a given monoid (in a monoidal (op-)fibred category). The category of modules is in general too complicated to deal with, and one therefore considers a subcategory of objects that are trivial in some sense (for instance free), or close to it. In the cases of commutative rings, of ringed spaces and of spaces, we observed that the modules considered satisfy a condition of local triviality with respect to a covering function on the base category (which is the category of monoids or a category equipped with a functor into it). We do not know, though, if finitely generated projective modules over non commutative rings can be characterized this way, neither do we for the categories of ``semi-finite'' modules over $S$-algebras mentioned in \cite{EKMM97}.

Thus, here is our general framework, that we state here only in the language of fibred categories (and not of opfibred ones) for the sake of simplicity. We start with a category \C, a monoidal fibration $P$ and a functor $F\colon\C\to\Mon(P)$. Taking the pullback of $\Mod(P)\to\Mon(P)$ along $F$, we end up with a fibration over \C that we denote $$\Mod_{\C}\to\C.$$ Moreover, this fibration is supposed to come with a structure of \emph{fibred site with trivial objects} (see \thref{def:fibsitetriv}). We then obtain a subfibration $\Loc\colon\Loc_{t}\to\Loc_{b}$ of locally trivial objects. We now define the $K$-theory of a locally trivial object $B\in\Loc_{b}$ in the base by applying the $K$-theory functor the fibre $(\Loc_{t})_{B}$ over $B$ in the fibration \Loc.

At this point, we would like to stress that we do not intend to be dogmatic on the framework of $K$-theory. Indeed, there might be relevant examples that do not fit in it. For instance, one might run into objects with meaningful $K$-theory whose associated fibration is not coming from a fibration of modules and monoids, but where the relevant fibre over an object is still that of locally trivial objects in the fibration.

Now, is our framework ready for $K$-theory? Not quite yet. We explain here the few more steps that have to be taken before applying $K$-theory functor.

First of all, the category to which the $K$-theory functor is applied must be essentially small. In addition, examples such as the $K$-theory of spaces and of rings show that one must only consider modules that are ``finitely generated'' or ``finitely presented'' in some sense, in order to avoid a trivial $K$-theory (this is due to the so called ``Eilenberg swindle'' \cite{WeiKbook}). In the examples of spaces and rings, the restriction to vector bundles whose fibres are finite dimensional vector spaces and to finitely generated modules solve both the size and the swindle problems. We will see in the article in preparation that these results can be generalized, under some conditions, when the fibred site is a stack. There, locally trivial objects inherits the desired smallness properties of the trivial ones.

We have already mentioned that the $K$-theory functor has its domain in Quillen-exact, or more generally Waldhausen categories%
\footnote{There is also a $K$-theory functor that is defined on symmetric monoidal categories. However, this functor seems to us to appear only accidentally in the questions we study. For instance, one can use this functor in order to calculate the $K$-theory of a ring $R$ because the category of finitely generated projective $R$-modules is split exact. However, in all cases where the category of modules is not split exact, it is the $K$-theory of the Quillen-exact category that is chosen, not the $K$-theory of the Cartesian monoidal category.}.
There is also a $K$-theory functor that takes into account a possible enrichment in the Cartesian closed category of $k$-spaces. The framework we have constructed up to now does not provide such a structure, however. We must therefore impose more structure on the monoidal fibred category $\E\to\B$ with which we work if one needs this form of $K$-theory.

Let us first attempt to obtain a Quillen-exact structure on the category of locally trivial modules over some monoid. Suppose that the fibre $\E_B$ is abelian and that its functors $-\otimes E$ are additive and preserve cokernels, for each $B\in\B$. The category of modules over a monoid $(B,R)\in\Mon(\E)$ is then abelian \cite{Ard04}, so that the question we must answer is the following.

Let $(P,K,\Triv)$ be a fibred site with trivial objects. Suppose each fibre of $P$ is an abelian category. Suppose that the trivial objects form a Quillen-exact subcategory in each fibre and the inverse image functors over arrows of $K$-coverings are exact functors. Then, under some natural conditions, locally trivial objects form a Quillen-exact subcategory in each fibre. Note that it is important not to require that all inverse image functors are exacts (consider the opfibration of commutative rings: the direct image functor along a localization is exact, but this is not true for a general morphism of rings).

One could also seek a Waldhausen structure on the category of locally trivial modules over some monoid. In this case, one should start with a monoidal fibration whose fibres $\E_B$ are monoidal model categories. Under some natural conditions, the category of modules over a monoid is then also model \cite{SS00}. One should thus consider the situation of a fibration with trivial objects and a covering function on its base whose fibres are model. Do the locally trivial objects then form a Waldhausen subcategory of each fibre? We have not yet worked out this question. Since each fibre is supposed to be a model category, one should probably consider a weakened notion of locally trivial objects. Locally trivial objects could be locally over a covering weakly equivalent to a trivial one. When the base is also a model category, it would be interesting to consider also weakened versions of the axioms of pretopologies, as defined in \cite{TV05} for instance.

Finally, the structure of a \emph{topological} Quillen-exact or Waldhausen category plays a role in topological $K$-theory (of spaces or Banach algebras, e.g.). Indeed, recall that if $X$ is a compact space, then the category of real vector bundles over $X$ is equivalent to the category $\Mod_{C(X)}^{fg,p}$ of finitely generated projective $C(X)$-modules, where $C(X)$ is the ring of continuous functions $X\to\R$. So one can as well calculate the $K$-theory of $X$ from $\Mod_{C(X)}^{fg,p}$. It is the $K$-theory functor that takes in account the topological enrichment of $\Mod_{C(X)}^{fg,p}$ (due to the fact that $C(X)$ is a topological ring) that gives the $K$-theory of $X$, not the $K$-functor that only takes the Quillen-exact structure into account, which is the algebraic $K$-theory of the discrete ring $C(X)$ \cite{Pal96,Gra94,Gra05}. Our theory should therefore have topologically enriched version in order to capture these examples.

\paragraph{}

We conclude with some possible applications of our theory. It turned out to be quite a job to describe these categorical foundations for $K$-theory and, as we have explained, there are still some gaps to fill. Along the way, we have encountered very rich categorical structures that are interesting for their own sake. Moreover, examples of these categorical notions arise in various parts of mathematics, as we have shown for locally trivial objects, for instance. Nevertheless, our initial and main motivation was $K$-theory and we therefore have hopes of applications of our theory into this field. 

We have developed a unified, conceptual framework for a number of important examples of objects studied in $K$-theory, thus providing new tools for obtaining $K$-theoretic information. One can now modify the parameters of the well-known examples, with all the control that the abstract theory gives to us. For example, within our general framework, it is now possible to ask the following interesting question: how does the $K$-theory vary under changes of covering functions, or of trivial objects? In general, it would be also interesting to study functors in $K$-theory induced by morphisms of fibrations. On the other hand, one can try to apply $K$-theory to new situations by means of this framework. We will address some of these questions in the cited article in preparation, in the framework of Quillen-exact categories.


\begin{appendix}
\chapter{A note on foundations}\label{cha:Found}
This appendix is meant to shortly recall the set theoretical questions that arise in category theory and what options have been proposed to deal with them. We then turn to a little more thorough treatment of the option we have chosen for this thesis. A complete treatment of the question would of course require much more space and time. This appendix will therefore be quite sketchy and not self-contained. In particular, it requires some basic knowledge in set theory and logic: \emph{first order language \textnormal{and} theories}, \emph{extension by definitions} of a theory, \emph{interpretation of a theory} in another, \emph{interpretation of formulae \textnormal{and} of definitions} of a theory in an interpretation of the theory. One also supposes that the reader is accustomed with the usual set theories \ZFC and \NBG. We learnt these notions mainly in the following books that we recommend. The reference \cite{Sho67} is a very rigorous and clear exposition of first order theory notions and takes seriously the question of what structure the meta-language is supposed to have. In particular, when studying models of set theories, the author uses the internal notion of an interpretation, not the external notion of a model. See also \cite{Kun80}, that has the same spirit, but with fewer details. We used \cite{Men97} especially for its full treatment of \NBG. The existence of this book is very precious for category theory because it is the only modern exposition of \NBG (as far as I know). Finally, \cite{Hin05} is a recent reference that, after having presented the usual notions of first order logic, restarts doing it \emph{inside} the formal theory \ZFC! To recommend to anyone who once found bizarre that when doing formal set theory, the meta-language seems to already contains all the notions that one is going to define … \cite{Kun80} has also clarifying discussions of these questions.

\section{Foundational issues of category theory}
Surveys of the subject can be found in \cite{HS73,Shu08,Hat82}.

It is well known since the very beginning of the field with Eilenberg and Mac Lane \cite{EM45} that category theory deals with set theoretical concepts that can lead to paradoxes if not treated with care. Actually it is hardly surprising that category theory asks foundational questions since it is designed to describe the structure of mathematical theories themselves and so has a foundational nature. Let me recall some of the difficulties that arise.

Category theory is often applied to a whole class of models of a theory, for example the class of all groups or all topological spaces. Yet these ``collections'' are precisely those for which one needs axiomatization in order to avoid Russell-like paradoxes. If one chooses \ZFC as the set theory, then its language states facts about only sets. A class is a mere formula, that is, an object of the meta-language. It is a formula $\phi(x)$ of the language \ZFC with at most one free variable%
\footnote{One could also consider formulae $\phi(x_1,\ldots,x_n,x)$ with more free variables. It is ``collectivizing'' if $\vdash\exists y\forall x(x\in y\leftrightarrow \phi(x_1,\ldots,x_n,x))$. But this yields the definition of a \emph{function symbol}, whereas a formula with only one free variable leads to the definition of a \emph{constant symbol}. For example, the formula $x\in x_1\vee x\in x_2$ gives rise to the function symbol $\cup$.}. 
A class may be ``collectivizing'', that is there exists a set of all $x$ such that $\phi(x)$. This is the case when \ZFC proves the following formula: 
\begin{equation}\label{eq:setphi}
\exists y\forall x(x\in y\leftrightarrow \phi(x)).
\end{equation}
For example, it is an axiom of \ZFC that the formula $x\neq x$ is collectivizing in the sense that the former formula for this $\phi$ is an axiom (this is $\emptyset$). But it is a theorem of \ZFC that the formula $x=x$ is not (the class $V$ of all sets). One says then that it is a proper class.
It is possible to speak about classes in \ZFC by identifying one class with a formula that defines it, but this is a drift in the meta-language since a formula is a word of the language, not something the language talks about. For instance, \ZFC cannot quantify over classes. We will come back to this later.
\paragraph{}
People in the field actually prefer to adopt an axiomatic set theory whose language speaks about classes and sets together, like \NBG\index[not]{NBG@\NBG} (Von Neumann-Bernays-G\"odel) or $MK$ (Morse-Kelley). In these theories, sets are particular classes, those that are elements of other classes. \emph{Proper classes} are the classes that are not sets. It is a theorem of \NBG that the classes of all sets, all groups etc. are proper classes. In such axiomatizations of the universe of sets, one can write formulae and state theorems about classes. Practically, this means that it makes sense to quantify over classes. Moreover, there is an \emph{interpretation} (or \emph{relativization}) of \ZFC in \NBG. That is, one can define in \NBG an unary predicate symbol $Set(X)$ by $\exists Y(X\in Y)$ such that the interpretation of the non logical axioms of \ZFC in \NBG via the predicate symbol $Set$ are theorems of \NBG \cite{Men97}. This implies the meta-theorem that if \NBG is consistent, then \ZFC is. In fact, the reciprocal is true [ibid.] and so \NBG is not ``riskier'' than \ZFC.

Some natural constructions of category theory are still not available though. Indeed, no notion of ``collection'' of proper classes is defined in the language since the only variables of the languages are classes, and proper classes belong to no class by definition. As within \ZFC, one could define such collections by means of formulae that define them. But again this requires a jump in the meta-language. Nevertheless, in category theory, there are very natural situations where one would like to talk about such ``collections'' of proper classes.

Let us, for now, adopt \NBG for set theory. A category \C is a class of objects $Ob\C$, a class of morphisms $Mor\C$, two functions \emph{domain} and \emph{codomain}, and a \emph{composition} relation. These data satisfy the usual axioms of category theory. Given two categories $\C$ and $\D$, a functor from \C to \D is a particular kind of function $F\colon Mor\C\to Mor\D$, i.e. a subclass $F\subset Mor\C\times Mor\D$. When \D is not empty, $F$ is a proper class if and only if $Mor\C$ is\footnote{Use the fact that if $f\colon X\to Y$ is a function between two classes, then $X\cong f$, unless $Y$ is empty (recall that a function is a subclass of the Cartesian product $X\times Y$). Then the replacement axiom (see further down \thref{axiom:replacement}) implies that one is a set if and only if the other is \cite{Men97}.}. In the same manner, when \D is not empty, a natural transformation $\alpha\colon F\Rightarrow G\colon\C\to\D$ is a proper class if and only if $Ob\C$ is. It is easy to check that there is a composition of these natural transformations that verifies all the axioms of composition of arrows in a category. But for \C a proper class and \D not empty, \NBG proves there is no class containing any of these functors or natural transformations since this would imply them to be sets. So they cannot gather in a category. But, of course, functors between large categories and natural transformations between them satisfy the axioms 2-categories.

For the same reason, in \NBG, there is no (2-)category $CAT$ of all (possibly large) categories. The previous problem is in fact just a part of this one: the hom-categories of $CAT$ are not defined in general. As a consequence, it doesn't make sense for example to talk about pseudo-functors from a category \C to $CAT$. However, the correspondence which sends a ring $R$ to its category of modules $M_{R}$ (and extension of scalars for morphisms) is readily checked to fulfil all axioms of a pseudo-functor.  In all these examples, one sees that ``it is a 2-category'' or ``it is a pseudo-functor'' but one cannot state it for set theoretical reasons. 
\paragraph{}
A first answer to this problem is to simply forget all about size problems. This isn't a real option though, because of the risk of paradoxes and because of the importance of size hypotheses in fundamental theorems of category theory. Herrlich and Strecker give a nice illustration of such a paradox \cite{HS73}: if one considers that categories form a category $CAT$, then is its full subcategory of categories that are not objects of themselves an object of itself? On the other hand, Shulman gives several examples in \cite{Shu08} of category theory statements where size completely changes the result. For instance, it is a very mild assumption for a category to have all products. But, if a category has all products indexed by proper classes, then it is necessarily a preoder!

A second way of handling these issues is to consider that all statements of the kind ``categories and functors between them form a category $CAT$'' or ``$\C\to CAT$ is a pseudo-functor'' really are just an abbreviation for the longer statement that lists the conditions that are involved. Yet, this is not satisfactory, firstly because it is dangerous: one is tempted to forget that it is an abbreviation and then one may use set-related results of category theory that are not available there. Secondly, it is of great help to apply category theory to these ``super-large'' collections such as $CAT$. For instance it quite simplifies the matter to know that ``$CAT$ is a 2-category'' when studying categories or to have the Yoneda embedding $\C\to Set^{\C^{op}}$ for a locally small category $\C$. Having a ``category'' $CAT$ also permits to precisely express the fact that Grothendieck fibrations on \C are equivalent to pseudo-functors from $\C^{op}$ to $CAT$.
\paragraph{}
There are different ways to take size problems seriously. One has at first a conceptual choice between two options to make. The first option is radical: get rid of set theory and found mathematics on an axiomatization of a certain fundamental categorical notion, using arrow-theoretic concepts. One can axiomatize the notion of a category: there is a first-order theory of categories, like there is one of groups \cite{Hat82,Law66}. But this is not sufficient. Indeed, one just obtains general theorems concerning objects and arrows. This theory says nothing about particular categories like the category of sets or of groups\footnote{Unless one adds axioms that force the theory to resemble the category one has in mind. This is successful for the category of sets, with the theory ETCS.}. Moreover, the notions of functor and natural transformations are absent. These notions would only appear as morphisms, and morphisms of morphisms, of \emph{models} of the theory of categories in a theory of \emph{sets}. To solve this problem, Lawvere proposed an axiomatization of the \emph{category of categories} \cite{Law66}. Sets will there just pop up as discrete categories. This is an interesting option, but not fully developed until now (see also \cite{McL91}). B\'enabou also provided such an axiomatization, based on \emph{fibred categories} \cite{Ben85}. Finally, Lawvere has developed an intermediate step with his Elementary Theory of the Category of Sets (ETCS). It is an axiomatization of sets via their category, expressing the fact that it is a particular topos \cite{Law03}. A fully satisfactory axiomatization using categories should in fact encompass categories, functors and natural transformations, but also bicategories, lax functors, lax natural transformations and modifications, tricategories and so on. So, ideally, one should provide an axiomatization of the $\infty$-category of $\infty$-categories.
\paragraph{}
The second option takes both size problems and set-foundation of mathematics seriously. It forces one therefore to modify the axiomatization of set theory, whether it is \ZFC or \NBG. Remember that \NBG was introduced in order to include classes in the language. This process can actually be continued. One can define a first-order theory with three different types of ``collections'': \emph{sets}, \emph{classes} and \emph{conglomerates}. The variables of the theory are meant to describe conglomerates, whereas sets and classes are particular conglomerates. Such an axiom system is sketched in \cite{HS73} and studied in more detail in \cite{Osb00} (with some differences). Category theorists often informally take this solution, but technically derive it from an ``universe axiom''. The reason why this three-level set theory is not widely used in category theory is probably that set theorists haven't studied it. Actually, even \NBG seems to be almost anecdotal to set theorists.
\paragraph{}
There is theoretically a way of dealing within \NBG, or even \ZFC, with ``collections of classes'', but it is both unmanageable and insufficient as we briefly explain now. One can define by formulae of the language of \ZFC or \NBG what are a first-order language $l$, an atomic terms of $l$, $l$-formulae, and a structure of $l$. One can then define what is an $l$-theory $t$, a formal proof and a model of a theory. Moreover, given a theory $T$ in the usual sense, there is a way to express it as a theory $\dot T$ defined in \ZFC \cite{Hin05}. Finally, the formula $Con(t)$ is the formula asserting that there exists a model for $t$\footnote{The theory $t$ must not be finitely axiomatizable. In particular, it makes sense to talk about the theory $\dot\ZFC$ in \ZFC and there is a formula of \ZFC saying that a set is a model of $\dot\ZFC$.}. 

Recall now that one part of the set theoretical issues of category theory is the impossibility to gather proper classes, like large categories. Yet, there is a formula of \NBG that expresses the fact that a class is a category (see further down \autoref{ssec:CatNBG}). One could thus consider the ``collection'' of all categories as this formula, but we remarked that this is a drift in the meta-language. However, as one can formalize the meta-language, for instance by \ZFC, this formula is now a set of \ZFC. This shows that it is \emph{in principe} correct to deal with ``super-classes'' as formulae and to apply set-theoretical tools on them. Nevertheless, there are two problems with this solution. Firstly, really working in $\dot{\NBG}$ would be a nightmare full of sequences of natural numbers! Secondly, the ``super-classes'' are indeed sets of \ZFC, since a formula of $\dot\NBG$ is coded as a sequence of natural numbers in \ZFC. Yet, they have no relationship with the sets and classes of $\dot\NBG$, because $\dot\NBG$ is a theory that \emph{talks} about sets and classes, whereas its formulae, that are sets that the ambient \ZFC is talking about, are just a coding of the language of $\dot\NBG$. It is as if I talk to someone about my wife's wonderful cakes, and that at his next visit we invite him to eat the recipe. Well, jokes aside, the problem with this incompatibility is that it is too restrictive, because in category theory, one does not just want to talk about a ``super-large'' category like \CAT, one really wants to use it, like for an indexed category $\B^{op}\to\CAT$.
\paragraph{}
The previous impracticable solution leads however to a useful option, but which requires a strengthening of set theory. Indeed, as already mentioned, one can define in \ZFC what is a model of a theory $l$ itself defined in \ZFC. Let us consider a model of $\dot\NBG$%
\footnote{The notion of a model of $\dot\NBG$ in \ZFC is the formalization of the notion of a model of \NBG in the unformalized meta-theory, the latter being an intuitive set $M$ together with a relation $\in^M\subset M\times M$ such that the axioms of \NBG are verified for this notion of elementhood (see \cite{Hin05} for the formal notion). There is also the notion of a model of \NBG (without the dot, that is, the theory \NBG defined in the meta-language) in \ZFC. It is an interpretation of \NBG in \ZFC whose universe $M$ (a predicate symbol) is a set (in other words, the the formula \hyperref[eq:setphi]{(\ref*{eq:setphi})} with $\phi=M$ is a theorem of \ZFC). It is a beautiful result that the notions of models of $\dot\NBG$ and of \NBG in \ZFC coincide. This applies in more general situations, see \cite{Hin05,Kun80}.}. 
In the model, the ``virtual'' sets and classes that $\dot\NBG$ talks about are ``realized'' as plain sets of \ZFC. In this manner, the meta-language where one drifts when talking about collections of classes is itself not only axiomatized but also in direct relationship with the language. Moreover, the meta-language is the usual language of set theory. 

The proof of the existence of such a model is of course too strong a request; if it were possible, then \ZFC would be inconsistent by Gödel's Second Incompleteness Theorem. One should thus try to add an axiom to \ZFC such that there exists a model $(M,\in^M)$ of \NBG (for example, the axiom ``there exists a model of \NBG''!). If one had such a model, then one could do all ordinary mathematics in the model, and when there is a need to gather proper classes of the model, this can be done in the surrounding \ZFC. 

One would like some compatibility between the model of \NBG and the theory \ZFC where it lives. As we already said, it is not just the question of studying $M$ for its own sake from outside. One wants to consider for example functions from a class of the model to a set outside of the model (remember the pseudo-functors $\B^{op}\to CAT$). The basic conditions one requires on the model $(M,\in^M)$ are that of a \emph{supertransitive standard model} (or \emph{supercomplete model}) \cite{BZ03,Sh51,Sh52,Sho67}. They express formally the fact that the relation $\in^M$ is the restriction of $\in$ to $M$, and that $M$ is closed under elements and subsets:
\begin{equation}\label{eq:supercomplete}
\begin{array}{lll}
(i) &\in^M=\in & (\text{standard model}),\medskip\\
(ii) & \vdash y\in x \wedge x\in M \to y\in M & (\text{transitivity}),\medskip\\
(iii) &\vdash y\subset x \wedge x\in M \to y\in M & (\text{supertransitivity}).
\end{array}
\end{equation}
A model $M$ of \NBG in \ZFC (or \NBG) is supercomplete if and only if $M=V_{\theta+1}$ for a (strongly) inaccessible cardinal $\theta$ \cite{Sh51,Sh52}, where, given an ordinal $\alpha$, $V_{\alpha}$ is the set of stage $\alpha$ in the von Neumann hierarchy of sets%
\footnote{These sets are defined by \emph{transfinite recursion} by $V_0=0$, $V_{\alpha+1}=\power{V_{\alpha}}$ and $V_{\alpha}=\bigcup_{\beta<\alpha}V_{\beta}$ for a limit ordinal $\alpha$. Recall that $0=\emptyset$, $\alpha+1=\alpha\cup\{\alpha\}$ is the \emph{successor ordinal} of $\alpha$ and $\alpha$ is \emph{limit} if there is no ordinal such that $\alpha$ is its successor (like $\emptyset$ or the ordinal $\omega$ of all natural numbers).}. 
Therefore, the sets of the model are the elements of $V_{\theta}$ and the classes are its subsets. Moreover, a set $\mathcal U$ satisfies the axioms of a \emph{universe} in the sense of Grothendieck if and only if $\mathcal U=V_{\theta}$ for a (strongly) inaccessible cardinal $\theta$\footnote{See \cite[Exposé I, Appendice]{SGA4}. For models of \ZF, this characterization is not true. All supercomplete models $M$ of \ZF in \ZF are of the form $M=V_{\alpha}$ for some ordinal $\alpha$ (which one calls \emph{natural models}). However, $\alpha$ need not to be inaccessible \cite{MV59}. A characterization of the possible ordinals $\alpha$ and in terms of universes is given in \cite{BZ03}.}.

We must notice that all axioms don't have the same importance with respect to consistency. For example, the Axiom of Choice $AC$ is independent of the theory $ZF$ in the sense that \ZF proves neither $AC$, nor its negative $\neg AC$. This implies that if \ZF is consistent, then both \ZFC and $ZF+\neg AC$ are consistent. Therefore, with respect to consistency, there is no change in adding $AC$ (or its negative) to the axioms of the theory. Adding to \ZFC the axiom $I$ of the existence of an inaccessible cardinal is of different nature (let us write $ZFI$ for $ZFC+I$). Indeed, one can show (in \ZFC) that, if \ZFC is consistent, then $ZFC+\neg I$ is consistent, or equivalently, that \ZFC does not prove $I$. On the contrary, one cannot prove (neither in \ZFC, nor in \ZFI) that the consistency of \ZFC implies the consistency of \ZFI. This is because \ZFI proves the consistency of \ZFC, and therefore, would prove its own consistency, violating the second incompleteness theorem of Gödel. Thus, adding the axiom $I$ really adds a risk. The axiom of infinity has the same status (with regards to consistency) as the axiom $I$ in the theory \ZF.

Before going to the concrete solution we adopt for this work, I should mention that there is still another approach to these set-theoretical problems, \emph{Feferman set theory}. See \cite{Shu08,Shu09} and references therein. 

\section{The strategy for this work}

We now turn to the solution adopted in this text. Our approach is similar to the usual methods proposed in the literature, since we also rely upon some ``Grothendieck universe'' trick. However, it provides both conceptual and technical new ingredients, as we explain now. 

\subsection{General approach}
We quite often use a three-level set theory, and we need once in while a fourth level. In this situation, people are used to adding an axiom of existence of \emph{two} inaccessible cardinals (e.g., \cite{DHKS04}) or even with a whole class of them (Grothendieck adds the axiom that every set belongs to an inaccessible cardinal \cite{SGA4}). We prefer to work in the theory $NBG$ together with the axiom $I$ of existence of an inaccessible cardinal, which we denote $NBG_I$\index[not]{NBGI@\NBI} (we haven't seen this axiomatic choice elsewhere). This theory allows indeed four levels of ``sets'' and is, in our point of view, much simpler. In particular, the two first levels respect the same axioms as the third and fourth (the latter being sets and classes in \NBG and the former, sets and classes in a model of it). Moreover, it is a weaker assumption than the existence of two inaccessible cardinals in \ZFC, since \NBG is equiconsistent with \ZFC, and \ZFI is strictly stronger than \NBG.

Authors also tend to make some inaccessible cardinal assumption at the beginning and then work all the time with the different levels it provides. One can then follow two different strategies. One either works at the two first levels when possible in order to make set theory the least visible, going to further levels only when one is forced to (e.g., \cite{HS73}). Or one works at the highest level possible, in order to attain the greatest generality (e.g., \cite{DHKS04}). We have in this text another point of view from these authors, but close to the former. We want to bother the reader the least possible with set theory intricacies, but remain correct and honest on this subject. That means we don't want to completely avoid the problem, neither to pretend to work with small categories when our categories of interest are obviously large. We also want to add the least axioms possible to \NBG, as it was noted that adding an inaccessible cardinal yields a strictly stronger set theory. So, unless it is necessary, we work in plain \NBG, with two kinds of categories: small categories and large categories. When needed, we turn to \NBI and then have four kinds of categories: small, large, extra large and extra extra large. All theorems proved in \NBG are also true in \NBI and in any model of \NBG (so in particular in a Grothendieck universe). This method has three advantages: it makes set theory the least visible possible, it makes the least additional set theory assumption and it has no loss of generality in comparison with the paradigm ``assume \NBI throughout the text + work at the highest level possible''. 

Now, we state precisely our set theory assumptions and define related basic category concepts.

\subsection{Some important facts about \textit{NBG}}\label{ssec:NBG}\index[not]{NBG@\NBG}

We recall some crucial definitions, axioms and propositions of \NBG that we use in the sequel. \emph{The} modern reference for \NBG set theory is the already cited (and very well written) book of Mendelson \cite{Men97} and we refer to his axiom system. Recall that in \NBG the variables represents objects that we call \emph{classes}\indexb{Class} and that a class is a \emph{set}\indexb{Set} \ssi it belongs to a class. A class is \emph{proper}\indexb{Class!Proper --} if it is not a set. Two classes are equal if and only if they have the same elements. In this part, we use capital letters for class variables and lowercase letters for set variables. One defines an unary predicate symbol $Set$ by $Set(X)\iff\exists Y(X\in Y)$, saying ``$X$ is a set''.

Let us first recall an important theorem of \NBG, also called \emph{Class Comprehension Schema}. 
\begin{Thm}[Class Existence Theorem]
Let $\phi(X_1,\ldots,X_n,Y_1,\ldots,Y_m)$ be a formula in \NBG whose variables occur among $X_1$, $\ldots,$ $X_n,$ $Y_1,$ $\ldots,$ $Y_m$ and in which only set variables are quantified. Such a formula is called \emph{predicative}. Then,
$$
\vdash_{NBG} \exists Z\ \forall x_1\ldots\forall x_n(\langle x_1,\ldots,x_n\rangle\in Z\iff\phi(x_1,\ldots,x_n,Y_1,\ldots,Y_m)).
$$	
\end{Thm}
In other words, \NBG proves that there exists a class of all $n$-tuples $\langle x_1,\ldots,x_n\rangle$ such that $$\phi(x_1,\ldots,x_n,Y_1,\ldots,Y_m).$$ The definition of equality of classes implies furthermore that this class is unique. One can, for instance, define the Cartesian product of two classes by means of this scheme. Another important feature is the \emph{power class} $\power X$ of a class $X$\indexb{Class!Power --}. It is defined, via the Class Existence Theorem, by the formula 
$(\forall x)(x\in\power X\iff x\subset X)$,
that is, $\power X$ is the class of all \emph{subsets} of $X$. It is then an axiom, the \emph{Power set axiom}, that the power class of a set is as set. 

From now on, we use informal language, but it should always be easy to return to the formal language. In particular, we denote $n$-tuples by usual parentheses. We sometimes return to formal language, when it seems useful for the purpose.

Another important feature of \NBG is its \emph{axiom of replacement}. We need some preliminary definitions. The \emph{domain} $\mathcal D(X)$ of a class $X$ is the class of all sets $v$ such that there exists a set $w$ with $(v,w)$ belonging to $X$. One defines similarly the \emph{image} $\mathcal I(X)$ of $X$. A class $X$ is \emph{univocal} if for all pairs $(x,y)$ and $(x,z)$ belonging to $X$, $y=z$. A class $R$ is a \emph{relation} if all its elements are pairs.  A class $F$ is a \emph{function from a class $X$ to a class $Y$} if it is a univocal relation whose domain is $X$ and image is a subclass of $Y$.

\begin{Axiom}[Axiom of Replacement]\label{axiom:replacement}
Let $Y$ be a univocal class and $x$ a set. Then the class of second components of pairs of $Y$ whose first components are in $x$ is a set.
\end{Axiom}
In particular, when $F$ is a function from $X$ to $Y$, then the image of the restriction of $F$ to a set is a set. Here are some useful consequences of the axiom of replacement.

\begin{Prop}\label{prop:NBG}
	\begin{enumerate}[(i)]
\item
Any subclass of a set is a set. Moreover, the domain and the image of a set are sets.
\item
The Cartesian product of two sets is a set.
\item\label{prop:NBG:functionset}
A relation is a set \ssi its domain and its image are sets. In particular, a function is a set \ssi its domain is a set.
\item\label{prop:NBG:functionclass}
There exists a class of all functions from a set $x$ to a class $Y$, denoted $Y^x$. Moreover, if $y$ is a set, then $y^x$ is a set. If $Y$ is a proper class and $x\neq\emptyset$, then $Y^x$ is a proper class.
\item\label{prop:NBG:InjClass}
If $X$ injects in $Y$, then if $X$ is a proper class, so is $Y$ (in other words, if $Y$ is a set, so is $X$).
\end{enumerate}
\end{Prop}

\subsection{Category theory in \textit{NBG}}\label{ssec:CatNBG}\index{Category}
Unless otherwise indicated by the vocabulary used or by an explicit assumption, we work in this thesis in the set theory \NBG with the axiom of choice for classes. We discuss now the definition of a \emph{category}. Notice that, since a proper class $X$ cannot be an element of a class, it cannot be a component of an ordered pair $(X,Y)$\footnote{At least in a non-trivial meaning. In fact, one defines $(X,Y)$ in \NBG by $\{\{X\},\{X,Y\}\}$, where $\{X,Y\}$ is defined, if both $X$ and $Y$ are sets, to be the set that contains both $X$ and $Y$ and, if not, to be $\emptyset$.\label{fnote:pairofclasses}}. Therefore, one cannot state the definition ``a category is a quintuple $(Ob\C,Mor\C,\circ,dom,cod)$ such that \ldots'' or ``a category is a pair $(Mor\C,\circ)$ such that \ldots''\footnote{We actually sometimes use this kind of sentences for large categories, like in ``a monoidal category is a sextuple\ldots'', but this should be read an informal sextuple, that is a sextuple in the meta-language.}. 

It can be defined by a sentence in the meta-language of the kind ``it consists of the following data''. More precisely, there is a first-order theory of categories \cite{Law66,Hat82}. Its signature consists of one ternary predicate symbol $C$ ($C(x,y,z)$ is interpreted as $z$ is the composite of $x$ with $y$) and two unary function symbols $dom$ and $cod$ (interpreted as domain and codomain). A category \C is an interpretation of the first-order theory of categories in \NBG whose universe is a class\footnote{In a general interpretation, the universe is just a formula $U(X)$ with at most one free variable $X$. This formula is not necessarily predicative. Moreover, it is not restrictive to consider only classes that can be defined by a formula, because we won't be able to collect large categories in an extra large category in this context. We can just talk, in the meta-language, about some specific categories we have defined.}. In other words, it consists of a class $Mor\C$ (determined by a formula with at most one free variable), a ternary relation $\circ$ on $Mor\C$, two functions $dom$ and $cod$ from $Mor\C$ to $Mor\C$, such that the interpretation of the axioms of category theory is verified. For example, there is a category of sets, written $Set$, whose class of objects, denoted $V$ (the \emph{universal class}), is proper. It is defined by the class of all $x$ such that $x\notin\emptyset$.

One can also define a category \C as the data of a (defined) class $Mor\C$ and a ternary relation $\circ$ on $Mor\C$ satisfying axioms. Finally a third way is to defined a category as the data of a (defined) class $Mor\C$, an unary relation $Ob\C$ on $Mor\C$, a ternary relation $\circ$ and two unary functions $dom$ and $cod$ on $Mor\C$. The two first definitions are equivalent, whereas the last one is ``essentially equivalent'', meaning that if the data ``$(Ob\C,Mor\C,\circ,dom,cod)$'' satisfies the latter definition, then the data ``$(Mor\C,\circ)$'' satisfies the former, and there is a unique such ``$(X,Mor\C,\circ,Y,Z)$'' up to isomorphism \cite{HS73}.

There also exists a formula $CAT(X)$ in the language of \NBG saying ``$X$ is a category'', but it is not very intuitive. Indeed, in order to avoid pairs of classes, one has to define a category to be its ternary relation of composition $C$. All morphisms appear as the composite of themselves with an identity. So ``$X$ is a category'' states for ``$X$ is a univocal (binary) relation such that $\mathcal D(X)\subset\mathcal I(X)\times \mathcal{I}(X)$ and verifying the usual finite set of axioms for the arrow-like definition of categories with $Mor\C=\mathcal{I}(X)$''. Then, given a ternary relation $X$ on a class $Y$, it is a category in the sense of the preceding paragraph \ssi $\vdash CAT(\bar X)$, where $\bar X$ is the binary relation corresponding to $X$ with domain a subset of $Y\times Y$.

A category is said \emph{large}\indexb{Category!Large --} if it is not small, that is, if its class of morphisms is proper (this is the case in particular if its class of objects is proper, by Proposition \hyperref[prop:NBG:InjClass]{\ref*{prop:NBG}(\ref*{prop:NBG:InjClass})}). We do not assume in general that the class of morphisms with fixed domain and codomain is a set. If needed, we specify this assumption by calling such a category \emph{locally small}\indexb{Category!Locally small --}.

Given two categories $\C$ and $\D$, there is no problem in defining functors from $\C$ to $\D$ in $NBG$. It is a function from $Mor\C$ to $Mor\D$ satisfying axioms. Similarly, for natural transformations.

We won't go into details for bicategories, but we want to remark that there is again an injection of objects into morphisms via identities of objects, and an injection of morphisms into 2-cells via identities of morphisms. 

\subsubsection{The 2-category \textit{Cat} of small categories}\label{sssec:Cat}

When restricting to sets, it is a much simpler matter to define categories. These categories are called \emph{small categories}\indexb{Category!Small --}. A formula $Cat(X)$ in \NBG saying ``$X$ is a small category'' can be written in the form ``$X$ is a quintuple $(Ob\C,Mor\C,\circ,dom,cod)$ satisfying axioms'' (indeed, the theory of categories is finitely axiomatizable). If $X$ is a proper class, then $Cat(X)$ is false, because a quintuple is always a set (see the convention when a component of a tuple is a proper class in \autoref{fnote:pairofclasses} of this Appendix). Notice that, by \thref{prop:NBG}, if $\C$ is a category and $Mor\C$ is a set, then the composition relation, the codomain and domain functions and the object class are also sets. Thus, the formula $Cat(X)$ is ``equivalent'' to the formula $CAT(X)\wedge Set(X)$ in the sense that \NBG proves that
$$Cat((x,y,z,v,w))\to CAT(z)\wedge \Set(z)$$ 
and that if $CAT(X)\wedge Set(X)$, then there is a quintuple $(Ob(X),\mathcal I(X),X,dom, cod)$ and
$$CAT(X)\wedge Set(X)\to Cat(Ob(X),\mathcal I(X),X,dom, cod).$$
It also useful to observe that if $\C$ is a category such that $Ob\C$ is a set and each hom-class $\C(A,B)=\{f\in Mor\C\mid dom (f)=A\wedge cod (f)=B\}$ is a set, then \C is a small category. Indeed, there is then a set $Z$ containing all the hom-sets as members and $Mor\C=\bigcup Z$ (it is an axiom of \NBG that the sum class $\bigcup x$ of a set $x$ is a set).

Now $Cat(X)$ is a predicative formula and therefore, there is a class $Ob(Cat)$ of all sets $x$ such that $\Cat(x)$, i.e. the class of all small categories. This class is proper because the proper class of all sets injects in it.

We define a predicative formula $Fun(X;Y,Z)$ of \NBG  with two parameters $Y$ and $Z$ asserting that $X$ is a functor from $Y$ to $Z$. It basically says that $X$ is a triple $(Y,F,Z)$ such that $Cat(Y)\wedge Cat(Z)$ (so in particular $Y$ and $Z$ must be sets) and $F\colon Mor Y\to Mor Z$ is a function (and so it must also be a set) satisfying the usual axioms. Therefore, there is a set of all functors from $Y$ to $Z$, which is empty if $Y$ or $Z$ is a proper class\footnote{Recall that we are defining the 2-category of \emph{small} categories. There is a non trivial category of functors from a small category to a large category, by Proposition \hyperref[prop:NBG:functionclass]{\ref*{prop:NBG}(\ref*{prop:NBG:functionclass})}.}. 

Similarly, there is a predicative formula $Nat(X;A,B,C,D)$ saying $X$ is a natural transformation from $A$ to $B$, $A$ and $B$ being functors from $C$ to $D$. It affirms that $X$ is a triple $(A,\alpha,B)$ such that $Fun(A;C,D)\wedge Fun(B;C,D)$ and that $\alpha\colon Ob(C)\to Mor(D)$ satisfying axioms. In particular, $X$ must be set. Thus, there is a set of all natural transformations between two given functors.

\subsection{Some important facts about \textit{NBGI}}\index[not]{NBGI@\NBI}

Recall that \NBI is the axiomatic set theory \NBG plus an axiom of existence of a strongly inaccessible cardinal $\theta $, or equivalently, of a supercomplete model $M=V_{\theta+1}$ of \NBG, or of a Grothendieck universe $\U=V_{\theta}$. We use preferably the Grothendieck notation $\U$, but warn the reader that she should remember that it is $\power\U$ that is a model of \NBG, not \U. So the axioms \hyperref[eq:supercomplete]{(\ref*{eq:supercomplete})} of a supercomplete model apply to $\power\U$. Recall moreover that $\U$ is a set of \NBI and that $\power\U$ is a model of \NBG for the elementhood predicate $\in$ of \NBI. Classes of the model run over all subsets of \U (and consequently are all sets of \NBI). A class $X$ of the model is a set of the model if there exists $Y\in\power\U$ such that $X\in Y$, that is \ssi $X\in\U$ (see further down for an explanation on how to translate formulae of \NBG into formulae of \NBI interpreted in the model). Note that since $\power\U$ is transitive \hyperref[eq:supercomplete]{(\ref*{eq:supercomplete})}, \emph{all} elements of \U are also subsets of \U. Indeed, transitivity implies $\vdash x\in\U\wedge\U\in\power\U\to x\in\power\U$. Thus, classes of the model are all subsets of $\U$, whereas sets of the model are all elements of \U. 

Now the discussion in \autoref{ssec:NBG} about \NBG remains true in \NBI. One has two different notions of objects, sets and classes, and there exists a proper class of all sets. The inner concepts of classes and sets of the model are called $\U$-classes and \U-sets for now, but we will change the terminology later on. We keep the distinction between capital and lowercase letters for formulae of \NBG and their counterpart in the model, and use only capital letters for general formulae of \NBI. In other words, a lowercase letter in a formula of \NBI always indicates that it is an element of \U.

Translating a formula $\phi$ of \NBG in its interpretation $\phi^{(\U)}$ in \NBI by the model amounts to the following procedure. Replace all parts of it of the form $(\exists X \psi)$ by $$\exists X(X\subset\U\wedge \psi),$$ and therefore all parts of it of the form $\forall X \psi$ by $$\forall X(X\subset\U\to \psi).$$ Let us write the formula thus obtained by $\phi_{\U}$. The formula $\phi$ might not be closed and therefore one finally needs to make the following change: $\phi^{(\U)}$ is the formula $$(X_1\subset\U\wedge\ldots\wedge X_n\subset\U)\to\phi_{\U},$$ where $X_1,\ldots,X_n$ are the free variables of $\phi$. We have already seen that the translation of the formula asserting that a class is a set amounts to the formula $X\in\U$. Therefore, the translation just described can be completed to formulae containing lowercase letter by exactly the same process, but using the element symbol instead of the subset symbol.

As $\power\U$ is a model of \NBG, all translations of theorems of \NBG are true in $\power\U$. One can also translate all the definitions of \NBG, like the definitions of the constant $\emptyset$, the binary predicate symbol $\subset$ or the function symbol $\power$ (these symbols are not originally part of the language, but added by definition). This is where the particular properties of our model start to play a role. Indeed, they insure compatibilities between notions of \NBG defined in \NBI and their translation in the model. The transitivity of the model provides lots of coherence already between the model concepts and the surrounding theory concept, but supertransitivity is needed for extended coherence\footnote{Note that we have already implicitly used the supercompleteness of the model to show that it is of the form ``powerset of its universe of sets''.}. 

Firstly, note that transitivity of \U-classes induces transitivity of \U-sets. Now, let us consider the very beginning of the theory, the definition of equality of classes:
\begin{equation*}
X=Y\leftrightarrow\forall z(z\in X\leftrightarrow z\in Y).
\end{equation*}
Let us first consider any transitive model $M$. The interpretation of the definition of equality is $X=^{M}Y\leftrightarrow \forall z\in M(z\in X\leftrightarrow z\in Y)$. This is a priori a much weaker condition, as one has only to check that they have the same $M$-sets. But, by transitivity, $M$-classes only have $M$-sets and therefore
\begin{equation}\label{eq:equalitymodel}X\in M\wedge Y\in M\to(X=^MY\leftrightarrow X=Y).\end{equation}
A predicate symbol that has the same property as the equality symbol in \hyperref[eq:equalitymodel]{(\ref*{eq:equalitymodel})}, that is, its translation in $M$ when restricted to variables in $M$ is equivalent to the original notion in \NBG, is called \emph{\op$M$\fp-absolute}, and similarly for function symbols. One sees that absoluteness is very desirable for function and predicate symbols if one wants mathematics inside and outside the model to ``look the same''. Another simple example is the constant $\emptyset$. The null set axiom of \NBG must be true in a model $M$ and therefore, the following formula is a theorem of $\NBI$: $$(\exists! x\in M) (\forall y\in M)(y\notin x).$$ The defining axiom of \NBI for the interpretation $\emptyset^M$ in $M$ of $\emptyset$ is the formula $$[(\forall x\in M)(x\notin\emptyset^M)]\wedge\emptyset^M\in M.$$ Now this notion is absolute in a transitive model, that is the inner empty set $\emptyset^M$ of the model is equal to the empty set $\emptyset$ of the ambient \NBI. Let us prove that $\emptyset^M$ does not contain any set of \NBI, and therefore is equal to $\emptyset$. If $X\in\emptyset^M$, then $X\in M$ by transitivity, but then $X\notin\emptyset^M$ by definition, a contradiction. Thus $\emptyset^M=\emptyset$.

The inclusion predicate $\subset$ has also the same meaning in the model and in the ambient set theory. Indeed, the definition of the subset predicate symbol in \NBG (and \NBI) is
$$
X\subset Y\iff\forall Z(Z\in X\to Z\in Y).
$$
Its translation $\subset^M$ in the model $M$ is 
$$
X\subset^M Y\iff\forall Z(Z\in M \to (Z\in X\to Z\in Y)).
$$
Yet, by transitivity, if $X,Y\in M$ and $Z\notin M$, then $Z\notin X$ by transitivity. Thus, the predicate $\subset$ and $\subset^M$ have the same meaning when restricted to $M$-classes. In fact, one can say more: by supertransitivity \hyperref[eq:supercomplete]{(\ref*{eq:supercomplete})}, if $Y\in M$, then $X\in M\wedge X\subset^M Y\iff X\subset Y$. The subclasses of an $M$-class are the sub-$M$-classes. Such a predicate is called \emph{complete}: $X\subset Y\wedge Y\in M\to X\in M$. This also a desirable property, but it is not reasonable to require for all notions. Indeed, the notion of subset is not complete, as we shall see now.

The power class of a class\index{Class!Power --} is an example of a concept that is not \U-absolute, and this non-absoluteness appears in our text (see Definition \hyperref[def:site:covfunc]{\ref*{def:site}(\ref*{def:site:covfunc})}). Indeed, if $X$ is a proper \U-class, then $\power^{\U}X=\{z\in\U\mid z\subset X\}\neq\{Z\mid Z\subset X\}=\power X$ since the latter contains $X$ itself (it is a set of \NBI!), which is not an element of \U by definition of a proper \U-class (recall our convention to keep lowercase letters only for \U-sets). But this notion is \emph{set-absolute}, in the sense that the model and global notions of power class coincide when applied to \U-sets. Indeed, it is a theorem of \NBG that a subclass of a set is a set, and this translate into a theorem of \NBI about the model $\power \U$ that any sub-\U-class of a \U-set is a \U-set.

One goes similarly through all the predicate and function symbols defined in \NBG. This is done in detail in \cite{Sho67} for models of \ZFC in \ZFC, and one can adapt the proofs for our case, with some care though. Indeed, $\U$ is a model of \ZFC, not $\power\U$ and therefore, the proofs in this reference give a priori set-absoluteness, not general absoluteness. The power class is an example of a set-absolute notion that it not absolute. Another simple example is the pair $\langle X,Y\rangle$ of two classes $X$ and $Y$. This notion is set-absolute, but not absolute (recall \autoref{fnote:pairofclasses}, p. \pageref{fnote:pairofclasses}). Yet, set-absoluteness of pairs is enough for absoluteness of the Cartesian product of two classes, that is, if $X,Y\in\power\U$, then $X\times^\U Y=X\times Y$. Moreover, the ternary \emph{function predicate symbol} $Fnc(X,Y,Z)$, asserting that a class $Z$ is a function from $X$ to $Y$ is absolute and complete. Thus, for $X,Y\subset\U$, \begin{equation}\label{eq:funabcomp}Fnc(Z,X,Y)\leftrightarrow (Z\subset\U\wedge Fnc^\U(Z,X,Y)).\end{equation}
Moreover, if $x,y\in\U$ are $\U$-sets and $Fnc(Z,x,y)$, then $Z\subset\U$ by the preceding equation. Since it is a theorem of \NBG that a function between two sets is a set, $Z$ must be a \U-set, by the Interpretation Theorem (that says that all theorems of a theory are true in any of its interpretations in an other theory).

From now on, we change the terminology for the different kind of sets and classes. We will call the \U-sets and \U-classes, ``sets'' and ``classes'', and provide new names for the general sets and classes of \NBI. Indeed, in this work, we stick to plain \NBG as long as there is no need of further levels and in this axiomatic system, use the terminology of sets and classes. When working in \NBI, the philosophy is that all the sets we considered in \NBG can be recovered inside $\U$ because one cannot get out of it by usual set constructions, and therefore all classes one can think of will stay in $\power \U$. Another justification is that one needs to axiomatize the meta-language in such a way that our usual set and class theory \NBG admits a model in it. Thus, the mathematics we did before now takes place in the model.

\begin{Term}
	We use the following terminology in \NBI:
	\begin{itemize}
\item
\emph{Sets}\indexb{Set}: Elements of \U,
\item
\emph{Classes}\indexb{Class}: Subsets of \U,
\item
\emph{Proper classes}: Classes that are not sets,
\item
\emph{Conglomerates}\indexb{Conglomerate}: Sets of \NBI,
\item
\emph{Proper conglomerates}: Conglomerates that are not classes,
\item
\emph{Cartels}\indexb{Cartel}: Classes of \NBI,
\item
\emph{Proper Cartels}: Cartels that are not conglomerates. 
\end{itemize}
\end{Term}

\subsection{Category theory in \textit{NBGI}}
The discussion about the definition of a category in \NBG applies to \NBI. The difference now is that there are, among the conglomerates (that is, the sets of \NBI), the conglomerates called classes and sets that will play a special role.

\begin{Term}
	We use the following terminology in \NBI. On the right hand side, the word ``category'' means ``category in \NBI''.
	\begin{itemize}
\item
\emph{Small category}\indexb{Category!Small --}: Category whose morphism cartel is a set.
\item
\emph{Category}\indexb{Category}: Category whose morphism cartel is a class.
\item
\emph{Large Category}\indexb{Category!Large --}: Category whose morphism cartel is a proper class.
\item
\emph{Extra large category \textnormal{or} XL-category}\indexb{Category!Extra large --}\indexb{Category!XL-{}--}: Category whose morphism cartel is a proper conglomerate.
\item
\emph{Extra extra large category \textnormal{or} XXL-category}\indexb{Category!Extra extra large --}\indexb{Category!XXL-{}--}: Category whose morphism cartel is a proper cartel.
\end{itemize}
\end{Term}
We denote by \XCAT\index[not]{CAT2@\XCAT} the 2-XXL-category of categories in \NBI whose morphism cartel is a conglomerate (thus it is the ``large category of small categories'' in \NBI). Recall that it is defined, via the Class Existence Theorem, as the class of all conglomerates $X$ satisfying some formula (of \NBG) that we had denoted $\Cat(X)$. We change the notation, and write it $\XCAT(X)$. It is a set-absolute notion, that is an $\U$-set $x$ satisfies $\XCAT^\U(x)$ if and only if it satisfies $\XCAT(x)$. Thus, the interpretation of the notion of a small category in the model $\power\U$ is equal to the notion of small category in \NBI. In particular, all theorems about small categories proved in \NBG apply to small categories of \NBI. The formulae of \NBG saying that a class is a functor between two small categories, or a natural transformation between two such functors are also set-absolute and moreover set-complete, in the sense that a conglomerate (in fact even a cartel) that is a functor between two small categories is automatically an \U-set (see \hyperref[eq:funabcomp]{(\ref*{eq:funabcomp})}) and the discussion that follows). We denote \Cat\index[not]{Cat0@\Cat} the large full sub-2-category of \XCAT consisting of small categories in \NBI.

The formula $\XCAT(X)$ is of course not absolute, because in \NBG, pairs of proper classes have a trivial definition, and thus, the interpretation of the definition of pairs in the model is also trivial for proper classes. Again, proper classes are conglomerates, and therefore the notion of pair of proper classes is non trivial in \NBI. Nevertheless, recall that there is a definition of (large) categories by means of a formula $\CAT(X)$ of \NBG. This notion is absolute. Thus, theorems of \NBG about (large) categories apply to (large) categories in \NBI. But now, unlike \NBG, \NBI has a conglomerate of all classes satisfying $\CAT(X)$. Moreover, a functor in \NBI between two large categories is a class (the notion is complete) and a functor in the model sense (the notion is absolute). One can say the same of natural transformation. There is therefore a full sub-2-XL-category of \XCAT consisting of large categories, which we denote \CAT\index[not]{CAT1@\CAT}. Note that \Cat is thus a full sub-2-category of \CAT. 

In conclusion, one has the following 2-categories in \NBI. \XCAT is the 2-XXL-category of categories in \NBI whose morphism cartel is a conglomerate (thus it is the ``large category of small categories'' in \NBI). \CAT is its full sub 2-XL-category of categories. \Cat is its full sub large 2-category of small categories. In the following, each inclusion is thus full
$$
\Cat\subset\CAT\subset\XCAT.
$$
Moreover, all theorems of \NBG about small and large categories apply to small and large categories of \NBI.

\end{appendix}

\cleardoublepage
\phantomsection

\clearpage
\phantomsection
\addcontentsline{toc}{chapter}{Glossary of notations}
\printindex[not]

\clearpage
\phantomsection
\addcontentsline{toc}{chapter}{Index}
\printindex

\end{document}